\pdfoutput=1
\documentclass[12pt, a4paper, twoside,leqno]{amsart}
\usepackage[centering, totalwidth = 435pt, totalheight = 640pt]{geometry}
\usepackage{amssymb, amsmath, amsthm, enumerate, microtype, stmaryrd, url,mathtools,tikz,rotating}
\usepackage[latin1]{inputenc}
\usepackage[arrow, matrix, tips, curve,graph]{xy}
\usepackage{color}
\definecolor{darkgreen}{rgb}{0,0.45,0}
\usepackage[pagebackref,colorlinks,citecolor=darkgreen,linkcolor=darkgreen]{hyperref}
\usepackage{newclude}
\usetikzlibrary{positioning}
\SelectTips{cm}{10}

 \DeclareMathOperator{\ob}{ob}

\newcommand{\cat}[1]{\mathbf{#1}}
\newcommand{\op}{\mathrm{op}}
\newcommand{\id}{\mathrm{id}}
\newcommand{\thg}{{\mathord{\text{--}}}}

\newcommand{\spn}[1]{{\langle{#1}\rangle}}

\newcommand{\defeq}{\mathrel{\mathop:}=}

\newcommand{\cd}[2][]{\vcenter{\hbox{\xymatrix#1{#2}}}}

\newcommand{\botimes}{\mathrel{\bar\otimes}}

\renewcommand{\phi}{\varphi}
\newcommand{\A}{{\mathcal A}}
\newcommand{\B}{{\mathcal B}}
\newcommand{\C}{{\mathcal C}}
\newcommand{\D}{{\mathcal D}}
\newcommand{\E}{{\mathcal E}}
\newcommand{\F}{{\mathcal F}}

\newcommand{\I}{{\mathcal I}}

\newcommand{\K}{{\mathcal K}}
\renewcommand{\L}{{\mathcal L}}
\newcommand{\M}{{\mathcal M}}

\renewcommand{\O}{{\mathcal O}}

\newcommand{\U}{{\mathcal U}}
\newcommand{\V}{{\mathcal V}}
\newcommand{\W}{{\mathcal W}}

\newcommand{\ti}{_\tau}
\newcommand{\lo}{_\lambda}

\newcommand{\f}[1]{\mathsf{#1}}
\makeatletter
\newcommand{\fcat}[1]{\@fcat#1\end@fcat}
\def\@fcat#1#2\end@fcat{\mathbb{#1}\mathbf{#2}}

\newcommand{\xtor}[1]{\cdl[@1]{{} \ar[r]|-{\object@{|}}^{#1} & {}}}
\newcommand{\tor}{\ensuremath{\relbar\joinrel\mapstochar\joinrel\rightarrow}}

\def\hookleftarrowfill@{\arrowfill@\leftarrow\relbar{\relbar\joinrel\rhook}}
\def\twoheadleftarrowfill@{\arrowfill@\twoheadleftarrow\relbar\relbar}
\def\leftbararrowfill@{\arrowdoublefill@{\leftarrow\mkern-5mu}\relbar\mapstochar\relbar\relbar}
\def\Leftbararrowfill@{\arrowdoublefill@{\Leftarrow\mkern-2mu}\Relbar\Mapstochar\Relbar\Relbar}
\def\leftringarrowfill@{\arrowdoublefill@{\leftarrow\mkern-3mu}\relbar{\mkern-3mu\circ\mkern-2mu}\relbar\relbar}
\def\lefttriarrowfill@{\arrowfill@{\mathrel\triangleleft\mkern0.5mu\joinrel\relbar}\relbar\relbar}
\def\Lefttriarrowfill@{\arrowfill@{\mathrel\triangleleft\mkern1mu\joinrel\Relbar}\Relbar\Relbar}

\def\hookrightarrowfill@{\arrowfill@{\lhook\joinrel\relbar}\relbar\rightarrow}
\def\twoheadrightarrowfill@{\arrowfill@\relbar\relbar\twoheadrightarrow}
\def\rightbararrowfill@{\arrowdoublefill@{\relbar\mkern-0.5mu}\relbar\mapstochar\relbar\rightarrow}
\def\Rightbararrowfill@{\arrowdoublefill@{\Relbar\mkern-2mu}\Relbar\Mapstochar\Relbar\Rightarrow}
\def\rightringarrowfill@{\arrowdoublefill@\relbar\relbar{\mkern-2mu\circ\mkern-3mu}\relbar{\mkern-3mu\rightarrow}}
\def\righttriarrowfill@{\arrowfill@\relbar\relbar{\relbar\joinrel\mkern0.5mu\mathrel\triangleright}}
\def\Righttriarrowfill@{\arrowfill@\Relbar\Relbar{\Relbar\joinrel\mkern1mu\mathrel\triangleright}}

\def\leftrightarrowfill@{\arrowfill@\leftarrow\relbar\rightarrow}
\def\mapstofill@{\arrowfill@{\mapstochar\relbar}\relbar\rightarrow}

\renewcommand*\xleftarrow[2][]{\ext@arrow 20{20}0\leftarrowfill@{#1}{#2}}
\providecommand*\xLeftarrow[2][]{\ext@arrow 60{22}0{\Leftarrowfill@}{#1}{#2}}
\providecommand*\xhookleftarrow[2][]{\ext@arrow 10{20}0\hookleftarrowfill@{#1}{#2}}
\providecommand*\xtwoheadleftarrow[2][]{\ext@arrow 60{20}0\twoheadleftarrowfill@{#1}{#2}}
\providecommand*\xleftbararrow[2][]{\ext@arrow 10{22}0\leftbararrowfill@{#1}{#2}}
\providecommand*\xLeftbararrow[2][]{\ext@arrow 50{24}0\Leftbararrowfill@{#1}{#2}}
\providecommand*\xleftringarrow[2][]{\ext@arrow 10{26}0\leftringarrowfill@{#1}{#2}}
\providecommand*\xlefttriarrow[2][]{\ext@arrow 80{24}0\lefttriarrowfill@{#1}{#2}}
\providecommand*\xLefttriarrow[2][]{\ext@arrow 80{24}0\Lefttriarrowfill@{#1}{#2}}

\renewcommand*\xrightarrow[2][]{\ext@arrow 01{20}0\rightarrowfill@{#1}{#2}}
\providecommand*\xRightarrow[2][]{\ext@arrow 04{22}0{\Rightarrowfill@}{#1}{#2}}
\providecommand*\xhookrightarrow[2][]{\ext@arrow 00{20}0\hookrightarrowfill@{#1}{#2}}
\providecommand*\xtwoheadrightarrow[2][]{\ext@arrow 03{20}0\twoheadrightarrowfill@{#1}{#2}}
\providecommand*\xrightbararrow[2][]{\ext@arrow 01{22}0\rightbararrowfill@{#1}{#2}}
\providecommand*\xRightbararrow[2][]{\ext@arrow 04{24}0\Rightbararrowfill@{#1}{#2}}
\providecommand*\xrightringarrow[2][]{\ext@arrow 01{26}0\rightringarrowfill@{#1}{#2}}
\providecommand*\xrighttriarrow[2][]{\ext@arrow 07{24}0\righttriarrowfill@{#1}{#2}}
\providecommand*\xRighttriarrow[2][]{\ext@arrow 07{24}0\Righttriarrowfill@{#1}{#2}}

\providecommand*\xmapsto[2][]{\ext@arrow 01{20}0\mapstofill@{#1}{#2}}
\providecommand*\xleftrightarrow[2][]{\ext@arrow 10{22}0\leftrightarrowfill@{#1}{#2}}
\providecommand*\xLeftrightarrow[2][]{\ext@arrow 10{27}0{\Leftrightarrowfill@}{#1}{#2}}

\makeatother
\makeatletter
\def\matrixobject@{%
 \edef \next@{={\DirectionfromtheDirection@ }}%
 \expandafter \toks@ \next@ \plainxy@
 \let\xy@@ix@=\xyq@@toksix@
 \xyFN@ \OBJECT@}
\let\xy@entry@@norm=\entry@@norm
\def\entry@@norm@patched{%
 \let\object@=\matrixobject@
 \xy@entry@@norm }
\AtBeginDocument{\let\entry@@norm\entry@@norm@patched}
\makeatother

\newcommand{\twocong}[2][0.5]{\ar@{}[#2] \save ?(#1)*{\cong}\restore}
\newcommand{\twoeq}[2][0.5]{\ar@{}[#2] \save ?(#1)*{=}\restore}
\newcommand{\rtwocell}[3][0.5]{\ar@{}[#2] \ar@{=>}?(#1)+/l 0.175cm/;?(#1)+/r 0.175cm/^{#3}}
\newcommand{\rtwocello}[3][0.5]{\ar@{}[#2] \ar@{=>}?(#1)+/l 0.175cm/;?(#1)+/r 0.175cm/_(0.35){#3}}
\newcommand{\ltwocell}[3][0.5]{\ar@{}[#2] \ar@{=>}?(#1)+/r 0.175cm/;?(#1)+/l 0.175cm/^{#3}}
\newcommand{\ltwocello}[3][0.5]{\ar@{}[#2] \ar@{=>}?(#1)+/r 0.175cm/;?(#1)+/l 0.175cm/_{#3}}
\newcommand{\dtwocell}[3][0.5]{\ar@{}[#2] \ar@{=>}?(#1)+/u  0.175cm/;?(#1)+/d 0.175cm/^{#3}}
\newcommand{\dtwocello}[3][0.5]{\ar@{}[#2] \ar@{=>}?(#1)+/u  0.175cm/;?(#1)+/d 0.175cm/_{#3}}
\newcommand{\dltwocell}[3][0.5]{\ar@{}[#2] \ar@{=>}?(#1)+/ur  0.175cm/;?(#1)+/dl 0.175cm/^{#3}}
\newcommand{\drtwocell}[3][0.5]{\ar@{}[#2] \ar@{=>}?(#1)+/ul  0.175cm/;?(#1)+/dr 0.175cm/^{#3}}
\newcommand{\dthreecell}[3][0.5]{\ar@{}[#2] \ar@3{->}?(#1)+/u  0.175cm/;?(#1)+/d 0.175cm/^{#3}}
\newcommand{\utwocell}[3][0.5]{\ar@{}[#2] \ar@{=>}?(#1)+/d 0.175cm/;?(#1)+/u 0.175cm/_{#3}}
\newcommand{\urtwocell}[3][0.5]{\ar@{}[#2] \ar@{=>}?(#1)+/dl 0.175cm/;?(#1)+/ur 0.175cm/_{#3}}
\newcommand{\dtwocelltarg}[3][0.5]{\ar@{}#2 \ar@{=>}?(#1)+/u  0.175cm/;?(#1)+/d 0.175cm/^{#3}}
\newcommand{\utwocelltarg}[3][0.5]{\ar@{}#2 \ar@{=>}?(#1)+/d  0.175cm/;?(#1)+/u 0.175cm/_{#3}}

\newdir{(}{{}*!<0em,-.14em>-\cir<.14em>{l^r}}
\newdir{ (}{{}*!/-5pt/\dir{(}}
\newdir{ >}{{}*!/-5pt/\dir{>}}
\newcommand{\sh}[2]{**{!/#1 -#2/}}

\theoremstyle{definition}

\swapnumbers
\theoremstyle{plain}
\newtheorem{Thm}[subsection]{Theorem}
\newtheorem{Prop}[subsection]{Proposition}
\newtheorem{Cor}[subsection]{Corollary}
\newtheorem{Lemma}[subsection]{Lemma}

\numberwithin{equation}{section}

\theoremstyle{definition}

\newtheorem{Ex}[subsection]{Example}

\newtheorem{Rk}[subsection]{Remark}

\newcommand{\Lan}{\mathrm{Lan}}

\newcommand{\lu}{\sigma}
\newcommand{\ru}{\tau}
\newcommand{\ass}{\alpha}
\newcommand{\fu}{\iota}
\newcommand{\fm}{\mu}
\newcommand{\h}[2]{\spn{#1,#2}}
\renewcommand{\t}{}
\newcommand{\laxu}{\mathfrak{i}}
\newcommand{\laxc}{\mathfrak{x}}

\title[Enriched categories as a free cocompletion]{Enriched categories as a free cocompletion}
\author{Richard Garner}
\address{Department of Computing, Macquarie University, North Ryde, NSW 2109, Australia}
\email{richard.garner@mq.edu.au}
\author{Michael Shulman}
\address{Department of Mathematics, Institute for Advanced Study, Princeton, NJ 08540, USA}
\email{mshulman@ias.edu}
\date{\today}
\thanks{The first author acknowledges the support of Australian Research Council Discovery Project grants DP110102360 and and DP130101969. The second author acknowledges the support of a United States National Science Foundation Postdoctoral Fellowship and a grant under agreement No. DMS-1128155.  Any opinions, findings, and conclusions or recommendations expressed in this material are those of the authors and do not necessarily reflect the views of the National Science Foundation.}

\begin{document}
 \leftmargini=2em

\begin{abstract}
This paper has two objectives. The first is to develop the theory of bicategories enriched in a monoidal bicategory---categorifying the classical theory of categories enriched in a monoidal category---up to a description of the free cocompletion of an enriched bicategory under a class of weighted bicolimits. The second objective is to describe a universal property of the process assigning to a monoidal category $\V$ the equipment of $\V$-enriched categories, functors, transformations, and modules; we do so by considering, more generally, the assignation sending an equipment $\C$ to the equipment of $\C$-enriched categories, functors, transformations, and modules, and exhibiting this as the free cocompletion of a certain kind of enriched bicategory under a certain class of weighted bicolimits.
\end{abstract}

\maketitle
\tableofcontents

\section{Introduction}
\label{sec:introduction}

The classical theory of categories enriched in a monoidal category~\cite{kelly:enriched} has many applications throughout mathematics.
The more general notion of a category enriched in a bicategory is less well-known, but it allows one to capture also \emph{internal} categories and \emph{indexed} categories through enrichment, and has been used in the study of sheaves and stacks~\cite{walters:sheaves-cauchy-1,bcsw:variation-enr,Street2005Enriched}.
More generally still, we can enrich categories in a double category or a proarrow equipment~\cite{wood:proarrows-i,leinster:gen-enr-cats}; the advantage of this over bicategory-enrichment is a better notion of enriched functor (see~\cite{shulman:eicats,CG2012Restriction} for some examples).

In this paper we do two things:
\begin{enumerate}
\item We categorify the theory of enriched categories to a theory of \emph{bicategories} enriched in a monoidal bicategory, or more generally in a tricategory.\label{item:enrbicat}\vskip0.25\baselineskip
\item We show that the construction ``$\C \mapsto$ categories enriched in $\C$'', for a bicategory or equipment $\C$, has a \emph{universal property}.\label{item:univprop}
\end{enumerate}
While these objectives are perhaps seemingly unrelated, in fact the former is necessary for the latter: the universal property of enriched categories is expressed as a free cocompletion of a certain kind of enriched bicategory.
This can be regarded as an instance of what Baez and Dolan~\cite{Baez1998Higher-Dimensional} term the \emph{microcosm principle}: the proper context in which to consider the theory of enriched categories is a categorified version of itself.

We now discuss~(\ref{item:enrbicat}) and~(\ref{item:univprop}) separately in somewhat more detail, beginning with~(\ref{item:enrbicat}).

\subsection{Enriched bicategories}
\label{sec:enrbicat}

Recall that if $\V$ is a monoidal category, then  a \emph{$\V$-enriched category} (or \emph{$\V$-category}) $\C$ comprises a collection of objects; for every pair of objects,  a hom-object $\C(x,y) \in \V$; and morphisms $\C(y,z) \otimes \C(x,y) \to \C(x,z)$   and $I \to \C(x,x)$ giving composition and identities. The associativity and unitality of composition in $\C$ is expressed through the commutativity of certain familiar diagrams in $\V$. By taking $\V$ to be $\cat{Set}$, $k\text-\cat{Vect}$, $\cat{SSet}$, or $\cat{DG}\text-R\text-\cat{Mod}$, for example, we recapture the notions of (locally small) category, $k$-linear category, simplicial category, and dg-category, respectively.

Certain kinds of higher-categorical structures can also be described using enriched categories: thus (locally small) $2$-categories are $\cat{Cat}$-categories, while semi-strict $3$-categories---to which every tricategory is equivalent---are precisely $\cat{Gray}$-categories, where as in~\cite{Gordon1995Coherence}, $\cat{Gray}$ is the category of $2$-categories equipped with the Gray tensor product. However, B{\'{e}nabou's \emph{bicategories}~\cite{benabou:bicats} cannot be described in this manner, since composition in an enriched category is always strictly associative and unital, whereas that in a bicategory is associative and unital only up to coherent isomorphism.

This last observation suggests the existence of a more general theory of enrichment, which stands in the same relation to the notion of bicategory as does the theory of $\V$-categories to the notion of ordinary category. An early development of such a theory was given by Bozapalid\`es in a series of papers deriving from his Ph.\,D.~thesis~\cite{Bozapalides1976Theorie}. However, at the time of his writing, the most appropriate kind of base $\V$ over which such a bicategory should be enriched had not yet been developed. This was achieved in~\cite{Gordon1995Coherence}: a \emph{monoidal bicategory} is a bicategory equipped with a tensor product~$\otimes$ that is associative and unital only up to \emph{equivalence}, together with higher cell data witnessing the coherence of these equivalences.
The first steps in the theory of bicategories enriched in a monoidal bicategory $\V$ were given in the two theses~\cite{Carmody1995Cobordism} and~\cite{Lack1995The-algebra}. A $\V$-enriched bicategory involves, as before, a collection of objects, a collection of hom-objects from $\V$, and composition and identity $1$-cells, but unlike before, the diagrams expressing the associativity and unitality of composition are no longer required to commute on the nose, but only up to coherent invertible $2$-cells of~$\V$. 

For example, taking $\V$ to be the cartesian monoidal bicategory $\cat{Cat}$, we recapture the notion of (locally small) bicategory. Less trivially, we may consider various $2$-dimensional analogues of the notion of $\cat{CMon}$- or $\cat{Ab}$-enriched category. For example, we can consider, as in~\cite{Lack1995The-algebra}, bicategories each of whose hom-categories admits finite coproducts preserved by composition with $1$-cells on either side. These are precisely $\V$-enriched bicategories when $\V$ is taken to be the  bicategory of categories with finite coproducts and finite-coproduct-preserving functors, under the tensor product that classifies functors preserving finite coproducts in each variable separately. Another possible generalisation involves bicategories each of whose hom-categories is a \emph{symmetric $2$-group}---a compact closed symmetric monoidal groupoid---and each of whose whiskering functors is coherently monoidal. These are $\V$-enriched bicategories for $\V$ the monoidal bicategory of symmetric $2$-groups, as defined in~\cite{Dupont2008Abelian}, for example. We shall see many more examples of enriched bicategories throughout this paper.

As we have said, the very basic definitions in the theory of enriched bicategories (amounting to most of our Section~\ref{sec:vbicats} and some of Section~\ref{sec:tricat-vbicats}) were given in~\cite{Carmody1995Cobordism} and~\cite{Lack1995The-algebra}. However, for the application of the theory that constitutes the second main objective of this article, we will need somewhat more than this: we must define and characterise the free cocompletion of an enriched bicategory under a class of colimits. Doing so will involve, along the way, developing a theory of \emph{modules} (a.k.a.\ profunctors or distributors) between enriched bicategories, including a construction of the tensor product and internal hom of such; appropriate versions of the Yoneda lemma; and some results concerning iterated colimits and left Kan extensions.

\subsection{Categories as monads}
\label{sec:cats-as-monads}
We now turn to a discussion of the second main objective of this paper, which is to use the theory of enriched bicategories to exhibit a universal property of the assignation ``$\C \mapsto$ categories enriched in $\C$'', for a bicategory or equipment $\C$. The universal property that we will present is the culmination of a long sequence of advances by many people. It begins with B\'enabou~\cite{benabou:bicats}, who observed that small categories can be identified with \emph{monads} in the bicategory of spans of sets.

Recall that a monad in a bicategory $\B$ consists of an object $A\in\B$, a morphism $t \colon A\to A$, and $2$-cells $tt\Rightarrow t$ and $1_A\Rightarrow t$ satisfying the usual associativity and unitality laws.
In the bicategory $\cat{Cat}$, this is a monad in the ordinary sense;
but in the bicategory $\cat{Span}$, a monad consists of a set $A_0$, a span $A_0 \leftarrow A_1 \to A_0$, and functions $A_1 \times_{A_0} A_1 \to A_1$ and $A_0 \to A_1$, satisfying axioms that state precisely that it is a small category.
More generally, if $\cat E$ is a category with pullbacks, then monads in $\cat{Span}(\cat E)$ are categories internal to $\cat E$.

The notion of monad in a bicategory is nice and general and has a good formal theory~\cite{street:ftm}.
A monad in $\B$ is equivalently a lax functor $1\to \B$, and the lax limit and lax colimit of this functor give abstract versions of the classical Eilenberg-Moore and Kleisli categories, referred to as \emph{Eilenberg-Moore objects} and \emph{Kleisli objects} respectively.
Thus, it is pleasing to see categories themselves arise as instances of monads---or it would be, if the identification did not break down at higher levels.

There is an obvious notion of \emph{lax morphism} between monads $(A,t)$ and $(B,s)$ in a bicategory, consisting of a morphism $f\colon A\to B$ and a $2$-cell $sf\Rightarrow ft$ satisfying some axioms.
These are the lax transformations between lax functors $1\to\B$, and induce morphisms between Eilenberg-Moore objects.
Dually we have \emph{colax morphisms}, involving a $2$-cell $ft \Rightarrow sf$; these can be identified with colax transformations, and induce morphisms of Kleisli objects.

However, neither sort of morphism of monads in $\cat{Span}$ gives what we expect as a morphism of categories.
A colax morphism (the more likely-looking one) from $A_1 \rightrightarrows A_0$ to $B_1 \rightrightarrows B_0$ involves a span $A_0 \leftarrow F_0 \to B_0$ and a function between pullbacks:
\begin{equation}\label{eq:mealy-intro}
  \vcenter{\hbox{\begin{tikzpicture}[->,scale=1.4]
    \node (A0a) at (0,2) {$A_0$};
    \node (A1) at (1,2) {$A_1$};
    \node (A0b) at (2,2) {$A_0$};
    \node (B0a) at (0,0) {$B_0$};
    \node (B1) at (1,0) {$B_1$};
    \node (B0b) at (2,0) {$B_0$\rlap{ .}};
    \draw (A1) -- (A0a); \draw (A1) -- (A0b);
    \draw (B1) -- (B0a); \draw (B1) -- (B0b);
    \node (F0a) at (0,1) {$F_0$};
    \node (F0b) at (2,1) {$F_0$};
    \draw (F0a) -- (B0a); \draw (F0b) -- (B0b);
    \node (pbb) at (1.3,1.3) {$\bullet$};
    \draw (F0a) -- (A0a); \draw (F0b) -- (A0b);
    \node (pba) at (.7,.7) {$\bullet$};
    \draw (pba) -- (F0a); \draw (pba) -- (B1);
    \draw[-] (.5,.6) -- (.5,.5) -- (.6,.5);
    \draw (pbb) -- (F0b);
    \draw[-] (1.5,1.4) -- (1.5,1.5) -- (1.4,1.5);
    \draw (pbb) -- (A1);
    \draw (pbb) -- (pba);
  \end{tikzpicture}}}
\end{equation}
If the first leg $F_0 \to A_0$ of the span is an identity, then these data reduce to functions $A_0 \to B_0$ and $A_1\to B_1$, which the axioms then assert to be a functor.
However, the general case does not correspond to any well-known notion of morphism between categories (but see Example~\ref{ex:mealy} and Remark~\ref{remark:mealy}).

The situation is even worse for $2$-cells.
A monad $2$-cell (which can be identified with a modification) consists of a $2$-cell $f\Rightarrow g$ satisfying some axioms.
But if $f$ and $g$ are functors, regarded as particular colax monad morphisms in $\cat{Span}$, then this is a function $F_0 \to G_0$ commuting with the identities $F_0 = A_0$ and $G_0=A_0$, hence itself merely an identity.
Thus we do not see the natural transformations at all.

This latter problem was rectified by Lack and Street~\cite{ls:ftm2} by considering the \emph{free cocompletion} of a bicategory $\C$ under Kleisli objects, $\cat{KL}(\C)$. (In fact, Lack and Street define $\cat{KL}(\C)$ only for $\C$ a $2$-category; but the description they give adapts without trouble to the bicategorical case.)
Up to equivalence, the objects of $\cat{KL}(\C)$ can be identified with monads in $\C$, and it turns out that its morphisms are simply colax monad morphisms.
However, a $2$-cell $f\Rightarrow g$ in $\cat{KL}(\C)$ consists instead of a $2$-cell $f\Rightarrow s g$ in $\C$ satisfying certain axioms.
In $\cat{Span}$, this datum is a function:
\begin{center}
  \begin{tikzpicture}[->,scale=1.4]
    \node (A0) at (1,2) {$A_0$};
    \node (F0) at (2,1) {$F_0$};
    \node (G0) at (0,1) {$G_0$};
    \node (B0a) at (0,0) {$B_0$};
    \node (B0b) at (2,0) {$B_0$\rlap{ .}};
    \node (B1) at (1,0) {$B_1$};
    \draw (F0) -- (A0); \draw (G0) -- (A0);
    \draw (F0) -- (B0b); \draw (G0) -- (B0a);
    \draw (B1) -- (B0a); \draw (B1) -- (B0b);
    \node (pb) at (.7,.7) {$\bullet$};
    \draw (pb) -- (G0); \draw (pb) -- (B1);
    \draw (F0) -- (pb);
    \draw[-] (.5,.6) -- (.5,.5) -- (.6,.5);
  \end{tikzpicture}
\end{center}
If $F$ and $G$ are functors, so that $F_0 \to A_0$ and $G_0 \to A_0$ are identities, then this is simply a function $A_0 \to B_1$, and the axioms assert indeed that this is a natural transformation.
Thus, the $2$-category $\cat{Cat}$ is a locally full sub-$2$-category of $\cat{KL} (\cat{Span})$.

There is still the problem of identifying the functors.
Towards this end, note that there is a functor $\cat{Set}\to\cat{Span}$ that is the identity on objects and that sends a function $A\to B$ to the span $1_A \colon A \leftarrow A \to B \colon f$.
This functor is moreover locally fully faithful, and each span of the form $(1_A, f)$ has a right adjoint in $\cat{Span}$, namely $(f, 1_A)$.
Thus, $\cat{Set}\to\cat{Span}$ is a \emph{proarrow equipment} in the sense of Wood~\cite{wood:proarrows-i}.

This at least gives us an abstract version of the construction: starting from a proarrow equipment $\K\to\M$, we can form $\cat{KL}(\M)$, and then its locally full sub-bicategory $\cat{KL}_\K(\M)$ on the morphisms whose underlying $\M$-morphism lies in $\K$.
It also suggests an improvement we might hope for, since $\cat{Cat}$ is itself part of a proarrow equipment $\cat{Cat}\to\cat{Mod}$.
Here $\cat{Mod}$ is the bicategory whose objects are small categories and whose morphisms are modules.
Indeed, this is the archetypical proarrow equipment that Wood sought to generalize.

Thus, we might hope for a general construction on proarrow equipments that when applied to $\cat{Set}\to\cat{Span}$ produces $\cat{Cat}\to\cat{Mod}$.
Such a construction is easy to write down: from $\K\to\M$ we produce $\cat{KL}_\K(\M) \to \cat{Mod}_1(\M)$, where the objects of $\cat{Mod}_1(\M)$ are monads as before, but its morphisms are a suitable notion of module.
(The subscript $1$ will be explained below.) 

Specifically, a \emph{module} from a monad $(A,t)$ to $(B,s)$ in a bicategory is a morphism $h \colon A\to B$ together with $2$-cells $ht\Rightarrow h$ and $sh\Rightarrow h$  
giving a compatible right action of $t$ and left action of $s$ on $h$.
In $\cat{Span}$, these are precisely the usual sort of bimodules between categories.
If we assume that $\M$ has \emph{local reflexive coequalisers}---that is, its hom-categories have reflexive coequalisers that are preserved by composition in each variable---then we can compose modules with a ``tensor product'', obtaining a bicategory $\cat{Mod}_1(\M)$ and a proarrow equipment $\cat{KL}_\K(\M) \to \cat{Mod}_1(\M)$.

The construction of $\cat{Mod}_1(\M)$ from $\M$ was studied abstractly by Street~\cite{street:cauchy-enr} and Carboni, Kasangian, and Walters~\cite{ckw:axiom-mod}.
They showed that it is \emph{idempotent}: $\cat{Mod}_1(\cat{Mod}_1(\M))\simeq \cat{Mod}_1(\M)$.
Moreover, a bicategory $\M$ is of the form $\cat{Mod}_1(\C)$ for some $\C$ (which can then be taken to be $\M$ itself) if and only if it has local reflexive coequalisers and Kleisli objects.
These two facts suggest that $\cat{Mod}_1(\M)$ is also some kind of completion of $\M$ under Kleisli objects, but in a different sense than $\cat{KL}(\M)$.

In this paper we unify these various threads, by observing that:
\begin{enumerate}[(a)]
\item the property of having local reflexive coequalisers, and\label{item:lc}
\item the structure enhancing a bicategory to a proarrow equipment\label{item:pe}
\end{enumerate}
can both be described as \emph{enrichments} of a bicategory in particular monoidal bicategories.
For~(\ref{item:lc}), this is much as in Section~\ref{sec:enrbicat} above: we enrich in the monoidal bicategory $\cat{Colim}_1$ of categories with reflexive coequalisers, with a tensor product that represents functors preserving reflexive coequalisers in each variable (again, the naming of this $2$-category and the subscript $1$ will be explained below).

We obtain~(\ref{item:pe}) by enriching in a monoidal bicategory $\F$ whose objects are pairs of categories together with a fully faithful functor between them.
Each hom-object of an $\F$-bicategory $\C$ is such a functor $\C_\tau(x,y) \to \C_\lambda(x,y)$, and the domains and codomains of these functors fit together into two bicategories $\C_\tau$ and $\C_\lambda$ with a locally fully faithful, identity-on-objects functor between them. This generalises constructions given for $1$-categories by Power~\cite{Power2002Premonoidal}, and for $2$-categories by Lack and Shulman~\cite{ls:limlax}.\footnote{In fact, proarrow equipments in the sense of~\cite{wood:proarrows-i} are rather special kinds of $\F$-bicategories: those for which every $1$-cell in the image of the functor $\C_\tau \to \C_\lambda$ is a \emph{map}: that is, admits a right adjoint. Our construction actually operates on the more general $\F$-bicategories, but has the property of sending proarrow equipments to proarrow equipments; see Lemma~\ref{lemma:fmod-preserves-equip}. In the body of this paper, we shall use ``equipment'' to mean ``$\F$-bicategory'', and ``map equipment'' to mean ``proarrow equipment'' in the sense of Wood.}

Finally, we can combine both of these structures, by enriching in a monoidal bicategory whose objects are fully faithful functors whose codomains have reflexive coequalisers; call this bicategory $\F_1$.
Then we will prove:

\begin{Thm}\label{thm:intro1}
  For any proarrow equipment $\K\to\M$, where $\M$ has local reflexive coequalisers, the proarrow equipment $\cat{KL}_\K(\M) \to \cat{Mod}_1(\M)$ is its free cocompletion, as an $\F_1$-enriched bicategory, under a class of $\F_1$-enriched colimits called \emph{tight Kleisli objects}.
\end{Thm}

In particular, $\cat{Cat}\to\cat{Mod}$ is obtained by freely cocompleting $\cat{Set}\to\cat{Span}$ in this manner.
We also have the following simpler version:

\begin{Thm}\label{thm:intro2}
  For any bicategory $\M$ with local reflexive coequalisers, the bicategory $\cat{Mod}_1(\M)$ is its free cocompletion under Kleisli objects as a $\cat{Colim}_1$-enriched bicategory.
  Moreover, Kleisli objects are an \emph{absolute colimit} for $\cat{Colim}_1$.
\end{Thm}

The latter fact explains the idempotence of $\cat{Mod}_1$, since any cocompletion under an absolute type of colimit (e.g.\ splitting of idempotents, or biproducts in additive categories) is idempotent.
Tight Kleisli objects, however, are not absolute for $\F_1$.

Finally, everything we have said so far has a ``many-object'' version.
Already in~\cite{benabou:bicats}, B\'enabou considered what he called \emph{polyads}, and which later authors have come to call \emph{categories enriched in a bicategory}.\footnote{These should not be confused with the \emph{bicategories} enriched in a \emph{monoidal} bicategory that we discussed in Section~\ref{sec:enrbicat}.
  Categories enriched in a bicategory are still a $1$-categorical notion; the corresponding two-dimensional notion is that of a bicategory enriched in a tricategory, which has as a special case the bicategories enriched in a monoidal bicategory discussed in this paper.
}
A category $\f A$ enriched in a bicategory $\M$ has a set of objects $x,y,\dots$, where to each object $x$ is assigned some $0$-cell $\epsilon x$ in $\M$ (called its \emph{extent}), together with $1$-cells $\f{A}(x,y)\colon \epsilon y \to \epsilon x$, and composition and unit $2$-cells satisfying the usual axioms.

If $\f A$ has exactly one object, then it is simply a monad in $\M$.
On the other hand, if $\M$ is a monoidal category, regarded as a one-object bicategory, then $\f A$ reduces to the usual sort of category enriched in a monoidal category (hence the name).

One can directly define functors, transformations, and modules between categories enriched in a bicategory.
However, the notion of functors considered by most authors is too limited, in that it requires them to preserve extents strictly.
In some cases, this can be circumvented with weak completeness conditions on the enriched categories, as in the situations of~\cite{walters:sheaves-cauchy-1,bcsw:variation-enr,shulman:eicats}.

However, a better solution is to consider instead categories enriched in a proarrow equipment $\K\to\M$ (more generally, an $\F$-bicategory), where the action of functors on objects is mediated by morphisms in $\K$.
When $\M$ is a monoidal category, we can take $\K$ to contain only the identity, so that this still reduces to the usual notion of functor in that case.
But in more general situations, it yields a more appropriate notion of functor; see~\cite{shulman:eicats,CG2012Restriction} for recent examples.

If $\M$ is locally cocomplete, we can compose modules between small $\M$-enriched categories.
Thus, we have a bicategory $\cat{Mod}_\infty(\M)$, and a proarrow equipment $\cat{Cat}_\K(\M) \to \cat{Mod}_\infty(\M)$.
The above two theorems generalize immediately.
Let $\cat{Colim}_\infty$ denote the $2$-category of cocomplete categories, with an appropriate tensor product, and $\F_\infty$ the $2$-category of fully faithful functors with cocomplete codomain.
The generalization of Kleisli objects from monads to enriched categories is called a \emph{collage}.

\begin{Thm}\label{thm:intro3}
  For any proarrow equipment $\K\to\M$, where $\M$ is locally cocomplete, the proarrow equipment $\cat{Cat}_\K(\M) \to \cat{Mod}_\infty(\M)$ is its free cocompletion, as an $\F_\infty$-enriched bicategory, under a class of $\F_\infty$-enriched colimits called \emph{tight collages}.
\end{Thm}

\begin{Thm}\label{thm:intro4}
  For any locally cocomplete bicategory $\M$, the bicategory $\cat{Mod}_\infty(\M)$ is its free cocompletion under collages, as a $\cat{Colim}_\infty$-enriched bicategory.
  Moreover, collages are an \emph{absolute colimit} for $\cat{Colim}_\infty$.
\end{Thm}

(The universal property described in the second of these theorems was exhibited for \emph{locally partially ordered} bicategories---ones whose every hom-category is  a partial order---in~\cite[Section~8]{Stubbe2005Categorical}.)

In fact, the $1$-case and the $\infty$-case are merely opposite ends of a spectrum; intermediate cases are parametrized by regular cardinals $\kappa$, and exhibit the totality of $\kappa$-small enriched categories, functors and modules as a free cocompletion in the world of ``locally $\kappa$-cocomplete bicategories''---those enriched over the monoidal bicategory $\cat{Colim}_\kappa$ of $\kappa$-cocomplete categories and $\kappa$-cocontinuous functors.

These theorems draw together the various descriptions of categories in a pleasing and abstractly well-behaved way.
They also emphasize the importance of considering categories enriched in equipments, rather than merely in bicategories, especially in order to obtain the right notion of functor.
One potential application is to a theory of exact completion for $2$-categories, since the locally posetal case of the above theorems was a basic ingredient in the general treatment of exact completion for 1-sites in~\cite{shulman:exsite}.

\subsection{Overview of the paper}
\label{sec:overview}

We conclude this introduction with a brief overview of the contents of this paper. The first part, Sections~\ref{sec:prelims}--\ref{sec:constructions}, develops  the theory of enriched bicategories up to the free cocompletion of an enriched bicategory under a class of colimits. The second part, Sections~\ref{sec:collages} and~\ref{sec:tight-collages}, applies this theory to prove the universal property of enriched categories. Although the first part is longer and more technically involved, there are no real surprises; so the reader primarily interested in the universal property of enriched categories could perfectly well skip directly to the second part, referring back to definitions and results from the first as necessary.

Section~\ref{sec:prelims} establishes notation and conventions and recalls some preliminary material. Section~\ref{sec:vbicats} defines bicategories enriched in a monoidal bicategory $\V$ and the various kinds of higher cell between them, while Section~\ref{sec:tricat-vbicats} describes the compositional structure of these cells, showing that they form a tricategory $\V\text-\cat{Bicat}$. In Section~\ref{sec:modules}, we define two-sided modules between $\V$-bicategories, and as a special case, the one-sided modules that correspond to covariant or contravariant presheaves over a $\V$-bicategory. In Section~\ref{sec:tensorproductofmodules} we define, in the $\V$-bicategorical world, the tensor product of a left $\A$-, right $\B$-module with a left $\B$-, right $\C$-module, while in Section~\ref{sec:internal-hom} we discuss the corresponding ``internal hom'' of $\V$-modules. All of this generalises the corresponding constructions on bimodules (a.k.a.\ profunctors) between ordinary categories.
We stop short of constructing the tricategory (and the resulting ``triequipment'') of $\V$-categories and modules, but this would be the natural next step.

In Section~\ref{sec:yoneda}, we state and prove a Yoneda lemma for $\V$-bicategories, and use it to prove some useful auxiliary results; then in Section~\ref{sec:vcatsofmodules}, we show that (under suitable size restrictions), the right modules over a $\V$-bicategory $\C$ themselves form a $\V$-bicategory $\M \C$ equipped with a fully faithful embedding $\C \to \M \C$.  Then in Section~\ref{sec:colimits}, we define weighted colimits in $\V$-bicategories, discuss their functoriality, and consider the closely related notion of left Kan extension, while in Section~\ref{sec:moderate}, we show that every $\V$-bicategory of the form $\M \C$ is cocomplete and prove some results related to the taking of iterated colimits. We draw together these strands in Section~\ref{sec:free-cocompletions}, by showing that the free cocompletion of a $\V$-bicategory $\C$ under a class of colimits may be constructed by closing $\C$ in $\M \C$ under colimits from that class.

Sections~\ref{sec:monoidal-adjunctions} and~\ref{sec:constructions} gather some further results relevant to the theory of enriched bicategories. Section~\ref{sec:monoidal-adjunctions}  considers the ``change of base'' operation $\V\text-\cat{Bicat} \to \W\text-\cat{Bicat}$ induced by a monoidal functor $L \colon \V \to \W$.
And in Section~\ref{sec:constructions}, we describe two ways of constructing new monoidal bicategories from old, via comma bicategories and via reflective sub-bicategories, and consider how these interact with the corresponding notions of enriched bicategory.

Finally, in Sections~\ref{sec:collages} and~\ref{sec:tight-collages}, we apply the theory developed throughout the rest of the paper to prove the universal property of enriched categories.
In Section~\ref{sec:collages}, we recall the definition of a category enriched in a bicategory, and of the modules between them; we define \emph{collages}---the kinds of colimit relevant for our free cocompletion results; we construct the monoidal bicategories $\cat{Colim}_1$ and $\cat{Colim}_\infty$;  and we prove Theorems~\ref{thm:intro2} and~\ref{thm:intro4}.
In fact, as anticipated above, we subsume these both into a more general statement, parametrized by a regular cardinal $\kappa$. Finally, in Section~\ref{sec:tight-collages}, we define equipments, and describe the equipment of categories, functors and modules enriched in an equipment. We then define tight collages, the relevant kind of colimit for our free cocompletion result; construct the monoidal bicategories $\F_1$ and $\F_\infty$; and prove Theorems~\ref{thm:intro1} and~\ref{thm:intro3}, again, by way of a more general statement parametric in a regular cardinal $\kappa$.

\section{Preliminaries}\label{sec:prelims}
We now begin our development of the theory of bicategories enriched over a monoidal bicategory. As noted in the introduction, some very basic aspects of this theory were developed in~\cite{Carmody1995Cobordism,Lack1995The-algebra}, but we will need to go significantly further; and since the two references just cited are not widely available, we have arranged to make our account self-contained. The material we describe is, of course, a two-dimensional generalisation of enriched category theory in the sense of~\cite{kelly:enriched}; however, it is not this that we will follow in our development, but rather~\cite{Street2005Enriched}. The key point is that we will not assume any kind of symmetry in the monoidal bicategory over which we are enriching. Although this restricts the range of constructions available to us---we cannot form the opposite of an enriched bicategory, or the tensor product or internal hom of two enriched bicategories---we still have enough flexibility to define enriched presheaf categories, limits and colimits, and the free completion under a class of weights. Working in the non-symmetric setting means that the theory we develop generalises without difficulty to the case of bicategories enriched in a tricategory; for the sake of simplicity, we have not given that generalisation here, but the reader should be able to make the relevant adaptations without fuss.

We will assume that the reader is familiar with the basic theory of bicategories, as set out in~\cite{Street1980Fibrations}, for example. We refer to homomorphisms of bicategories simply as \emph{functors}, and refer to pseudonatural transformations simply as \emph{transformations}. Throughout this article, $\V$ will be a monoidal bicategory in the sense of~\cite{Gordon1995Coherence}; our notational conventions will be those of~\cite{Gurski2006An-algebraic} which for the sake of self-containedness we now spell out. We write
$\otimes \colon \V \times \V \rightarrow \V$ and $I \colon 1 \rightarrow \V$ for the binary and nullary tensor product functors, and write $I$ also for the unit object picked out by the nullary tensor. The associativity and unit equivalence transformations of $\V$ we write as
\begin{equation*}
  \mathfrak a \colon \mathord{\otimes} \circ (\mathord{\otimes} \times 1) \Rightarrow \mathord{\otimes} \circ (1 \times \mathord{\otimes}) \qquad \mathfrak l \colon \mathord{\otimes} \circ (I \times 1) \Rightarrow \id \qquad \mathfrak r \colon \mathord{\otimes} \circ (1 \times I) \Rightarrow \id\rlap{ ,}
\end{equation*}
thus with $1$-cell components $\mathfrak a_{ABC} \colon (A \otimes B) \otimes C \rightarrow A \otimes (B \otimes C)$, $\mathfrak l_A \colon I \otimes A \rightarrow A$ and $\mathfrak r_A \colon A \otimes I \rightarrow A$; and we write $\mathfrak a^\centerdot$, $\mathfrak l^\centerdot$ and $\mathfrak r^\centerdot$ for specified choices of adjoint pseudoinverse. Finally, we write $\pi$,  $\lambda$, $\rho$ and $\nu$ for the invertible modifications with  components
\begin{equation*}
 \xybox{<0.1cm,0cm>:
    \POS (0,15)*+{(A  B)  (C  D)}="0", 
    (-14,5)*+{((A  B)  C)  D}="1", 
    (-9,-12)*+!/r 1em/{(A  (B  C))  D}="2", 
    (9,-12)*+!/l 1em/{A  ((B  C)  D)}="3", 
    (14,5)*+{A  (B  (C  D))}="4",
    (0,0.5)*+{}="C"
    \POS"1" \ar "0"^(0.45){\labelstyle {\mathfrak a}}|{}="01"
    \POS"1" \ar "2"_{\labelstyle {\mathfrak a 1}}|{}="12"
    \POS"2" \ar "3"_{\labelstyle {\mathfrak a}}|{}="23"
    \POS"3" \ar "4"_{\labelstyle {1 \mathfrak a}}|{}="34"
    \POS"0" \ar "4"^(0.55){\labelstyle {\mathfrak a}}|{}="04"
    \POS"C"+/u 0.2cm/ \ar@2 "C"+/d 0.2cm/^{\labelstyle \pi}
  }  
  \cd[@C-2.3em]{
    {(I  A)  B} \ar[rr]^-{\mathfrak a} \ar[dr]_-{\mathfrak l  1} & \rtwocell{d}{\lambda} &
    {I  (A  B)} \ar[dl]^-{\mathfrak l} &
    & {A  B} \ar[dl]_-{\mathfrak r^\centerdot} \ar[dr]^-{1  \mathfrak r^\centerdot} \ltwocell{d}{\rho} \\ &
    {A  B} & & 
    {(A  B)  I} \ar[rr]_-{\mathfrak a} & &
    {A  (B  I)}
  }  
\cd[@!@-1.7em]{
    \sh{l}{0.4em}{(A  I)  B} \ar[r]^-{\mathfrak a} \ar@{<-}[d]_{\mathfrak r^\centerdot  1} \dtwocell{dr}{\nu} &
    \sh{r}{0.4em}{A  (I  B)} \ar[d]^{1  \mathfrak l} \\
    {A  B} \ar[r]_-{1} &
    {A  B}\rlap{ ,}
  }
\end{equation*}
where for conciseness we write the tensor product $\otimes$ as mere juxtaposition (note that our $\nu$ was called $\mu$ in~\cite{Gordon1995Coherence} and \cite{Gurski2006An-algebraic}). 

We will generally use string diagrams rather than pasting diagrams to define compound $2$-cells in $\V$, 
with objects represented by regions, 1-cells by strings, and generating 2-cells by vertices. To avoid clutter, we will omit the symbol $\otimes$ in string diagrams, and will not explicitly label regions with objects of $\V$; the appropriate labels can always be recovered from the $1$-cell labels on strings.  For example, with these conventions the coherence $2$-cells displayed above, which will always be notated explicitly in our string diagrams, are given by:
\begin{equation*}
\vcenter{\hbox{\begin{tikzpicture}[y=0.80pt, x=0.80pt, yscale=-1.000000, xscale=1.000000, inner sep=0pt, outer sep=0pt, every text node part/.style={font=\tiny} ]
\path[draw=black,line join=miter,line cap=butt,even odd rule,line width=0.650pt] (110.0000,932.3622) .. controls (150.0000,932.3622) and (150.0000,952.3622) ..  node[above right=0.12cm,at start] {$\mathfrak a 1$}(150.0000,952.3622);
\path[draw=black,line join=miter,line cap=butt,even odd rule,line width=0.650pt] (110.0000,952.3622) -- node[above right=0.12cm,at start] {$\mathfrak a $}(150.0000,952.3622);
\path[draw=black,line join=miter,line cap=butt,even odd rule,line width=0.650pt] (110.0000,972.3622) .. controls (150.0000,972.3622) and (150.0000,952.3622) .. node[above right=0.12cm,at start] {$1 \mathfrak a$}(150.0000,952.3622) .. controls (150.0000,952.3622) and (160.0000,942.3622) .. node[above left=0.12cm,at end] {$\mathfrak a$} (190.0000,942.3622);
\path[draw=black,line join=miter,line cap=butt,even odd rule,line width=0.650pt] (150.0000,952.3622) .. controls (150.0000,952.3622) and (160.0000,962.3622) .. node[above left=0.12cm,at end] {$\mathfrak a$}(190.0000,962.3622);
\path[fill=black,line join=miter,line cap=butt,line width=0.650pt] (150.0000,952.3622) node[circle, draw, line width=0.65pt, minimum width=5mm, fill=white, inner sep=0.25mm] (text4023) {$\pi$};
\end{tikzpicture}}} \qquad \quad 
\vcenter{\hbox{\begin{tikzpicture}[y=0.80pt, x=0.80pt, yscale=-1.000000, xscale=1.000000, inner sep=0pt, outer sep=0pt, every text node part/.style={font=\tiny} ]
\path[draw=black,line join=miter,line cap=butt,even odd rule,line width=0.650pt] (0.0000,1002.3622) -- node[above right=0.12cm,at start] {$\mathfrak l 1$}(40.1100,1002.3622);
\path[draw=black,line join=miter,line cap=butt,even odd rule,line width=0.650pt] (80.0000,992.3622) .. controls (40.0000,992.3622) and (40.0000,1002.3622) .. node[above left=0.12cm,at start] {$\mathfrak a$}(40.0000,1002.3622) .. controls (40.0000,1002.3622) and (40.0000,1011.1315) .. node[above left=0.12cm,at end] {$\mathfrak l$}(80.0000,1012.3622);
\path[fill=black,line join=miter,line cap=butt,line width=0.650pt] (40.0000,1002.3622) node[circle, draw, line width=0.65pt, minimum width=5mm, fill=white, inner sep=0.25mm] (text3975) {$\lambda$};
\end{tikzpicture}}} \qquad \quad 
\vcenter{\hbox{\begin{tikzpicture}[y=0.80pt, x=0.80pt, yscale=-1.000000, xscale=1.000000, inner sep=0pt, outer sep=0pt, every text node part/.style={font=\tiny} ]
\path[draw=black,line join=miter,line cap=butt,even odd rule,line width=0.650pt] (0.0000,1002.3622) -- node[above right=0.12cm,at start] {$1 \mathfrak r^\centerdot$}(40.1100,1002.3622);
\path[draw=black,line join=miter,line cap=butt,even odd rule,line width=0.650pt] (80.0000,992.3622) .. controls (40.0000,992.3622) and (40.0000,1002.3622) .. node[above left=0.12cm,at start] {$\mathfrak r^\centerdot$}(40.0000,1002.3622) .. controls (40.0000,1002.3622) and (40.0000,1011.1315) .. node[above left=0.12cm,at end] {$\mathfrak a$}(80.0000,1012.3622);
\path[fill=black,line join=miter,line cap=butt,line width=0.650pt] (40.0000,1002.3622) node[circle, draw, line width=0.65pt, minimum width=5mm, fill=white, inner sep=0.25mm] (text3975) {$\rho$};
\end{tikzpicture}}}
\qquad \text{and} \qquad 
\vcenter{\hbox{\begin{tikzpicture}[y=0.80pt, x=0.80pt, yscale=-1.000000, xscale=1.000000, inner sep=0pt, outer sep=0pt, every text node part/.style={font=\tiny} ]
\path[draw=black,line join=miter,line cap=butt,even odd rule,line width=0.650pt] (160.0000,1002.3622) -- node[above right=0.12cm,at start] {$\mathfrak a$}(200.0000,1002.3622) .. controls (200.0000,1002.3622) and (199.8619,982.3622) .. node[above right=0.12cm,at end] {$\mathfrak r^\centerdot 1$}(160.0000,982.3622);
\path[draw=black,line join=miter,line cap=butt,even odd rule,line width=0.650pt] (200.0000,1002.3622) .. controls (200.0000,1002.3622) and (198.6890,1022.3622) .. node[above right=0.12cm,at end] {$1 \mathfrak l$}(160.0000,1022.3622);
\path[fill=black,line join=miter,line cap=butt,line width=0.650pt] (200.0000,1002.3622) node[circle, draw, line width=0.65pt, minimum width=5mm, fill=white, inner sep=0.25mm] (text4011) {$\nu$};
\end{tikzpicture}}}\rlap{\ \ \  .}
\end{equation*}

The following additional conventions will prove useful. First, if $\xi \colon A \rightarrow B$ is a $1$-cell in $\V$ with specified adjoint pseudoinverse $\xi^\centerdot \colon B \rightarrow A$, then we depict the unit and counit $2$-cells of the adjoint equivalence in string diagrams as simple cups and caps:
\begin{equation*}
\vcenter{\hbox{\begin{tikzpicture}[y=0.90pt, x=0.80pt, yscale=1.000000, xscale=-1.000000, inner sep=0pt, outer sep=0pt, every text node part/.style={font=\tiny} ]
  \path[draw=black,line join=miter,line cap=butt,even odd rule,line width=0.650pt] (30.0000,852.3622) .. controls (50.0000,852.3622) and (80.0000,852.3622) .. node[above left=0.12cm,at start] {$\xi$} (80.0000,862.3622) .. controls (80.0000,872.3622) and (50.0000,872.3622) .. node[above left=0.12cm,at end] {$\xi^\centerdot$}(30.0000,872.3622);
\end{tikzpicture}}} \qquad \qquad \qquad 
\vcenter{\hbox{\begin{tikzpicture}[y=0.90pt, x=0.80pt, yscale=-1.000000, xscale=1.000000, inner sep=0pt, outer sep=0pt, every text node part/.style={font=\tiny} ]
  \path[draw=black,line join=miter,line cap=butt,even odd rule,line width=0.650pt] (30.0000,852.3622) .. controls (50.0000,852.3622) and (80.0000,852.3622) .. node[above right=0.12cm,at start] {$\xi$} (80.0000,862.3622) .. controls (80.0000,872.3622) and (50.0000,872.3622) .. node[above right=0.12cm,at end] {$\xi^\centerdot$}(30.0000,872.3622);
\end{tikzpicture}}} \quad \rlap{ .}
\end{equation*}
Note that, in this situation, $\xi^\centerdot$ also has specified adjoint pseudoinverse $\xi$, so that we may without ambiguity exchange the position of $\xi$ and $\xi^\centerdot$ in these cups and caps.
Our next convention concerns the pseudonaturality constraint $2$-cells of $\mathfrak a$, $\mathfrak l$ and $\mathfrak r$:
\begin{equation*}
  \cd{
    {(A \otimes B) \otimes C} \ar[r]^-{\mathfrak a_{ABC}} \ar[d]_{(f \otimes g) \otimes h} \rtwocell{dr}{\mathfrak a_{fgh}} &
    {A \otimes (B \otimes C)} \ar[d]^{f\otimes  (g \otimes h)} \\
    {(A'\otimes  B')\otimes  C'} \ar[r]_-{\mathfrak a_{A'B'C'}} &
    {A' \otimes (B' \otimes C')}
  } \qquad 
  \cd{
    {I \otimes A} \ar[r]^-{\mathfrak l_{A}} \ar[d]_{1  \otimes f} \rtwocell{dr}{\mathfrak l_{f}} &
    {A} \ar[d]^{f} \\
    {I \otimes A'} \ar[r]_-{\mathfrak l_{A'}} &
    {A'}
  } \qquad 
  \cd{
    {A \otimes I} \ar[r]^-{\mathfrak r_{A}} \ar[d]_{f  \otimes 1} \rtwocell{dr}{\mathfrak r_{f}} &
    {A} \ar[d]^{f} \\
    {A' \otimes I} \ar[r]_-{\mathfrak r_{A'}} &
    {A'}\rlap{ .}
  }
\end{equation*}
We will draw these and their inverses as string crossings, with the convention that the string labelled by $\mathfrak a$, $\mathfrak l$ or $\mathfrak r$ should remain uppermost; so $\mathfrak a_{fgh}$ and $\mathfrak a_{fgh}^{-1}$ are drawn as:
\begin{equation*}
\vcenter{\hbox{\begin{tikzpicture}[y=0.90pt, x=1pt, yscale=-1.000000, xscale=1.000000, inner sep=0pt, outer sep=0pt, every text node part/.style={font=\tiny} ]
\path[draw=black,line join=miter,line cap=butt,even odd rule,line width=0.650pt] 
(0.0000,922.3622) .. controls (40.0000,922.3622) and (40.0000,902.3622) .. 
 node[above right=0.12cm,at start] {$\mathfrak a$}  node[above left=0.12cm,at end] {$\mathfrak a$}(80.0000,902.3622)
(0.0000,902.3622) .. controls (18.1366,902.3622) and (28.0498,906.4739) .. 
 node[above right=0.12cm,at start] {$\!(fg)h$}(37.1968,910.9687)
(42.8032,913.7557) .. controls (51.9502,918.2505) and (61.8634,922.3622) .. 
 node[above left=0.12cm,at end] {$f(gh)\!\!$}(80.0000,922.3622);
\end{tikzpicture}}} \quad \qquad \text{and} \quad \qquad \vcenter{\hbox{\begin{tikzpicture}[y=0.90pt, x=1pt, yscale=-1.000000, xscale=-1.000000, inner sep=0pt, outer sep=0pt, every text node part/.style={font=\tiny} ]
\path[draw=black,line join=miter,line cap=butt,even odd rule,line width=0.650pt] 
(0.0000,922.3622) .. controls (40.0000,922.3622) and (40.0000,902.3622) .. 
 node[above left=0.12cm,at start] {$\mathfrak a$}  node[above right=0.12cm,at end] {$\mathfrak a$}(80.0000,902.3622)
(0.0000,902.3622) .. controls (18.1366,902.3622) and (28.0498,906.4739) .. 
 node[above left=0.12cm,at start] {$(fg)h\!$}(37.1968,910.9687)
(42.8032,913.7557) .. controls (51.9502,918.2505) and (61.8634,922.3622) .. 
 node[above right=0.12cm,at end] {$\!f(gh)$}(80.0000,922.3622);
\end{tikzpicture}}}\quad \rlap{ ,}
\end{equation*}
and correspondingly for $\mathfrak l$ and $\mathfrak r$. We also allow ourselves to apply this convention to the pseudonaturality constraints of the pseudoinverse transformations $\mathfrak a^\centerdot$, $\mathfrak l^\centerdot$ and $\mathfrak r^\centerdot$.

Our final convention concerns the pseudofunctoriality of $\otimes$. Given $1$-cells $f \colon A \rightarrow A'$ and $g \colon B \rightarrow B'$ in $\V$, this pseudofunctoriality gives canonical invertible $2$-cells 
\begin{equation*}
  \cd{
    {A \otimes B} \ar[r]^-{f \otimes 1} \ar[d]_{1 \otimes g} \ar[dr]|-{f \otimes g}&
    {A' \otimes B} \ar[d]^{1 \otimes g} \rtwocell[0.3]{dl}{} \rtwocell[0.7]{dl}{} \\
    {A \otimes B'} \ar[r]_-{f \otimes 1} &
    {A' \otimes B'}\rlap{ .}
  }
\end{equation*}
We notate instances of these two $2$-cells and their inverses by string splittings and joinings:
\begin{equation*}
\vcenter{\hbox{\begin{tikzpicture}[y=0.90pt, x=0.80pt, yscale=-1.000000, xscale=-1.000000, inner sep=0pt, outer sep=0pt, every text node part/.style={font=\tiny} ]
  \path[draw=black,line join=miter,line cap=butt,even odd rule,line width=0.650pt] (0.0000,1002.3622) --  node[above left=0.12cm,at start] {$fg$}(40.1100,1002.3622);
\path[draw=black,line join=miter,line cap=butt,even odd rule,line width=0.650pt] (80.0000,992.3622) .. controls (40.0000,992.3622) and (40.0000,1002.3622) .. node[above right=0.12cm,at start] {$1 g$}(40.0000,1002.3622) .. controls (40.0000,1002.3622) and (40.0000,1011.1315) .. node[above right=0.12cm,at end] {$f 1 $}(80.0000,1012.3622);
\path[xscale=-1.000,yscale=1.000,draw=black,fill=black] (-40.0000,1002.3622) circle (0.0353cm);
\end{tikzpicture}}} \qquad \quad 
\vcenter{\hbox{\begin{tikzpicture}[y=0.90pt, x=0.80pt, yscale=-1.000000, xscale=1.000000, inner sep=0pt, outer sep=0pt, every text node part/.style={font=\tiny} ]
  \path[draw=black,line join=miter,line cap=butt,even odd rule,line width=0.650pt] (0.0000,1002.3622) --  node[above right=0.12cm,at start] {$fg$}(40.1100,1002.3622);
\path[draw=black,line join=miter,line cap=butt,even odd rule,line width=0.650pt] (80.0000,992.3622) .. controls (40.0000,992.3622) and (40.0000,1002.3622) .. node[above left=0.12cm,at start] {$f 1$}(40.0000,1002.3622) .. controls (40.0000,1002.3622) and (40.0000,1011.1315) .. node[above left=0.12cm,at end] {$1 g$}(80.0000,1012.3622);
\path[xscale=-1.000,yscale=1.000,draw=black,fill=black] (-40.0000,1002.3622) circle (0.0353cm);
\end{tikzpicture}}} \qquad \quad 
\vcenter{\hbox{\begin{tikzpicture}[y=0.90pt, x=0.80pt, yscale=-1.000000, xscale=-1.000000, inner sep=0pt, outer sep=0pt, every text node part/.style={font=\tiny} ]
  \path[draw=black,line join=miter,line cap=butt,even odd rule,line width=0.650pt] (0.0000,1002.3622) --  node[above left=0.12cm,at start] {$fg$}(40.1100,1002.3622);
\path[draw=black,line join=miter,line cap=butt,even odd rule,line width=0.650pt] (80.0000,992.3622) .. controls (40.0000,992.3622) and (40.0000,1002.3622) .. node[above right=0.12cm,at start] {$f1$}(40.0000,1002.3622) .. controls (40.0000,1002.3622) and (40.0000,1011.1315) .. node[above right=0.12cm,at end] {$1g $}(80.0000,1012.3622);
\path[xscale=-1.000,yscale=1.000,draw=black,fill=black] (-40.0000,1002.3622) circle (0.0353cm);
\end{tikzpicture}}} \qquad \quad 
\vcenter{\hbox{\begin{tikzpicture}[y=0.90pt, x=0.80pt, yscale=-1.000000, xscale=1.000000, inner sep=0pt, outer sep=0pt, every text node part/.style={font=\tiny} ]
  \path[draw=black,line join=miter,line cap=butt,even odd rule,line width=0.650pt] (0.0000,1002.3622) --  node[above right=0.12cm,at start] {$fg$}(40.1100,1002.3622);
\path[draw=black,line join=miter,line cap=butt,even odd rule,line width=0.650pt] (80.0000,992.3622) .. controls (40.0000,992.3622) and (40.0000,1002.3622) .. node[above left=0.12cm,at start] {$1g$}(40.0000,1002.3622) .. controls (40.0000,1002.3622) and (40.0000,1011.1315) .. node[above left=0.12cm,at end] {$f 1 $}(80.0000,1012.3622);
\path[xscale=-1.000,yscale=1.000,draw=black,fill=black] (-40.0000,1002.3622) circle (0.0353cm);
\end{tikzpicture}}}\qquad \rlap{ .}
\end{equation*}
To avoid clutter, where joinings and splittings occur to the extreme left or right of a diagram, we may choose to omit them; this means that the displayed $2$-cell will have an invalid domain or codomain $1$-cell, but no confusion should arise, since the correct diagram may be readily constructed by reappending the omitted joinings or splittings. Furthermore, we notate instances of the composite $2$-cell $(f \otimes 1) \circ (1 \otimes g) \Rightarrow (1 \otimes g) \circ (f \otimes 1)$ and its inverse again by string crossings, where we now decide (arbitrarily) to keep the string labelled by $f \otimes 1$ uppermost, as in:
\begin{equation*}
\vcenter{\hbox{\begin{tikzpicture}[y=0.90pt, x=1pt, yscale=-1.000000, xscale=1.000000, inner sep=0pt, outer sep=0pt, every text node part/.style={font=\tiny} ]
\path[draw=black,line join=miter,line cap=butt,even odd rule,line width=0.650pt] 
(0.0000,922.3622) .. controls (40.0000,922.3622) and (40.0000,902.3622) .. 
 node[above right=0.12cm,at start] {$\!f1$}  node[above left=0.12cm,at end] {$f1$}(80.0000,902.3622)
(0.0000,902.3622) .. controls (18.1366,902.3622) and (28.0498,906.4739) .. 
 node[above right=0.12cm,at start] {$\!1g$}(37.1968,910.9687)
(42.8032,913.7557) .. controls (51.9502,918.2505) and (61.8634,922.3622) .. 
 node[above left=0.12cm,at end] {$1g$}(80.0000,922.3622);
\end{tikzpicture}}} \quad \qquad \text{and} \quad \qquad \vcenter{\hbox{\begin{tikzpicture}[y=0.90pt, x=1pt, yscale=-1.000000, xscale=-1.000000, inner sep=0pt, outer sep=0pt, every text node part/.style={font=\tiny} ]
\path[draw=black,line join=miter,line cap=butt,even odd rule,line width=0.650pt] 
(0.0000,922.3622) .. controls (40.0000,922.3622) and (40.0000,902.3622) .. 
 node[above left=0.12cm,at start] {$f1\!$}  node[above right=0.12cm,at end] {$\!f1$}(80.0000,902.3622)
(0.0000,902.3622) .. controls (18.1366,902.3622) and (28.0498,906.4739) .. 
 node[above left=0.12cm,at start] {$1g\!$}(37.1968,910.9687)
(42.8032,913.7557) .. controls (51.9502,918.2505) and (61.8634,922.3622) .. 
 node[above right=0.12cm,at end] {$\!1g$}(80.0000,922.3622);
\end{tikzpicture}}}\quad \rlap{ .}
\end{equation*}

As a first application of our diagrammatic conventions, we use them to prove the following lemma, which categorifies the well-known result that, in a monoidal category, we have $\lambda_I = \rho_I \colon I \otimes I \to I$.

\begin{Lemma}\label{lem:II}
There is an invertible $2$-cell $\theta \colon \mathfrak l_I \Rightarrow \mathfrak r_I \colon I \otimes I \to I$ in $\V$.
\end{Lemma}
\begin{proof}
First, for any $A \in \V$ there are invertible $2$-cells 
\begin{equation*}
 \mathfrak r^\centerdot_{A \otimes I} \Rightarrow \mathfrak r^\centerdot_A \otimes 1 \colon A \otimes I \rightarrow (A \otimes I) \otimes I \quad \text{and} \quad 
\mathfrak l_{I \otimes A} \Rightarrow 1 \otimes \mathfrak l_A \colon I \otimes (I \otimes A) \rightarrow I \otimes A
\end{equation*}
given by the string diagrams
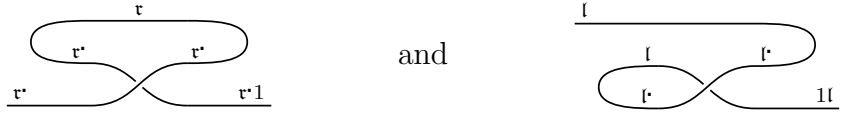
\begin{equation*}
\vcenter{\hbox{\begin{tikzpicture}[y=0.80pt, x=0.9pt, yscale=-1.000000, xscale=0.500000, inner sep=0pt, outer sep=0pt, every text node part/.style={font=\tiny} ]
\path[draw=black,line join=miter,line cap=butt,even odd rule,line width=0.650pt] 
(180.0000,782.3622) .. controls (200.0000,782.3622) and (230.0000,782.3622) .. 
(230.0000,792.3622) .. controls (230.0000,802.3622) and (200.0000,802.3622) .. 
node[above=0.07cm,pos=0.85] {$\mathfrak r^\centerdot$}(180.0000,802.3622);
\path[draw=black,line join=miter,line cap=butt,even odd rule,line width=0.650pt] 
(100.0000,822.3622) .. controls (140.0000,822.3622) and (140.0000,802.3622) .. 
(180.0000,802.3622)
(100.0000,802.3622) .. controls (118.1366,802.3622) and (128.0498,806.4739) .. 
(137.1968,810.9687)
(142.8032,813.7557) .. controls (151.9502,818.2505) and (161.8634,822.3622) .. 
(180.0000,822.3622);
\path[draw=black,line join=miter,line cap=butt,even odd rule,line width=0.650pt] 
(100.0000,782.3622) .. controls (80.0000,782.3622) and (50.0000,782.3622) .. 
(50.0000,792.3622) .. controls (50.0000,802.3622) and (80.0000,802.3622) .. 
node[above=0.07cm,pos=0.85]{$\mathfrak r^\centerdot$}(100.0000,802.3622);
\path[draw=black,line join=miter,line cap=butt,even odd rule,line width=0.650pt] 
(100.0000,782.3622) -- node[above=0.07cm,pos=0.5] {$\mathfrak r$}(180.0000,782.3622);
\path[draw=black,line join=miter,line cap=butt,even odd rule,line width=0.650pt] 
(180.0000,822.3622) -- node[above left=0.12cm,at end] {$\mathfrak r^\centerdot 1$}(250.0000,822.3622);
\path[draw=black,line join=miter,line cap=butt,even odd rule,line width=0.650pt] 
(30.0000,822.3622) -- node[above right=0.12cm,at start] {$\mathfrak r^\centerdot$}(100.0000,822.3622);
\end{tikzpicture}}}\qquad \qquad \text{and} \qquad \qquad 
\vcenter{\hbox{\begin{tikzpicture}[y=0.80pt, x=0.9pt, yscale=-1.000000, xscale=0.500000, inner sep=0pt, outer sep=0pt, every text node part/.style={font=\tiny} ]
\path[draw=black,line join=miter,line cap=butt,even odd rule,line width=0.650pt]
(180.0000,722.3622) .. controls (200.0000,722.3622) and (230.0000,722.3622) .. 
(230.0000,732.3622) .. controls (230.0000,742.3622) and (200.0000,742.3622) .. 
(180.0000,742.3622) node[above=0.07cm,pos=0.85] {$\mathfrak l^\centerdot$};
\path[draw=black,line join=miter,line cap=butt,even odd rule,line width=0.650pt] 
(100.0000,762.3622) .. controls (140.0000,762.3622) and (140.0000,742.3622) .. (180.0000,742.3622)
(100.0000,742.3622) .. controls (118.1366,742.3622) and (128.0498,746.4739) .. (137.1968,750.9687)
(142.8032,753.7557) .. controls (151.9502,758.2505) and (161.8634,762.3622) .. (180.0000,762.3622);
\path[draw=black,line join=miter,line cap=butt,even odd rule,line width=0.650pt] 
(100.0000,742.3622) .. controls (80.0000,742.3622) and (50.0000,742.3622) .. node[above=0.07cm,pos=0.15] {$\mathfrak l$}
(50.0000,752.3622) .. controls (50.0000,762.3622) and (80.0000,762.3622) .. node[above=0.07cm,pos=0.85] {$\mathfrak l^\centerdot$}
(100.0000,762.3622);
\path[draw=black,line join=miter,line cap=butt,even odd rule,line width=0.650pt]
(30.0000,722.3622) -- node[above right=0.12cm,at start] {$\mathfrak l$} (180.0000,722.3622);
\path[draw=black,line join=miter,line cap=butt,even odd rule,line width=0.650pt]
(180.0000,762.3622) -- node[above left=0.12cm,at end] {$1\mathfrak l$} (250.0000,762.3622);
\end{tikzpicture}}}
\end{equation*}
which given as pastings are the following composites of coherence $2$-cells:
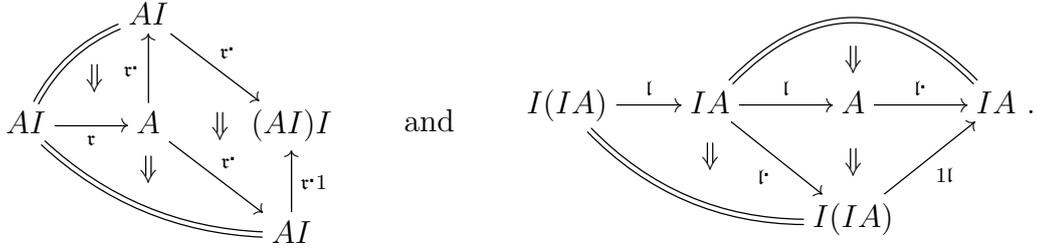
\begin{equation*}
  \cd{
    \dtwocell[0.56]{dr}{} & A I \ar@/_8pt/@{=}[dl] \ar[dr]^-{\mathfrak r^\centerdot} \\
    A I \ar@/_12pt/@{=}[drr] \ar[r]_-{\mathfrak r} & A \dtwocell{r}{} \dtwocell[0.4]{d}{} \ar[u]^-{\mathfrak r^\centerdot} \ar[dr]^-{\mathfrak r^\centerdot} & (AI)I \\
    & & AI \ar[u]_-{\mathfrak r^\centerdot 1}
  } \qquad \text{and} \qquad 
  \cd{ & {\ }\\
    I(IA) \ar@{=}@/_16pt/[drr] \ar[r]^-{\mathfrak l} & IA \dtwocell[0.45]{d}{} \ar@{=}@/^32pt/[rr] \ar[dr]_-{\mathfrak l^\centerdot} \ar[r]^-{\mathfrak l} & A \dtwocell{d}{} \dtwocell[-0.4]{d}{} \ar[r]^-{\mathfrak l^\centerdot} & IA\rlap{ .} \\
    & & I(IA) \ar[ur]_-{1 \mathfrak l}
  }
\end{equation*}
Taking $A = I$, we may now incorporate these two $2$-cells into the following
more complex string diagram specifying the desired $\theta$:
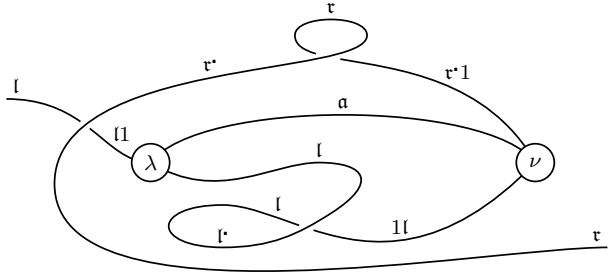
\begin{equation*}
\begin{tikzpicture}[y=0.8pt, x=0.9pt,yscale=-1, inner sep=0pt, outer sep=0pt, every text node part/.style={font=\tiny} ]
\path[use as bounding box] (-10, 885) rectangle (230,1015);
\path[draw=black,line join=miter,line cap=butt,line width=0.650pt] 
(40.0000,962.3622) .. controls (60.0000,932.3622) and (180.0000,932.3622) .. 
node[above=0.09cm,pos=0.5] {$\mathfrak a$}(200.0000,962.3622);
\path[draw=black,line join=miter,line cap=butt,even odd rule,line width=0.650pt] 
(40.0000,962.3622) .. controls (67.8089,981.8274) and (90.5607,962.3622) .. 
node[above=0.08cm,at end]{$\mathfrak l$}(110.0000,962.3622) .. controls (130.0000,962.3622) and (134.5151,972.0851) .. 
(115.2778,986.0659) .. controls (101.1787,996.3124) and (90.0000,1002.3622) .. 
(70.0000,1002.3622) .. controls (40.0000,1002.3622) and (40.0000,982.3622) ..
node[above right=0.08cm,pos=0.05]{$\mathfrak l^\centerdot$}(70.0000,982.3622) .. controls (79.4864,982.3622) and (89.9728,987.3618) ..
node[above right=0.08cm,pos=0.65]{$\mathfrak l$}(101.7752,991.9857)
(107.7718,994.2359) .. controls (132.2224,1002.9204) and (162.1268,1006.5476) ..
node[above left=0.08cm,pos=0.5]{$1 \mathfrak l$}(200.0000,962.3622);
\path[draw=black,line join=miter,line cap=butt,even odd rule,line width=0.650pt] 
(40.0000,962.3622) .. controls (29.9226,962.1914) and (23.4739,954.0321) .. 
node[above right=0.07cm,pos=0.67]{$\mathfrak l$1}(14.6820,946.2668)
(10.8060,943.0287) .. controls (3.3823,937.2214) and (-5.9576,932.3622) .. 
node[above right=0.12cm,at end]{$\mathfrak l$}(-20.0000,932.3622)
(200.0000,962.3622) .. controls (180.0268,922.6353) and (141.8842,918.1984) .. 
node[above right=0.08cm,pos=0.5]{$\mathfrak r^\centerdot 1$}(118.9143,913.5606)
(108.7270,911.0198) .. controls (103.2816,909.1753) and (100.0000,906.6861) .. 
(100.0000,902.3622) .. controls (100.0000,892.3622) and (130.0000,892.3622) .. 
node[above=0.09cm,pos=0.5]{$\mathfrak r$}(130.0000,902.3622) .. controls (130.0000,922.3622) and (0.0000,912.3622) ..
node[above=0.09cm,pos=0.5]{$\mathfrak r^\centerdot$}(0.0000,972.3622) .. controls (0.0000,1041.8558) and (180.0000,1002.3622) ..
node[above left=0.12cm,at end]{$\mathfrak r$}(230.0000,1002.3622);
\path[fill=black] (40.0000,962.3622) node[circle, draw, line width=0.65pt, minimum width=5mm, fill=white, inner sep=0.25mm] (text2987) {$\lambda$};
\path[fill=black] (200.0000,962.3622) node[circle, draw, line width=0.65pt, minimum width=5mm, fill=white, inner sep=0.25mm] (text2991) {$\nu$};
\end{tikzpicture}
\quad \text. \qedhere
\end{equation*}
\end{proof}

\section{$\V$-bicategories and their morphisms}\label{sec:vbicats}
We now describe the notion of bicategory enriched in a monoidal bicategory $\V$, together with the various kinds of higher cells between these: enriched functors, transformations and modifications, and additionally the enriched analogue of the \emph{icons} of~\cite{lack:icons}. We also describe the construction assigning to every $\V$-enriched bicategory its underlying ordinary bicategory, and correspondingly for the cells between them; ordinary bicategories in fact arise as $\V$-enriched bicategories on taking $\V = \cat{Cat}$.

\subsection{$\V$-bicategories}
A \emph{$\V$-bicategory} $\B$, or \emph{bicategory enriched in $\V$},  is given by:
\begin{itemize}
\item A set of objects $\ob \B$;
\item For each $x, y \in \B$ a hom-object $\B(x,y) \in \V$;
\item For each $x \in \B$ a morphism $j_x \colon I \to \B(x,x)$ in $\V$;
\item For each $x,y,z \in \B$ a morphism $m_{xyz} \colon \B(y,z) \otimes \B(x,y) \to \B(x,z)$ in $\V$;
\item For each $x,y \in \B$, invertible $2$-cells
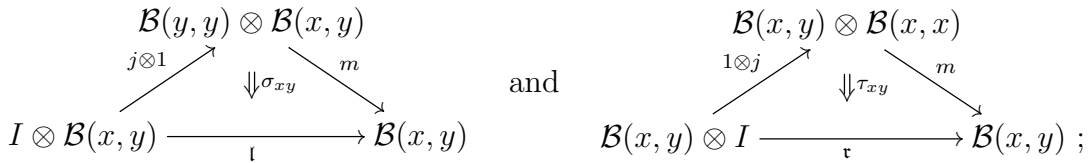
\begin{equation*}
\cd[@!C@C-4.5em]{
 & \B(y, y) \otimes \B(x, y) \ar[dr]^{m} \dtwocell[0.55]{d}{\lu_{xy}} \\
 I \otimes \B(x, y) \ar[rr]_{\mathfrak l} \ar[ur]^{j \otimes 1} & &
 \B(x, y)
}\quad \text{and} \quad
\cd[@!C@C-4.5em]{
 & \B(x, y) \otimes \B(x, x) \ar[dr]^{m} \dtwocell[0.55]{d}{\ru_{xy}} \\
 \B(x, y) \otimes I \ar[rr]_{\mathfrak r} \ar[ur]^{1 \otimes j} & &
 \B(x, y)\rlap{ ;}
}
\end{equation*}
\vskip0.5\baselineskip
\item For each $w,x,y,z \in \B$, an invertible $2$-cell
\begin{equation*}
\cd{
 (\B(y, z) \otimes \B(x, y)) \otimes \B(w, x) \ar[rr]^{\mathfrak a} \ar[d]_{m \otimes 1}  & {} \rtwocell[0.55]{d}{{\alpha}_{wxyz}} &
 \B(y, z) \otimes (\B(x, y) \otimes \B(w, x)) \ar[d]^{1 \otimes m} \\
 \B(x, z) \otimes \B(w, x) \ar[r]_-m & \B(w, z) & \B(y, z) \otimes \B(w, y) \ar[l]^-{\ m}
}
\end{equation*}
\end{itemize}
 subject to the axioms that, for each $x,y,z \in \B$, we have
\begin{equation*}
\vcenter{\hbox{
\begin{tikzpicture}[y=0.80pt, x=0.7pt,yscale=-1, inner sep=0pt, outer sep=0pt, every text node part/.style={font=\tiny} ]
  \path[draw=black,line join=miter,line cap=butt,line width=0.650pt]
    (8.2143,947.7193) .. controls (43.2143,947.0050) and (140.3571,944.8622) ..
    node[above right=0.15cm,at start] {$\mathfrak r^\centerdot \t 1$}(157.5000,962.0050);
  \path[draw=black,line join=miter,line cap=butt,line width=0.650pt]
    (8.2143,968.7908) .. controls (23.3018,968.7908) and (45.8868,971.1394) ..
    node[above right=0.15cm,at start] {$\!\!(1 \t j) \t 1$}(66.7018,975.4731)
    (74.3039,977.1662) .. controls (91.8110,981.3341) and (107.1382,986.9540) ..
    node[above right=0.06cm,pos=0.45]{$1 \t (j \t 1)$} (114.2857,993.7908)
    (38.9286,997.7193) .. controls (55.7143,987.3622) and (80.7143,955.2193) ..
    node[above left=0.12cm,pos=0.72] {$\mathfrak a$}(157.5000,964.8622);
  \path[draw=black,line join=miter,line cap=butt,line width=0.650pt]
    (38.5714,1000.9337) .. controls (58.2143,1006.2907) and (102.5000,1004.8622)
    .. node[above=0.13cm,pos=0.45]{$1 \t m$} (114.6429,997.3622);
  \path[draw=black,line join=miter,line cap=butt,line width=0.650pt]
    (39.2857,1003.7908) .. controls (45.2089,1021.2793) and (140.6406,1018.4336)
    .. node[above left=0.12cm,at end]{$m$} (166.2143,1018.4336);
  \path[draw=black,line join=miter,line cap=butt,line width=0.650pt]
    (116.4286,997.0050) .. controls (135.7143,997.7193) and (136.0714,979.8622) ..
    node[below right=0.07cm,pos=0.55]{$1 \t \mathfrak l$}(158.2143,967.7193);
  \path[draw=black,line join=miter,line cap=butt,line width=0.650pt]
    (8.2143,989.1479) .. controls (33.6008,989.0186) and (35.7143,998.0765) ..
    node[above right=0.15cm,at start] {$m \t 1$}(35.7143,998.0765);
  \path[draw=black,line join=miter,line cap=butt,line width=0.650pt]
    (8.2143,1010.9336) .. controls (33.2143,1011.2907) and (36.0714,1002.3622) ..
    node[above right=0.15cm,at start] {$m$}(36.0714,1002.3622);
  \path[fill=black] (158.57143,965.93359) node[circle, draw, line width=0.65pt,
    minimum width=5mm, fill=white, inner sep=0.25mm] (text3445) {$\nu$    };
  \path[fill=black] (114.1893,996.59198) node[circle, draw, line width=0.65pt,
    minimum width=5mm, fill=white, inner sep=0.25mm] (text3449) {$1 \t \lu$    };
  \path[fill=black] (36.42857,1001.2907) node[circle, draw, line width=0.65pt,
    minimum width=5mm, fill=white, inner sep=0.25mm] (text4000) {${\alpha}$    };
\end{tikzpicture}
}}
\quad \ =\quad\ \
\vcenter{\hbox{
\begin{tikzpicture}[y=0.75pt, x=0.59pt,yscale=-1, inner sep=0pt, outer sep=0pt, every text node part/.style={font=\tiny} ]
  \path[draw=black,line join=miter,line cap=butt,line width=0.650pt]
    (0.7143,990.5765) .. controls (148.2143,986.2907) and (118.2143,1013.0765) ..
    node[above right=0.15cm,at start] {$\mathfrak r^\centerdot \t 1$} (98.2143,1013.0765);
  \path[draw=black,line join=miter,line cap=butt,line width=0.650pt]
    (0.7143,1011.6479) .. controls (35.8741,1011.6379) and (70.9174,1008.3155) ..
    node[above right=0.15cm,at start] {$\!\!(1 \t j) \t 1$}(94.2857,1012.3622);
  \path[draw=black,line join=miter,line cap=butt,line width=0.650pt]
    (0.7143,1032.7193) .. controls (29.6429,1033.0764) and (82.8571,1026.6479) ..
    node[above right=0.15cm,at start] {$m \t 1$}(94.2857,1015.2193);
  \path[draw=black,line join=miter,line cap=butt,line width=0.650pt]
    (0.7143,1052.0051) .. controls (56.4286,1052.3622) and (91.5714,1031.2907) ..
    node[above right=0.15cm,at start] {$m$}(115.5000,1031.6479);
  \path[fill=black] (95.714287,1013.7908) node[circle, draw, line width=0.65pt,
    minimum width=5mm, fill=white, inner sep=0.25mm] (text3445) {$\ru \t 1$    };
\end{tikzpicture}
}}
\end{equation*}
in $\V(\B(y, z)  \otimes \B(x,y), \B(x,z))$, and that for each $v,w,x,y,z \in \B$, we have
\begin{equation*}
\vcenter{\hbox{\begin{tikzpicture}[y=0.9pt, x=0.8pt,yscale=-1, inner sep=0pt, outer sep=0pt, every text node part/.style={font=\tiny} ]
  \path[draw=black,line join=miter,line cap=butt,line width=0.650pt]
    (0.5761,996.5573) .. controls (25.8299,996.0522) and (35.7165,994.4461) ..
    node[above right=0.15cm,at start] {$\!\!(m \t 1) \t 1$}(48.0259,1000.9550);
  \path[draw=black,line join=miter,line cap=butt,line width=0.650pt]
    (0.5761,1020.4270) .. controls (32.7994,1020.6795) and (38.4294,1009.3612) ..
    node[above right=0.17cm,at start] {$m \t 1$}(48.0259,1003.0477);
  \path[draw=black,line join=miter,line cap=butt,line width=0.650pt]
    (49.8969,1000.5044) .. controls (58.2307,980.6554) and (151.0215,979.4784) ..
    node[above right=0.25cm,pos=0.4] {$\mathfrak a \t 1$}(166.1738,987.8122);
  \path[draw=black,line join=miter,line cap=butt,line width=0.650pt]
    (49.6135,1003.1578) .. controls (60.9777,1006.6933) and (72.9109,1017.9473) ..
    node[below left=0.07cm,pos=0.55]{$m \t 1$} (80.4466,1027.1123);
  \path[draw=black,line join=miter,line cap=butt,line width=0.650pt]
    (0.5761,1042.0698) .. controls (14.4657,1042.3223) and (75.2981,1040.0698) ..
    node[above right=0.15cm,at start] {$m$} (80.8539,1029.7127);
  \path[draw=black,line join=miter,line cap=butt,line width=0.650pt]
    (82.2825,1030.8845) .. controls (90.1821,1054.0933) and (173.6677,1051.1703)
    .. node[above left=0.12cm,at end] {$m$}(214.7867,1051.3122);
  \path[draw=black,line join=miter,line cap=butt,line width=0.650pt]
    (141.5590,1019.2256) .. controls (143.0743,1011.6494) and (166.0258,992.5776)
    .. node[below right=0.02cm] {$1 \t \mathfrak a$}(166.0258,992.5776);
  \path[draw=black,line join=miter,line cap=butt,line width=0.650pt]
    (141.6203,1021.3685) .. controls (141.6203,1021.3685) and (163.4571,1014.0062)
    .. node[above left=0.12cm,at end] {$1 \t (1 \t m)\!\!$}(214.7868,1014.7638);
  \path[draw=black,line join=miter,line cap=butt,line width=0.650pt]
    (141.4111,1023.0061) .. controls (141.4111,1023.0061) and (142.5493,1033.5038)
    .. node[above left=0.12cm,at end] {$1 \t m$}(214.7868,1033.7564);
  \path[draw=black,line join=miter,line cap=butt,line width=0.650pt]
    (167.5590,988.4831) .. controls (167.5590,988.4831) and (177.6142,980.4754) ..
    node[above left=0.12cm,at end] {$\mathfrak a$} (214.7868,979.9703);
  \path[draw=black,line join=miter,line cap=butt,line width=0.650pt]
    (167.9162,991.3403) .. controls (172.8729,995.3375) and (185.9479,996.5546) ..
    node[above left=0.12cm,at end] {$\mathfrak a$}(214.7867,996.5546);
  \path[draw=black,line join=miter,line cap=butt,line width=0.650pt]
    (82.5536,1028.3546) .. controls (96.1906,1037.3648) and (132.5520,1035.8496)
    .. node[below right=0.12cm,pos=0.42] {$1 \t m$} (139.8756,1022.0412);
  \path[draw=black,line join=miter,line cap=butt,line width=0.650pt]
    (82.7084,1025.8212) .. controls (84.1609,1008.5230) and (138.1030,990.5459) ..
    node[above left=0.07cm,pos=0.6] {$\mathfrak a$}(166.0786,990.2620)
    (49.8145,1001.8228) .. controls (60.2802,1001.0663) and (79.7259,1003.0374) ..
    node[above right=0.08cm,pos=0.4,rotate=-9] {$(1 \t m) \t 1$}(98.0839,1006.4270)
    (105.9744,1007.9878) .. controls (121.2679,1011.2244) and (134.6866,1015.3648) ..
    node[below left=0.05cm,pos=0.68,rotate=-18] {$1 \t (m \t 1)$}(139.8779,1019.5827);
  \path[fill=black] (166.83624,990.0094) node[circle, draw, line width=0.65pt,
    minimum width=5mm, fill=white, inner sep=0.25mm] (text3313) {$\pi$    };
  \path[fill=black] (48.501251,1001.512) node[circle, draw, line width=0.65pt,
    minimum width=5mm, fill=white, inner sep=0.25mm] (text3305) {$\ass \t 1$    };
  \path[fill=black] (81.600487,1027.7435) node[circle, draw, line width=0.65pt,
    minimum width=5mm, fill=white, inner sep=0.25mm] (text3309) {$\ass$    };
  \path[fill=black] (140.38065,1021.0309) node[circle, draw, line width=0.65pt,
    minimum width=5mm, fill=white, inner sep=0.25mm] (text3317) {$1 \t \ass$  };
\end{tikzpicture}}}
\quad\ \  = \quad\ \
\vcenter{\hbox{\begin{tikzpicture}[y=0.9pt, x=0.7pt,yscale=-1, inner sep=0pt, outer sep=0pt, every text node part/.style={font=\tiny} ]
  \path[draw=black,line join=miter,line cap=butt,line width=0.650pt]
    (0.5761,1010.9221) .. controls (22.7994,1011.1746) and (33.5330,1012.5294) ..
    node[above right=0.15cm,at start] {$m \t 1$}(49.1780,1027.2420);
  \path[draw=black,line join=miter,line cap=butt,line width=0.650pt]
    (0.5761,1042.0698) .. controls (14.4657,1042.3223) and (43.6151,1040.3578) ..
    node[above right=0.15cm,at start] {$m$}(49.1709,1030.0007);
  \path[draw=black,line join=miter,line cap=butt,line width=0.650pt]
    (51.1755,1030.5965) .. controls (67.4280,1056.6856) and (126.7968,1041.6654)
    .. node[below right=0.14cm] {$ m$} (139.8511,1033.1665);
  \path[draw=black,line join=miter,line cap=butt,line width=0.650pt]
    (142.8500,1030.5854) .. controls (142.8500,1030.5854) and (184.5612,1028.3914)
    .. node[above left=0.12cm,at end] {$1 \t m$}(225.8909,1029.1490);
  \path[draw=black,line join=miter,line cap=butt,line width=0.650pt]
    (142.6407,1033.9511) .. controls (142.6407,1033.9511) and (172.9964,1047.0248)
    .. node[above left=0.12cm,at end] {$m$}(225.8909,1046.9894);
  \path[draw=black,line join=miter,line cap=butt,line width=0.650pt]
    (51.8895,1026.1092) .. controls (64.3077,997.3630) and (142.1189,975.3128) ..
    node[above left=0.12cm,at end] {$\mathfrak a$}(225.8909,975.6050)
    (52.8868,1028.6427) .. controls (64.4224,1023.3508) and (87.4860,1019.3533) ..
    node[below right=0.1cm,pos=0.53] {$1 \t m$} (113.0660,1016.4776)
    (121.4833,1015.5795) .. controls (135.8648,1014.1237) and (150.7159,1013.0048) ..
    node[above=0.1cm] {$1 \t m$} (164.5150,1012.1936)
    (173.4370,1011.7035) .. controls (191.3809,1010.7870) and (206.7880,1010.4220) ..
    node[above left=0.12cm,at end] {$1 \t (1 \t m)\!\!$} (225.8909,1010.5363)
    (143.3647,1028.4425) .. controls (173.9707,1013.3776) and (183.0558,992.0340) ..
    node[above left=0.12cm,at end] {$\mathfrak a$}(225.8909,992.6100)
    (1.1521,981.0038) .. controls (30.9106,980.4472) and (58.1432,987.9204) ..
    node[above right=0.15cm,at start] {$\!\!(m \t 1) \t 1$}(81.9942,998.1170)
    (88.1242,1000.8286) .. controls (108.0234,1009.9228) and (125.3927,1020.6737) ..
    node[above right=0.1cm,pos=0.22] {$m \t 1$}
    (139.6961,1029.7577);
  \path[fill=black] (49.917496,1028.8956) node[circle, draw, line width=0.65pt,
    minimum width=5mm, fill=white, inner sep=0.25mm] (text3309) {$\ass$    };
  \path[fill=black] (139.03424,1031.976) node[circle, draw, line width=0.65pt,
    minimum width=5mm, fill=white, inner sep=0.25mm] (text3317) {$\ass$  };
\end{tikzpicture}}}
\end{equation*}
in  $\V(\,((\B(y,z) \otimes  \B(x,y)) \otimes \B(w,x)) \otimes \B(v,w),\, \B(v,z)\,)$.

\begin{Ex}\label{ex:unit-vcat}
  The \emph{unit $\V$-bicategory} $\I$ has one object, $\star$, with $\I(\star,\star)=I$, $j_\star = 1_I$, $m_{\star\star\star} = \mathfrak l$, and with $\sigma$, $\tau$, and $\alpha$ constructed from coherence cells for $\V$ (with $\tau$ involving $\theta$).
\end{Ex}

\subsection{Underlying bicategory}
To any $\V$-bicategory $\B$ we may associate an ordinary bicategory $\B_0$ with the same objects as $\B$, and with hom-categories $\B_0(x,y) = \V(I, \B(x,y))$. The identity morphism at $x \in \B_0$ is  $j_x \in \V(I, \B(x,x))$, while composition is mediated by the functors:
\begin{equation*}
\V(I, \B(y,z)) \times \V(I, \B(x,y)) \xrightarrow{\otimes} \V(I \otimes I, \B(y, z) \otimes \B(x,y)) \xrightarrow{\V(\mathfrak l^\centerdot, m)} \V(I,\B(x,z))\rlap{ .}
\end{equation*}
The left and right unit constraint $2$-cells for $\B_0$ are constructed from $\lu$ and $\ru$ respectively, with the right constraint also involving $\theta$; the associativity constraint $2$-cells are built from $\ass$ and $\lambda$.

\subsection{Hom-functors}\label{subsec:hom-functors}
For any  $\V$-bicategory $\B$ and any $x \in \B$, there is a functor $\B(x, \thg) \colon \B_0 \to \V$ that on objects sends $y$ to $\B(x,y)$, on $1$-cells sends $f \colon y \to z$ to the composite
\begin{equation*}
\B(1, f) \defeq \B(x,y) \xrightarrow{\mathfrak l^\centerdot} I \otimes \B(x,y) \xrightarrow{f \otimes 1} \B(y,z) \otimes \B(x,y) \xrightarrow m \B(x,z)\rlap{ ,}
\end{equation*}
and on $2$-cells sends $\gamma \colon f \Rightarrow g$ to $m \circ (\gamma \otimes 1) \circ \mathfrak l^\centerdot$. The nullary and binary functoriality constraints for $\B(x, \thg)$ 
are built from $\lu$ and from $\ass$ and $\lambda$, respectively.
Similarly, we obtain for each $y \in \B$ a functor $\B(\thg, y) \colon \B_0^\op \to \V$ that sends $x$ to $\B(x,y)$, sends $g \colon w \to x$ to the composite
\begin{equation*}
\B(g, 1) \defeq \B(y,x) \xrightarrow{\mathfrak r^\centerdot}  \B(x,y) \otimes I \xrightarrow{1 \otimes g} \B(x,y) \otimes \B(w,x) \xrightarrow m \B(w,y)\rlap{ ,}
\end{equation*}
and sends a $2$-cell $\gamma$ to $m \circ (1 \otimes \gamma) \circ \mathfrak r^\centerdot$; its binary and nullary functoriality constraints built from $\tau$ and from $\alpha$ and $\rho$, respectively.
 
For any maps $f \colon y \to z$ and $g \colon w \to x$ in $\B_0$, there is an interchange isomorphism $\B(g,1) \circ \B(1, f) \Rightarrow \B(1,f) \circ \B(g,1) \colon \B(x,y) \to \B(w,z)$ built from $\ass$ and $\lambda$, and
using these, we may assemble together the functors defined above into a functor $\B(\thg, \thg) \colon \B_0^\op \times \B_0 \to \V$. With respect to this functor structure, the unit and composition maps $j_x \colon I \to \B(x,x)$ and $m_{xyz} \colon \B(y,z) \otimes \B(x,y) \to \B(x,z)$ now become pseudonatural in each variable, in that we have invertible $2$-cells
\begin{equation}\label{eq:pseudonat-of-composition}
\begin{gathered}
\cd{
  I \ar[r]^-{j_x} \ar[d]_{j_y} \rtwocell{dr}{\upsilon_f} & \B(x,x) \ar[d]^{\B(1,f)} \\
  \B(y,y) \ar[r]_-{\B(f,1)} & \B(x,y)
} \qquad
\cd[@C+1em]{
  \B(y,z) \otimes \B(x,y) \rtwocell{dr}{\ass_{fyz}} \ar[d]_{m_{xyz}} \ar[r]^-{1 \otimes \B(f,1)} &
  \B(y,z) \otimes \B(w,y) \ar[d]^{m_{wyz}} \\
  \B(x,z) \ar[r]_-{\B(f,1)}  & \B(w,z)
}\\
\cd[@C+0.5em]{
  \B(y,z) \otimes \B(w,x) \rtwocell{dr}{\ass_{wfz}} \ar[d]_{\B(f,1) \otimes 1} \ar[r]^-{1 \otimes \B(1,f)} &
  \B(y,z) \otimes \B(w,y) \ar[d]^{m_{wyz}} \\
  \B(x,z) \otimes \B(w,x) \ar[r]_-{m_{wxz}} &
  \B(w,z)
}\quad \text{and} \quad
\cd{
  \B(x,y) \otimes \B(w,x) \rtwocell{dr}{\ass_{wxf}} \ar[r]^-{\vphantom{\B}m_{wxy}} \ar[d]_{\B(1,f) \otimes 1} & \B(w,y) \ar[d]^{\B(1,f)} \\
  \B(x,z) \otimes \B(w,x) \ar[r]_-{m_{wxz}} & \B(w,z)
}
\end{gathered}\end{equation}satisfying coherence axioms. The $2$-cells $\upsilon$ are built using $\lu$ and $\ru$ for $\B$, together with two instances of $\theta$; while those of the remaining three kinds are built from instances of $\ass$ for $\B$, together with $\rho$, $\nu$ and $\lambda$ respectively.

\subsection{$\V$-functors}
\label{subsec:vfunctors}

If $\B$ and $\C$ are $\V$-bicategories, then a \emph{$\V$-functor} $\B \to \C$ is given by:
\begin{itemize}
\item A function $F \colon \ob \B \to \ob \C$;
\item For each $x, y \in \B$, a morphism $F_{xy} \colon \B(x,y) \to \C(Fx,Fy)$ in $\V$;
\item For all $x \in \B$ and for all $x,y,z \in \B$, invertible $2$-cells
\begin{equation*}
\cd[@!C@C-2.5em]{
  I \ar[rr]^j \ar[dr]_{j} & & \C(Fx,Fx) \\
  & \B(x, x) \ar[ur]_F \dtwocell[0.55]{u}{\fu_x}
} \quad \text{and} \quad
\cd{
\B(y, z) \otimes \B(x, y)  \ar[d]_m \ar[r]^-{F \otimes F} \dtwocell{dr}{\fm_{xyz}} &
\C(Fy,Fz) \otimes \C(Fx,Fy) \ar[d]^m \\
\B(x, z) \ar[r]_-F  & \C(Fx,Fz)\rlap{ ,}
}
\end{equation*}
\end{itemize}
 subject to the axioms that for all $x, y \in \B$, we have:
\begin{equation*}
\vcenter{\hbox{\begin{tikzpicture}[y=0.80pt, x=0.75pt,yscale=-1, inner sep=0pt, outer sep=0pt, every text node part/.style={font=\tiny} ]
\path[draw=black,line join=miter,line cap=butt,line width=0.650pt]
  (15.0000,1007.3622) .. controls (27.1429,1007.0051) and (35.0000,1005.9335) ..
  node[above right=0.15cm,at start] {$1 \t j$}(44.6429,1015.2193);
\path[draw=black,line join=miter,line cap=butt,line width=0.650pt]
  (15.0000,1028.0765) .. controls (29.6429,1027.7193) and (36.4286,1027.7194) ..
  node[above right=0.15cm,at start] {$m$}(44.6429,1017.7194);
\path[draw=black,line join=miter,line cap=butt,line width=0.650pt]
  (46.4286,1016.2907) .. controls (71.7857,1016.2907) and (96.0714,988.0765) ..
  node[above left=0.12cm,at end] {$\mathfrak r$}(120.0000,987.7193)
  (15.3571,987.7193) .. controls (42.0379,987.4444) and
  (60.6951,995.4383) .. node[above right=0.15cm,at start] {$F \t 1$}(77.3631,1002.9319)
  (83.0016,1005.4652) .. controls
  (95.1081,1010.8639) and (106.4966,1015.4336) ..
  node[above left=0.14cm,at end] {$F\!$}(120, 1015.5765);
\path[fill=black] (45,1016.2908) node[circle, draw, line width=0.65pt, minimum
  width=5mm, fill=white, inner sep=0.25mm] (text4302) {$\ru$  };
\end{tikzpicture}}}
\qquad = \qquad
\vcenter{\hbox{\begin{tikzpicture}[y=0.80pt, x=0.7pt,yscale=-1, inner sep=0pt, outer sep=0pt, every text node part/.style={font=\tiny} ]
  \path[draw=black,line join=miter,line cap=butt,line width=0.650pt]
    (0.0000,1029.7675) .. controls (12.1429,1029.4104) and (14.7626,1029.0768) ..
    node[above right=0.15cm,at start] {$1 \t j$}(30.1198,1029.7915);
  \path[draw=black,line join=miter,line cap=butt,line width=0.650pt]
    (0.0000,1052.3622) .. controls (14.6429,1052.0050) and (135.1479,1052.3952) ..
    node[above right=0.15cm,at start] {$m$}(144.1198,1036.8389);
  \path[draw=black,line join=miter,line cap=butt,line width=0.650pt]
    (31.9048,1026.9337) .. controls (34.8704,1021.0641) and (44.5979,1016.0173) ..
    node[below right=0.05cm,pos=0.58] {$1 \t j$}(57.8107,1012.0740)
    (70.7201,1008.7946) .. controls (109.8433,1000.3402) and
    (165.8453,999.9424) ..
    node[above=0.08cm,pos=0.45] {$1 \t j$} (180.9055,1012.5532)
    (0.0000,1007.4579) .. controls
    (63.8853,1003.6139) and (136.9531,1024.0533) ..
    node[above right=0.15cm,at start] {$F \t 1$}
    node[above right=0.07cm,pos=0.72] {$F \t F$}(144.4531,1034.0532);
  \path[draw=black,line join=miter,line cap=butt,line width=0.650pt]
    (31.9293,1031.2198) .. controls (41.2150,1039.4341) and (94.4774,1030.8420) ..
    node[below right=0.05cm,pos=0.58] {$1 \t F$}(111.7478,1019.1597);
  \path[draw=black,line join=miter,line cap=butt,line width=0.650pt]
    (146.9127,1035.8324) .. controls (160.4841,1035.4753) and (168.7150,1023.6008)
    .. node[below right=0.07cm] {$m$} (181.9293,1015.3865);
  \path[draw=black,line join=miter,line cap=butt,line width=0.650pt]
    (145.5721,1036.8151) .. controls (145.5721,1036.8151) and (180.5721,1050.7197)
    .. node[above left=0.12cm,at end] {$F\!$} (222.8579,1050.0293);
  \path[draw=black,line join=miter,line cap=butt,line width=0.650pt]
    (182.6436,1013.6008) .. controls (193.6912,1012.5532) and (193.2150,1021.8388)
    .. node[above left=0.12cm,at end] {$\mathfrak r$} (222.8579,1021.4579);
  \path[shift={(34.179628,88.617761)},draw=black,fill=black]
    (77.8571,931.1122)arc(0.000:180.000:1.250)arc(-180.000:0.000:1.250) -- cycle;
  \path[fill=black] (31.215008,1030.1722) node[circle, draw, line width=0.65pt,
    minimum width=5mm, fill=white, inner sep=0.25mm] (text4302) {$1 \t \fu$    };
  \path[fill=black] (144.50072,1036.1007) node[circle, draw, line width=0.65pt,
    minimum width=5mm, fill=white, inner sep=0.25mm] (text4373) {$\fm$    };
  \path[fill=black] (181.57216,1013.2437) node[circle, draw, line width=0.65pt,
    minimum width=5mm, fill=white, inner sep=0.25mm] (text4383) {$\ru$  };
\end{tikzpicture}}}
\end{equation*}
in $\V(\B(x,y) \otimes I, \C(Fx,Fy))$, and
\begin{equation*}
\vcenter{\hbox{\begin{tikzpicture}[y=0.80pt, x=0.75pt,yscale=-1, inner sep=0pt, outer sep=0pt, every text node part/.style={font=\tiny} ]
\path[draw=black,line join=miter,line cap=butt,line width=0.650pt]
  (15.0000,1007.3622) .. controls (27.1429,1007.0051) and (35.0000,1005.9335) ..
  node[above right=0.15cm,at start] {$j \t 1$}(44.6429,1015.2193);
\path[draw=black,line join=miter,line cap=butt,line width=0.650pt]
  (15.0000,1028.0765) .. controls (29.6429,1027.7193) and (36.4286,1027.7194) ..
  node[above right=0.15cm,at start] {$m$}(44.6429,1017.7194);
\path[draw=black,line join=miter,line cap=butt,line width=0.650pt]
  (46.4286,1016.2907) .. controls (71.7857,1016.2907) and (96.0714,988.0765) ..
  node[above left=0.12cm,at end] {$\mathfrak l$}(120.0000,987.7193)
  (15.3571,987.7193) .. controls (42.0379,987.4444) and
  (60.6951,995.4383) .. node[above right=0.15cm,at start] {$1 \t F$}(77.3631,1002.9319)
  (83.0016,1005.4652) .. controls
  (95.1081,1010.8639) and (106.4966,1015.4336) ..
  node[above left=0.14cm,at end] {$F\!$}(120, 1015.5765);
\path[fill=black] (45,1016.2908) node[circle, draw, line width=0.65pt, minimum
  width=5mm, fill=white, inner sep=0.25mm] (text4302) {$\lu$  };
\end{tikzpicture}}}
\qquad = \qquad
\vcenter{\hbox{\begin{tikzpicture}[y=0.80pt, x=0.7pt,yscale=-1, inner sep=0pt, outer sep=0pt, every text node part/.style={font=\tiny} ]
  \path[draw=black,line join=miter,line cap=butt,line width=0.650pt]
    (0.0000,1029.7675) .. controls (12.1429,1029.4104) and (14.7626,1029.0768) ..
    node[above right=0.15cm,at start] {$j \t 1$}(30.1198,1029.7915);
  \path[draw=black,line join=miter,line cap=butt,line width=0.650pt]
    (0.0000,1052.3622) .. controls (14.6429,1052.0050) and (135.1479,1052.3952) ..
    node[above right=0.15cm,at start] {$m$}(144.1198,1036.8389);
  \path[draw=black,line join=miter,line cap=butt,line width=0.650pt]
    (31.9048,1026.9337) .. controls (34.8704,1021.0641) and (44.5979,1016.0173) ..
    node[below right=0.05cm,pos=0.58] {$j\t 1$}(57.8107,1012.0740)
    (70.7201,1008.7946) .. controls (109.8433,1000.3402) and
    (165.8453,999.9424) ..
    node[above=0.08cm,pos=0.45] {$j\t 1$} (180.9055,1012.5532)
    (0.0000,1007.4579) .. controls
    (63.8853,1003.6139) and (136.9531,1024.0533) ..
    node[above right=0.15cm,at start] {$1\t F$}
    node[above right=0.07cm,pos=0.72] {$F \t F$}(144.4531,1034.0532);
  \path[draw=black,line join=miter,line cap=butt,line width=0.650pt]
    (31.9293,1031.2198) .. controls (41.2150,1039.4341) and (94.4774,1030.8420) ..
    node[below right=0.05cm,pos=0.58] {$F \t 1$}(111.7478,1019.1597);
  \path[draw=black,line join=miter,line cap=butt,line width=0.650pt]
    (146.9127,1035.8324) .. controls (160.4841,1035.4753) and (168.7150,1023.6008)
    .. node[below right=0.07cm] {$m$} (181.9293,1015.3865);
  \path[draw=black,line join=miter,line cap=butt,line width=0.650pt]
    (145.5721,1036.8151) .. controls (145.5721,1036.8151) and (180.5721,1050.7197)
    .. node[above left=0.12cm,at end] {$F\!$} (222.8579,1050.0293);
  \path[draw=black,line join=miter,line cap=butt,line width=0.650pt]
    (182.6436,1013.6008) .. controls (193.6912,1012.5532) and (193.2150,1021.8388)
    .. node[above left=0.12cm,at end] {$\mathfrak l$} (222.8579,1021.4579);
  \path[shift={(34.179628,88.617761)},draw=black,fill=black]
    (77.8571,931.1122)arc(0.000:180.000:1.250)arc(-180.000:0.000:1.250) -- cycle;
  \path[fill=black] (31.215008,1030.1722) node[circle, draw, line width=0.65pt,
    minimum width=5mm, fill=white, inner sep=0.25mm] (text4302) {$\fu \t 1$    };
  \path[fill=black] (144.50072,1036.1007) node[circle, draw, line width=0.65pt,
    minimum width=5mm, fill=white, inner sep=0.25mm] (text4373) {$\fm$    };
  \path[fill=black] (181.57216,1013.2437) node[circle, draw, line width=0.65pt,
    minimum width=5mm, fill=white, inner sep=0.25mm] (text4383) {$\lu$  };
\end{tikzpicture}}}
\end{equation*}
in $\V(I \otimes \B(x,y), \C(Fx,Fy))$, and that for all $w,x,y,z \in \B$, we have
\begin{equation*}
\vcenter{\hbox{\begin{tikzpicture}[y=0.80pt, x=0.8pt,yscale=-1, inner sep=0pt, outer sep=0pt, every text node part/.style={font=\tiny} ]
   \path[use as bounding box] (-11.7, 978) rectangle (239, 1064);
\path[draw=black,line join=miter,line cap=butt,line width=0.650pt]
    (-11.6628,992.2856) .. controls (21.8267,992.5650) and (77.1166,991.3917) ..
    node[above right=0.15cm,at start] {$\!(F \t F) \t F$}
    node[above right=0.15cm,pos=0.5] {$\!(F \t F) \t 1$}(91.9540,1001.3431);
  \path[draw=black,line join=miter,line cap=butt,line width=0.650pt]
    (-11.6628,1052.2090) .. controls (-11.6628,1052.2090) and (148.9540,1054.3431)
    .. node[above right=0.15cm,at start] {$m$} (162.2874,1036.6764);
  \path[draw=black,line join=miter,line cap=butt,line width=0.650pt]
    (-11.6628,1024.1630) .. controls (28.3372,1024.1630) and (82.8544,1026.4580)
    .. node[above right=0.15cm,at start] {$m \t 1$}(91.9540,1003.3430)
    (32.6743,993.0519) .. controls (36.0759,1004.3545) and
    (46.0440,1013.0991) .. node[above right=0.05cm,pos=0.6] {$1 \t F$}(58.7983,1018.9604)
    (67.3929,1022.3443) .. controls
    (91.7588,1030.4851) and (122.2897,1029.3916) ..
    node[below left=0.12cm,pos=0.6] {$1 \t F$}(138.0843,1017.2664);
  \path[draw=black,line join=miter,line cap=butt,line width=0.650pt]
    (92.9540,1001.3431) .. controls (99.9540,993.6764) and (179.2874,982.6764) ..
    node[above=0.1cm] {$m \t 1$} (192.9540,997.6764);
  \path[draw=black,line join=miter,line cap=butt,line width=0.650pt]
    (93.6207,1003.3431) .. controls (105.6207,1003.3431) and (155.6207,1020.6764)
    .. node[above right=0.08cm,pos=0.32] {$F \t 1$}
    node[above right=0.08cm,pos=0.7] {$F \t F$}(162.2874,1033.6764);
  \path[draw=black,line join=miter,line cap=butt,line width=0.650pt]
    (164.6207,1035.0097) .. controls (181.6207,1024.3431) and (174.2874,1006.0097)
    .. node[below right=0.1cm,pos=0.5] {$m$} (192.6207,999.6764);
  \path[draw=black,line join=miter,line cap=butt,line width=0.650pt]
    (164.9540,1037.3431) .. controls (174.6207,1056.3431) and (214.9540,1043.3431)
    .. node[above left=0.12cm,at end] {$F$}(238.2874,1044.3431);
  \path[draw=black,line join=miter,line cap=butt,line width=0.650pt]
    (195.2874,997.6764) .. controls (195.2874,997.6764) and (209.6207,983.3430) ..
    node[above left=0.12cm,at end] {$\mathfrak a$}(239.6207,983.3430);
  \path[draw=black,line join=miter,line cap=butt,line width=0.650pt]
    (195.2874,999.6764) .. controls (195.2874,999.6764) and (214.6207,1005.3431)
    .. node[above left=0.12cm,at end] {$1 \t m$}(239.6207,1005.3431);
  \path[draw=black,line join=miter,line cap=butt,line width=0.650pt]
    (194.2874,1001.3431) .. controls (194.2874,1001.3431) and (203.6207,1026.0097)
    .. node[above left=0.12cm,at end] {$m$}(238.9540,1025.6764);
  \path[fill=black] (100,964.02887) node[above right] (text3136) {     };
  \path[fill=black] (92.287354,1002.3431) node[circle, draw, line width=0.65pt,
    minimum width=5mm, fill=white, inner sep=0.25mm] (text3140) {$\fm \t 1$    };
  \path[fill=black] (163.28735,1035.6764) node[circle, draw, line width=0.65pt,
    minimum width=5mm, fill=white, inner sep=0.25mm] (text3144) {$\fm$    };
  \path[fill=black] (193.95403,999.34308) node[circle, draw, line width=0.65pt,
    minimum width=5mm, fill=white, inner sep=0.25mm] (text3148) {$\ass$    };
  \path[shift={(-43.656948,61.096743)},draw=black,fill=black]
    (77.8571,931.1122)arc(0.000:180.000:1.250)arc(-180.000:0.000:1.250) -- cycle;
  \path[shift={(60.940753,86.000958)},draw=black,fill=black]
    (77.8571,931.1122)arc(0.000:180.000:1.250)arc(-180.000:0.000:1.250) -- cycle;
\end{tikzpicture}}} \ \ = \ \
\vcenter{\hbox{\begin{tikzpicture}[y=0.80pt, x=0.8pt,yscale=-1, inner sep=0pt, outer sep=0pt, every text node part/.style={font=\tiny} ]
  \path[use as bounding box] (0, 990) rectangle (260, 1076);
  \path[draw=black,line join=miter,line cap=butt,line width=0.650pt]
    (0.5977,1052.2090) .. controls (16.3065,1052.5921) and (29.2069,1052.4274) ..
    node[above right=0.15cm,at start] {$m$}(38.3257,1041.2741);
  \path[draw=black,line join=miter,line cap=butt,line width=0.650pt]
    (0.5977,1024.1630) .. controls (15.5977,1023.1630) and (28.6245,1025.6917) ..
    node[above right=0.15cm,at start] {$m \t 1$}(38.8736,1038.5921);
  \path[draw=black,line join=miter,line cap=butt,line width=0.650pt]
    (147.3793,1020.4503) .. controls (164.3793,1009.7837) and (232.5249,1009.4579)
    .. node[above left=0.12cm,at end] {$1 \t m$} (258.9042,1009.6381);
  \path[draw=black,line join=miter,line cap=butt,line width=0.650pt]
    (146.9464,1024.3163) .. controls (156.6130,1043.3163) and (193.4981,1024.9523)
    .. node[above right=0.1cm,pos=0.26] {$1 \t F$}
    node[above right=0.1cm,pos=0.68] {$F \t F$}(215.2989,1040.8948);
  \path[draw=black,line join=miter,line cap=butt,line width=0.650pt]
    (41.0575,1038.6725) .. controls (41.4271,1001.7354) and (72.7867,992.0723) ..
    (117.5280,991.6088) .. controls (150.7190,991.2648) and (206.3195,990.5078) ..
    node[above left=0.12cm,at end] {$\mathfrak a$}(258.3946,990.2396)
    (0.5977,997.5889) .. controls (24.5112,997.3873) and
    (36.9043,1002.3448) .. node[above right=0.15cm,at start] {$\!(F \t F) \t F$} (49.9130,1007.7561)
    (54.9332,1009.8339) .. controls
    (72.1492,1016.8567) and (92.9504,1023.6391) ..
    node[above right=0.05cm,pos=0.2] {$\!\!F \t (F \t F)$}
    node[above right=0.07cm,pos=0.66] {$\!\!1 \t (F \t F)$}(143.2950,1020.1170);
  \path[draw=black,line join=miter,line cap=butt,line width=0.650pt]
    (41.0575,1040.6725) .. controls (41.0575,1040.6725) and (68.7834,1043.1403) ..
    node[above=0.08cm,pos=0.7] {$1 \t m$}(95.5968,1041.1633)
    (102.2951,1040.5693) .. controls (121.1463,1038.5961) and
    (138.3175,1034.0731) .. node[below=0.12cm,pos=0.62] {$1 \t m$} (143.4521,1024.5006)
    (88.7634,1019.0597) .. controls
    (91.8979,1061.3165) and (159.4440,1061.3657) ..
    node[below=0.1cm,pos=0.62] {$F \t 1$} (186.0954,1033.7415);
  \path[draw=black,line join=miter,line cap=butt,line width=0.650pt]
    (40.0575,1042.3393) .. controls (97.5640,1090.1537) and (184.1058,1075.3070)
    .. node[below=0.1cm,pos=0.42] {$m$}  (215.5824,1042.5346);
  \path[fill=black] (144.01149,1022.6496) node[circle, draw, line width=0.65pt,
    minimum width=5mm, fill=white, inner sep=0.25mm] (text3140) {$1 \t \fm$    };
  \path[fill=black] (41.724148,1040.3392) node[circle, draw, line width=0.65pt,
    minimum width=5mm, fill=white, inner sep=0.25mm] (text3148) {$\ass$    };
  \path[shift={(12.568792,88.26273)},draw=black,fill=black]
    (77.8571,931.1122)arc(0.000:180.000:1.250)arc(-180.000:0.000:1.250) -- cycle;
  \path[shift={(109.52888,102.7749)},draw=black,fill=black]
    (77.8571,931.1122)arc(0.000:180.000:1.250)arc(-180.000:0.000:1.250) -- cycle;
  \path[draw=black,line join=miter,line cap=butt,line width=0.650pt]
    (216.8582,1040.8679) .. controls (224.5211,1028.6074) and (242.1456,1028.6075)
    .. node[above left=0.12cm,at end] {$m$} (259.0038,1028.2243);
  \path[draw=black,line join=miter,line cap=butt,line width=0.650pt]
    (217.6245,1043.5499) .. controls (224.9042,1049.2970) and (237.9310,1049.2971)
    .. node[above left=0.12cm,at end] {$F$} (259.7701,1048.9139);
  \path[fill=black] (216.13411,1042.0825) node[circle, draw, line width=0.65pt,
    minimum width=5mm, fill=white, inner sep=0.25mm] (text3144) {$\fm$  };
\end{tikzpicture}}}
\end{equation*}
in $\V(\,(\B(y,z) \otimes \B(x,y)) \otimes \B(w,x) ,\, \C(Fw, Fz)\,)$.

We call a $\V$-functor $F \colon \B \to \C$ \emph{fully faithful} if each $F_{xy}$ is an equivalence in $\V$.

\subsection{Underlying ordinary functor}
Given a $\V$-functor $F \colon \B \to \C$, there is an ordinary functor $F_0 \colon \B_0 \to \C_0$ whose action on objects is that of $F$, and whose action on hom-categories is given by $\V(1, F_{xy}) \colon \V(I, \B(x,y)) \to \V(I, \C(Fx, Fy))$. The nullary and binary functoriality constraints of $F_0$ are obtained using $\fu$ and $\fm$ respectively. We may without ambiguity write the action of $F_0$ on a $1$-cell $f$ or $2$-cell $\alpha$ of $\B_0$ as $Ff$ or $F\alpha$, respectively.
The interaction of the underlying functor $F_0$ of $F$ with the hom-functors of Section~\ref{subsec:hom-functors} is expressed by the existence of invertible
$2$-cells
\begin{equation}\label{eq:functoriality-ordinary-2cells}
\cd{
 \B(x,y) \ar[r]^-{F_{xy}} \ar[d]_{\B(1,f)} \dtwocell{dr}{\fm_{xf}} & \C(Fx,Fy) \ar[d]^{\C(1, Ff)} \\ 
 \B(x,z) \ar[r]_-{F_{xz}} & \C(Fx,Fz)
} \qquad \text{and} \qquad
\cd{
 \B(x,y) \ar[r]^-{F_{xy}} \ar[d]_{\B(f,1)} \dtwocell{dr}{\fm_{fy}} & \C(Fx,Fy) \ar[d]^{\C(Ff,1)} \\ 
 \B(w,y) \ar[r]_-{F_{wy}} & \C(Fw,Fy)
}\end{equation}
built from the binary constraint cells $\fm$ for $F$.

\subsection{$\V$-transformations}
Let $F, G \colon \B \to \C$ be $\V$-functors. A \emph{$\V$-transformation} $\gamma \colon F \Rightarrow G$ is given by:
\begin{itemize}
\item For each $x \in \B$, a morphism $\gamma_x \colon Fx \to Gx$ in $\C_0$;
\item For each $x,y \in \B$, an invertible $2$-cell
\begin{equation*}
\cd[@+0.7em]{
\B(x, y) \dtwocell{dr}{{\bar \gamma}_{xy}} \ar[d]_{F} \ar[r]^-{G} &
\C(Gx, Gy) \ar[d]^{\C(\gamma_x,1)} \\
\C(Fx,Fy) \ar[r]_-{\C(1,\gamma_y)}
& \C(Fx,Gy) \rlap{ ;}
}
\end{equation*}
\end{itemize}
subject to the axioms that for all $x \in \B$, we have
\begin{equation*}
\vcenter{\hbox{\begin{tikzpicture}[y=0.80pt, x=0.8pt,yscale=-1, inner sep=0pt, outer sep=0pt, every text node part/.style={font=\tiny} ]
  \path[draw=black,line join=miter,line cap=butt,line width=0.650pt]
  (0.0000,1012.3622) .. controls (41.7187,1012.3622) and (68.9821,1012.8078) ..
  node[above right=0.12cm,at start] {\!$\C(\gamma,1)$}(79.9759,1012.8078);
  \path[draw=black,line join=miter,line cap=butt,line width=0.650pt]
  (0.0000,992.3622) --
  node[above right=0.12cm,at start] {\!$j$}(40.0000,992.3622);
  \path[draw=black,line join=miter,line cap=butt,line width=0.650pt]
  (41.8397,993.9571) .. controls (49.7762,1004.2745) and (73.3379,999.4464) ..
  node[above=0.07cm] {$G$}(79.5451,1010.2312);
  \path[draw=black,line join=miter,line cap=butt,line width=0.650pt]
  (41.6994,989.9888) .. controls (48.4567,981.3801) and (98.3921,982.3622) ..
  node[above left=0.12cm,at end] {$j$\!}(130.0000,982.3622);
  \path[draw=black,line join=miter,line cap=butt,line width=0.650pt]
  (81.4612,1011.1949) .. controls (83.3095,1007.9935) and (93.6846,1002.3622) ..
  node[above left=0.12cm,at end] {$F$\!}(130.0000,1002.3622);
  \path[draw=black,line join=miter,line cap=butt,line width=0.650pt]
  (81.4201,1013.7020) .. controls (83.4526,1017.2224) and (93.1184,1022.3622) ..
  node[above left=0.12cm,at end] {$\C(1, \gamma)$\!}(130.0000,1022.3622);
  \path[fill=black] (41.068821,992.23395) node[circle, draw, line width=0.65pt, minimum width=5mm, fill=white, inner sep=0.25mm] (text3765) {$\fu$     };
  \path[fill=black] (81.122116,1012.3558) node[circle, draw, line width=0.65pt, minimum width=5mm, fill=white, inner sep=0.25mm] (text3769) {$\bar \gamma$   };
\end{tikzpicture}}}\quad = \quad
\vcenter{\hbox{\begin{tikzpicture}[y=0.80pt, x=0.7pt,yscale=-1, inner sep=0pt, outer sep=0pt, every text node part/.style={font=\tiny} ]
  \path[draw=black,line join=miter,line cap=butt,line width=0.650pt]
  (5.0000,1022.3622) .. controls (38.8112,1022.3622) and (49.9918,1018.2842) ..
  node[above right=0.12cm,at start] {\!$\C(\gamma,1)$}(52.8689,1014.0693);
  \path[draw=black,line join=miter,line cap=butt,line width=0.650pt]
  (5.0000,1002.6428) .. controls (29.2639,1002.6428) and (49.9178,1004.8735) ..
  node[above right=0.12cm,at start] {\!$j$}(53.0044,1010.7661);
  \path[draw=black,line join=miter,line cap=butt,line width=0.650pt]
  (93.6728,1004.3392) .. controls (98.6331,1011.0852) and (110.0000,1012.3622) ..
  node[above left=0.12cm,at end] {$F$\!}(130.0000,1012.3622);
  \path[draw=black,line join=miter,line cap=butt,line width=0.650pt]
  (93.5325,1000.3709) .. controls (98.8157,993.6365) and (110.5747,992.3622) ..
  node[above left=0.12cm,at end] {$j$\!}(130.0000,992.3622);
  \path[draw=black,line join=miter,line cap=butt,line width=0.650pt]
  (55.8064,1010.5997) .. controls (61.8252,1005.6071) and (76.5035,1002.3622) ..
  node[above=0.07cm] {$j$}(92.0000,1002.3622);
  \path[draw=black,line join=miter,line cap=butt,line width=0.650pt]
  (56.0317,1014.3463) .. controls (63.5306,1027.3347) and (84.3372,1032.3622) ..
  node[above left=0.12cm,at end] {$\C(1, \gamma)$\!}(130.0000,1032.3622);
  \path[fill=black] (92.901894,1002.616) node[circle, draw, line width=0.65pt, minimum width=5mm, fill=white, inner sep=0.25mm] (text3765) {$\fu$     };
  \path[fill=black] (54.376083,1012.3557) node[circle, draw, line width=0.65pt, minimum width=5mm, fill=white, inner sep=0.25mm] (text3769) {$\upsilon$   };
\end{tikzpicture}}}
\end{equation*}
in $\V(I, \C(Fx, Gx))$, and that for all $x, y, z \in \B$, we have
\begin{equation*}
\vcenter{\hbox{
\begin{tikzpicture}[y=0.8pt, x=0.8pt,yscale=-1, inner sep=0pt, outer sep=0pt, every text node part/.style={font=\tiny} ]
  \path[draw=black,line join=miter,line cap=butt,line width=0.650pt]
  (271.0231,1000.8587) .. controls (274.1128,997.7689) and (285.6349,992.3622) ..
  node[above left=0.12cm,at end] {$m$\!}(305.0000,992.3622);
  \path[draw=black,line join=miter,line cap=butt,line width=0.650pt]
  (271.2555,1003.6365) .. controls (273.8520,1006.2330) and (285.2381,1012.3622) ..
  node[above left=0.12cm,at end] {$F$\!}(305.0000,1012.3622);
  \path[draw=black,line join=miter,line cap=butt,line width=0.650pt]
  (169.3556,1006.0635) .. controls (178.0542,996.2426) and (261.1887,993.9978) ..
  node[above=0.07cm,pos=0.7] {$F \t F$}
  node[above=0.07cm,pos=0.3] {$F \t 1$}  (267.6425,1000.7322);
  \path[draw=black,line join=miter,line cap=butt,line width=0.650pt]
  (25.6025,993.4366) .. controls (39.1351,979.3201) and (60.6553,973.0472) ..
  (79.9702,973.0776) .. controls (122.3556,973.1444) and (163.1678,999.2120) ..
  node[above=0.08cm,at start] {$G \t 1$}(166.8303,1006.3441)
  (61.5912,1000.4516) .. controls (76.0857,989.5768) and (91.7953,982.8004) ..
  node[above left=0.07cm,pos=0.4] {$1 \t F$}(107.5282,978.9538)
  (118.3906,976.7476) .. controls (171.2588,968.1205) and (221.2706,989.7375) ..
  node[above=0.09cm,pos=0.46] {$1 \t F$}(222.7468,996.8038);
  \path[draw=black,line join=miter,line cap=butt,line width=0.650pt]
  (230.6035,1021.7770) .. controls (236.7767,1016.1650) and (263.7141,1009.7113) ..
  node[above left=0.07cm,pos=0.53] {$m$}(267.9231,1004.6605);
  \path[draw=black,line join=miter,line cap=butt,line width=0.650pt]
  (230.0423,1024.5830) .. controls (237.3578,1031.8985) and (296.2250,1032.3622) ..
  node[above left=0.12cm,at end] {$\C(1, \gamma)$\!}(305.0000,1032.3622);
  \path[draw=black,line join=miter,line cap=butt,line width=0.650pt]
  (168.7945,1008.0277) .. controls (173.5646,1013.3591) and (223.0274,1016.7262) ..
  node[above right=0.045cm,pos=0.3,rotate=-11] {$\C(1, \gamma) \t 1$}(227.5170,1022.6188);
  \path[draw=black,line join=miter,line cap=butt,line width=0.650pt]
  (111.4749,1020.9352) .. controls (114.2809,1015.8844) and (161.7795,1013.6397) ..
  node[above left=0.06cm,pos=0.74,rotate=8] {$\C(\gamma,1 ) \t 1$}(167.1109,1008.8695);
  \path[draw=black,line join=miter,line cap=butt,line width=0.650pt]
  (111.4749,1024.5830) .. controls (121.8571,1031.0368) and (220.5020,1031.3173) ..
  node[above=0.09cm,pos=0.46] {$m$}(228.0782,1024.3024);
  \path[draw=black,line join=miter,line cap=butt,line width=0.650pt]
  (61.7315,1004.3799) .. controls (66.7823,1011.1143) and (104.7406,1014.2008) ..
  node[above right=0.045cm,pos=0.22,rotate=-11] {$1 \t \C(1, \gamma)$}(108.1077,1020.6546);
  \path[draw=black,line join=miter,line cap=butt,line width=0.650pt]
  (0.0000,1032.3622) .. controls (31.7088,1032.3622) and (103.1106,1029.2756) ..
  node[above right=0.12cm,at start] {\!$m$}(107.8807,1024.2249);
  \path[draw=black,line join=miter,line cap=butt,line width=0.650pt]
  (0.0000,1012.3622) .. controls (31.1577,1012.3622) and (54.4896,1007.8727) ..
  node[above right=0.12cm,at start] {\!$1 \t \C(\gamma,1)$}(58.1373,1004.2249);
  \path[draw=black,line join=miter,line cap=butt,line width=0.650pt]
  (0.0000,992.3622) .. controls (30.0358,992.3622) and (54.4896,997.2726) ..
  node[above right=0.12cm,at start] {\!$G \t G$}
  node[above right=0.08cm,pos=0.46] {$1 \t G$}(58.1373,1000.9204);
  \path[fill=black] (58.925564,1003.5381) node[circle, draw, line width=0.65pt, minimum width=5mm, fill=white, inner sep=0.25mm] (text4450) {$1 \t \bar \gamma$     };
  \path[fill=black] (109.43319,1022.8994) node[circle, draw, line width=0.65pt, minimum width=5mm, fill=white, inner sep=0.25mm] (text4454) {$\ass$     };
  \path[fill=black] (167.79756,1007.4665) node[circle, draw, line width=0.65pt, minimum width=5mm, fill=white, inner sep=0.25mm] (text4458) {$\bar \gamma \t 1$     };
  \path[fill=black] (229.2485,1023.4606) node[circle, draw, line width=0.65pt, minimum width=5mm, fill=white, inner sep=0.25mm] (text4462) {$\ass$     };
  \path[fill=black] (268.81281,1001.8546) node[circle, draw, line width=0.65pt, minimum width=5mm, fill=white, inner sep=0.25mm] (text4466) {$\fm$     };
  \path[shift={(-50.758029,62.653701)},draw=black,fill=black] (77.8600,931.1000)arc(0.000:180.000:1.250)arc(-180.000:0.000:1.250) -- cycle;
  \path[shift={(146.00905,65.805858)},draw=black,fill=black] (77.8600,931.1000)arc(0.000:180.000:1.250)arc(-180.000:0.000:1.250) -- cycle;
\end{tikzpicture}}}\quad=\quad
\vcenter{\hbox{
\begin{tikzpicture}[y=0.80pt, x=0.8pt,yscale=-1, inner sep=0pt, outer sep=0pt, every text node part/.style={font=\tiny} ]
  \path[draw=black,line join=miter,line cap=butt,line width=0.650pt]
  (5.0000,962.3622) .. controls (35.8287,962.3622) and (83.6718,964.6070) ..
  node[above right=0.12cm,at start] {\!$G \t G$}(88.7225,970.7801);
  \path[draw=black,line join=miter,line cap=butt,line width=0.650pt]
  (5.0000,982.3622) .. controls (30.2050,982.3622) and (46.0758,987.0349) ..
  node[above right=0.12cm,at start] {\!$1 \t \C(\gamma,1)$}(48.3205,990.9633);
  \path[draw=black,line join=miter,line cap=butt,line width=0.650pt]
  (5.0000,1002.3622) .. controls (29.9402,1002.3622) and (45.7030,997.4882) ..
  node[above right=0.12cm,at start] {\!$m$}(48.2284,994.1210);
  \path[draw=black,line join=miter,line cap=butt,line width=0.650pt]
  (53.4418,994.3652) .. controls (57.6507,998.8547) and (125.1523,998.8547) ..
  node[above=0.07cm] {$\C(\gamma,1)$}(128.5195,993.8040);
  \path[draw=black,line join=miter,line cap=butt,line width=0.650pt]
  (53.4418,990.7801) .. controls (56.8089,984.6070) and (83.1253,979.1353) ..
  node[above left=0.07cm,pos=0.57] {$m$}(88.4567,974.0845);
  \path[draw=black,line join=miter,line cap=butt,line width=0.650pt]
  (91.8627,970.6398) .. controls (95.7910,965.4488) and (152.5839,962.3622) ..
  node[above left=0.12cm,at end] {$m$\!}(176.0000,962.3622);
  \path[draw=black,line join=miter,line cap=butt,line width=0.650pt]
  (91.8627,974.2248) .. controls (95.2298,979.1353) and (126.3522,986.8518) ..
  node[above right=0.07cm,pos=0.44] {$G$}(128.4567,990.9204);
  \path[draw=black,line join=miter,line cap=butt,line width=0.650pt]
  (131.6004,990.6306) .. controls (136.6512,985.0187) and (145.0546,981.6589) ..
  node[above left=0.12cm,at end] {$F$\!}(176.0000,982.3622);
  \path[draw=black,line join=miter,line cap=butt,line width=0.650pt]
  (131.6004,993.9978) .. controls (136.0900,998.7680) and (144.2540,1002.3622) ..
  node[above left=0.12cm,at end] {$\C(1, \gamma)$\!}(176.0000,1002.3622);
  \path[fill=black] (49.603172,993.23523) node[circle, draw, line width=0.65pt, minimum width=5mm, fill=white, inner sep=0.25mm] (text4238) {$\ass^{-1}$     };
  \path[fill=black] (89.682541,972.20343) node[circle, draw, line width=0.65pt, minimum width=5mm, fill=white, inner sep=0.25mm] (text4242) {$\fm$     };
  \path[fill=black] (129.36508,992.83838) node[circle, draw, line width=0.65pt, minimum width=5mm, fill=white, inner sep=0.25mm] (text4246) {$\bar \gamma$   };
\end{tikzpicture}}}
\end{equation*}
in $\V(\,(\B(y,z) \otimes \B(x,y)) \otimes I,\, \C(Fx,Gz)\,)$. If in the preceding definition, we remove the requirement that $\bar\gamma$ be invertible, we obtain a notion of \emph{lax $\V$-transformation}, while if we allow it to go in the other direction, we obtain \emph{oplax $\V$-transformations}. If necessary, we will call $\V$-transformations \emph{pseudo} to differentiate them from the lax and oplax variants.

If $\alpha \colon F \Rightarrow G$ is a pseudo, lax or oplax $\V$-transformation, then we obtain a pseudo, lax or oplax transformation $\alpha_0 \colon F_0 \Rightarrow G_0$ whose $1$-cell components are those of $\alpha$, and whose $2$-cell component at a map $f \colon x \to y$ of $\B_0$ is obtained by whiskering the $2$-cell $\bar \gamma_{xy}$ with $f \colon I \to \B(x,y)$. 

\subsection{$\V$-icons}
\label{sec:vicons}
 Following~\cite{Lack20082-nerves,lack:icons}, a $\V$-icon is essentially an oplax $\V$-trans\-formation whose $1$-cell components $\gamma_x$ are identities; as in the unenriched case, we can formulate an equivalent notion more simply.
Namely, if $F,G\colon \B\to\C$ are $\V$-functors that agree on objects, then a \emph{$\V$-icon} $\gamma\colon F\Rightarrow G$ is given by a collection of $2$-cells
\[\bar \gamma_{xy} \colon F_{xy} \Rightarrow G_{xy} \colon \B(x,y) \to \C(Fx,Fy) = \C(Gx,Gy)\] in $\V$ for all $x,y \in \B$,
subject to the axioms that for all $x\in \B$, the equality below on the left holds in $\V(I, \C(Gx,Gx))$, and that for all $x,y,z \in \B$, the equality below on the right holds in $\V(\B(y,z) \otimes \B(x,y),\C(Gx,Gz))$:
\begin{equation*}
\vcenter{\hbox{\begin{tikzpicture}[y=0.80pt, x=0.7pt,yscale=-1, inner sep=0pt, outer sep=0pt, every text node part/.style={font=\tiny} ]
  \path[use as bounding box] (60, 956) rectangle (125, 1003);
  \path[draw=black,line join=miter,line cap=butt,line width=0.650pt] (125.0000,972.3622) .. controls (110.0000,972.3622) and (100.0000,972.3622) .. 
  node[above left=0.12cm,at start] {$j$\!} (91.7678,979.4332);
  \path[draw=black,line join=miter,line cap=butt,line width=0.650pt] (125.0000,992.3622) .. controls (110.0000,992.3622) and (100.0000,991.6046) .. 
  node[above left=0.12cm,at start] {$G$\!} (91.6415,984.2810);
  \path[draw=black,line join=miter,line cap=butt,line width=0.650pt] (90.0000,982.3622) -- 
    node[above right=0.12cm,at end] {\!$j$} (60.0000,982.3622);
  \path[fill=black] (89.651039,982.15656) node[circle, draw, line width=0.65pt, minimum width=5mm, fill=white, inner sep=0.25mm] (text6017) {$\iota$     };
\end{tikzpicture}}}\quad = \quad
\vcenter{\hbox{\begin{tikzpicture}[y=0.80pt, x=0.7pt,yscale=-1, inner sep=0pt, outer sep=0pt, every text node part/.style={font=\tiny} ]
  \path[use as bounding box] (60, 956) rectangle (160, 1005);
  \path[draw=black,line join=miter,line cap=butt,line width=0.650pt] (160.0000,972.3622) .. controls (130.0000,972.3622) and (101.3156,969.8854) .. 
  node[above left=0.12cm,at start] {$j$\!}(91.7678,979.4332);
  \path[draw=black,line join=miter,line cap=butt,line width=0.650pt] (126.8185,992.3622) .. controls (99.6148,992.3622) and (95.1766,987.3784) .. 
  node[above right=0.08cm,pos=0.4] {$F$}(91.6415,984.2810);
  \path[draw=black,line join=miter,line cap=butt,line width=0.650pt] (90.0000,982.3622) -- 
  node[above right=0.12cm,at end] {\!$j$}(60.0000,982.3622);
  \path[draw=black,line join=miter,line cap=butt,line width=0.650pt] (160.0000,992.3622) -- 
  node[above left=0.12cm,at start] {$G$\!}(129.2424,992.3622);
  \path[fill=black] (89.651039,982.15656) node[circle, draw, line width=0.65pt, minimum width=5mm, fill=white, inner sep=0.25mm] (text6017) {$\iota$     };
  \path[fill=black] (129.50302,992.22101) node[circle, draw, line width=0.65pt, minimum width=5mm, fill=white, inner sep=0.25mm] (text6017-9) {$\bar \gamma$
   };
\end{tikzpicture}}}
\qquad\text{and}\qquad
\vcenter{\hbox{\begin{tikzpicture}[y=0.80pt, x=0.7pt,xscale=-1, inner sep=0pt, outer sep=0pt, every text node part/.style={font=\tiny} ]
  \path[use as bounding box] (50, 959) rectangle (160, 1012);
  \path[draw=black,line join=miter,line cap=butt,line width=0.650pt] (160.0000,972.3622) .. controls (110.0000,972.3622) and (100.0000,972.3622) .. 
  node[above right=0.12cm,at start] {\!$m$}(91.7678,979.4332);
  \path[draw=black,line join=miter,line cap=butt,line width=0.650pt] (122.5254,992.3622) .. controls (112.5254,992.3622) and (100.0000,992.3622) .. 
  node[above right=0.12cm,pos=0.35] {$GG$}(91.6415,984.2810);
  \path[draw=black,line join=miter,line cap=butt,line width=0.650pt] (87.6307,979.8837) .. controls (80.0000,972.3622) and (69.7239,972.3622) .. 
  node[above left=0.12cm,at end] {$G$\!}(50.0000,972.3622);
  \path[draw=black,line join=miter,line cap=butt,line width=0.650pt] (50.0000,992.3622) .. controls (70.0000,992.3622) and (80.0000,992.3622) .. 
  node[above left=0.12cm,at start] {$m$\!}(87.8833,983.6718);
  \path[draw=black,line join=miter,line cap=butt,line width=0.650pt] (160.0000,992.3622) -- 
  node[above right=0.12cm,at start] {\!$FF$}(124.5013,992.3622);
  \path[fill=black] (89.651039,982.15656) node[circle, draw, line width=0.65pt, minimum width=5mm, fill=white, inner sep=0.25mm] (text6017) {$\fm$
     };
  \path[fill=black] (124.53553,992.36218) node[circle, draw, line width=0.65pt, minimum width=5mm, fill=white, inner sep=0.25mm] (text6017-9) {$\bar \gamma \bar \gamma$
   };
  \end{tikzpicture}}}
  \quad=\quad
\vcenter{\hbox{\begin{tikzpicture}[y=0.80pt, x=0.7pt,yscale=-1, inner sep=0pt, outer sep=0pt, every text node part/.style={font=\tiny} ]
  \path[use as bounding box] (50, 954) rectangle (160, 1004);
  \path[draw=black,line join=miter,line cap=butt,line width=0.650pt] (160.0000,972.3622) .. controls (110.0000,972.3622) and (100.0000,972.3622) .. 
  node[above left=0.12cm,at start] {$m$\!}(91.7678,979.4332);
  \path[draw=black,line join=miter,line cap=butt,line width=0.650pt] (122.5254,992.3622) .. controls (112.5254,992.3622) and (100.0000,992.3622) .. 
  node[above right=0.08cm,pos=0.6] {$F$}(91.6415,984.2810);
  \path[draw=black,line join=miter,line cap=butt,line width=0.650pt] (87.6307,979.8837) .. controls (80.0000,972.3622) and (69.7239,972.3622) .. 
  node[above right=0.12cm,at end] {\!$F\t F$}(50.0000,972.3622);
  \path[draw=black,line join=miter,line cap=butt,line width=0.650pt] (50.0000,992.3622) .. controls (70.0000,992.3622) and (80.0000,992.3622) .. 
  node[above right=0.12cm,at start] {\!$m$}(87.8833,983.6718);
  \path[draw=black,line join=miter,line cap=butt,line width=0.650pt] (160.0000,992.3622) -- 
  node[above left=0.12cm,at start] {$G$\!}(124.5013,992.3622);
  \path[fill=black] (89.651039,982.15656) node[circle, draw, line width=0.65pt, minimum width=5mm, fill=white, inner sep=0.25mm] (text6017) {$\fm$
     };
  \path[fill=black] (130.53553,992.36218) node[circle, draw, line width=0.65pt, minimum width=5mm, fill=white, inner sep=0.25mm] (text6017-9) {$\bar \gamma$
   };
  \end{tikzpicture}\rlap{ .}}}
\end{equation*}

Given a $\V$-icon $\gamma \colon F \Rightarrow G$, let $\gamma_x$ be the identity morphism in $\C_0(Fx,Gx)$ for each $x \in \B$; now the morphisms $\C(\gamma_x,1)$ and $\C(1,\gamma_x)$ are isomorphic to identities in $\V$, and so there is a bijection between $2$-cells $\C(1,\gamma_y)\circ F_{xy} \Rightarrow \C(\gamma_x,1) \circ G_{xy}$ and $2$-cells $F_{xy} \Rightarrow G_{xy}$. Under this correspondence, the two $\V$-transformation axioms corresponds to the two $\V$-icon axioms; and thus we have:

\begin{Prop}\label{thm:icons}
  There is a bijection between $\V$-icons $\gamma\colon  F\Rightarrow G$ and oplax $\V$-transformations $\gamma\colon  F\Rightarrow G$ whose $1$-cell components are identities, under which invertible icons correspond to pseudo $\V$-transformations.
\end{Prop}

\subsection{$\V$-modifications}
\label{subsec:vmodif}

If $\gamma, \delta \colon F \Rightarrow G$ are $\V$-transformations, then a \emph{$\V$-modification} $\Gamma \colon \gamma \Rrightarrow \delta$ comprises $2$-cells $\Gamma_x \colon \gamma_x \Rightarrow \delta_x$ in $\C_0$ for each $x \in \B$,
subject to the axiom that for all $x,y \in \B$, we have $(\C(1,\Gamma_y) \otimes 1) \circ \bar \gamma_{xy} = \bar \delta_{xy} \circ (\C(\Gamma_x,1) \otimes 1)$.

Any $\V$-modification $\Gamma \colon \alpha \Rrightarrow \beta$ has an underlying ordinary modification $\Gamma_0 \colon \alpha_0 \Rrightarrow \beta_0$ with the same components as $\Gamma$.

\section{The tricategory of $\V$-bicategories}\label{sec:tricat-vbicats}
\subsection{Local structure}
For any $\V$-bicategories $\B$ and $\C$, we may define an (ordinary) bicategory $\V\text-\cat{Bicat}(\B, \C)$ whose objects, $1$-cells and $2$-cells are $\V$-functors, $\V$-transformations and $\V$-modifications from $\B$ to $\C$. 
 The identity $\V$-transformation $1_F \colon F \Rightarrow F$ has $1$-cell components $(1_F)_x = 1_{Fx}$ and $2$-cell component $(\overline{1_F})_{xy}$ coming from the nullary functoriality constraints of $\C(Fx,\thg)$ and of $\C(\thg, Fy)$. 
The composite of $\V$-transformations $\gamma \colon F \Rightarrow G$ and $\delta \colon G \Rightarrow H$ has $1$-cell components $(\delta\gamma)_x = \delta_x \circ \gamma_x$,
and $2$-cell component $(\overline{\delta\gamma})_{xy}$ built from $\bar \delta_{xy}$ and $\bar \gamma_{xy}$ together with the the binary functoriality constraints of $\C(\thg, Hy)$ and $\C(Fx,\thg)$.
The remaining data for this bicategory are obtained componentwise from the corresponding data in $\C_0$. 

\subsection{Composition of $\V$-functors}
Given $\V$-functors $F \colon \A \to \B$ and $G \colon \B \to \C$, we define the composite $G F \colon \A \to \C$ to have action on objects $(G F)(X) = G(FX)$, action on homs $(GF)_{xy} = G_{Fx,Fy} \circ F_{xy}$, and functoriality coherence $2$-cells \begin{equation*}
\vcenter{\hbox{\begin{tikzpicture}[y=0.80pt, x=0.7pt,yscale=-1, inner sep=0pt, outer sep=0pt, every text node part/.style={font=\tiny} ]
  \path[draw=black,line join=miter,line cap=butt,line width=0.650pt]
  (10.0000,1002.3622) --
  node[above right=0.12cm,at start]{$j$}(40.0000,1002.3622);
  \path[draw=black,line join=miter,line cap=butt,line width=0.650pt]
  (41.2195,1000.7362) .. controls (48.9837,990.7362) and (70.7932,992.3622) ..
  node[above left=0.07cm,pos=0.55]{$j$}(78.5168,992.3622);
  \path[draw=black,line join=miter,line cap=butt,line width=0.650pt]
  (80.7724,990.8581) .. controls (85.5693,982.5496) and (114.0917,982.3622) ..
  node[above left=0.12cm,at end]{$j$}(120.0000,982.3622);
  \path[draw=black,line join=miter,line cap=butt,line width=0.650pt]
  (80.5691,993.7849) .. controls (85.2966,1001.9732) and (114.9187,1002.3622) ..
  node[above left=0.12cm,at end]{$F$}(120.0000,1002.3622);
  \path[draw=black,line join=miter,line cap=butt,line width=0.650pt]
  (41.1382,1003.6630) .. controls (44.6705,1017.7236) and (87.0321,1022.3622) ..
  node[above left=0.12cm,at end]{$G$}(120.0000,1022.3622);
  \path[fill=black] (40.650406,1002.7687) node[circle, draw, line width=0.65pt, minimum width=5mm, fill=white, inner sep=0.25mm] (text4799) {$\fu$};
  \path[fill=black] (79.10569,992.19958) node[circle, draw, line width=0.65pt, minimum width=5mm, fill=white, inner sep=0.25mm] (text4803) {$\fu$};
\end{tikzpicture}}}
\qquad \text{and} \qquad
\vcenter{\hbox{\begin{tikzpicture}[y=0.80pt, x=0.8pt,yscale=-1, inner sep=0pt, outer sep=0pt, every text node part/.style={font=\tiny} ]
  \path[draw=black,line join=miter,line cap=butt,line width=0.650pt]
  (-13.0000,1002.3622) .. controls (47.3532,1001.8099) and (102.4701,989.3219) ..
  node[above right=0.12cm,at start]{\!\!$(G F)\t (G F)$}
  node[above=0.07cm,pos=0.55]{$F\t F$}(113.7516,1001.2125);
  \path[draw=black,line join=miter,line cap=butt,line width=0.650pt]
  (-13.0000,1022.3622) .. controls (40.0000,1022.3622) and (66.1195,1032.9371) ..
  node[above right=0.12cm,at start]{$m$}(78.7670,1023.4515);
  \path[draw=black,line join=miter,line cap=butt,line width=0.650pt]
  (51.0209,999.4729) .. controls (55.0451,1013.2701) and (71.1498,1010.9249) ..
  node[below left=0.06cm,pos=0.5]{$G \t G$}(78.6233,1020.9854);
  \path[draw=black,line join=miter,line cap=butt,line width=0.650pt]
  (81.3460,1021.1747) .. controls (85.2265,1014.7073) and (109.3111,1011.5455) ..
  node[below right=0.07cm,pos=0.45]{$m$}(114.0539,1003.4971);
  \path[draw=black,line join=miter,line cap=butt,line width=0.650pt]
  (115.4911,1001.4850) .. controls (122.5335,994.2990) and (130.0000,992.3622) ..
  node[above left=0.12cm,at end]{$m$}(150.0000,992.3622);
  \path[draw=black,line join=miter,line cap=butt,line width=0.650pt]
  (115.6884,1003.5120) .. controls (118.3307,1008.5933) and (130.0000,1012.3622) ..
  node[above left=0.12cm,at end]{$F$}(150.0000,1012.3622);
  \path[draw=black,line join=miter,line cap=butt,line width=0.650pt]
  (81.0730,1023.3799) .. controls (85.1381,1029.0710) and (116.9919,1032.3622) ..
  node[above left=0.12cm,at end]{$G$}(150.0000,1032.3622);
  \path[shift={(-25.521879,68.238085)},draw=black,fill=black] (77.8600,931.1000)arc(0.000:180.000:1.250)arc(-180.000:0.000:1.250) -- cycle;
  \path[fill=black] (79.878052,1022.2809) node[circle, draw, line width=0.65pt, minimum width=5mm, fill=white, inner sep=0.25mm] (text5170) {$\fm$
     };
  \path[fill=black] (114.26829,1002.3622) node[circle, draw, line width=0.65pt, minimum width=5mm, fill=white, inner sep=0.25mm] (text5174) {$\fm$
   };
\end{tikzpicture}\ \ \rlap{ .}}} 
\end{equation*}

Now the assignation $G \mapsto G F$ provides the action on objects of a ``whiskering'' functor {$(\thg) F \colon \V\text-\cat{Bicat}(\B,\C) \to \V\text-\cat{Bicat}(\A,\C)$} that on morphisms, sends $\gamma \colon G \Rightarrow H$ to the $\V$-transformation with $1$-cell components $(\gamma F)_x = \gamma_{Fx}$ and $2$-cell components
  obtained by whiskering those of $\gamma$ with $F$ and applying associativity constraints; and on $2$-cells,
sends  $\Gamma \colon \gamma \Rrightarrow \delta$ to the $\V$-modification with components $(\Gamma F)_x = \Gamma_{Fx}$. It is easy to see that $(\thg) F$ in fact \emph{strictly} preserves identities and composition.

On the other hand, the assignation $F \mapsto G  F$ is the action on objects of a  functor {$G (\thg) \colon \V\text-\cat{Bicat}(\A,\B) \to \V\text-\cat{Bicat}(\A,\C)$} that on $1$-cells sends $\gamma \colon F \Rightarrow K$ to the $\V$-transformation with $1$-cell components $(G \gamma)_x = G(\gamma_x)$ and $2$-cell components
\begin{equation*}
\vcenter{\hbox{
\begin{tikzpicture}[y=0.80pt, x=0.8pt,yscale=-1, inner sep=0pt, outer sep=0pt, every text node part/.style={font=\tiny} ]
  \path[draw=black,line join=miter,line cap=butt,line width=0.650pt]
  (-15.0000,962.3622) .. controls (35.8287,962.3622) and (83.6718,964.6070) ..
  node[above right=0.12cm,at start] {\!$K$}(88.7225,970.7801);
  \path[draw=black,line join=miter,line cap=butt,line width=0.650pt]
  (-15.0000,982.3622) .. controls (30.2050,982.3622) and (46.0758,987.0349) ..
  node[above right=0.12cm,at start] {\!$G$}(48.3205,990.9633);
  \path[draw=black,line join=miter,line cap=butt,line width=0.650pt]
  (-15.0000,1002.3622) .. controls (29.9402,1002.3622) and (45.7030,997.4882) ..
  node[above right=0.12cm,at start] {\!$\C(G\gamma,1)$}(48.2284,994.1210);
  \path[draw=black,line join=miter,line cap=butt,line width=0.650pt]
  (53.4418,994.3652) .. controls (57.6507,998.8547) and (125.1523,998.8547) ..
  node[above=0.07cm] {$G$}(128.5195,993.8040);
  \path[draw=black,line join=miter,line cap=butt,line width=0.650pt]
  (53.4418,990.7801) .. controls (56.8089,984.6070) and (83.1253,979.1353) ..
  node[above left=0.07cm,pos=0.57] {$\B(\gamma,1)$}(88.4567,974.0845);
  \path[draw=black,line join=miter,line cap=butt,line width=0.650pt]
  (91.8627,970.6398) .. controls (95.7910,965.4488) and (152.5839,962.3622) ..
  node[above left=0.12cm,at end] {$F$\!}(186.0000,962.3622);
  \path[draw=black,line join=miter,line cap=butt,line width=0.650pt]
  (91.8627,974.2248) .. controls (95.2298,979.1353) and (126.3522,986.8518) ..
  node[above right=0.07cm,pos=0.44] {$\B(1,\gamma)$}(128.4567,990.9204);
  \path[draw=black,line join=miter,line cap=butt,line width=0.650pt]
  (131.6004,990.6306) .. controls (136.6512,985.0187) and (145.0546,981.6589) ..
  node[above left=0.12cm,at end] {$G$\!}(186.0000,982.3622);
  \path[draw=black,line join=miter,line cap=butt,line width=0.650pt]
  (131.6004,993.9978) .. controls (136.0900,998.7680) and (144.2540,1002.3622) ..
  node[above left=0.12cm,at end] {$\C(1, G\gamma)$\!}(186.0000,1002.3622);
  \path[fill=black] (49.603172,993.23523) node[circle, draw, line width=0.65pt, minimum width=5mm, fill=white, inner sep=0.25mm] (text4238) {$\fm$     };
  \path[fill=black] (89.682541,972.20343) node[circle, draw, line width=0.65pt, minimum width=5mm, fill=white, inner sep=0.25mm] (text4242) {$\bar \gamma$     };
  \path[fill=black] (129.36508,992.83838) node[circle, draw, line width=0.65pt, minimum width=5mm, fill=white, inner sep=0.25mm] (text4246) {$\fm^{-1}$   };
\end{tikzpicture}\rlap{\quad;}}}
\end{equation*}
and on $2$-cells sends $\Gamma \colon \gamma \Rrightarrow \delta$ to the $\V$-modification with components $(G \circ \Gamma)_x = G(\Gamma_x)$. The functoriality constraint $2$-cells of $G (\thg)$ are obtained pointwise from those of $G_0 \colon \B_0 \to \C_0$. The operations just described give rise to a functor
\[\ast \colon \V\text-\cat{Bicat}(\B, \C) \times \V\text-\cat{Bicat}(\A,\B) \to \V\text-\cat{Bicat}(\A,\C)\rlap{ ,}\]
whose action on $\V$-transformations $\gamma \colon F \Rightarrow G \colon \A \to \B$ and $\delta \colon H \Rightarrow K \colon \B \to \C$ is given by $\delta \ast \gamma \defeq (K \gamma) \circ (\delta F) \colon HF \Rightarrow KG$, and  on $\V$-modifications $\Gamma \colon \gamma \Rightarrow \delta$ and $\Delta \colon \epsilon \Rightarrow \nu$ by $\Delta \ast \Gamma \defeq (K \Gamma) \circ (\Delta F)$. The functoriality constraints of $\ast$ are obtained from those of the left and right whiskering functors together with 
the coherent interchange $\V$-modifications $(K \gamma) \circ (\delta F) \Rrightarrow (\delta G) \circ (H \gamma) \colon HF \Rightarrow KG$ whose component at $x$ is given by whiskering $\bar \delta_{Fx,Gx}$ with $\gamma_x \colon I \to \B(Fx,Gx)$. 

For any $\V$-bicategory $\B$, there is an identity $\V$-functor $1_\B \colon \B \to \B$ that is the identity on objects, has action on homs $(1_\B)_{xy} = 1_{\B(x,y)}$, and functoriality $2$-cells obtained from the unitality constraints of the bicategory $\V$.

\subsection{The tricategory $\V\text-\cat{Bicat}$}\label{sec:tricat-vbicat}
In order to make out of the structures defined above a tricategory
$\V\text-\cat{Bicat}$, it remains only to describe the coherent constraint cells associated to composition along $0$-cell boundaries. Given $\V$-functors $F \colon \A \to \B$, $G \colon \B \to \C$ and $H \colon \C \to \D$, we observe that the composites $(HG)F$ and $(HG)F$ agree on objects, and differ on hom-objects only up to associativity constraints in the bicategory $\V$;  these constraints form the components of an invertible $\V$-icon $(HG)F \Rightarrow H(GF)$. The $\V$-transformations corresponding to these $\V$-icons under Proposition~\ref{thm:icons} now provide 
the $1$-cell components of the associativity pseudonatural equivalence
$\mathfrak a_{\A\B\C\D}$ for $\V\text-\cat{Bicat}$; the corresponding $2$-cell components are constructed from unit constraints in $\D_0$. Similarly, for any $\V$-functor $F \colon \C \to \D$, we have invertible $\V$-icons $F1_\C \Rightarrow F$ and $1_\D F \Rightarrow F$ giving rise to the unit constraints $\mathfrak l$ and $\mathfrak r$ for $\V\text-\cat{Bicat}$.

The remaining data needed to make $\V\text-\cat{Bicat}$ into a tricategory are the invertible modifications $\pi$, $\lambda$, $\rho$ and $\nu$ witnessing the satisfaction up to coherent isomorphism of the pentagon axiom and three unit axioms. The components of these modifications are $\V$-modifications between $\V$-transformations built from 
$\mathfrak a$, $\mathfrak l$ and $\mathfrak r$; since these $\V$-transformations have identity $1$-cell components, the components of the required $\V$-modifications are built from coherence constraints in the monoidal bicategory $\V$. The coherence theorem for monoidal bicategories now ensures that the $\pi$, $\lambda$, $\rho$ and $\nu$ so defined satisfy the tricategory axioms.

\begin{Rk}\label{rk:tricat-vbicat}
The tricategory structure of $\V\text-\cat{Bicat}$ can be derived from a different kind of structure that captures more faithfully the constraint data involved in its compositions. Observe first that, since the associativity and unit constraints for composition of $\V$-functors are given by invertible $\V$-icons, we have a bicategory $\V\text-\cat{Bicat}_2$
 of $\V$-bicategories, $\V$-functors and $\V$-icons (note that this bicategory will not be a $2$-category unless $\V$ is too). $\V\text-\cat{Bicat}_2$ captures the up-to-isomorphism associativity of $\V$-functor composition, but does not contain the general $\V$-transformations or the $\V$-modifications. However, we may retain both of these by combining the structures of $\V\text-\cat{Bicat}$ and $\V\text-\cat{Bicat}_2$ into a \emph{locally cubical bicategory}  in the sense of~\cite{gg:ldstr-tricat}. A locally cubical bicategory is just a $\cat{DblCat}$-bicategory, where $\cat{DblCat}$ is the monoidal 2-category of pseudo double categories, pseudo double functors, and vertical transformations in the sense of~\cite{Grandis1999Limits}; in our case, we obtain such a structure $\V\text-\underline{\cat{Bicat}}$ whose hom-pseudo double category from $\B$ to $\C$ has $\V\text-\cat{Bicat}(\B,\C)$ as its horizontal bicategory, $\V$-icons as its vertical maps, and as cells, \emph{cubical modifications}
 \begin{equation*}
  \vcenter{\xymatrix{
      F\ar@{=>}[r]|-@{|}^\alpha \ar@{=>}[d]_\gamma \dthreecell{dr}{\Gamma} &
      G\ar@{=>}[d]^\delta\\
      H\ar@{=>}[r]|-@{|}_\beta &
      K
      }}
\end{equation*}
that consist of $2$-cells $\Gamma_x\colon \alpha_x \Rightarrow \beta_x$ in $\C_0$ such that for all $x,y\in\B$ we have

\begin{equation*}
\vcenter{\hbox{
\begin{tikzpicture}[y=1pt, x=0.85pt,yscale=-1, inner sep=0pt, outer sep=0pt, every text node part/.style={font=\tiny} ]
  \path[draw=black,line join=miter,line cap=butt,line width=0.650pt] (124.7857,972.3622) .. controls (110.0000,972.3622) and (100.0000,972.3622) .. 
  node[above=0.1cm,pos=0.35] {$F$}(91.6415,980.5697);
  \path[draw=black,line join=miter,line cap=butt,line width=0.650pt] (125.5254,992.3622) .. controls (110.0000,992.3622) and (100.0000,992.3622) .. 
  node[above=0.1cm,pos=0.22] {$\C(1, \alpha)$}(91.6415,984.2810);
  \path[draw=black,line join=miter,line cap=butt,line width=0.650pt] (88.0095,980.7676) .. controls (80.0000,972.3622) and (69.7239,972.3622) .. 
  node[above right=0.12cm,at end] {\!$G$}(40.0000,972.3622);
  \path[draw=black,line join=miter,line cap=butt,line width=0.650pt] (40.0000,992.3622) .. controls (70.0000,992.3622) and (80.0000,992.3622) .. 
  node[above right=0.12cm,at start] {\!$\C(\alpha,1)$}(88.1358,983.9244);
  \path[draw=black,line join=miter,line cap=butt,line width=0.650pt] (180.0000,992.3622) -- 
  node[above left=0.12cm,at start] {$\C(1, \beta)$\!\!}(124.5013,992.3622);
  \path[draw=black,line join=miter,line cap=butt,line width=0.650pt] (180.0000,972.3622) -- 
  node[above left=0.12cm,at start] {$H$\!}(124.6429,972.3622);
  \path[fill=black] (140.53553,992.36218) node[circle, draw, line width=0.65pt, minimum width=5mm, fill=white, inner sep=0.25mm] (text6017-9) {$\Gamma$     };
  \path[fill=black] (89.651039,982.15656) node[circle, draw, line width=0.65pt, minimum width=5mm, fill=white, inner sep=0.25mm] (text6017) {$\bar \alpha$     };
  \path[fill=black] (140.21429,972.18359) node[circle, draw, line width=0.65pt, minimum width=5mm, fill=white, inner sep=0.25mm] (text6017-9-1) {$\bar \gamma$   };
  \end{tikzpicture}}}
  \quad=\quad
\vcenter{\hbox{
\begin{tikzpicture}[y=1pt, x=0.85pt,yscale=-1, xscale=-1,inner sep=0pt, outer sep=0pt, every text node part/.style={font=\tiny} ]
  \path[draw=black,line join=miter,line cap=butt,line width=0.650pt] (124.7857,972.3622) .. controls (110.0000,972.3622) and (100.0000,972.3622) .. 
  node[above=0.1cm,pos=0.35] {$K$}(91.6415,980.5697);
  \path[draw=black,line join=miter,line cap=butt,line width=0.650pt] (125.5254,992.3622) .. controls (110.0000,992.3622) and (100.0000,992.3622) .. 
  node[above=0.1cm,pos=0.22] {$\C(\beta, 1)$}(91.6415,984.2810);
  \path[draw=black,line join=miter,line cap=butt,line width=0.650pt] (88.0095,980.7676) .. controls (80.0000,972.3622) and (69.7239,972.3622) .. 
  node[above left=0.12cm,at end] {$H$\!}(40.0000,972.3622);
  \path[draw=black,line join=miter,line cap=butt,line width=0.650pt] (40.0000,992.3622) .. controls (70.0000,992.3622) and (80.0000,992.3622) .. 
  node[above left=0.12cm,at start] {$\C(1, \beta)$\!}(88.1358,983.9244);
  \path[draw=black,line join=miter,line cap=butt,line width=0.650pt] (180.0000,992.3622) -- 
  node[above right=0.12cm,at start] {\!$\C(\alpha,1)$}(124.5013,992.3622);
  \path[draw=black,line join=miter,line cap=butt,line width=0.650pt] (180.0000,972.3622) -- 
  node[above right=0.12cm,at start] {\!$G$}(124.6429,972.3622);
  \path[fill=black] (140.53553,992.36218) node[circle, draw, line width=0.65pt, minimum width=5mm, fill=white, inner sep=0.25mm] (text6017-9) {$\Gamma$     };
  \path[fill=black] (89.651039,982.15656) node[circle, draw, line width=0.65pt, minimum width=5mm, fill=white, inner sep=0.25mm] (text6017) {$\bar \beta$     };
  \path[fill=black] (140.21429,972.18359) node[circle, draw, line width=0.65pt, minimum width=5mm, fill=white, inner sep=0.25mm] (text6017-9-1) {$\bar \delta$   };
  \end{tikzpicture}}}
\end{equation*}
in $\V(\B(x,y),\C(Fx,Gy))$. (In these two diagrams, we write $\Gamma$ as an abbreviation for the $2$-cells $\C(\Gamma, 1)$ and $\C(1, \Gamma)$.)

For any invertible $\V$-icon $\gamma\colon F\Rightarrow G$, we can find---using Proposition~\ref{thm:icons}---a $\V$-transformation $\hat{\gamma} \colon  F\Rightarrow G$ and an invertible cubical modification
\begin{equation*}
  \vcenter{\xymatrix{
      F\ar@{=>}[r]|-@{|}^{\hat\gamma} \ar@{=>}[d]_\gamma \dthreecell{dr}{\epsilon} &
      G\ar@{=>}[d]^{1_G}\\
      G\ar@{=>}[r]|-@{|}_{1_G} &
      G\rlap{ .}
      }}
\end{equation*}
This means that the locally cubical bicategory $\V\text-\underline{\cat{Bicat}}$ is \emph{locally fibrant} in the sense of~\cite{gg:ldstr-tricat}; whence, by Theorem 24 of that paper, its ``locally horizontal'' part can be equipped with the structure of a tricategory. This is precisely the tricategory structure on $\V\text-\cat{Bicat}$ described above.
\end{Rk}

\subsection{The trihomomorphism $(\thg)_0$}\label{subsec:ordinary-trifunctor}
Taking underlying ordinary structures now yields a trifunctor $(\thg)_0 \colon \V\text-\cat{Bicat} \to \cat{Bicat}$; it is in fact a \emph{strict} functor on hom-bicategories $\V\text-\cat{Bicat}(\B,\C) \to \cat{Bicat}(\B_0, \C_0)$, and preserves binary and nullary composition of $0$-cells up to an invertible icon; in fact, we may see $(\thg)_0$ as having been induced from a morphism of locally cubical bicategories $\V\text-\underline{\cat{Bicat}} \to \underline{\cat{Bicat}}$. We will see in Example~\ref{ex:underlying-change-of-base} below that this trifunctor is an instance of \emph{change of base}, here along the monoidal functor $\V(I, \thg) \colon \V \to \cat{Cat}$. Alternatively, we may obtain it as the hom-functor $\V\text-\cat{Bicat}(\I, \thg)$, where $\I$ is the unit $\V$-bicategory of Example~\ref{ex:unit-vcat}.

Finally in this section, we note the following result, which says that the trihomomorphism $(\thg)_0$ is ``locally conservative''.
\begin{Prop}\label{prop:equivifcomponentsequiv}
A $\V$-transformation $\gamma \colon F \Rightarrow G \colon \B \to \C$ is an equivalence in $\V\text-\cat{Bicat}(\B,\C)$ if and only if each of its components is an equivalence in $\C_0$.
\end{Prop}
\begin{proof}
The ``only if'' direction is immediate, since composition in $\V\text-\cat{Bicat}(\B,\C)$ is lifted componentwise from $\C_0$. Conversely, if each $\gamma_x$ has equivalence pseudoinverse $\gamma^\centerdot_x$ in $\C_0$, then we obtain 
a $\V$-transformation $\gamma^\centerdot \colon G \Rightarrow F$ with $1$-cell components $\gamma^\centerdot_x$ and $2$-cell components
\begin{equation*}
\begin{tikzpicture}[y=0.8pt, x=1pt,yscale=-1, inner sep=0pt, outer sep=0pt, every text node part/.style={font=\tiny} ]
  \path[draw=black,line join=miter,line cap=butt,line width=0.650pt] 
  (151.8683,992.3622) -- 
  node[above left=0.12cm,at end] {$G$\!}(216.0000,992.3622);
  \path[draw=black,line join=miter,line cap=butt,line width=0.650pt] 
  (84.0000,992.3622) -- 
  node[above right=0.12cm,at start] {\!$F$}(148.1317,992.3622);
  \path[draw=black,line join=miter,line cap=butt,line width=0.650pt] 
  (147.9621,996.4002) .. controls (143.0915,996.4002) and (113.8822,997.9754) .. 
  (113.8822,1008.5811) .. controls (113.8822,1013.9475) and (125.8353,1017.2287) .. 
  node[above right=0.01cm,pos=0.35,rotate=-10] {$\C(1, \gamma^\centerdot)$}(142.0411,1019.2331)
  (157.2427,1020.6905) .. controls (178.1546,1022.2246) and (201.7381,1022.3622) .. 
  node[above left=0.12cm,at end] {$\C(1, \gamma^\centerdot)$\!\!}(216.0000,1022.3622)
  (151.7956,996.4523) .. controls (156.8860,996.4523) and (185.8755,997.8461) .. 
  (185.8755,1008.6332) .. controls (185.8755,1021.6211) and (113.8298,1022.3622) .. 
  node[above right=0.12cm,at end] {\!\!$\C(\gamma^\centerdot,1)$}(84.0000,1022.3622);
  \path[fill=black] (149.55309,993.92542) node[circle, draw, line width=0.65pt, minimum width=5mm, fill=white, inner sep=0.25mm] (text3650) {$\bar \gamma^{-1}$   };
\end{tikzpicture}\rlap{ .}
\end{equation*}
The $2$-cells witnessing each $\gamma^\centerdot_x$ as pseudoinverse to $\gamma_x$ now form the components of invertible $\V$-modifications witnessing $\gamma^\centerdot$ as pseudoinverse to $\gamma$.
\end{proof}

\section{Modules}\label{sec:modules}
In this section, we define left, right, and two-sided modules between $\V$-bicategories. In the presence of suitable extra structure on $\V$, we could define left $\A$-modules as $\V$-functors $\A \to \V$, right $\B$-modules as $\V$-functors $\B^\op \to \V$, and $\A$-$\B$-modules as $\V$-functors $\B^\op \otimes \A \to \V$; the extra structure required on $\V$ would be right closedness (in order to view $\V$ itself as a $\V$-bicategory) and some form of symmetry (to form the opposite $\V$-bicategory of $\B$, and the tensor product $\B^\op \otimes \A$). However, we prefer to introduce modules as a separate notion, rather than reducing them to $\V$-functors; this allows us to present them in a manner that both is simpler to work with, and also valid even in the absence of the above-mentioned extra structure on $\V$; in particular, we need not require any form of symmetry.

Throughout this section, $\A$ and $\B$ will be $\V$-bicategories.
\subsection{Modules and bimodules}
A \emph{right $\B$-module} $W$, also written as $W \colon \bullet \tor \B$, is given by the following data:
\begin{itemize}
\item For each $x \in \B$, an object $W x$ of $\V$;
\item For each $x, y \in \B$, a morphism $m_{xy} \colon W y \otimes \B(x,y) \to W x$ of $\V$;
\item For each $x,y,z \in \B$, invertible $2$-cells
\begin{equation*}
\cd[@!C@C-5.2em]{
 & W  x \otimes \B(x, x) \ar[dr]^{m} \dtwocell[0.55]{d}{\ru_{x}} \\
 W  x \otimes I  \ar[rr]_{\mathfrak r} \ar[ur]^{1 \otimes j} & &
 W  x
}\  \text{and} \ \ \
\cd[@C-2.5em]{
 (W  z \otimes \B(y, z)) \otimes \B(x, y) \ar[rr]^{\mathfrak a} \ar[d]_{m \otimes 1}  & {} \rtwocell[0.55]{d}{\ass_{xyz}} &
 W  z \otimes (\B(y, z) \otimes \B(x, y)) \ar[d]^{1 \otimes m} \\
 W  y \otimes \B(x, y) \ar[r]_-m & W  x & W  z \otimes \B(x, z) \ar[l]^-{m}
}
\end{equation*}
\end{itemize}
subject to two axioms that as string diagrams coincide with the axioms for a $\V$-bicategory, but now interpreted in the respective categories $\V(W y \otimes \B(x,y), W x)$ and $\V(\,((W z \otimes \B(y,z)) \otimes \B(x,y)) \otimes \B(v,x),\, W v\,)$. 

Dually, a \emph{left $\A$-module} $W$ (written $W \colon \A \tor \bullet$) is given by objects $Wx$ in $\V$ for each $x \in \A$, morphisms $m_{xy} \colon \A(x,y) \otimes W x \to W y$ for each $x,y \in \A$, and for each $x,y,z \in \A$, invertible $2$-cells
\begin{equation*}
\cd[@!C@C-5.2em]{
 & \A(x, x) \otimes W  x \ar[dr]^{m} \dtwocell[0.55]{d}{\lu_{x}} \\
 I \otimes W  x  \ar[rr]_{\mathfrak l} \ar[ur]^{j \otimes 1} & &
 W  x
}\ \text{and} \ \ \ \cd[@C-2.5em]{
 (\A(y, z) \otimes \A(x, y)) \otimes W  x \ar[rr]^{\mathfrak a} \ar[d]_{m \otimes 1}  & {} \rtwocell[0.55]{d}{\ass_{xyz}} &
 \A(y, z) \otimes (\A(x, y) \otimes W  x) \ar[d]^{1 \otimes m} \\
 \A(x, z) \otimes W  x \ar[r]_-m & W  z & \A(y, z) \otimes W  y \ar[l]^-{\ \, m}
}
\end{equation*}
subject to two axioms that as string diagrams again coincide with the two axioms for a $\V$-bicategory, interpreted in appropriate hom-categories.

Finally, an \emph{$\A$-$\B$-bimodule} $M \colon \A \tor \B$ is given by objects $M(b,a) \in \V$ for each $a \in \A$ and $b \in \B$, together with left $\A$-module structures on each $M(b, \thg)$, right $\B$-module structures on each $M (\thg, a)$, and for each $a,a' \in \A$ and $b,b' \in \B$, invertible $2$-cells
\begin{equation*}
\cd{
 (\A(a, a') \otimes M (b, a)) \otimes \B(b', b) \ar[rr]^{\mathfrak a} \ar[d]_{m \otimes 1}  & {} \rtwocell[0.55]{d}{\ass_{aa'bb'}} &
 \A(a, a') \otimes (M  (b, a) \otimes \B(b', b)) \ar[d]^{1 \otimes m} \\
 M (b, a') \otimes \B(b', b) \ar[r]_-m & M (b', a') & \A(a, a') \otimes M (b', a) \ar[l]^-{m}
}
\end{equation*}
satisfying two axioms that as string diagrams both coincide with the second $\V$-bicategory axiom, but interpreted in appropriate hom-categories.

In what follows, we will give a number of constructions and definitions that operate equally well on left modules, right modules and bimodules. We shall typically describe those constructions in the most involved case, that of bimodules, and leave it to the reader to derive the corresponding one-sided versions in which an indexing $\V$-bicategory has been replaced by $\bullet$. On occasions, such replacement will lead us to the consideration of ``bimodules'' $\bullet \tor \bullet$; by definition, these are simply objects of $\V$, with the morphisms and transformations between them being the $1$- and $2$-cells of $\V$.

\begin{Ex}\label{ex:hom-module}
For each $\V$-bicategory $\C$, there is a bimodule $\C \colon \C \tor \C$, the \emph{hom-module} of $\C$, whose components are the hom-objects $\C(b,c)$, whose left and right action morphisms are given by composition in $\C$, and whose coherence $2$-cells are unit and associativity constraints for $\C$.
\end{Ex}
\begin{Ex}\label{ex:reindexing}
Given $F \colon \B \to \C$ and $M \colon \A \tor \C$, there is an $\A$-$\B$-bimodule $M(F,1)$ whose component at $(b,a)$ is $M(Fb,a)$, whose left actions are those of $M$, and whose right actions are given by
\begin{equation*}
M(Fb, a) \otimes \B(b',b) \xrightarrow{1 \otimes F} M(Fb, a) \otimes \C(Fb', Fb) \xrightarrow m M(Fb',a)
\end{equation*}
with unit and associativity constraint $2$-cells given by\begin{equation*}
\vcenter
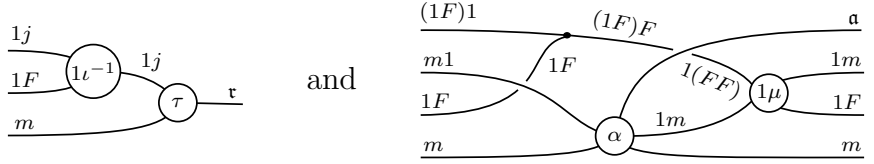
{\hbox{\begin{tikzpicture}[y=0.80pt, x=0.8pt,yscale=-1, inner sep=0pt, outer sep=0pt, every text node part/.style={font=\tiny} ]
  \path[draw=black,line join=miter,line cap=butt,line width=0.650pt]
  (0.0000,1022.3622) .. controls (20.0000,1022.3622) and (33.1989,1021.0637) ..
  node[above right=0.12cm,at start] {$\!1 \t F$}(38.5075,1013.9418);
  \path[draw=black,line join=miter,line cap=butt,line width=0.650pt]
  (0.0000,1002.3622) .. controls (20.0000,1002.3622) and (33.9263,1005.0313) ..
  node[above right=0.12cm,at start] {$\!1 \t j$}(38.7575,1010.8373);
  \path[draw=black,line join=miter,line cap=butt,line width=0.650pt]
  (79.3505,1025.8977) .. controls (72.5970,1018.8876) and (63.0000,1011.3622) ..
  node[above right=0.09cm,pos=0.65] {$1 \t j$}(40.5000,1012.3622);
  \path[draw=black,line join=miter,line cap=butt,line width=0.650pt]
  (0.0000,1042.3622) .. controls (61.7748,1042.3622) and (73.2500,1037.3548) ..
  node[above right=0.12cm,at start] {$m$}(79.5000,1028.1048);
  \path[draw=black,line join=miter,line cap=butt,line width=0.650pt]
  (80.2566,1026.9063) --
  node[above left=0.12cm,at end] {$\mathfrak r$}(109.9551,1027.6134);
  \path[fill=black] (40,1012.1122) node[circle, draw, line width=0.65pt, minimum width=5mm, fill=white, inner sep=0.25mm] (text3086) {$1 \t \fu^{-1}$     };
  \path[fill=black] (79.549515,1026.9062) node[circle, draw, line width=0.65pt, minimum width=5mm, fill=white, inner sep=0.25mm] (text3099) {$\ru$};
\end{tikzpicture}}}
\qquad
\text{and}
\qquad
\vcenter{\hbox{\begin{tikzpicture}[y=0.80pt, x=1pt,yscale=-1, inner sep=0pt, outer sep=0pt, every text node part/.style={font=\tiny} ]
  \path[draw=black,line join=miter,line cap=butt,line width=0.650pt]
  (166.0000,1022.3622) .. controls (146.0000,1022.3622) and (132.8011,1021.0637) ..
  node[above left=0.12cm,at start] {$1 \t F$\!}(127.4925,1013.9418);
  \path[draw=black,line join=miter,line cap=butt,line width=0.650pt]
  (166.0000,1002.3622) .. controls (146.0000,1002.3622) and (132.0737,1005.0313) ..
  node[above left=0.12cm,at start] {$1 \t m$\!}(127.2425,1010.8373);
  \path[draw=black,line join=miter,line cap=butt,line width=0.650pt]
  (74.1495,1032.3475) .. controls (89.4814,1032.2412) and (113.1281,1031.0184) ..
  node[above=0.1cm,pos=0.37] {$1 \t m$}(125.3750,1013.9872);
  \path[draw=black,line join=miter,line cap=butt,line width=0.650pt]
  (166.0000,1042.3622) .. controls (104.2252,1042.3622) and (79.5000,1043.1048) ..
  node[above left=0.12cm,at start] {$m$\!}(73.2500,1033.8548);
  \path[draw=black,line join=miter,line cap=butt,line width=0.650pt]
  (71.0792,1034.2274) .. controls (63.5792,1040.7274) and (50.0000,1042.3622) ..
  node[above right=0.12cm,at end] {$\!m$}(0.0000,1042.3622);
  \path[draw=black,line join=miter,line cap=butt,line width=0.650pt]
  (73.2353,1029.8600) .. controls (78.5230,1001.7179) and (89.5585,982.3622) ..
  node[above left=0.12cm,at end] {$\mathfrak a$\!}(165.0000,982.3622)
  (125.4878,1010.5393) .. controls (122.3177,1003.1709) and (113.3145,997.7190) ..
  node[below left=0.07cm,pos=0.22,rotate=-23] {$1 \t (F \t F)$}(101.4916,993.6875)
  (94.5053,991.5502) .. controls (56.9607,981.2329) and (0.0000,982.3622) ..
  node[above right=0.12cm,at end] {$\!\!(1 \t F) \t 1$}
  node[above right=0.12cm,pos=0.25,rotate=-12] {$\!\!(1 \t F) \t F$}(0.0000,982.3622);
  \path[draw=black,line join=miter,line cap=butt,line width=0.650pt]
  (55.3347,985.2122) .. controls (45.6642,987.3922) and (45.2719,997.3438) ..
  node[below right=0.08cm,pos=0.5] {$1 \t F$}(39.5209,1006.4000)
  (36.4306,1010.4871) .. controls (30.6218,1016.9162) and (20.5344,1022.0920) ..
  node[above right=0.14cm,at end] {$\!\!1 \t F$}(0.0000,1022.3622)
  (71.2812,1032.1551) .. controls (49.8912,1025.2608) and (54.0000,1001.8622) ..
  node[above right=0.12cm,at end] {$\!m \t 1$}(0.0000,1002.3622);
  \path[cm={{-1.0,0.0,0.0,1.0,(131.38279,53.645108)}},draw=black,fill=black] (77.8571,931.1122)arc(0.000:180.000:1.250)arc(-180.000:0.000:1.250) -- cycle;
  \path[fill=black] (72.605469,1032.2372) node[circle, draw, line width=0.65pt, minimum width=5mm, fill=white, inner sep=0.25mm] (text3099) {$\ass$
     };
  \path[fill=black] (130.34198,1012.0086) node[circle, draw, line width=0.65pt, minimum width=5mm, fill=white, inner sep=0.25mm] (text3086) {$1 \t \fm$
   };
\end{tikzpicture}\rlap{ .}}}
\end{equation*}
The bimodule constraint cells are given by
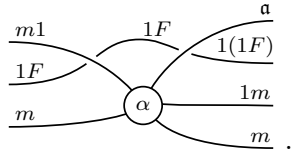
\begin{equation*}
\begin{tikzpicture}[y=0.80pt, x=0.9pt,yscale=-1, inner sep=0pt, outer sep=0pt, every text node part/.style={font=\tiny} ]
  \path[draw=black,line join=miter,line cap=butt,line width=0.650pt]
  (127.8878,970.8146) .. controls (139.8191,972.8759) and (148.6704,971.6046) ..
  node[above left=0.12cm,at end] {$1 \t m$\!}(180.0000,972.3622);
  \path[draw=black,line join=miter,line cap=butt,line width=0.650pt]
  (127.6786,973.3558) .. controls (131.2137,983.1147) and (146.4071,992.3622) ..
  node[above left=0.12cm,at end] {$m$\!}(180.0000,992.3622);
  \path[draw=black,line join=miter,line cap=butt,line width=0.650pt]
  (70.0000,982.3622) .. controls (94.7819,982.3622) and (118.6172,978.5843) ..
  node[above right=0.12cm,at start] {$m$\!}(125.9408,972.7759);
  \path[draw=black,line join=miter,line cap=butt,line width=0.650pt]
  (126.7723,969.6175) .. controls (132.3342,954.5869) and (152.8816,932.3622) ..
  node[above left=0.12cm,at end] {$\mathfrak a$\!}(180.0000,932.3622)
  (70.0000,962.3622) .. controls (87.3281,962.3622) and (95.9169,957.5941) ..
  node[above right=0.12cm,at start] {$1 \t F$\!}(102.3037,952.5793)
  (106.4006,949.1692) .. controls (111.4451,944.8565) and (115.9351,941.1519) ..
  (124.3477,941.1519) .. controls (131.0779,941.1519) and (134.8266,943.4620) ..
  node[above=0.08cm,pos=0.42] {$1 \t F$}(140.6994,945.9848)
  (145.9739,948.0176) .. controls (152.8718,950.3424) and (162.8225,952.3622) ..
  node[above left=0.12cm,at end] {$1 \t (1 \t F)$\!\!}(180.0000,952.3622)
  (70.0000,942.3622) .. controls (100.6008,942.3622) and (116.2826,958.0172) ..
  node[above right=0.12cm,at start] {$m \t 1$\!}(125.0912,969.9746);
  \path[fill=black] (125.89511,972.28424) node[circle, draw, line width=0.65pt, minimum width=5mm, fill=white, inner sep=0.25mm] (text3317) {$\ass$};
\end{tikzpicture}\rlap{ .}
\end{equation*}
Dualising the above, from each $G \colon \A \to \B$ and $M \colon \B \tor \C$, we obtain an $\A$-$\C$-bimodule $M(1, G)$ with $M(1, G)(c,a) = M(c, Ga)$, and remaining data obtained analogously to before. Combining the two constructions, we obtain from any $F \colon \B \to \D$, $G \colon \A \to \C$ and $M \colon \C \tor \D$, a bimodule $M(F,G) \colon \A \tor \B$.
\end{Ex}

\subsection{Module morphisms}
If $V, W \colon \bullet \tor \B$, then a \emph{right module morphism} $\phi \colon V \mathop{\rightarrow} W $ is given by:
\begin{itemize}
\item For each $x \in \B$, a morphism $\phi_x \colon Vx \to W x$ in $\V$;
\item For each $x, y \in \B$, invertible $2$-cells
\begin{equation*}
\cd[@C+2.5em]{
 Vy \otimes \B(x,y) \ar[r]^-{m} \ar[d]_{\phi \otimes 1} \rtwocell{dr}{\bar \phi_{xy}} &
 Vx \ar[d]^{\phi} \\
 W y \otimes \B(x,y) \ar[r]_-{m} & W x
}
\end{equation*}
\end{itemize}
subject to the axioms that for all $x \in \B$, we have
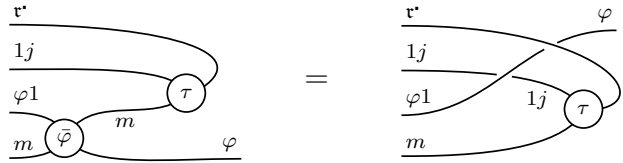
\begin{equation*}
\vcenter{\hbox{
\begin{tikzpicture}[y=0.80pt, x=0.73pt,yscale=-1, inner sep=0pt, outer sep=0pt, every text node part/.style={font=\tiny} ]
\path[use as bounding box] (15, 950) rectangle (135,1020);
  \path[draw=black,line join=miter,line cap=butt,line width=0.650pt]
  (15.0000,980.2193) .. controls (45.0000,980.2193) and (98.5424,978.1914) ..
  node[above right=0.15cm,at start] {$\!1 \t j$}(104.7179,990.3701);
  \path[draw=black,line join=miter,line cap=butt,line width=0.650pt]
  (46.0643,1009.1803) .. controls (60.0503,992.6170) and (88.4477,1005.2557) ..
  node[below=0.1cm,pos=0.5] {$m$}(104.3679,993.2345);
  \path[draw=black,line join=miter,line cap=butt,line width=0.650pt]
  (45.7179,1014.1587) .. controls (57.5251,1026.6048) and (97.6482,1022.3622) ..
  node[above left=0.15cm,at end] {$\phi\!$}(134.0000,1022.3622);
  \path[draw=black,line join=miter,line cap=butt,line width=0.650pt]
  (107.2143,992.4336) .. controls (147.7330,973.5876) and (98.7783,959.1212) ..
  node[above right=0.15cm,at end] {$\!\mathfrak r^\centerdot$}(15.0000,959.1327);
  \path[draw=black,line join=miter,line cap=butt,line width=0.650pt]
  (15.0000,1000.5765) .. controls (40.0294,1000.4472) and (42.1429,1009.5051) ..
  node[above right=0.15cm,at start] {$\!\phi \t 1$}(42.1429,1009.5051);
  \path[draw=black,line join=miter,line cap=butt,line width=0.650pt]
  (15.0000,1022.0086) .. controls (40.0000,1022.3657) and (42.8571,1013.7908) ..
  node[above right=0.15cm,at start] {$\!m$}(42.8571,1013.7908);
  \path[fill=black] (43.214283,1012.7192) node[circle, draw, line width=0.65pt, minimum width=5mm, fill=white, inner sep=0.25mm] (text4050) {$\bar \phi$   };
  \path[fill=black] (105.77222,991.40356) node[circle, draw, line width=0.65pt, minimum width=5mm, fill=white, inner sep=0.25mm] (text3449) {$\ru$     };
\end{tikzpicture}
}}
\quad \  =\quad\ \
\vcenter{\hbox{
\begin{tikzpicture}[y=0.80pt, x=0.73pt,yscale=-1, inner sep=0pt, outer sep=0pt, every text node part/.style={font=\tiny} ]
  \path[draw=black,line join=miter,line cap=butt,line width=0.650pt]
  (0.7143,1032.7193) .. controls (33.4402,1032.9653) and (50.7053,1013.5463) ..
  node[above right=0.15cm,at start] {$\!\phi \t 1$}(73.9303,1001.6762)
  (80.6246,998.5644) .. controls (89.4024,994.9073) and (99.2003,992.5432) ..
  node[above left=0.15cm,at end] {$\phi\!$}(111.0000,992.7193)
  (0.7143,1011.6479) .. controls (18.2407,1011.6479) and (35.1220,1012.0913) ..
  node[above right=0.15cm,at start] {$\!1 \t j$}(49.5738,1013.5828)
  (57.3368,1014.5222) .. controls (74.8965,1016.9929) and (87.9017,1021.4249) ..
  node[below left=0.08cm,pos=0.45] {$1 \t j$}(93.2143,1029.1479)
  (0.7143,990.5765) .. controls (102.7520,989.9005) and (134.2353,1027.0446) ..
  node[above right=0.15cm,at start] {$\!\mathfrak r^\centerdot$}(98.7857,1030.2194);
  \path[draw=black,line join=miter,line cap=butt,line width=0.650pt]
  (0.7143,1052.0051) .. controls (56.4286,1052.3622) and (84.4286,1044.8621) ..
  node[above right=0.15cm,at start] {$\!m$}(93.7143,1031.6479);
  \path[fill=black] (94.071426,1029.8622) node[circle, draw, line width=0.65pt, minimum width=5mm, fill=white, inner sep=0.25mm] (text3445) {$\ru$};
\end{tikzpicture}
}}
\end{equation*}
in $\V(Vx, W x)$, and that for all $x,y,z \in \B$, we have
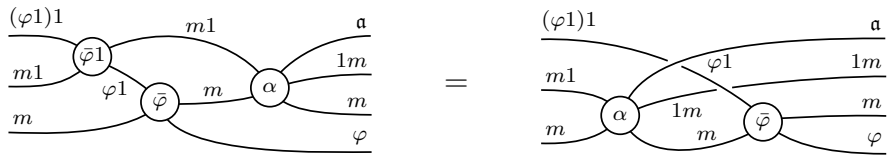
\begin{equation*}
\vcenter{\hbox{
\begin{tikzpicture}[y=0.80pt, x=0.7pt,yscale=-1, inner sep=0pt, outer sep=0pt, every text node part/.style={font=\tiny} ]
  \path[draw=black,line join=miter,line cap=butt,line width=0.650pt]
  (0.5761,996.5573) .. controls (18.8212,996.0522) and (33.5263,996.1982) ..
  node[above right=0.15cm,at start] {$\!\!(\phi \t 1) \t 1$}(44.0835,1005.7735);
  \path[draw=black,line join=miter,line cap=butt,line width=0.650pt]
  (0.5761,1020.4270) .. controls (22.7994,1020.6795) and (31.8588,1019.4362) ..
  node[above right=0.15cm,at start] {$\!m \t 1$}(44.9596,1008.7423);
  \path[draw=black,line join=miter,line cap=butt,line width=0.650pt]
  (46.9852,1007.9763) .. controls (58.3494,1011.5118) and (75.9772,1020.5756) ..
  node[below left=0.06cm,pos=0.5] {$\phi \t 1$}(80.4466,1027.1123);
  \path[draw=black,line join=miter,line cap=butt,line width=0.650pt]
  (0.5761,1042.0698) .. controls (14.4657,1042.3223) and (59.5285,1044.8883) ..
  node[above right=0.15cm,at start] {$\!m$}(80.8539,1029.7127);
  \path[draw=black,line join=miter,line cap=butt,line width=0.650pt]
  (82.2825,1030.8845) .. controls (90.1821,1054.0933) and (173.6677,1051.1703) ..
  node[above left=0.15cm,at end] {$\phi\!$}(194.7867,1051.3122);
  \path[draw=black,line join=miter,line cap=butt,line width=0.650pt]
  (141.5590,1019.2256) .. controls (145.5182,1011.6494) and (166.2283,996.5535) ..
  node[above left=0.15cm,at end] {$\mathfrak a\!$}(194.3349,996.5535);
  \path[draw=black,line join=miter,line cap=butt,line width=0.650pt]
  (141.6203,1021.3685) .. controls (141.6203,1021.3685) and (163.4571,1014.0062) ..
  node[above left=0.15cm,at end] {$1 \t m\!$}(194.7868,1014.7638);
  \path[draw=black,line join=miter,line cap=butt,line width=0.650pt]
  (141.4111,1023.0061) .. controls (149.4759,1031.0709) and (162.5493,1033.5038) ..
  node[above left=0.15cm,at end] {$m\!$}(194.7868,1033.7564);
  \path[draw=black,line join=miter,line cap=butt,line width=0.650pt]
  (82.5536,1028.3546) .. controls (96.1906,1029.3648) and (132.5520,1027.8496) ..
  node[above=0.09cm,pos=0.45] {$m$}(139.8756,1022.0412);
  \path[draw=black,line join=miter,line cap=butt,line width=0.650pt]
  (46.7482,1005.3271) .. controls (86.1067,987.1663) and (125.5659,1000.5348) ..
  node[above right=0.1cm,pos=0.4] {$m \t 1$}(139.8779,1019.5827);
  \path[fill=black] (45.434944,1006.3305) node[circle, draw, line width=0.65pt, minimum width=5mm, fill=white, inner sep=0.25mm] (text3305) {$\bar \phi \t 1$     };
  \path[fill=black] (81.600487,1027.7435) node[circle, draw, line width=0.65pt, minimum width=5mm, fill=white, inner sep=0.25mm] (text3309) {$\bar \phi$     };
  \path[fill=black] (140.38065,1021.0309) node[circle, draw, line width=0.65pt, minimum width=5mm, fill=white, inner sep=0.25mm] (text3317) {$\ass$   };
\end{tikzpicture}
}}
\quad\ \  = \quad\ \
\vcenter{\hbox{
\begin{tikzpicture}[y=0.80pt, x=0.6pt,yscale=-1, inner sep=0pt, outer sep=0pt, every text node part/.style={font=\tiny} ]
\path[use as bounding box] (0, 980) rectangle (216,1047);
  \path[draw=black,line join=miter,line cap=butt,line width=0.650pt]
  (0.5761,1016.9221) .. controls (22.7994,1017.1746) and (40.1997,1018.2637) ..
  node[above right=0.15cm,at start] {$\!m \t 1$}(49.1780,1027.2420);
  \path[draw=black,line join=miter,line cap=butt,line width=0.650pt]
  (0.5761,1042.0698) .. controls (14.4657,1042.3223) and (37.4465,1041.7251) ..
  node[above right=0.15cm,at start] {$\!m$}(49.1709,1030.0007);
  \path[draw=black,line join=miter,line cap=butt,line width=0.650pt]
  (51.1755,1030.5965) .. controls (67.4280,1056.6856) and (126.7968,1041.6654) ..
  node[above left=0.08cm,pos=0.65] {$m$}(139.8511,1033.1665);
  \path[draw=black,line join=miter,line cap=butt,line width=0.650pt]
  (142.8500,1030.5854) .. controls (142.8500,1030.5854) and (184.5612,1028.3914) ..
  node[above left=0.15cm,at end] {$m \!$}(215.8909,1029.1490);
  \path[draw=black,line join=miter,line cap=butt,line width=0.650pt]
  (142.6407,1033.9511) .. controls (154.1168,1040.7959) and (172.9964,1047.0248) ..
  node[above left=0.15cm,at end] {$\phi \!$}(215.8909,1046.9894);
  \path[draw=black,line join=miter,line cap=butt,line width=0.650pt]
  (0.5761,993.0038) .. controls (29.9053,992.4660) and (56.3004,996.4552) ..
  node[above right=0.15cm,at start] {$\!\!(\phi \t 1) \t 1$}(79.5654,1003.1390)
  (88.0853,1005.7499) .. controls (108.0007,1012.2403) and (125.3833,1020.6678) ..
  node[above right=0.08cm,pos=0.23] {$\phi \t 1$}(139.6961,1029.7577)
  (52.8868,1028.6427) .. controls (63.9796,1023.5539) and (85.7325,1019.6621) ..
  node[below right=0.08cm,pos=0.54] {$1 \t m$}(110.1322,1016.8137)
  (119.7394,1015.7580) .. controls (157.2114,1011.8826) and (198.3095,1010.3156) ..
  node[above left=0.15cm,at end] {$1 \t m\!$}(215.8909,1010.5363)
  (51.8895,1026.1092) .. controls (64.3077,997.3630) and (141.5428,992.3100) ..
  node[above left=0.15cm,at end] {$\mathfrak a \!$}(215.3148,992.6021);
  \path[fill=black] (49.917496,1028.8956) node[circle, draw, line width=0.65pt, minimum width=5mm, fill=white, inner sep=0.25mm] (text3309) {$\ass$     };
  \path[fill=black] (139.03424,1031.976) node[circle, draw, line width=0.65pt, minimum width=5mm, fill=white, inner sep=0.25mm] (text3317) {$\bar \phi$   };
\end{tikzpicture}
}}
\end{equation*}
in  $\V(\,(Vz \otimes \B(y,z)) \otimes \B(x,y),\, W x\,)$.

Dually, if $V, W \colon \A \tor \bullet$, then a left module morphism $\phi \colon V \to W$ is given by $1$-cells $\phi_x \colon Vx \to Wx$ and invertible $2$-cells
\begin{equation*}
\cd[@C+2.5em]{
 \A(x,y) \otimes Vx \ar[d]_{m} \ar[r]^-{1 \otimes \phi} \rtwocell{dr}{\bar \phi_{xy}} &
 \A(x,y) \otimes Wx \ar[d]^{m} \\
 Vy \ar[r]_-{\phi} & W y
}
\end{equation*}
subject to the axioms that for all $x \in \A$, we have
\begin{equation*}
\vcenter{\hbox{
\begin{tikzpicture}[y=0.80pt, x=0.8pt,yscale=-1, inner sep=0pt, outer sep=0pt, every text node part/.style={font=\tiny} ]
\path[use as bounding box] (8, 960) rectangle (85,1039);
  \path[draw=black,line join=miter,line cap=butt,line width=0.650pt]
  (10.0000,1012.3622) .. controls (31.9384,1012.3622) and (44.6754,1009.3804) ..
  node[above right=0.12cm,at start] {\!$m$}(48.2961,1005.5467);
  \path[draw=black,line join=miter,line cap=butt,line width=0.650pt]
  (10.0000,992.3622) .. controls (30.6605,992.3622) and (45.3144,994.7050) ..
  node[above right=0.12cm,at start] {\!$j \t 1$}(48.7221,999.6036);
  \path[draw=black,line join=miter,line cap=butt,line width=0.650pt]
  (10.0000,972.3622) .. controls (61.5891,972.3622) and (99.1240,1001.2238) ..
  node[above right=0.12cm,at start] {\!$\mathfrak l^\centerdot$}(51.2779,1002.6059);
  \path[draw=black,line join=miter,line cap=butt,line width=0.650pt]
  (10.0000,1032.3622) --
  node[above right=0.12cm,at start] {\!$\phi$}(79.9398,1032.3622);
  \path[fill=black] (49.096386,1002.6634) node[circle, draw, line width=0.65pt, minimum width=5mm, fill=white, inner sep=0.25mm] (text6241) {$\lu$   };
\end{tikzpicture}
}}
 =\quad\ \
\vcenter{\hbox{
\begin{tikzpicture}[y=0.80pt, x=0.8pt,yscale=-1, inner sep=0pt, outer sep=0pt, every text node part/.style={font=\tiny} ]
  \path[draw=black,line join=miter,line cap=butt,line width=0.650pt]
  (0.0000,982.3622) .. controls (102.0377,981.6862) and (110.8067,1017.2589) ..
  node[above right=0.12cm,at start] {\!$\mathfrak l^\centerdot$}(80.9286,1021.6480)
  (0.0000,1002.3622) .. controls (18.1994,1002.3622) and (32.9656,1002.8214) ..
  node[above right=0.12cm,at start] {\!$j \t 1$}(44.7561,1004.3771)
  (50.6372,1005.2912) .. controls (64.2974,1007.7724) and (73.2978,1012.1959) ..
  node[above right=0.08cm,pos=0.4] {$j \t 1$}(78.5714,1019.8622)
  (31.4286,1029.5050) .. controls (35.2452,1023.7186) and (43.6454,1003.1626) ..
  node[below right=0.08cm,pos=0.48] {$1 \t \phi$}(64.1786,991.1874)
  (70.7370,987.8640) .. controls (78.8039,984.3576) and (88.4503,982.1898) ..
  node[above left=0.12cm,at end] {$\phi$\!}(100.0000,982.3622);
  \path[draw=black,line join=miter,line cap=butt,line width=0.650pt]
  (32.1429,1032.7194) .. controls (52.3214,1037.8980) and (67.2857,1034.8621) ..
  node[below=0.07cm,pos=0.48] {$m$}(78.3571,1023.4336);
  \path[draw=black,line join=miter,line cap=butt,line width=0.650pt]
  (29.2857,1029.5050) .. controls (20.0000,1022.3622) and (15.3571,1022.3622) ..
  node[above right=0.12cm,at end] {\!$m$}(0.0000,1022.3622);
  \path[draw=black,line join=miter,line cap=butt,line width=0.650pt]
  (29.6429,1033.7908) .. controls (20.0000,1042.3622) and (14.8214,1042.3622) ..
  node[above right=0.12cm,at end] {\!$\phi$}(0.0000,1042.3622);
  \path[fill=black] (79.428574,1021.2908) node[circle, draw, line width=0.65pt, minimum width=5mm, fill=white, inner sep=0.25mm] (text3445) {$\lu$     };
  \path[fill=black] (30.357143,1031.2908) node[circle, draw, line width=0.65pt, minimum width=5mm, fill=white, inner sep=0.25mm] (text5100) {$\bar \phi$   };
\end{tikzpicture}
}}
\end{equation*}
in $\V(Vx, W x)$, and that for all $x,y,z \in \A$, we have
\begin{equation*}
\vcenter{\hbox{
\begin{tikzpicture}[y=0.80pt, x=0.75pt,yscale=-1, inner sep=0pt, outer sep=0pt, every text node part/.style={font=\tiny} ]
  \path[draw=black,line join=miter,line cap=butt,line width=0.650pt]
  (9.7893,987.6072) .. controls (28.0344,987.1021) and (38.7396,987.2481) ..
  node[above right=0.12cm,at start] {\!$m \t 1$}(49.2967,996.8234);
  \path[draw=black,line join=miter,line cap=butt,line width=0.650pt]
  (9.7893,1011.4770) .. controls (32.0127,1011.7295) and (37.0720,1010.4862) ..
  node[above right=0.12cm,at start] {\!$m$}(50.1728,999.7923);
  \path[draw=black,line join=miter,line cap=butt,line width=0.650pt]
  (52.1984,999.0263) .. controls (63.5627,1002.5618) and (81.1905,1011.6256) ..
  node[below left=0.07cm,pos=0.535] {$m$}(85.6598,1018.1623);
  \path[draw=black,line join=miter,line cap=butt,line width=0.650pt]
  (9.7893,1033.1198) .. controls (23.6789,1033.3723) and (64.7418,1035.9383) ..
  node[above right=0.12cm,at start] {\!$\phi$}(86.0672,1020.7627);
  \path[draw=black,line join=miter,line cap=butt,line width=0.650pt]
  (87.4957,1021.9345) .. controls (95.3954,1045.1433) and (198.8809,1042.2203) ..
  node[above left=0.12cm,at end] {$m$\!}(220.0000,1042.3622);
  \path[draw=black,line join=miter,line cap=butt,line width=0.650pt]
  (146.6243,1014.0561) .. controls (154.6891,1022.1209) and (187.7626,1022.1096) ..
  node[above left=0.12cm,at end] {$1 \t m$\!}(220.0000,1022.3622);
  \path[draw=black,line join=miter,line cap=butt,line width=0.650pt]
  (87.7668,1019.4046) .. controls (101.4039,1020.4148) and (137.7652,1018.8996) ..
  node[above=0.07cm,pos=0.45] {$1 \t \phi$}(145.0888,1013.0912);
  \path[draw=black,line join=miter,line cap=butt,line width=0.650pt]
  (51.9614,996.3770) .. controls (98.0288,990.6100) and (130.7791,991.5847) ..
  node[above=0.07cm,pos=0.45] {$1 \t m$}(145.0912,1010.6327);
  \path[draw=black,line join=miter,line cap=butt,line width=0.650pt]
  (146.7723,1010.2756) .. controls (148.8731,1006.2555) and (161.4487,998.6424) ..
  node[below right=0.03cm,pos=0.57] {$1 \t (1 \t \phi)$}(176.9905,992.3813)
  (183.7964,989.8003) .. controls (195.8431,985.5150) and (208.9335,982.3622) ..
  node[above left=0.12cm,at end] {$1 \t \phi$\!}(220.0000,982.3622)
  (51.6892,994.3214) .. controls (64.4562,979.3255) and (90.3328,975.4064) ..
  (125.4967,977.3769) .. controls (172.4952,980.0106) and (191.8918,1002.3622) ..
  node[above left=0.12cm,at end] {$\mathfrak a$\!}(220.0000,1002.3622);
  \path[fill=black] (50.648197,997.38043) node[circle, draw, line width=0.65pt, minimum width=5mm, fill=white, inner sep=0.25mm] (text3305) {$\ass$     };
  \path[fill=black] (86.813736,1018.7935) node[circle, draw, line width=0.65pt, minimum width=5mm, fill=white, inner sep=0.25mm] (text3309) {$\bar \phi$     };
  \path[fill=black] (145.5939,1012.0807) node[circle, draw, line width=0.65pt, minimum width=5mm, fill=white, inner sep=0.25mm] (text3317) {$1 \t \bar \phi$     };
\end{tikzpicture}
}}
\quad\ \  = \quad\ \
\vcenter{\hbox{
\begin{tikzpicture}[y=0.80pt, x=0.7pt,yscale=-1, inner sep=0pt, outer sep=0pt, every text node part/.style={font=\tiny} ]
  \path[draw=black,line join=miter,line cap=butt,line width=0.650pt]
  (13.2613,1006.6822) .. controls (35.4846,1006.9347) and (44.8849,1008.0238) ..
  node[above right=0.12cm,at start] {\!$m$}(53.8632,1017.0021);
  \path[draw=black,line join=miter,line cap=butt,line width=0.650pt]
  (13.2613,1031.8299) .. controls (27.1509,1032.0824) and (42.1317,1031.4852) ..
  node[above right=0.12cm,at start] {\!$\phi$}(53.8561,1019.7608);
  \path[draw=black,line join=miter,line cap=butt,line width=0.650pt]
  (55.8608,1020.3566) .. controls (74.9645,1037.4844) and (131.4820,1031.4255) ..
  node[above=0.08cm] {$m$}(144.5363,1022.9266);
  \path[draw=black,line join=miter,line cap=butt,line width=0.650pt]
  (147.9672,1022.2177) .. controls (147.9672,1022.2177) and (158.6704,1021.6046) ..
  node[above left=0.12cm,at end] {$1 \t m$\!}(190.0000,1022.3622);
  \path[draw=black,line join=miter,line cap=butt,line width=0.650pt]
  (147.8398,1024.9544) .. controls (152.2923,1029.3859) and (160.1636,1042.3622) ..
  node[above left=0.12cm,at end] {$m$\!}(190.0000,1042.3622);
  \path[draw=black,line join=miter,line cap=butt,line width=0.650pt]
  (57.5720,1018.4028) .. controls (68.0253,1002.9750) and (120.0000,982.3622) ..
  node[above left=0.12cm,at end] {$1 \t \phi$\!}(190.0000,982.3622)
  (13.8373,982.7639) .. controls (44.9228,982.1825) and (70.9382,986.8922) ..
  node[above right=0.12cm,at start] {\!$m \t 1$}(93.3970,994.5778)
  (100.9856,997.3359) .. controls (117.1189,1003.5453) and (131.3758,1011.2581) ..
  node[above right=0.08cm,pos=0.45] {$m \t 1$}(144.3813,1019.5178);
  \path[draw=black,line join=miter,line cap=butt,line width=0.650pt]
  (147.8821,1018.8904) .. controls (153.8384,1012.8048) and (160.0000,1002.3622) ..
  node[above left=0.12cm,at end] {$\mathfrak a$\!}(190.0000,1002.3622);
  \path[fill=black] (54.602711,1018.6556) node[circle, draw, line width=0.65pt, minimum width=5mm, fill=white, inner sep=0.25mm] (text3309) {$\bar \phi$     };
  \path[fill=black] (146.16768,1021.5919) node[circle, draw, line width=0.65pt, minimum width=5mm, fill=white, inner sep=0.25mm] (text3317) {$\ass$   };
\end{tikzpicture}
}}
\end{equation*}
in  $\V(\,(\A(y,z) \otimes \A(x,y)) \otimes Vx,\, W z\,)$.

Finally, if $M, N \colon \A \tor \B$, then a bimodule morphism $\phi \colon M \to N$ is given by morphisms $\phi_{ab} \colon M(b,a) \to N(b,a)$ together with $2$-cells making each $\phi_{\thg b}$ into a left $\A$-module morphism and each $\phi_{a \thg}$ into a right $\B$-module morphism, and such that for all $a,a' \in \A$ and $b, b' \in \B$, we have
\begin{equation*}
\vcenter{\hbox{
\begin{tikzpicture}[y=0.80pt, x=0.75pt,yscale=-1, inner sep=0pt, outer sep=0pt, every text node part/.style={font=\tiny} ]
  \path[draw=black,line join=miter,line cap=butt,line width=0.650pt]
  (9.7893,987.6072) .. controls (28.0344,987.1021) and (38.7396,987.2481) ..
  node[above right=0.12cm,at start] {\!$m \t 1$}(49.2967,996.8234);
  \path[draw=black,line join=miter,line cap=butt,line width=0.650pt]
  (9.7893,1011.4770) .. controls (32.0127,1011.7295) and (37.0720,1010.4862) ..
  node[above right=0.12cm,at start] {\!$m$}(50.1728,999.7923);
  \path[draw=black,line join=miter,line cap=butt,line width=0.650pt]
  (52.1984,999.0263) .. controls (63.5627,1002.5618) and (81.1905,1011.6256) ..
  node[below left=0.07cm,pos=0.535] {$m$}(85.6598,1018.1623);
  \path[draw=black,line join=miter,line cap=butt,line width=0.650pt]
  (9.7893,1033.1198) .. controls (23.6789,1033.3723) and (64.7418,1035.9383) ..
  node[above right=0.12cm,at start] {\!$\phi$}(86.0672,1020.7627);
  \path[draw=black,line join=miter,line cap=butt,line width=0.650pt]
  (87.4957,1021.9345) .. controls (95.3954,1045.1433) and (198.8809,1042.2203) ..
  node[above left=0.12cm,at end] {$m$\!}(220.0000,1042.3622);
  \path[draw=black,line join=miter,line cap=butt,line width=0.650pt]
  (146.6243,1014.0561) .. controls (154.6891,1022.1209) and (187.7626,1022.1096) ..
  node[above left=0.12cm,at end] {$1 \t m$\!}(220.0000,1022.3622);
  \path[draw=black,line join=miter,line cap=butt,line width=0.650pt]
  (87.7668,1019.4046) .. controls (101.4039,1020.4148) and (137.7652,1018.8996) ..
  node[above=0.07cm,pos=0.45] {$1 \t \phi$}(145.0888,1013.0912);
  \path[draw=black,line join=miter,line cap=butt,line width=0.650pt]
  (51.9614,996.3770) .. controls (98.0288,990.6100) and (130.7791,991.5847) ..
  node[above=0.07cm,pos=0.45] {$1 \t m$}(145.0912,1010.6327);
  \path[draw=black,line join=miter,line cap=butt,line width=0.650pt]
  (146.7723,1010.2756) .. controls (148.8731,1006.2555) and (161.4487,998.6424) ..
  node[below right=0.03cm,pos=0.65] {$1 \t (\phi \t 1)$}(176.9905,992.3813)
  (183.7964,989.8003) .. controls (195.8431,985.5150) and (208.9335,982.3622) ..
  node[above left=0.12cm,at end] {$(1 \t \phi) \t 1$\!}(220.0000,982.3622)
  (51.6892,994.3214) .. controls (64.4562,979.3255) and (90.3328,975.4064) ..
  (125.4967,977.3769) .. controls (172.4952,980.0106) and (191.8918,1002.3622) ..
  node[above left=0.12cm,at end] {$\mathfrak a$\!}(220.0000,1002.3622);
  \path[fill=black] (50.648197,997.38043) node[circle, draw, line width=0.65pt, minimum width=5mm, fill=white, inner sep=0.25mm] (text3305) {$\ass$     };
  \path[fill=black] (86.813736,1018.7935) node[circle, draw, line width=0.65pt, minimum width=5mm, fill=white, inner sep=0.25mm] (text3309) {$\bar \phi$     };
  \path[fill=black] (145.5939,1012.0807) node[circle, draw, line width=0.65pt, minimum width=5mm, fill=white, inner sep=0.25mm] (text3317) {$1 \t \bar \phi^{-1}$     };
\end{tikzpicture}
}}
\quad\ \  = \quad\ \
\vcenter{\hbox{
\begin{tikzpicture}[y=0.80pt, x=0.7pt, inner sep=0pt, outer sep=0pt, every text node part/.style={font=\tiny} ]
  \path[draw=black,line join=miter,line cap=butt,line width=0.650pt]
  (0.5761,996.5573) .. controls (18.8212,996.0522) and (33.5263,996.1982) ..
  node[above right=0.15cm,at start] {$\!\phi$}(44.0835,1005.7735);
  \path[draw=black,line join=miter,line cap=butt,line width=0.650pt]
  (0.5761,1020.4270) .. controls (22.7994,1020.6795) and (31.8588,1019.4362) ..
  node[above right=0.15cm,at start] {$\!m$}(44.9596,1008.7423);
  \path[draw=black,line join=miter,line cap=butt,line width=0.650pt]
  (46.9852,1007.9763) .. controls (58.3494,1011.5118) and (75.9772,1020.5756) ..
  node[above left=0.06cm,pos=0.54] {$\phi \t 1$}(80.4466,1027.1123);
  \path[draw=black,line join=miter,line cap=butt,line width=0.650pt]
  (0.5761,1042.0698) .. controls (14.4657,1042.3223) and (59.5285,1044.8883) ..
  node[above right=0.15cm,at start] {$\!m \t 1$}(80.8539,1029.7127);
  \path[draw=black,line join=miter,line cap=butt,line width=0.650pt]
  (82.2825,1030.8845) .. controls (90.1821,1054.0933) and (173.6677,1051.1703) ..
  node[above left=0.15cm,at end] {$(1 \t \phi) \t 1\!$}(194.7867,1051.3122);
  \path[draw=black,line join=miter,line cap=butt,line width=0.650pt]
  (141.5590,1019.2256) .. controls (145.5182,1011.6494) and (166.2283,996.5535) ..
  node[above left=0.15cm,at end] {$m\!$}(194.3349,996.5535);
  \path[draw=black,line join=miter,line cap=butt,line width=0.650pt]
  (141.6203,1021.3685) .. controls (141.6203,1021.3685) and (163.4571,1014.0062) ..
  node[above left=0.15cm,at end] {$1 \t m\!$}(194.7868,1014.7638);
  \path[draw=black,line join=miter,line cap=butt,line width=0.650pt]
  (141.4111,1023.0061) .. controls (149.4759,1031.0709) and (162.5493,1033.5038) ..
  node[above left=0.15cm,at end] {$\mathfrak a\!$}(194.7868,1033.7564);
  \path[draw=black,line join=miter,line cap=butt,line width=0.650pt]
  (82.5536,1028.3546) .. controls (96.1906,1029.3648) and (132.5520,1027.8496) ..
  node[above=0.09cm,pos=0.45] {$m \t 1$}(139.8756,1022.0412);
  \path[draw=black,line join=miter,line cap=butt,line width=0.650pt]
  (46.7482,1005.3271) .. controls (86.1067,987.1663) and (125.5659,1000.5348) ..
  node[below right=0.1cm,pos=0.4] {$m$}(139.8779,1019.5827);
  \path[fill=black] (45.434944,1006.3305) node[circle, draw, line width=0.65pt, minimum width=5mm, fill=white, inner sep=0.25mm] (text3305) {$\bar \phi^{-1}$     };
  \path[fill=black] (81.600487,1027.7435) node[circle, draw, line width=0.65pt, minimum width=5mm, fill=white, inner sep=0.25mm] (text3309) {$\bar \phi \t 1$     };
  \path[fill=black] (140.38065,1021.0309) node[circle, draw, line width=0.65pt, minimum width=5mm, fill=white, inner sep=0.25mm] (text3317) {$\ass$   };
\end{tikzpicture}
}}
\end{equation*}
in  $\V(\,(\A(a,a') \otimes M(b,a)) \otimes \B(b',b),\, N(b',a')\,)$.

\begin{Ex}\label{ex:pseudoinversemodulemorphism}
Let $\phi \colon M \to N$ be a bimodule morphism all of whose $1$-cell components are equivalences. Then we obtain a bimodule morphism $\phi^\centerdot \colon N \to M$ whose $1$-cell components are adjoint inverses of the $1$-cell components of $\phi$, and whose $2$-cell components are the mates under adjunction of the inverses of the $2$-cell components of $\phi$.
In fact, $\phi^\centerdot$ is pseudoinverse to $\phi$ in the bicategory of bimodules as defined below.
\end{Ex}

\begin{Ex}\label{ex:modulemorphisminducedbyfunctor}
Given a $\V$-functor $F \colon \B \to \C$, there is a morphism of $\B$-$\B$-bimodules $\hat F \colon \B \to \C(F,F)$ with $1$-cell components $F_{b,b'}$, and $2$-cell components obtained from binary functoriality constraints $\fm$ for $F$ and pseudofunctoriality of $\otimes$.
\end{Ex}
\begin{Ex}\label{ex:module-morphism-restriction-along-functor}
Given a map of $\C$-$\D$-bimodules $\phi \colon M \to N$ together with $\V$-functors $F \colon \A \to \C$ and $G \colon \B \to \D$, we induce a map $\phi(G,F) \colon M(G,F) \to N(G,F)$ of $\A$-$\B$-bimodules, whose $1$-cell component at $(x,y)$ is $\phi_{Gx,Fy}$ and whose $2$-cell components are obtained from those of $\phi$ and interchange isomorphisms for $\otimes$.
\end{Ex}

\subsection{Module transformations}
If $\phi, \psi \colon V \to W $ are right $\B$-module morphisms, then a \emph{module transformation} $\Gamma \colon \phi \Rightarrow \psi$ is given by $2$-cells $\Gamma_x \colon \phi_x \Rightarrow \psi_x$ for each $x \in \B$, subject to the axiom that for all $x, y \in \B$, we have
\begin{equation*}
\vcenter{\hbox{\begin{tikzpicture}[y=0.80pt, x=0.7pt,yscale=-1, inner sep=0pt, outer sep=0pt, every text node part/.style={font=\tiny} ]
  \path[draw=black,line join=miter,line cap=butt,line width=0.650pt]
  (0.0000,1022.3622) .. controls (20.0000,1022.3622) and (33.4489,1022.0637) ..
  node[above right=0.12cm,at start] {$\!m$}(38.7575,1014.9418);
  \path[draw=black,line join=miter,line cap=butt,line width=0.650pt]
  (0.0000,1002.3622) .. controls (20.0000,1002.3622) and (33.6763,1003.7813) ..
  node[above right=0.12cm,at start] {$\!\phi \t 1$}(38.5075,1009.5873);
  \path[draw=black,line join=miter,line cap=butt,line width=0.650pt]
  (88.9541,1022.3622) .. controls (71.2005,1022.6021) and (47.2500,1020.3622) ..
  node[above right=0.09cm,pos=0.55] {$\phi$}(40.5000,1013.8622);
  \path[draw=black,line join=miter,line cap=butt,line width=0.650pt]
  (130.0000,1002.3622) .. controls (116.5000,1002.3622) and (50.6321,999.4801) ..
  node[above left=0.12cm,at start] {$m\!$}(40.7500,1009.3622);
  \path[draw=black,line join=miter,line cap=butt,line width=0.650pt]
  (90.4541,1022.3622) --
  node[above left=0.12cm,at end] {$\psi\!$}(130.0000,1022.3622);
  \path[fill=black] (40,1012.1122) node[circle, draw, line width=0.65pt, minimum width=5mm, fill=white, inner sep=0.25mm] (text3086) {$\bar \phi$};
  \path[xscale=1.005,yscale=0.995,fill=black] (89.022888,1027.5282) node[circle, draw, line width=0.65pt, minimum width=5mm, fill=white, inner sep=0.25mm] (text3092) {$\Gamma$};
\end{tikzpicture}}} \quad \ = \quad \
\vcenter{\hbox{\begin{tikzpicture}[y=0.80pt, x=0.7pt, xscale=-1,inner sep=0pt, outer sep=0pt, every text node part/.style={font=\tiny} ]
  \path[draw=black,line join=miter,line cap=butt,line width=0.650pt]
  (0.0000,1022.3622) .. controls (20.0000,1022.3622) and (33.4489,1022.0637) ..
  node[above left=0.12cm,at start] {$m\!$}(38.7575,1014.9418);
  \path[draw=black,line join=miter,line cap=butt,line width=0.650pt]
  (0.0000,1002.3622) .. controls (20.0000,1002.3622) and (33.6763,1003.7813) ..
  node[above left=0.12cm,at start] {$\psi\!$}(38.5075,1009.5873);
  \path[draw=black,line join=miter,line cap=butt,line width=0.650pt]
  (88.9541,1022.3622) .. controls (71.2005,1022.6021) and (47.2500,1020.3622) ..
  node[above right=0.09cm,pos=0.37] {$\psi \t 1$}(40.5000,1013.8622);
  \path[draw=black,line join=miter,line cap=butt,line width=0.650pt]
  (130.0000,1002.3622) .. controls (116.5000,1002.3622) and (50.6321,999.4801) ..
  node[above right=0.12cm,at start] {$\!m$}(40.7500,1009.3622);
  \path[draw=black,line join=miter,line cap=butt,line width=0.650pt]
  (90.4541,1022.3622) --
  node[above right=0.12cm,at end] {$\! \phi\t 1$}(130.0000,1022.3622);
  \path[fill=black] (40,1012.1122) node[circle, draw, line width=0.65pt, minimum width=5mm, fill=white, inner sep=0.25mm] (text3086) {$\bar \psi$};
  \path[xscale=1.005,yscale=0.995,fill=black] (89.022888,1027.5282) node[circle, draw, line width=0.65pt, minimum width=5mm, fill=white, inner sep=0.25mm] (text3092) {$\Gamma \t 1$};
\end{tikzpicture}}}\end{equation*}
in $\V(Vy \otimes \B(x,y), W x)$.
The notion of left $\B$-module transformation is dual; while if $\phi, \psi \colon M \to N$ are $\A$-$\B$-bimodule morphisms, then a bimodule transformation $\Gamma \colon \phi \Rightarrow \psi$ is given by $2$-cells $\Gamma_{ab} \colon \phi_{ab} \Rightarrow \psi_{ab}$ making each $\Gamma_{a \thg}$ and each $\Gamma_{\thg b}$ into a one-sided module transformation.

\subsection{Compositional structure}
The  $\A$-$\B$-bimodules, morphisms and transformations form a bicategory ${}_\A \cat{Mod}_\B$; composition of module morphisms is given by that in $\V$ at the level of $1$-cell data, and by pasting in $\V$ together with the pseudofunctoriality of $\otimes$ at the $2$-cell level. The remaining data for this bicategory---horizontal and vertical composition of $2$-cells, and associativity and unit coherence isomorphisms---is obtained pointwise from that in $\V$. Similarly, we obtain bicategories ${}_\A\cat{Mod}_\bullet$ of left $\A$-modules, and ${}_\bullet \cat{Mod}_\B$ of $\A$-$\B$-bimodules (and in accordance with our convention, we also have ${}_\bullet \cat{Mod}_\bullet \defeq \V$).

\subsection{Comparing $\bullet$ with $\I$}
\label{sec:bullet-unit}

There are evident forgetful functors ${}_\A \cat{Mod}_\B \to {}_\bullet \cat{Mod}_\B$, ${}_\A \cat{Mod}_\B \to {}_\A \cat{Mod}_\bullet$, and so on.
It is easy to see that with $\I$ the unit $\V$-bicategory from Example~\ref{ex:unit-vcat}, the functors ${}_\I \cat{Mod}_\B \to {}_\bullet \cat{Mod}_\B$, ${}_\A \cat{Mod}_\I \to {}_\A \cat{Mod}_\bullet$, and so on, are biequivalences.
By passing across these biequivalences, many statements about bimodules literally imply the corresponding statements for one-sided modules.
On the other hand, many definitions and proofs are convenient to perform for one-sided modules first, so the notion of one-sided module is still useful to have around.

\subsection{Ordinary functors induced by modules}
If $M$ is an $\A$-$\B$-bimodule, then the construction of Section~\ref{subsec:hom-functors} carries over, mutatis mutandis, to yield a functor $\B_0^\op \times \A_0 \to \V$ sending $(b,a)$ to $M(b,a)$, and so on. With respect to this functor structure, the action morphisms $M(b,a) \otimes \B(b',b) \to M(b',a)$ and $\A(a,a') \otimes M(b,a) \to M(b,a')$ now become pseudonatural in each variable, precisely as in Section~\ref{subsec:hom-functors}; as there, we shall uniformly denote the $2$-cells witnessing this pseudonaturality by $\alpha$.

\begin{Ex}\label{ex:homs-induced-by-modules}
Given $M \colon \A \tor \C$, there is for any $\V$-bicategory $\B$ a functor $M(\thg, 1) \colon \V\text-\cat{Bicat}(\B, \C)^\op \to {}_\A \cat{Mod}_\B$ that on $0$-cells sends $F$ to $M(F,1)$ as defined in Example~\ref{ex:reindexing}. On $1$-cells, it sends a $\V$-transformation $\gamma \colon F \Rightarrow G$ to the bimodule morphism $M(\gamma, 1) \colon M(G,1) \to M(F,1)$ whose $1$-cell components are $M(\gamma_b, a) \colon M(Gb, a) \to M(Fb,a)$, whose $2$-cell components for the left actions are the pseudonaturality morphisms $\alpha$ of~\eqref{eq:pseudonat-of-composition}, and whose $2$-cell components for the right actions are given by
\begin{equation*}
\vcenter{\hbox{\begin{tikzpicture}[y=0.8pt, x=0.95pt,yscale=-1, inner sep=0pt, outer sep=0pt, every text node part/.style={font=\tiny} ]
  \path[draw=black,line join=miter,line cap=butt,line width=0.650pt]
  (50.0000,1002.3622) .. controls (64.1580,1002.3622) and (74.9680,998.7367) ..
  node[above right=0.12cm,at start] {\!$1 \t F$}(84.6786,994.5064)
  (89.8515,992.1703) .. controls (100.2347,987.3581) and (109.7701,982.4744) ..
  (121.5012,981.6080) .. controls (134.6487,980.6370) and (141.0127,990.7542) ..
  node[above=0.07cm, at start] {$1 \t F$}(150.7225,990.7801)
  (50.0000,982.3622) .. controls (90.0895,982.6122) and (98.1688,1007.0350) ..
  node[above right=0.12cm,at start] {\!$\C(\gamma, 1) \t 1$}(102.3205,1010.9633);
  \path[draw=black,line join=miter,line cap=butt,line width=0.650pt]
  (50.0000,1022.3622) .. controls (69.9402,1022.3622) and (99.7030,1017.4882) ..
  node[above right=0.12cm,at start] {\!$m$}(102.2284,1014.1210);
  \path[draw=black,line join=miter,line cap=butt,line width=0.650pt]
  (107.4418,1014.3652) .. controls (111.6507,1018.8547) and (195.1523,1018.8547) ..
  node[above=0.07cm] {$m$}(198.5195,1013.8040);
  \path[draw=black,line join=miter,line cap=butt,line width=0.650pt]
  (107.4418,1010.7801) .. controls (110.8089,1004.6070) and (145.1254,999.1353) ..
  node[above left=0.06cm,rotate=18,pos=0.73] {$1 \t \C(1, \gamma)$}(150.4567,994.0845);
  \path[draw=black,line join=miter,line cap=butt,line width=0.650pt]
  (153.8627,990.6398) .. controls (157.7910,985.4488) and (222.5839,982.3622) ..
  node[above left=0.12cm,at end] {$1 \t G$\!}(246.0000,982.3622);
  \path[draw=black,line join=miter,line cap=butt,line width=0.650pt]
  (153.8627,994.2248) .. controls (157.2298,999.1353) and (196.3522,1006.8518) ..
  node[above right=0.04cm,rotate=-14,pos=0.27] {$1 \t \C(\gamma, 1)$}(198.4567,1010.9204);
  \path[draw=black,line join=miter,line cap=butt,line width=0.650pt]
  (201.6004,1010.6306) .. controls (206.6512,1005.0187) and (215.0546,1001.6589) ..
  node[above left=0.12cm,at end] {$m$\!}(246.0000,1002.3622);
  \path[draw=black,line join=miter,line cap=butt,line width=0.650pt]
  (201.6004,1013.9978) .. controls (206.0900,1018.7680) and (214.2540,1022.3622) ..
  node[above left=0.12cm,at end] {$\C(\gamma, 1)$\!}(246.0000,1022.3622);
  \path[fill=black] (103.60317,1013.2352) node[circle, draw, line width=0.65pt, minimum width=5mm, fill=white, inner sep=0.25mm] (text4238) {$\ass$     };
  \path[fill=black] (151.68254,992.20343) node[circle, draw, line width=0.65pt, minimum width=5mm, fill=white, inner sep=0.25mm] (text4242) {$1 \t\bar \gamma^{-1}$     };
  \path[fill=black] (202.36508,1012.8384) node[circle, draw, line width=0.65pt, minimum width=5mm, fill=white, inner sep=0.25mm] (text4246) {$\alpha^{-1}$   };
\end{tikzpicture}\rlap{ .}}}
\end{equation*}
On $2$-cells, $M(\thg, 1)$ sends a $\V$-modification $\Gamma \colon \gamma \Rrightarrow \delta$ to the bimodule transformation $M(\gamma,1) \Rightarrow M(\delta,1)$ with components $M(\Gamma_b, a)$. The functoriality constraints of $M(\thg, 1)$ are obtained pointwise from those of the functor $\B_0^\op \times \A_0 \to \V$ as defined in the preceding section. 

Dualising the above construction, we obtain a functor $M(1, \thg) \colon \V\text-\cat{Bicat}(\B, \A) \to {}_{\B}\cat{Mod}_\C$; and replacing the bicategory $\B$ by $\bullet$, we obtain two further variants: a functor $\C_0 \to {}_\A\cat{Mod}_\bullet$ sending $c$ to $M(c, \thg)$, and a functor $\A_0 \to {}_\bullet\cat{Mod}_\C$ sending $a$ to $M(\thg, a)$.
\end{Ex}

\subsection{Copowers of modules}\label{subsec:copowers-of-modules}
If $W$ is a right $\B$-module and $A \in \V$, then there is a right $\B$-module $A \otimes W$, called the \emph{copower of $W$ by $A$}, that has $(A \otimes W)(x) = A \otimes Wx$, action morphisms
\begin{equation*}
(A \otimes Wy) \otimes \B(x,y) \xrightarrow{\mathfrak a} A \otimes (Wy \otimes \B(x,y)) \xrightarrow{1 \otimes m} A \otimes Wx
\end{equation*}
and unit and associativity constraints given by the respective $2$-cells
\begin{equation*}
\vcenter{\hbox{
\begin{tikzpicture}[y=0.80pt, x=0.9pt,yscale=-1, inner sep=0pt, outer sep=0pt, every text node part/.style={font=\scriptsize} ]
  \path[draw=black,line join=miter,line cap=butt,line width=0.650pt]
  (10.0000,982.3622) .. controls (32.6405,982.3622) and (76.4458,985.2537) ..
  node[above right=0.12cm,at start] {\!$\mathfrak a$}(86.9880,977.1212)
  (10.0000,962.3622) .. controls (28.0025,962.3622) and (38.7468,970.7123) ..
  node[above right=0.12cm,at start] {\!$1 \t j$}(45.6613,980.3646)
  (48.8638,985.3504) .. controls (53.8461,993.9740) and (56.3591,1002.6793) ..
  node[left=0.08cm,pos=0.4] {$1 \t (1 \t j)$}(58.6747,1006.7959);
  \path[draw=black,line join=miter,line cap=butt,line width=0.650pt]
  (60.0000,1008.3622) .. controls (60.0000,1008.3622) and (140.2169,990.2538) ..
  node[below right=0.08cm,pos=0.5] {$1 \t \mathfrak r$}(96.3012,974.5308);
  \path[draw=black,line join=miter,line cap=butt,line width=0.650pt]
  (84.1446,971.0851) .. controls (63.2253,966.5670) and (71.4657,953.1504) ..
  (87.5104,957.6685) .. controls (95.9012,960.0313) and (118.2849,972.9646) ..
  node[above left=0.12cm,at end] {$\mathfrak r$\!}(149.8193,972.9646);
  \path[draw=black,line join=miter,line cap=butt,line width=0.650pt]
  (10.0000,1011.0730) .. controls (29.5783,1011.0730) and (53.9157,1011.3742) ..
  node[above right=0.12cm,at start] {\!$1 \t m$}(58.1325,1008.3622);
  \path[fill=black] (90.662651,975.25372) node[circle, draw, line width=0.65pt, minimum width=5mm, fill=white, inner sep=0.25mm] (text6454) {$\rho^{-1}$
     };
  \path[fill=black] (58.73494,1008.6634) node[circle, draw, line width=0.65pt, minimum width=5mm, fill=white, inner sep=0.25mm] (text6458) {$1 \t \tau$
     };
\end{tikzpicture}}}
\quad \text{and} \quad
\vcenter{\hbox{
\begin{tikzpicture}[y=0.90pt, x=1pt,yscale=-1, inner sep=0pt, outer sep=0pt, every text node part/.style={font=\tiny} ]
  \path[draw=black,line join=miter,line cap=butt,line width=0.650pt]
  (70.0000,922.3622) .. controls (101.2330,922.3622) and (143.6145,921.7598) ..
  node[above right=0.12cm,at start] {\!$\mathfrak a \t 1$}(151.3870,930.2041);
  \path[draw=black,line join=miter,line cap=butt,line width=0.650pt]
  (126.7723,969.6175) .. controls (128.2875,962.0413) and (140.0945,940.6924) ..
  node[above left=0.08cm,pos=0.47] {$1 \t \mathfrak a$}(151.2391,934.9696);
  \path[draw=black,line join=miter,line cap=butt,line width=0.650pt]
  (127.6786,973.3558) .. controls (136.1123,980.1329) and (156.4071,982.3622) ..
  node[above left=0.12cm,at end] {$1 \t m$\!}(190.0000,982.3622);
  \path[draw=black,line join=miter,line cap=butt,line width=0.650pt]
  (152.7723,930.8750) .. controls (155.6057,928.5814) and (169.8180,922.9558) ..
  node[above left=0.12cm,at end] {$\mathfrak a$\!}(190.0000,922.3622);
  \path[draw=black,line join=miter,line cap=butt,line width=0.650pt]
  (153.1294,933.7322) .. controls (155.4312,936.0339) and (161.1229,962.3622) ..
  node[above left=0.12cm,at end] {$\mathfrak a$\!}(190.0000,962.3622)
  (127.8878,970.8146) .. controls (136.2744,968.2835) and (145.9898,957.8154) ..
  node[below right=0.05cm,pos=0.35,rotate=35] {$1 \t (1 \t m)$}(159.5167,950.3078)
  (164.0504,947.9860) .. controls (171.4928,944.4990) and (180.0280,942.1210) ..
  node[above left=0.12cm,at end] {$1 \t m$\!}(190.0000,942.3622);
  \path[draw=black,line join=miter,line cap=butt,line width=0.650pt]
  (70.0000,982.3622) .. controls (94.7819,982.3622) and (118.6172,978.5843) ..
  node[above right=0.12cm,at start] {\!$1 \t m$}(125.9408,972.7759);
  \path[draw=black,line join=miter,line cap=butt,line width=0.650pt]
  (70.0000,962.3622) .. controls (104.0413,962.3622) and (95.4124,932.6539) ..
  node[above right=0.12cm,at start] {\!$\mathfrak a$}(151.2919,932.6539)
  (70.0000,942.3622) .. controls (81.3043,942.3622) and (90.5727,944.4985) ..
  node[above right=0.12cm,at start] {\!\!$(1 \t m) \t 1$}(98.2107,947.7957)
  (103.1302,950.1901) .. controls (113.3893,955.7771) and (120.2698,963.4297) ..
  node[below left=0.05cm,pos=0.74,rotate=-40] {$1 \t (m \t 1)$}(125.0912,969.9746);
  \path[fill=black] (152.0495,932.40131) node[circle, draw, line width=0.65pt, minimum width=5mm, fill=white, inner sep=0.25mm] (text3313) {$\pi$     };
  \path[fill=black] (125.89511,972.28424) node[circle, draw, line width=0.65pt, minimum width=5mm, fill=white, inner sep=0.25mm] (text3317) {$1 \t \ass$   };
\end{tikzpicture}\rlap{ .}
}}
\end{equation*}

The assignation $(A, W) \mapsto A \otimes W$ is the action on objects of a functor
\begin{equation}\label{eq:cotensorfunctor}
\mathord{\otimes} \colon \V \times {}_\bullet\cat{Mod}_\B \to {}_\bullet\cat{Mod}_\B\rlap{ ,}
\end{equation}
that on morphisms, sends $(f, \phi) \colon (A, V) \to (A',V')$ to the module morphism $f \otimes \phi$ with $1$-cell components $(f \otimes \phi)_x = f \otimes \phi_x$, and $2$-cell components
\begin{equation*}
\begin{tikzpicture}[y=0.85pt, x=0.85pt,yscale=-1, inner sep=0pt, outer sep=0pt, every text node part/.style={font=\tiny} ]
  \path[draw=black,line join=miter,line cap=butt,line width=0.650pt]
  (10.0000,1022.3622) .. controls (46.7914,1022.3622) and (111.8675,1019.1996) ..
  node[above right=0.12cm,at start] {\!$1 \t m$}(117.8916,1014.0791);
  \path[draw=black,line join=miter,line cap=butt,line width=0.650pt]
  (10.0000,982.3622) .. controls (21.7418,982.3622) and (35.7774,984.2742) ..
  node[above right=0.12cm,at start] {\!\!$(f \t \phi) \t 1$}(49.7809,987.1303)
  (56.3982,988.5490) .. controls (86.7372,995.3632) and (115.1847,1005.9278) ..
  node[below left=0.1cm,pos=0.28,rotate=-16] {$f \t (\phi \t 1)$}
  node[below left=0.1cm,pos=0.73,rotate=-20] {$1 \t (\phi \t 1)$}(117.7410,1010.2537)
  (10.0000,1002.3622) .. controls (55.0260,1002.3622) and (65.6012,965.5276) ..
  node[above right=0.12cm,at start] {\!$\mathfrak a$}
  (97.4433,963.3588) .. controls (132.6073,960.9638) and (164.4699,982.3622) ..
  node[above left=0.12cm,at end] {$\mathfrak a$\!}(200.0000,982.3622);
  \path[draw=black,line join=miter,line cap=butt,line width=0.650pt]
  (121.3554,1014.1092) .. controls (128.1650,1020.9188) and (177.3375,1022.3622) ..
  node[above left=0.08cm,pos=0.49] {$1 \t \phi$}
  node[above left=0.12cm,at end] {$f \t \phi$\!}(200.0000,1022.3622);
  \path[draw=black,line join=miter,line cap=butt,line width=0.650pt]
  (83.6339,995.7873) .. controls (138.7056,959.8477) and (170.3645,1011.2520) ..
  node[above right=0.06cm,pos=0.4] {$f \t 1$}(170.2147,1021.5661)
  (121.0542,1010.4947) .. controls (125.5992,1006.5178) and (140.8270,1004.4690) ..
  node[above left=0.1cm,pos=0.79] {$1 \t m$}(157.7638,1003.4232)
  (163.9675,1003.0862) .. controls (176.8918,1002.4731) and (190.1676,1002.3622) ..
  node[above left=0.12cm,at end] {$1 \t m$\!}(200.0000,1002.3622);
  \path[fill=black] (118.84506,1012.3212) node[circle, draw, line width=0.65pt, minimum width=5mm, fill=white, inner sep=0.25mm] (text5042) {$1 \t \bar \phi$};
  \path[shift={(7.3423566,64.856253)},draw=black,fill=black] (77.8600,931.1000)arc(0.000:180.000:1.250)arc(-180.000:0.000:1.250) -- cycle;
  \path[shift={(93.570722,90.539316)},draw=black,fill=black] (77.8600,931.1000)arc(0.000:180.000:1.250)arc(-180.000:0.000:1.250) -- cycle;
\end{tikzpicture}\rlap{ ;}
\end{equation*}
and which on $2$-cells, sends $(\gamma, \Gamma) \colon (f, \phi) \Rightarrow (g, \psi)$ to the module transformation $\gamma \otimes \Gamma$ with components $(\gamma \otimes \Gamma)_x = \gamma \otimes \Gamma_x$. The pseudofunctoriality constraints of $\otimes$ are obtained pointwise from the pseudofunctoriality of the tensor product on $\V$.

The functor~\eqref{eq:cotensorfunctor} in fact underlies an action of the monoidal bicategory $\V$ on the bicategory ${}_\bullet \cat{Mod}_\B$. Thus, for each right $\B$-module $W$, there is a module equivalence $\mathfrak l \colon I \otimes W \to W$ in ${}_\bullet \cat{Mod}_\B$ with $1$-cell components $\mathfrak l_{Wx} \colon I \otimes Wx \to Wx$ and $2$-cell components obtained using $\lambda$ and pseudonaturality of $\mathfrak l$. 
Similarly, for each $W \in {}_\bullet \cat{Mod}_\B$ and $A, B \in \V$, there is a module equivalence $\mathfrak a \colon (A \otimes B) \otimes W \to A \otimes (B \otimes W)$ with $1$-cell components $\mathfrak a_{A,B,Wx}$, and $2$-cell components obtained using $\pi$ and pseudonaturality of $\mathfrak a$.
Finally, there are invertible module modifications $\pi$, $\nu$ and $\lambda$ whose components are those of the corresponding coherence constraints for the monoidal bicategory $\V$; it follows that these satisfy axioms corresponding to the axioms for a monoidal bicategory.

In a completely analogous way, we can define the copower $W \otimes B$ of a left $\A$-module $W$ by an object $B$; it has components $(W \otimes B)(x) = Wx \otimes B$ and action morphisms
\[\A(x,y) \otimes (Wx \otimes B) \xrightarrow{\mathfrak a^\centerdot} (\A(x,y) \otimes Wx) \otimes B \xrightarrow{m \otimes 1} Wy \otimes B\rlap{ ,}\] and provides the assignation on objects of a right action $\mathord{\otimes} \colon {}_\A\cat{Mod}_\bullet \times \V \to {}_\A\cat{Mod}_\bullet$ of the
monoidal bicategory $\V$ on ${}_\A\cat{Mod}_\bullet$.

Finally, we may combine the above two constructions; given a left $\A$-module $V$ and a right $\B$-module $W$, there is an $\A$-$\B$-bimodule $V \otimes W$ with $(V \otimes W)(b,a) = Va \otimes Wb$, with right $\B$-actions $(1 \otimes m) \circ \mathfrak a \colon (Va \otimes Wb) \otimes \B(b',b) \to Va \otimes Wb'$ and with left $\A$-actions $(m \otimes 1) \circ \mathfrak a^\centerdot \colon \A(a,a') \otimes (Va \otimes Wb) \to Va' \otimes Wb$. This construction underlies a functor
\begin{equation}
\label{eq:tensor-left-right}
\mathord{\otimes} \colon {}_\A\cat{Mod}_\bullet \times {}_\bullet\cat{Mod}_\B \to {}_\A\cat{Mod}_\B\rlap{ ,}
\end{equation}
which is compatible with the actions of $\V$ on ${}_\A\cat{Mod}_\bullet$ and ${}_\bullet\cat{Mod}_\B$ in the sense that we have, for every $V \colon \A \tor \bullet$, $A \in \V$ and $W \colon \bullet \tor \B$, an equivalence $\mathfrak a \colon (V \otimes A) \otimes W \to V \otimes (A \otimes W)$ of $\A$-$\B$-bimodules with $1$-cell components $\mathfrak a_{Va, A, Wb}$, together with further invertible bimodule transformations $\pi$ and $\nu$ witnessing the coherence of these data.

\subsection{Bimodules via copowers}\label{subsec:bimod-via-copower}
We may use copowers to restate the definition of bimodule purely in terms of one-sided modules. Given a bimodule $M \colon \B \tor \C$, the left $\B$-action maps $m_{bb'c} \colon \B(b,b') \otimes M(c,b) \to M(c,b')$ may be seen as the $1$-cell components of a family of right $\C$-module morphisms $m_{bb'\thg} \colon \B(b,b') \otimes M(\thg, b) \to M(\thg, b')$, where in the domain we are forming the copower of $M(\thg, b)$ by $\B(b,b')$; the $2$-cell components of $m_{bb
\thg}$ are precisely the bimodule compatibility $2$-cells of $M$. In fact, it is easy to see that giving the data of a $\B$-$\C$-bimodule is precisely equivalent to giving a family of right $\C$-modules $M(\thg, b)$, a family of right $\C$-module morphisms $\B(b,b') \otimes M(\thg, b) \to M(\thg, b')$, and a family of right $\C$-module transformations expressing the unit and multiplication constraints required for us to have a left $\B$-module ``in the world of right $\C$-modules''. Of course, we may dually view a $\B$-$\C$-bimodule as a right $\C$-module ``in the world of left $\B$-modules''. An entirely similar argument allows us to restate the definition of bimodule morphism solely in terms of one-sided module morphisms.

\section{Tensor product of modules}\label{sec:tensorproductofmodules}
An important fact about modules between (small) enriched categories is that---under suitable cocompleteness assumptions on the base---they form a bicategory $\V\text-\cat{Mod}$ whose $1$-cell composition is given by \emph{tensor product} of bimodules. There is an analogue of this at the level of enriched bicategories; for a sufficiently cocomplete $\V$, there is a tricategory $\V\text-\cat{Mod}$ whose objects are small $\V$-bicategories, with $\V\text-\cat{Mod}(\B,\C) = {}_\B \cat{Mod}_\C$ and with composition of $1$-cells given by a suitable tensor product. We shall not construct $\V\text-\cat{Mod}$ in its entirety in this paper, but we will need, among other things, the tensor product giving its $1$-cell composition. The tensor product of two bimodules will be characterised as a classifier for ``bilinear maps''; we begin, therefore, by discussing the relevant notion of bilinearity.

\subsection{Module bimorphisms}\label{subsec:modulebimorphisms}
Given $V \colon \bullet \tor \B$, $W \colon \B \tor \bullet$, and $A \in \V$, a \emph{module bimorphism} $\phi \colon V,W \to A$ is given by:
\begin{itemize}
\item For each $x \in \B$, a morphism $\phi_x \colon Vx \otimes Wx \to A$ in $\V$;
\item For each $x, y \in \B$, invertible $2$-cells
\begin{equation*}
\cd[@C-0.5em]{
 (Vy \otimes \B(x,y)) \otimes Wx \ar[rr]^{\mathfrak a} \ar[d]_{m \otimes 1}  & {} \rtwocell[0.55]{d}{\bar \phi_{xy}} &
 Vy \otimes (\B(x,y) \otimes Wx) \ar[d]^{1 \otimes m} \\
 Vx \otimes Wx \ar[r]_-\phi & A & Vy \otimes Wy \ar[l]^-{\phi}
}
\end{equation*}
\end{itemize}
subject to two axioms that as string diagrams, once again take the same shape as the two axioms for a $\V$-bicategory, though now with $\phi$ and $\bar \phi$ replacing some instances of $m$ and $\alpha$.
Given $M \colon \A \tor \B$, $N \colon \B \tor \C$ and $P \colon \A \tor \C$, an \emph{$(\A,\B,\C)$-module bimorphism} $\phi \colon M, N \to P$ is given by $1$-cells $\phi_{abc} \colon M(b, a) \otimes N(c, b) \to P(c,a)$ together with $2$-cells making each $\phi_{ \thg b c}$ into a left $\A$-module morphism,  each $\phi_{a \thg c}$ into a module bimorphism over $\B$, and each $\phi_{ab \thg}$ into a right $\C$-module morphism. There are corresponding notions of module bimorphism in which any of $\A$, $\B$ or $\C$ have been replaced by $\bullet$.

\begin{Ex}\label{ex:yoneda-bimorphisms}
Given $M \colon \A \tor \B$, the action morphisms $M(b,a) \otimes \B(b',b) \to M(b',a)$ are the $1$-cell components of a module bimorphism $m \colon M, \B \to M$, whose $2$-cell components are obtained from the associativity constraints of the action. Similarly, the left action morphisms of $M$ are components of a module bimorphism $m \colon \A, M \to M$.
\end{Ex}
\begin{Ex}\label{ex:pullback-of-bimorphisms}
Let $\phi \colon M, N \to P$ be an ($\A$, $\B$, $\C$)-module bimorphism, and $F \colon \B' \to \B$ a $\V$-functor. There is a module bimorphism $M(F, 1), N(1, F) \to P$ whose $1$-cell components are $\phi_{a,Fb,c}$ and whose $2$-cell components are those of $\phi$ for the $\A$- and $\C$-actions, and are obtained from those of $\phi$ together with pseudonaturality of $\mathfrak a$ for the $\B'$-action. By abuse of notation, we refer to this induced bimorphism also as $\phi$.
\end{Ex}

\subsection{Bimorphism transformations}
Given  module bimorphisms
$\phi, \psi \colon V,W \to A$, with $V$, $W$ and $A$ as before, a \emph{transformation} $\Gamma \colon \phi \Rightarrow \psi$ is given by $2$-cells $\Gamma_x \colon \phi_x \Rightarrow \psi_x$ for each $x \in \A$, subject to the axiom that
\begin{equation*}
\vcenter{\hbox{\begin{tikzpicture}[y=0.80pt, x=0.8pt,yscale=-1, inner sep=0pt, outer sep=0pt, every text node part/.style={font=\tiny} ]
\path[use as bounding box] (0, 983) rectangle (130,1040);
  \path[draw=black,line join=miter,line cap=butt,line width=0.650pt]
  (0.0000,1022.3622) .. controls (20.0000,1022.3622) and (33.4489,1022.0637) ..
  node[above right=0.12cm,at start] {$\!\phi$}(38.7575,1014.9418);
  \path[draw=black,line join=miter,line cap=butt,line width=0.650pt]
  (0.0000,1002.3622) .. controls (20.0000,1002.3622) and (33.6763,1003.7813) ..
  node[above right=0.12cm,at start] {$\!m \t 1$}(38.5075,1009.5873);
  \path[draw=black,line join=miter,line cap=butt,line width=0.650pt]
  (89.1875,1032.3622) .. controls (62.9485,1032.3622) and (48.1250,1025.2372) ..
  node[above right=0.09cm,pos=0.55] {$\phi$}(40.8750,1014.8622);
  \path[draw=black,line join=miter,line cap=butt,line width=0.650pt]
  (130.0000,1012.3622) .. controls (100.5000,1012.3622) and (74.1700,1012.3622) ..
  node[above left=0.12cm,at start] {$1 \t m\!$}(42.2500,1012.3622);
  \path[draw=black,line join=miter,line cap=butt,line width=0.650pt]
  (90.4541,1032.3622) --
  node[above left=0.12cm,at end] {$\psi\!$}(130.0000,1032.3622);
  \path[draw=black,line join=miter,line cap=butt,line width=0.650pt]
  (41.1250,1009.6122) .. controls (53.1250,991.6122) and (72.9933,992.3622) ..
  node[above left=0.12cm,at end] {$\mathfrak a\!$}(130.0000,992.3622);
  \path[xscale=1.005,yscale=0.995,fill=black] (89.022888,1037.5812) node[circle, draw, line width=0.65pt, minimum width=5mm, fill=white, inner sep=0.25mm] (text3092) {$\Gamma$};
  \path[fill=black] (40,1012.1122) node[circle, draw, line width=0.65pt, minimum width=5mm, fill=white, inner sep=0.25mm] (text3086) {$\bar \phi$};
\end{tikzpicture}}}
\qquad \text= \qquad
\vcenter{\hbox{\begin{tikzpicture}[y=0.80pt, x=0.8pt, inner sep=0pt, outer sep=0pt, every text node part/.style={font=\tiny} ]
  \path[use as bounding box] (0, 963) rectangle (130,1020);
  \path[draw=black,line join=miter,line cap=butt,line width=0.650pt]
  (130.0000,992.3622) .. controls (110.0000,992.3622) and (106.0673,992.5498) ..
  node[above left=0.12cm,at start] {$1 \t m\!$}(92.9925,992.5498);
  \path[draw=black,line join=miter,line cap=butt,line width=0.650pt]
  (130.0000,1012.3622) .. controls (110.0000,1012.3622) and (98.1869,1004.3486) ..
  node[above left=0.12cm,at start] {$\mathfrak a\!$}(91.4925,994.9043);
  \path[draw=black,line join=miter,line cap=butt,line width=0.650pt]
  (41.0459,982.1294) .. controls (58.7995,981.8895) and (82.7500,984.1294) ..
  node[above right=0.09cm,pos=0.3] {$\psi$}(89.5000,990.6294);
  \path[draw=black,line join=miter,line cap=butt,line width=0.650pt]
  (0.0000,1002.1294) .. controls (30.0000,1002.1294) and (79.1524,1004.4770) ..
  node[above right=0.12cm,at start] {$\!m \t 1$}(89.2500,994.3794);
  \path[draw=black,line join=miter,line cap=butt,line width=0.650pt]
  (39.5459,982.1294) --
  node[above right=0.12cm,at end] {$\!\phi$}(0.0000,982.1294);
  \path[draw=black,line join=miter,line cap=butt,line width=0.650pt]
  (92.0000,990.1294) .. controls (98.2500,981.3622) and (110.5067,972.3622) ..
  node[above left=0.12cm,at end] {$\psi\!$}(130.0000,972.3622);
  \path[xscale=1.005,yscale=0.995,fill=black] (40.075989,987.28436) node[circle, draw, line width=0.65pt, minimum width=5mm, fill=white, inner sep=0.25mm] (text3092) {$\Gamma$};
  \path[fill=black] (90.003906,993.08643) node[circle, draw, line width=0.65pt, minimum width=5mm, fill=white, inner sep=0.25mm] (text3086) {$\bar \psi$};
\end{tikzpicture}}}
\end{equation*}
in $\V(\, (Vy \otimes \A(x,y)) \otimes Wx,\,A\,)$. More generally, given $M$, $N$ and $P$ as before, a bimodule transformation $\Gamma \colon \phi \Rightarrow \psi \colon M, N \to P$ is given by $2$-cells $\Gamma_{abc} \colon \phi_{abc} \Rightarrow \psi_{abc}$ such that each $\Gamma_{\thg bc}$, each $\Gamma_{a \thg c}$ and each $\Gamma_{ab \thg}$ is a module transformation of the appropriate kind.

\subsection{Compositional structure}
Given $M \colon \A \tor \B$, $N \colon \B \tor \C$ and $P \colon \A \tor \C$, the totality of the module bimorphisms and transformations $M,N \to P$ form a category $\cat{Bimor}_{\A\B\C}(M,N; P)$, and as $M$, $N$ and $P$ vary, these categories constitute the action on objects of a functor
\begin{equation}\label{eq:bimor-functor}
\cat{Bimor}_{\A\B\C}(\thg,\thg; \thg) \colon {}_\A\cat{Mod}_\B^\op \times {}_\B\cat{Mod}_\C^\op \times {}_\A\cat{Mod}_\C \to \cat{CAT}\rlap{ .}
\end{equation}
 At a morphism $(\gamma, \delta, \varepsilon) \colon (M,N,P) \to (M',N',P')$ of ${}_\A\cat{Mod}_\B^\op \times {}_\B\cat{Mod}_\C^\op \times {}_\A\cat{Mod}_\C$, the induced functor $\cat{Bimor}(\gamma,\delta; \varepsilon)$ sends $\phi \colon M, N  \to P$ to the bimorphism $M', N' \to P'$ with $1$-cell components $\varepsilon_{ac} \circ \phi_{abc} \circ (\gamma_{ab} \otimes \delta_{bc})$ and with
compatibility $2$-cells for the $\A$, $\B$ and $\C$-actions obtained from those of $\varphi$ in combination with those of $\gamma$ and $\varepsilon$ (for the $\A$-action), of $\delta$ and $\gamma$ (for the $\B$-action) and of $\delta$ and $\varepsilon$ (for the $\C$-action).
On morphisms, $\cat{Bimor}(\gamma,\delta; \varepsilon)$ sends $\Gamma \colon \phi \Rightarrow \phi'$ to the module transformation with components $\varepsilon_{ac} \circ \Gamma_{abc} \circ (\gamma_{ab} \otimes \delta_{bc})$. 
This defines~\eqref{eq:bimor-functor} on $1$-cells; at a $2$-cell $(\Gamma, \Delta, \Upsilon)$, the
 induced natural transformation $\cat{Bimor}(\Gamma,\Delta; \Upsilon)$ has its component at $\phi$ the bimodule transformation with components $\Upsilon_{ac} \circ \phi_{abc} \circ (\Gamma_{ab} \otimes \Delta_{bc})$. The functoriality constraint $2$-cells of~\eqref{eq:bimor-functor} are obtained from the bicategorical associativity and unit constraints in $\V$ and the functoriality constraints of $\otimes$.

\subsection{Tensor products}\label{subsec:construction-of-tensor-products}
Given $M \colon \A \tor \B$ and $N \colon \B \tor \C$, a \emph{tensor product}  of $M$ and $N$ over $\B$ is a birepresenting element $M, N \to M \otimes_\B N$ for the functor $\cat{Bimor}(M,N; \thg) \colon {}_\A\cat{Mod}_\C \to \cat{CAT}$. In this section, we consider circumstances under which such tensor products may be shown to exist.

Note  that there is one case that we can dispatch immediately: when $\B = \bullet$, the tensor product of a left $\A$-module $M$ and a right $\C$-module $N$ may be taken to be that defined by~\eqref{eq:tensor-left-right}. In particular, when $\A = \B = \bullet$, so that $M$ is simply an object of $\V$, and $N$ a right $\C$-module, a tensor product of $M$ with $N$ is given by the copower $M \otimes N$; and dually when $\B = \C = \bullet$.

We consider next the case of a general $\B$, but with $\A = \C = \bullet$: given a right $\B$-module $V$ and left $\B$-module $W$, we will describe a colimit in $\V$ that, if it exists, must underlie the tensor product $V \otimes_\B W$, and whose existence, conversely, is guaranteed by the existence of the tensor product.

Let $\D$ denote the category whose object set is $\ob \B + (\ob \B)^2 + (\ob \B)^3$, and whose morphisms are generated by arrows
\begin{align*}
i_x \colon  x &\to (x,x) & d_{xy} \colon  (x,y) &\to x & c_{xy} \colon (x,y) &\to y\\
p_{xyz} \colon  (x,y,z) &\to (x,y) & q_{xyz} \colon  (x,y,z) &\to (y,z) & n_{xyz} \colon (x,y,z) &\to (x,z)
\end{align*}
for all $x,y,z$ in $\B$, subject to the simplicial identities $di = 1$, $ci =1$, $cp = dq$, $dp = dm$ and $cm = cq$. We define a functor $F^{VW} \colon \D \to \V$ that on objects is given by
\begin{itemize}
\item $F^{VW}(x) = Vx \otimes Wx$;
\item $F^{VW}(x,y) = (Vy \otimes \B(x,y)) \otimes Wx$;
\item $F^{VW}(x,y,z) = ((Vz \otimes \B(y,z)) \otimes \B(x,y)) \otimes Wx$.
\end{itemize}
and on generating morphisms by
\begin{itemize}
\item $F^{VW}(i_x) = ((1 \otimes j) \otimes 1) \circ (\mathfrak r^\centerdot \otimes 1)$;
\item $F^{VW}(d_{xy}) = m \otimes 1$;
\item $F^{VW}(c_{xy}) = (1 \otimes m) \circ \mathfrak a$;
\item $F^{VW}(p_{xyz}) = ((m \otimes 1) \otimes 1)$;
\item $F^{VW}(n_{xyz}) = ((1 \otimes m) \otimes 1) \circ (\mathfrak a \otimes 1)$;
\item $F^{VW}(q_{xyz}) = (1 \otimes m) \circ \mathfrak a$.
\end{itemize}

To extend these assignations to a functor, it suffices to exhibit functoriality coherence cells with respect to the generating simplicial identities. These are given by the respective string diagrams:
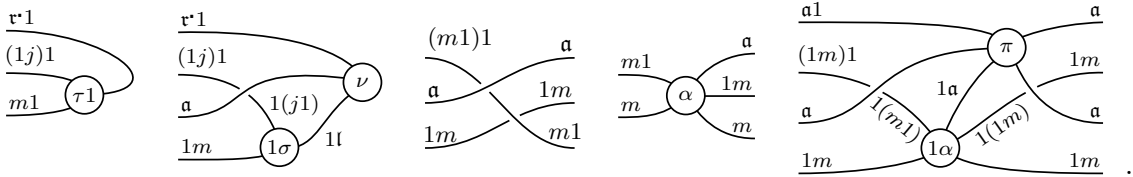
\begin{equation*}
\vcenter{\hbox{
\begin{tikzpicture}[y=0.80pt, x=0.75pt,yscale=-1, inner sep=0pt, outer sep=0pt, every text node part/.style={font=\tiny} ]
\path[use as bounding box] (8, 960) rectangle (85,1039);
  \path[draw=black,line join=miter,line cap=butt,line width=0.650pt]
  (10.0000,1012.3622) .. controls (31.9384,1012.3622) and (44.6754,1009.3804) ..
  node[above right=0.12cm,at start] {\!$m \t 1$}(48.2961,1005.5467);
  \path[draw=black,line join=miter,line cap=butt,line width=0.650pt]
  (10.0000,992.3622) .. controls (30.6605,992.3622) and (45.3144,994.7050) ..
  node[above right=0.12cm,at start] {\!\!$(1 \t j) \t 1$}(48.7221,999.6036);
  \path[draw=black,line join=miter,line cap=butt,line width=0.650pt]
  (10.0000,972.3622) .. controls (61.5891,972.3622) and (99.1240,1001.2238) ..
  node[above right=0.12cm,at start] {\!$\mathfrak r^\centerdot \t 1$}(51.2779,1002.6059);
  \path[fill=black] (49.096386,1002.6634) node[circle, draw, line width=0.65pt, minimum width=5mm, fill=white, inner sep=0.25mm] (text6241) {$\tau \t 1$   };
\end{tikzpicture}}} \
\vcenter{\hbox{
\begin{tikzpicture}[y=0.80pt, x=0.7pt,yscale=-1, inner sep=0pt, outer sep=0pt, every text node part/.style={font=\tiny} ]
  \path[draw=black,line join=miter,line cap=butt,line width=0.650pt]
  (60.0000,942.3622) .. controls (95.0000,941.6479) and (140.3571,944.8622) ..
  node[above right=0.12cm,at start] {\!$\mathfrak r^\centerdot \t 1$}(157.5000,962.0050);
  \path[draw=black,line join=miter,line cap=butt,line width=0.650pt]
  (60.0000,982.3622) .. controls (102.9365,982.3622) and (80.7143,955.2193) ..
  node[above right=0.12cm,at start] {\!$\mathfrak a$}(157.5000,964.8622)
  (60.0000,962.3622) .. controls (74.6167,962.3622) and (85.1196,964.9993) ..
  node[above right=0.12cm,at start] {\!\!$(1 \t j) \t 1$}(92.8564,969.0169)
  (97.6790,971.9277) .. controls (106.7104,978.2283) and (111.0730,986.6095) ..
  node[above right=0.07cm,pos=0.53] {\!$1 \t (j \t 1)$}(114.2857,993.7908);
  \path[draw=black,line join=miter,line cap=butt,line width=0.650pt]
  (60.0000,1002.3622) .. controls (82.6960,1002.3622) and (106.6071,1003.0765) ..
  node[above right=0.12cm,at start] {\!$1 \t m$}(114.6429,997.3622);
  \path[draw=black,line join=miter,line cap=butt,line width=0.650pt]
  (116.4286,997.0050) .. controls (135.7143,997.7193) and (136.0714,979.8622) ..
  node[below right=0.07cm,pos=0.5] {$1 \t \mathfrak l$}(158.2143,967.7193);
  \path[fill=black] (158.57143,965.93359) node[circle, draw, line width=0.65pt, minimum width=5mm, fill=white, inner sep=0.25mm] (text3445) {$\nu$};
  \path[fill=black] (114.1893,996.59198) node[circle, draw, line width=0.65pt, minimum width=5mm, fill=white, inner sep=0.25mm] (text3449) {$1 \t \lu$};
\end{tikzpicture}}}\ \ \
\vcenter{\hbox{
\begin{tikzpicture}[y=0.85pt, x=0.8pt,yscale=-1, inner sep=0pt, outer sep=0pt, every text node part/.style={font=\scriptsize} ]
  \path[draw=black,line join=miter,line cap=butt,line width=0.650pt]
  (40.0000,952.3622) .. controls (50.6250,952.3622) and (58.7659,957.3795) ..
  node[above right=0.12cm,at start] {\!$(m \t 1) \t 1$}(66.1736,963.8603)
  (70.1584,967.5252) .. controls (81.9745,978.8442) and (92.6741,992.3622) ..
  node[above left=0.16cm,at end] {$m \t 1$\!\!\!\!}(110.0000,992.3622)
  (40.0000,972.3622) .. controls (70.0000,972.3622) and (80.0000,952.3622) ..
  node[above right=0.12cm,at start] {\!$\mathfrak a$}
  node[above left=0.12cm,at end] {$\mathfrak a$\!}(110.0000,952.3622)
  (40.0000,992.3622) .. controls (57.2229,992.3622) and (67.8540,985.7705) ..
  node[above right=0.12cm,at start] {\!\!$1 \t m$}(79.4620,980.1556)
  (84.3780,977.8918) .. controls (91.6284,974.7583) and (99.6161,972.3622) ..
  node[above left=0.12cm,at end] {$1 \t m$\!}(110.0000,972.3622);
\end{tikzpicture}}}\ \ \
\vcenter{\hbox{
\begin{tikzpicture}[y=0.80pt, x=0.8pt,yscale=-1, inner sep=0pt, outer sep=0pt, every text node part/.style={font=\tiny} ]
  \path[draw=black,line join=miter,line cap=butt,line width=0.650pt]
  (122.0000,962.3622) .. controls (104.2362,962.3622) and (93.5355,974.6350) ..
  node[above left=0.12cm,at start] {$\mathfrak a$\!}(91.7678,979.4333);
  \path[draw=black,line join=miter,line cap=butt,line width=0.650pt]
  (122.0000,982.3622) --
  node[above left=0.12cm,at start] {$1 \t m$\!}(92.7779,982.3622);
  \path[draw=black,line join=miter,line cap=butt,line width=0.650pt]
  (122.0000,1002.3622) .. controls (103.4643,1002.3622) and (95.4296,990.3419) ..
  node[above left=0.12cm,at start] {$m$\!}(91.6415,985.0386);
  \path[draw=black,line join=miter,line cap=butt,line width=0.650pt]
  (58.0000,972.3622) .. controls (67.6795,972.3622) and (80.2525,972.1096) ..
  node[above right=0.12cm,at start] {\!$m \t 1$}(88.3338,980.4434);
  \path[draw=black,line join=miter,line cap=butt,line width=0.650pt]
  (58.0000,992.3622) .. controls (66.6752,992.3622) and (80.2525,992.6147) ..
  node[above right=0.12cm,at start] {\!$m$}(87.8287,983.5234);
  \path[fill=black] (89.651039,982.15656) node[circle, draw, line width=0.65pt, minimum width=5mm, fill=white, inner sep=0.25mm] (text6017) {$\ass$ };
\end{tikzpicture}}}\ \ \
\vcenter{\hbox{
\begin{tikzpicture}[y=0.95pt, x=0.95pt,yscale=-1, inner sep=0pt, outer sep=0pt, every text node part/.style={font=\tiny} ]
  \path[draw=black,line join=miter,line cap=butt,line width=0.650pt]
  (70.0000,922.3622) .. controls (101.2330,922.3622) and (143.6145,921.7598) ..
  node[above right=0.12cm,at start] {\!$\mathfrak a \t 1$}(151.3870,930.2041);
  \path[draw=black,line join=miter,line cap=butt,line width=0.650pt]
  (126.7723,969.6175) .. controls (128.2875,962.0413) and (140.0945,940.6924) ..
  node[above left=0.08cm,pos=0.47] {$1 \t \mathfrak a$}(151.2391,934.9696);
  \path[draw=black,line join=miter,line cap=butt,line width=0.650pt]
  (127.6786,973.3558) .. controls (136.1123,980.1329) and (156.4071,982.3622) ..
  node[above left=0.12cm,at end] {$1 \t m$\!}(190.0000,982.3622);
  \path[draw=black,line join=miter,line cap=butt,line width=0.650pt]
  (152.7723,930.8750) .. controls (155.6057,928.5814) and (169.8180,922.9558) ..
  node[above left=0.12cm,at end] {$\mathfrak a$\!}(190.0000,922.3622);
  \path[draw=black,line join=miter,line cap=butt,line width=0.650pt]
  (153.1294,933.7322) .. controls (155.4312,936.0339) and (161.1229,962.3622) ..
  node[above left=0.12cm,at end] {$\mathfrak a$\!}(190.0000,962.3622)
  (127.8878,970.8146) .. controls (136.2744,968.2835) and (145.9898,957.8154) ..
  node[below right=0.05cm,pos=0.35,rotate=35] {$1 \t (1 \t m)$}(159.5167,950.3078)
  (164.0504,947.9860) .. controls (171.4928,944.4990) and (180.0280,942.1210) ..
  node[above left=0.12cm,at end] {$1 \t m$\!}(190.0000,942.3622);
  \path[draw=black,line join=miter,line cap=butt,line width=0.650pt]
  (70.0000,982.3622) .. controls (94.7819,982.3622) and (118.6172,978.5843) ..
  node[above right=0.12cm,at start] {\!$1 \t m$}(125.9408,972.7759);
  \path[draw=black,line join=miter,line cap=butt,line width=0.650pt]
  (70.0000,962.3622) .. controls (104.0413,962.3622) and (95.4124,932.6539) ..
  node[above right=0.12cm,at start] {\!$\mathfrak a$}(151.2919,932.6539)
  (70.0000,942.3622) .. controls (81.3043,942.3622) and (90.5727,944.4985) ..
  node[above right=0.12cm,at start] {\!\!$(1 \t m) \t 1$}(98.2107,947.7957)
  (103.1302,950.1901) .. controls (113.3893,955.7771) and (120.2698,963.4297) ..
  node[below left=0.05cm,pos=0.74,rotate=-40] {$1 \t (m \t 1)$}(125.0912,969.9746);
  \path[fill=black] (152.0495,932.40131) node[circle, draw, line width=0.65pt, minimum width=5mm, fill=white, inner sep=0.25mm] (text3313) {$\pi$     };
  \path[fill=black] (125.89511,972.28424) node[circle, draw, line width=0.65pt, minimum width=5mm, fill=white, inner sep=0.25mm] (text3317) {$1 \t \ass$   };
\end{tikzpicture} \rlap{ .}}}
\end{equation*}

\begin{Prop}\label{prop:tensor-product-from-bicolimit}
With notation as above, the tensor product $V \otimes_\B W$ and the  conical bicolimit $\Delta 1 \star F^{VW}$ represent pseudonaturally equivalent functors $\V \to \cat{CAT}$; in particular, the one exists if and only if the other does.
\end{Prop}
\begin{proof}
For brevity, we write $F^{VW}$ simply as $F$ for the duration of the proof. We must exhibit a pseudonatural correspondence between transformations $\theta \colon \Delta 1 \to \V(F, A)$, and module bimorphisms $V,W \to A$. Given a transformation $\theta$, its components at $x \in \B$ and at $(x,y) \in \B^\mathbf 2$ pick out morphisms $\phi_x \colon Fx \to A$ and $\psi_{xy} \colon F(x,y) \to A$; while the pseudonaturality constraints at the maps $d_{xy}$ and $c_{xy}$ of $\D$ pick out invertible $2$-cells $\delta_{xy} \colon \phi_x \circ Fd_{xy} \Rightarrow \psi_{xy}$ and $\gamma_{xy} \colon \phi_x \circ Fc_{xy} \Rightarrow \psi_{xy}$.
Now by replacing $\psi_{xy}$ by $\phi_{xy} \circ Fd_{xy}$, replacing $\delta_{xy}$ by the identity $2$-cell, and replacing $\gamma_{xy}$ by $\delta_{xy}^{-1} \circ \gamma_{xy}$, we obtain $\theta' \colon \Delta 1 \to \V(F,A)$, isomorphic to $\theta$, with the property that its pseudonaturality component at $d_{xy}$ is an identity. A similar transport of structure argument shows that we may replace $\theta'$ with an isomorphic $\theta''$ whose pseudonaturality component at $p_{xyz}$ is also an identity.

Consider now the full subcategory $\E$ of $\cat{Hom}(\Delta 1, \V(F,A))$ spanned by objects of the form $\theta''$. By the preceding construction, the inclusion of this subcategory is an equivalence; and it is easy to see that any $\theta''$ is completely determined by the $1$-cells $\phi_x \colon Vx \otimes Wx \to A$ picked out by the component at $x \in \B$ together with the $2$-cells $\bar \phi_{xy} \colon \phi_x \circ (m \otimes 1) \Rightarrow \phi_x \circ (1 \otimes m) \circ \mathfrak a$ picked out by the pseudonaturality constraint at $c_{xy}$. We may now verify that such data represent an object of the subcategory just when they satisfy the bimorphism axioms. We may argue similarly to show that modifications $\theta'' \to \omega''$ correspond precisely to transformations between the associated bimorphisms. Consequently, the full subcategory $\E$ is  isomorphic to $\cat{Bimor}(V,W;A)$, and we therefore have, for each $A \in \V$, an injective equivalence of categories
\begin{equation}\label{eq:functor-bimodules-to-cocones}
\cat{Bimor}(V,W;A) \rightarrow \cat{Hom}(\Delta 1, \V(F,A))\ \text,
\end{equation}
and these are easily verified to be pseudonatural in $A$, as required.\end{proof}

\begin{Cor}
If $\B$ is a small $\V$-bicategory, and $\V$ is cocomplete, then the tensor product $V \otimes_\B W$ exists for every right $\B$-module $V$ and left $\B$-module $W$.
\end{Cor}

We now turn to the construction of the tensor product of a general pair of bimodules $M \colon \A \tor \B$ and $N \colon \B \tor \C$. Just as before, the existence of such may be reduced to the existence of certain bicolimits, though now in the bicategory ${}_\A\cat{Mod}_\C$. To this end, let $\D$ be the category considered before, and define $F^{MN} \colon \D \to {}_\A\cat{Mod}_\C$ by
\begin{itemize}
\item $F^{MN}(x) = M(x, \thg) \otimes N(\thg, x)$;
\item $F^{MN}(x,y) = (M(y, \thg) \otimes \B(x,y)) \otimes N(\thg, x)$;
\item $F^{MN}(x,y,z) = ((M(z, \thg) \otimes \B(y,z)) \otimes \B(x,y)) \otimes N(\thg, x)$
\end{itemize}
(where the tensor products used are now those of Section~\ref{subsec:copowers-of-modules})
and with remaining data as before. Bearing in mind the remarks of Section~\ref{subsec:bimod-via-copower}, the proof of Proposition~\ref{prop:tensor-product-from-bicolimit} now adapts immediately to yield:
\begin{Prop}
The tensor product $M \otimes_\B N$ and the bicolimit $\Delta 1 \star F^{MN}$ represent pseudonaturally equivalent functors ${}_\A\cat{Mod}_\C \to \cat{CAT}$; in particular, the one exists if and only if the other does.
\end{Prop}
\begin{Cor}\label{prop:tensorproductwithbimodule}
Let $\V$ be cocomplete, and let each functor $\thg \otimes X$ and $X \otimes \thg \colon \V \to \V$ preserve bicolimits. Given bimodules $M \colon \A \tor \B$ and $N \colon \B \tor \C$, with $\B$ small,
the tensor product $M \otimes_\B N$ exists, and is pointwise in the sense that $(M \otimes_\B N)(c,a) = M(\thg, a) \otimes_\B N(c, \thg)$.
\end{Cor}
\begin{proof}
Because each $\thg \otimes X$ and $X \otimes \thg$ preserves bicolimits, the bicategory ${}_\A\cat{Mod}_\C$ admits bicolimits created by the forgetful functor ${}_\A\cat{Mod}_\C \to \V^{\ob \A \times \ob \C}$.
\end{proof}

\subsection{Compatibility with copowers}
Let $V \colon \bullet \tor \B$, $M \colon \B \tor \C$ and $W \colon \bullet \tor \C$. For any
 bimorphism $\phi \colon V,M \to W$ and any $B \in \V$, we may construct a bimorphism $B \otimes \phi \colon B \otimes V,\, M \to B \otimes W$ with $1$-cell components $(1 \otimes \phi_{bc}) \circ \mathfrak a \colon (B \otimes Vb) \otimes M(c,b) \to B \otimes Wc$ and with $2$-cell components for the $\B$-actions and the $\C$-actions given by
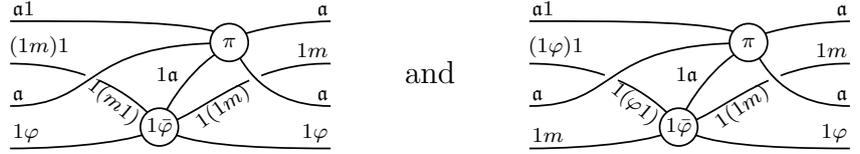
\begin{equation*}
\vcenter{\hbox{
\begin{tikzpicture}[y=0.80pt, x=1pt,yscale=-1, inner sep=0pt, outer sep=0pt, every text node part/.style={font=\tiny} ]
  \path[draw=black,line join=miter,line cap=butt,line width=0.650pt]
  (70.0000,922.3622) .. controls (101.2330,922.3622) and (143.6145,921.7598) ..
  node[above right=0.12cm,at start] {\!$\mathfrak a \t 1$}(151.3870,930.2041);
  \path[draw=black,line join=miter,line cap=butt,line width=0.650pt]
  (126.7723,969.6175) .. controls (128.2875,962.0413) and (140.0945,940.6924) ..
  node[above left=0.08cm,pos=0.47] {$1 \t \mathfrak a$}(151.2391,934.9696);
  \path[draw=black,line join=miter,line cap=butt,line width=0.650pt]
  (127.6786,973.3558) .. controls (136.1123,980.1329) and (156.4071,982.3622) ..
  node[above left=0.12cm,at end] {$1 \t \phi$\!}(190.0000,982.3622);
  \path[draw=black,line join=miter,line cap=butt,line width=0.650pt]
  (152.7723,930.8750) .. controls (155.6057,928.5814) and (169.8180,922.9558) ..
  node[above left=0.12cm,at end] {$\mathfrak a$\!}(190.0000,922.3622);
  \path[draw=black,line join=miter,line cap=butt,line width=0.650pt]
  (153.1294,933.7322) .. controls (155.4312,936.0339) and (161.1229,962.3622) ..
  node[above left=0.12cm,at end] {$\mathfrak a$\!}(190.0000,962.3622)
  (127.8878,970.8146) .. controls (136.2744,968.2835) and (145.9898,957.8154) ..
  node[below right=0.05cm,pos=0.35,rotate=35] {$1 \t (1 \t m)$}(159.5167,950.3078)
  (164.0504,947.9860) .. controls (171.4928,944.4990) and (180.0280,942.1210) ..
  node[above left=0.12cm,at end] {$1 \t m$\!}(190.0000,942.3622);
  \path[draw=black,line join=miter,line cap=butt,line width=0.650pt]
  (70.0000,982.3622) .. controls (94.7819,982.3622) and (118.6172,978.5843) ..
  node[above right=0.12cm,at start] {\!$1 \t \phi$}(125.9408,972.7759);
  \path[draw=black,line join=miter,line cap=butt,line width=0.650pt]
  (70.0000,962.3622) .. controls (104.0413,962.3622) and (95.4124,932.6539) ..
  node[above right=0.12cm,at start] {\!$\mathfrak a$}(151.2919,932.6539)
  (70.0000,942.3622) .. controls (81.3043,942.3622) and (90.5727,944.4985) ..
  node[above right=0.12cm,at start] {\!\!$(1 \t m) \t 1$}(98.2107,947.7957)
  (103.1302,950.1901) .. controls (113.3893,955.7771) and (120.2698,963.4297) ..
  node[below left=0.05cm,pos=0.74,rotate=-40] {$1 \t (m \t 1)$}(125.0912,969.9746);
  \path[fill=black] (152.0495,932.40131) node[circle, draw, line width=0.65pt, minimum width=5mm, fill=white, inner sep=0.25mm] (text3313) {$\pi$     };
  \path[fill=black] (125.89511,972.28424) node[circle, draw, line width=0.65pt, minimum width=5mm, fill=white, inner sep=0.25mm] (text3317) {$1 \t \bar \phi$   };
\end{tikzpicture}\
}}\qquad \text{and} \qquad
\vcenter{\hbox{
\begin{tikzpicture}[y=0.80pt, x=1pt,yscale=-1, inner sep=0pt, outer sep=0pt, every text node part/.style={font=\tiny} ]
  \path[draw=black,line join=miter,line cap=butt,line width=0.650pt]
  (70.0000,922.3622) .. controls (101.2330,922.3622) and (143.6145,921.7598) ..
  node[above right=0.12cm,at start] {\!$\mathfrak a \t 1$}(151.3870,930.2041);
  \path[draw=black,line join=miter,line cap=butt,line width=0.650pt]
  (126.7723,969.6175) .. controls (128.2875,962.0413) and (140.0945,940.6924) ..
  node[above left=0.08cm,pos=0.47] {$1 \t \mathfrak a$}(151.2391,934.9696);
  \path[draw=black,line join=miter,line cap=butt,line width=0.650pt]
  (127.6786,973.3558) .. controls (136.1123,980.1329) and (156.4071,982.3622) ..
  node[above left=0.12cm,at end] {$1 \t \phi$\!}(190.0000,982.3622);
  \path[draw=black,line join=miter,line cap=butt,line width=0.650pt]
  (152.7723,930.8750) .. controls (155.6057,928.5814) and (169.8180,922.9558) ..
  node[above left=0.12cm,at end] {$\mathfrak a$\!}(190.0000,922.3622);
  \path[draw=black,line join=miter,line cap=butt,line width=0.650pt]
  (153.1294,933.7322) .. controls (155.4312,936.0339) and (161.1229,962.3622) ..
  node[above left=0.12cm,at end] {$\mathfrak a$\!}(190.0000,962.3622)
  (127.8878,970.8146) .. controls (136.2744,968.2835) and (145.9898,957.8154) ..
  node[below right=0.05cm,pos=0.35,rotate=35] {$1 \t (1 \t m)$}(159.5167,950.3078)
  (164.0504,947.9860) .. controls (171.4928,944.4990) and (180.0280,942.1210) ..
  node[above left=0.12cm,at end] {$1 \t m$\!}(190.0000,942.3622);
  \path[draw=black,line join=miter,line cap=butt,line width=0.650pt]
  (70.0000,982.3622) .. controls (94.7819,982.3622) and (118.6172,978.5843) ..
  node[above right=0.12cm,at start] {\!$1 \t m$}(125.9408,972.7759);
  \path[draw=black,line join=miter,line cap=butt,line width=0.650pt]
  (70.0000,962.3622) .. controls (104.0413,962.3622) and (95.4124,932.6539) ..
  node[above right=0.12cm,at start] {\!$\mathfrak a$}(151.2919,932.6539)
  (70.0000,942.3622) .. controls (81.3043,942.3622) and (90.5727,944.4985) ..
  node[above right=0.12cm,at start] {\!\!$(1 \t \phi) \t 1$}(98.2107,947.7957)
  (103.1302,950.1901) .. controls (113.3893,955.7771) and (120.2698,963.4297) ..
  node[below left=0.05cm,pos=0.74,rotate=-40] {$1 \t (\phi \t 1)$}(125.0912,969.9746);
  \path[fill=black] (152.0495,932.40131) node[circle, draw, line width=0.65pt, minimum width=5mm, fill=white, inner sep=0.25mm] (text3313) {$\pi$     };
  \path[fill=black] (125.89511,972.28424) node[circle, draw, line width=0.65pt, minimum width=5mm, fill=white, inner sep=0.25mm] (text3317) {$1 \t \bar \phi$   };
\end{tikzpicture}\rlap{ .}
}}
\end{equation*}
\begin{Prop}\label{prop:copowers-of-universal}
Suppose that $B \otimes \thg$ and each $\thg \otimes X \colon \V \to \V$ preserve bicolimits. If $\phi \colon V, M \to W$ is a universal bimorphism, then so too is $B \otimes \phi$; which is equally to say that if $V \otimes_\A M$ exists with universal morphism $\phi$, then so does $(B \otimes V) \otimes_\A M$, and the canonical morphism $(B \otimes V) \otimes_\A M \to B \otimes (V \otimes_\A M)$ induced by $B \otimes \phi$ is an equivalence.
\end{Prop}
Observe that this result is really constructing a simple instance of an associativity constraint in the tricategory of $\V$-bimodules.
\begin{proof}
There is a pseudonatural equivalence $\gamma \colon F^{B \otimes V,M} \Rightarrow B \otimes F^{VM} \colon \D \to {}_\bullet \cat{Mod}_\C$ with $1$-cell components
\begin{align*}
\gamma_{x} &= \mathfrak a\ \text, & \gamma_{(x,y)} &= \mathfrak a \circ (\mathfrak a \otimes 1)\ \text,\ & \
\gamma_{(x,y,z)} &=\mathfrak a \circ (\mathfrak a \otimes 1) \circ ((\mathfrak a \otimes 1) \otimes 1)
\end{align*}
and with pseudonaturality $2$-cells given at the generating morphisms $i_x$, $d_{xy}$, $c_{xy}$, $p_{xyz}$, $m_{xyz}$ and $q_{xyz}$ of $\D$ by the respective composites
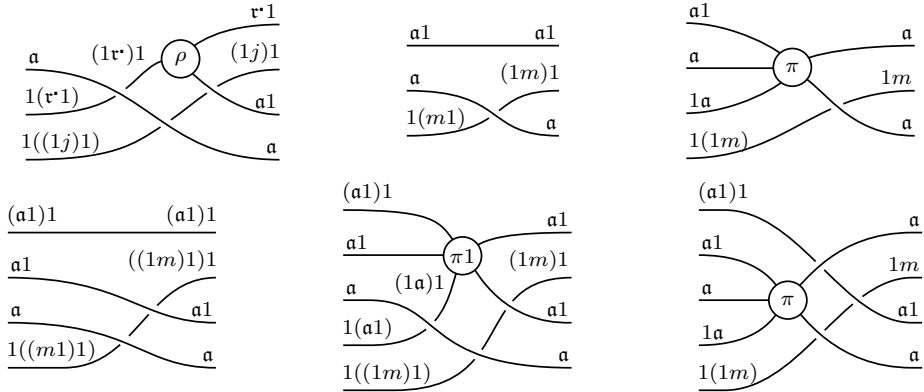
\begin{equation*}
\vcenter{\hbox{
\begin{tikzpicture}[y=0.85pt, x=0.95pt,yscale=-1, xscale=-1,inner sep=0pt, outer sep=0pt, every text node part/.style={font=\tiny} ]
  \path[draw=black,line join=miter,line cap=butt,line width=0.650pt]
  (10.0000,962.3622) .. controls (33.2390,962.3622) and (43.6865,970.6959) ..
  node[above left=0.12cm,at start] {$\mathfrak r^\centerdot \t 1$\!}(48.2322,976.2518);
  \path[draw=black,line join=miter,line cap=butt,line width=0.650pt]
  (10.0000,1002.3622) .. controls (27.1784,1002.3622) and (41.9188,987.2619) ..
  node[above left=0.12cm,at start] {$\mathfrak a \t 1$\!}(48.7373,979.1807)
  (110.0000,1002.3622) .. controls (90.3427,1002.3622) and (81.0837,998.6623) ..
  node[above right=0.12cm,at start] {\!\!$1 \t (\mathfrak r^\centerdot \t 1)$}(74.6653,994.1266)
  (70.5818,990.8430) .. controls (63.3802,984.3420) and (59.9087,977.5329) ..
  node[above left=0.05cm] {$(1 \t \mathfrak r^\centerdot) \t 1$}(51.1612,977.6655)
  (10.0000,982.3622) .. controls (20.4192,982.3622) and (27.8791,985.3848) ..
  node[above left=0.12cm,at start] {$(1 \t j) \t 1$\!\!}(34.4523,989.7682)
  (38.6409,992.8076) .. controls (43.3111,996.4463) and (47.7408,1000.6365) .. (52.7802,1004.6966)
  (56.8775,1007.8455) .. controls (67.6864,1015.7181) and (75.8091,1022.3622) ..
  node[above right=0.12cm,at end] {\!\!\!$1 \t ((1 \t j) \t 1)$}(110.0000,1022.3622)
  (10.0000,1022.3622) .. controls (59.0481,1022.3622) and (61.5490,982.3622) ..
  node[above left=0.12cm,at start] {$\mathfrak a$\!}
  node[above right=0.12cm,at end] {\!$\mathfrak a$}(110.0000,982.3622);
  \path[fill=black] (48.992397,977.6109) node[circle, draw, line width=0.65pt, minimum width=5mm, fill=white, inner sep=0.25mm] (text3023) {$\rho$
   };
\end{tikzpicture}}}\qquad\qquad
\vcenter{\hbox{\begin{tikzpicture}[y=0.85pt, x=0.95pt,yscale=-1, xscale=-1, inner sep=0pt, outer sep=0pt, every text node part/.style={font=\tiny} ]
  \path[draw=black,line join=miter,line cap=butt,line width=0.650pt]
  (10.0000,1022.3622) .. controls (40.0000,1022.3622) and (29.8985,1002.3622) ..
  node[above left=0.12cm,at start] {$\mathfrak a$\!}
  node[above right=0.12cm,at end] {\!$\mathfrak a $}(70.0000,1002.3622)
  (10.0000,1002.3622) .. controls (23.8456,1002.3622) and (28.4503,1006.6222) ..
  node[above left=0.12cm,at start] {$(1 \t m) \t 1$\!\!}(33.6595,1011.2101)
  (37.5839,1014.4657) .. controls (43.2082,1018.6743) and (51.3701,1022.3622) ..
  node[above right=0.12cm,at end] {\!\!$1 \t (m \t 1)$}(70.0000,1022.3622);
  \path[draw=black,line join=miter,line cap=butt,line width=0.650pt]
  (10.0000,982.3622) --
  node[above left=0.12cm,at start] {$\mathfrak a \t 1$\!}
  node[above right=0.12cm,at end] {\!$\mathfrak a \t 1$}(70.0000,982.3622);
\end{tikzpicture}}}\qquad\qquad
\vcenter{\hbox{\begin{tikzpicture}[y=0.85pt, x=0.95pt,yscale=-1, inner sep=0pt, outer sep=0pt, every text node part/.style={font=\tiny} ]
  \path[draw=black,line join=miter,line cap=butt,line width=0.650pt]
  (90.0000,992.3622) .. controls (69.2919,992.3622) and (47.3236,995.6858) ..
  node[above left=0.12cm,at start] {$\mathfrak a$\!}(42.7779,1000.9891);
  \path[draw=black,line join=miter,line cap=butt,line width=0.650pt]
  (90.0000,1012.3622) .. controls (78.6264,1012.3622) and (69.7021,1014.7510) ..
  node[above left=0.12cm,at start] {$1 \t m$\!}(61.8484,1018.1805)
  (57.0164,1020.4728) .. controls (39.6011,1029.3577) and (26.4331,1042.3622) ..
  node[above right=0.12cm,at end] {\!$1 \t (1 \t m)$}(0.0000,1042.3622)
  (90.0000,1032.3622) .. controls (62.4722,1032.3622) and (54.6472,1009.7264) ..
  node[above left=0.12cm,at start] {$\mathfrak a$\!}(42.5254,1003.9180);
  \path[draw=black,line join=miter,line cap=butt,line width=0.650pt]
  (0.0000,982.3622) .. controls (17.9320,982.3622) and (38.6085,994.5716) ..
  node[above right=0.12cm,at start] {\!$\mathfrak a \t 1$}(40.6288,999.8749);
  \path[draw=black,line join=miter,line cap=butt,line width=0.650pt]
  (0.0000,1022.3622) .. controls (18.9824,1022.3622) and (39.3687,1009.4580) ..
  node[above right=0.12cm,at start] {\!$1 \t \mathfrak a$}(40.6313,1004.9123);
  \path[draw=black,line join=miter,line cap=butt,line width=0.650pt]
  (40.0000,1002.3622) --
  node[above right=0.12cm,at end] {\!$\mathfrak a$}(0.0000,1002.3622);
  \path[xscale=-1.000,yscale=1.000,fill=black] (-41.785713,1002.005) node[circle, draw, line width=0.65pt, minimum width=5mm, fill=white, inner sep=0.25mm] (text3250) {$\pi$   };
\end{tikzpicture}}}
\end{equation*}
\begin{equation*}
\vcenter{\hbox{\begin{tikzpicture}[y=0.85pt, x=0.95pt,yscale=-1, xscale=-1, inner sep=0pt, outer sep=0pt, every text node part/.style={font=\tiny} ]
  \path[draw=black,line join=miter,line cap=butt,line width=0.650pt]
  (10.0000,972.3622) --
  node[above left=0.12cm,at start] {$(\mathfrak a \t 1) \t 1$\!\!}
  node[above right=0.12cm,at end] {\!\!$(\mathfrak a \t 1) \t 1$}(92.0000,972.3622);
  \path[draw=black,line join=miter,line cap=butt,line width=0.650pt]
  (10.0000,1032.3622) .. controls (40.0000,1032.3622) and (41.7857,1012.3622) ..
  node[above left=0.12cm,at start] {$\mathfrak a$\!}
  node[above right=0.12cm,at end] {\!$\mathfrak a$}(92.0000,1012.3622)
  (10.0000,1012.3622) .. controls (40.0000,1012.3622) and (45.3571,992.3622) ..
  node[above left=0.12cm,at start] {$\mathfrak a \t 1$\!}
  node[above right=0.12cm,at end] {\!$\mathfrak a \t 1$}(92.0000,992.3622)
  (10.0000,992.3622) .. controls (20.1300,992.3622) and (27.0170,997.8321) ..
  node[above left=0.12cm,at start] {$((1 \t m) \t 1) \t 1$\!\!}(33.1370,1004.7265)
  (36.4760,1008.6803) .. controls (39.0559,1011.8538) and (41.5772,1015.1555) .. (44.2450,1018.2506)
  (47.7488,1022.0668) .. controls (53.5777,1027.9552) and (60.3978,1032.3622) ..
  (70.0000,1032.3622) .. controls (72.7622,1032.3622) and (76.7927,1032.3622) ..
  node[above right=0.12cm,at end] {\!\!$1 \t ((m \t 1) \t 1)$}(92.0000,1032.3622);
\end{tikzpicture}}}\qquad\qquad
\vcenter{\hbox{\begin{tikzpicture}[y=0.85pt, x=0.95pt,yscale=-1, inner sep=0pt, outer sep=0pt, every text node part/.style={font=\tiny} ]
  \path[draw=black,line join=miter,line cap=butt,line width=0.650pt]
  (120.0000,972.3622) .. controls (97.8543,972.3622) and (84.1557,974.0792) ..
  node[above left=0.12cm,at start] {$\mathfrak a \t 1$\!}(79.1557,981.5792);
  \path[draw=black,line join=miter,line cap=butt,line width=0.650pt]
  (76.1940,980.2373) .. controls (68.2069,968.5921) and (70.0000,962.3622) ..
  node[above right=0.12cm,at end] {\!\!$(\mathfrak a \t 1) \t 1$}(30.0000,962.3622);
  \path[draw=black,line join=miter,line cap=butt,line width=0.650pt]
  (120.0000,1012.3622) .. controls (93.2143,1013.0765) and (83.1889,990.9695) ..
  node[above left=0.12cm,at start] {$\mathfrak a \t 1$\!}(79.2603,984.5409)
  (120.0000,1032.3622) .. controls (60.0000,1032.3622) and (70.0000,1002.3622) ..
  node[above left=0.12cm,at start] {$\mathfrak a$\!}
  (40.0000,1002.3622) .. controls (36.9697,1002.3622) and (33.4171,1002.3622) ..
  node[above right=0.12cm,at end] {\!$\mathfrak a$}(30.0000,1002.3622)
  (120.0000,992.3622) .. controls (108.1726,992.3622) and (101.5178,996.6735) ..
  node[above left=0.12cm,at start] {$(1 \t m) \t 1$\!}(96.7661,1002.7641)
  (93.9430,1006.8970) .. controls (90.6368,1012.3470) and (88.1179,1018.6068) .. (84.6550,1024.3356)
  (81.8130,1028.5658) .. controls (75.9156,1036.3860) and (67.0899,1042.3622) ..
  (50.0000,1042.3622) .. controls (48.2816,1042.3622) and (31.8506,1042.3622) ..
  node[above right=0.12cm,at end] {\!\!$1 \t ((1 \t m) \t 1)$}(30.0000,1042.3622)
  (76.4898,984.5410) .. controls (73.6966,989.8482) and (73.1193,1002.3872) ..
  node[above left=0.06cm,pos=0.65] {$(1 \t \mathfrak a) \t 1$}(66.0703,1011.5793)
  (62.6702,1015.2453) .. controls (57.7971,1019.5713) and (50.6923,1022.4939) ..
  (40.0000,1022.3622) .. controls (37.8866,1022.3362) and (32.3339,1022.3622) ..
  node[above right=0.12cm,at end] {\!\!$1 \t (\mathfrak a \t 1)$}(30.0000,1022.3622);
  \path[draw=black,line join=miter,line cap=butt,line width=0.650pt]
  (75.0508,982.3622) --
  node[above right=0.12cm,at end] {\!$\mathfrak a \t 1$}(30.0000,982.3622);
  \path[fill=black] (77.024132,982.66168) node[circle, draw, line width=0.65pt, minimum width=5mm, fill=white, inner sep=0.25mm] (text3440) {$\pi \t 1$   };
\end{tikzpicture}}}\qquad\qquad
\vcenter{\hbox{\begin{tikzpicture}[y=0.85pt, x=0.95pt,yscale=-1, inner sep=0pt, outer sep=0pt, every text node part/.style={font=\tiny} ]
  \path[draw=black,line join=miter,line cap=butt,line width=0.650pt]
  (16.0000,982.3622) .. controls (40.0000,982.3622) and (45.5357,969.8622) ..
  node[above right=0.12cm,at start] {\!$1 \t \mathfrak a$}(48.9286,965.2193);
  \path[draw=black,line join=miter,line cap=butt,line width=0.650pt]
  (16.0000,962.3622) --
  node[above right=0.12cm,at start] {\!$\mathfrak a$}(50.0000,962.3622);
  \path[draw=black,line join=miter,line cap=butt,line width=0.650pt]
  (16.0000,942.3622) .. controls (40.0000,942.3622) and (46.2500,955.2193) ..
  node[above right=0.12cm,at start] {\!$\mathfrak a \t 1$}(48.9286,959.5050);
  \path[draw=black,line join=miter,line cap=butt,line width=0.650pt]
  (16.0000,1002.3622) .. controls (19.6372,1002.3622) and (22.6962,1002.3622) ..
  node[above right=0.12cm,at start] {\!\!$1 \t (1 \t m)$}
  (26.0000,1002.3622) .. controls (40.1871,1002.3622) and (51.4373,991.1803) .. (62.2251,979.3925)
  (65.5914,975.6929) .. controls (69.2893,971.6220) and (72.9630,967.6104) .. (76.7139,964.0925)
  (80.4917,960.7574) .. controls (86.5743,955.7479) and (92.9396,952.3622) ..
  (100.0000,952.3622) .. controls (101.0482,952.3622) and (103.0242,952.3622) ..
  node[above left=0.12cm,at end] {$1 \t m$\!}(104.0000,952.3622)
  (16.0000,922.3622) .. controls (19.4080,922.3622) and (22.8932,922.3622) ..
  node[above right=0.12cm,at start] {\!\!$(\mathfrak a \t 1) \t 1$}
  (26.0000,922.3622) .. controls (40.1341,922.3622) and (51.3122,933.4606) .. (62.0390,945.1998)
  (65.4020,948.9024) .. controls (76.2500,960.8685) and (86.9027,972.3622) ..
  (100.0000,972.3622) .. controls (101.0125,972.3622) and (103.0548,972.3622) ..
  node[above left=0.12cm,at end] {$\mathfrak a \t 1$\!}(104.0000,972.3622)
  (52.1429,959.8622) .. controls (58.1342,954.3844) and (70.0000,932.3622) ..
  (100.0000,932.3622) .. controls (100.8352,932.3622) and (103.1722,932.3622) ..
  node[above left=0.12cm,at end] {$\mathfrak a$\!}(104.0000,932.3622)
  (52.1429,963.7908) .. controls (56.9060,968.8941) and (70.0000,992.3622) ..
  (100.0000,992.3622) .. controls (100.9634,992.3622) and (103.0524,992.3622) ..
  node[above left=0.12cm,at end] {$\mathfrak a$\!}(104.0000,992.3622);
  \path[fill=black] (51.07143,962.36218) node[circle, draw, line width=0.65pt, minimum width=5mm, fill=white, inner sep=0.25mm] (text3381) {$\pi$};
\end{tikzpicture}\rlap{\ \  .}}}
\end{equation*}

Now if $\phi \colon V, M \to W$ is a universal bimorphism, then applying  (the analogue for bimodules of) the functor~\eqref{eq:functor-bimodules-to-cocones} yields a colimiting cylinder $\tilde \phi \colon \Delta 1 \to {}_\bullet \cat{Mod}_\C(F^{VM}, W)$. Because $B \otimes \thg \colon \V \to \V$ preserves bicolimits,
and the forgetful ${}_\bullet \cat{Mod}_\C \to \V^{\ob \C}$ creates them, it follows that $B \otimes \thg \colon {}_\bullet \cat{Mod}_\C \to {}_\bullet \cat{Mod}_\C$ also preserves bicolimits. Thus the composite cylinder
\begin{equation*}
\Delta 1 \xrightarrow{\tilde \phi} {}_\bullet \cat{Mod}_\C(F^{VM}, W) \xrightarrow{B \otimes \thg} {}_\bullet \cat{Mod}_\C(B \otimes F^{VM}, B \otimes W) \xrightarrow{(\thg) \circ \gamma} {}_\bullet \cat{Mod}_\C(F^{B \otimes V, M}, B \otimes W)
\end{equation*}
is also colimiting. Applying the pseudoinverse of~\eqref{eq:functor-bimodules-to-cocones} to this  cylinder, we obtain a universal bimorphism $B \otimes V, M \to B \otimes W$ which, by tracing through the construction, we see to be isomorphic to $B \otimes \phi$.
\end{proof}

\section{Internal hom of modules}\label{sec:internal-hom}
\subsection{Left and right homs}
Given $M \colon \A \tor \B$, $N \colon \B \tor \C$ and $P \colon \A \tor \C$, we can consider each of the partial functors:
\begin{align*}
\cat{Bimor}_{\A\B\C}(M, N; \thg) \colon {}_\A\cat{Mod}_\C & \to \cat{CAT}\\
\cat{Bimor}_{\A\B\C}(\thg, N; P) \colon {}_\A\cat{Mod}_\B^\op & \to \cat{CAT}\\
\cat{Bimor}_{\A\B\C}(M, \thg; P) \colon {}_\B\cat{Mod}_\C^\op & \to \cat{CAT}\rlap{ .}
\end{align*}
We have already defined the \emph{tensor product} $M \otimes_\B N$ of $M$ and $N$ to be a birepresentation for the first of these; we now define the \emph{right hom} $\h N P$ of $N$ and $P$ as a birepresentation for the second, and the \emph{left hom} $\h M P_\ell$ of $M$ and $P$ as a birepresentation for the third. In terms of the tricategory $\V\text-\cat{Mod}$ of $\V$-bimodules, the existence of left or right homs amounts to the existence of (tricategorical) \emph{right extensions} and \emph{right liftings}.
In this section, we discuss the construction of left and right homs between modules; in fact, we shall concentrate on the case of right homs, since that is what we will need for our further development. The arguments for the left case are entirely analogous.

\subsection{Right closed bicategories}\label{subsec:right-closed-bicat}
In order to assure the existence of right homs, we shall assume that our base monoidal bicategory $\V$ is \emph{right closed}, meaning that
each functor $\thg \otimes B \colon \V \to \V$ admits a right biadjoint $[B, \thg]$. (If we were interested in the construction of left homs, we would require instead left closedness, involving the existence of right biadjoints to each $B \otimes \thg$). Thus we have unit and counit $1$-cells
$\mathfrak u_{AB} \colon A \to [B, A \otimes B]$ and $\mathfrak e_{BC} \colon [B,C] \otimes B \to C$, inducing adjoint equivalences of categories
\begin{equation}\label{eq:exptransp}
\V(A \otimes B, C) \simeq \V(A, [B,C])
\end{equation}
pseudonatural in $A$ and $C$. We call the process of applying either direction of this equivalence \emph{exponential transpose}, and use a bar to denote its action; thus, given $f \colon A \otimes B \to C$ and $g \colon A \to [B,C]$, we have
\begin{equation*}
\bar f = A \xrightarrow{\mathfrak u} [B, A \otimes B] \xrightarrow{[1, f]} [B,C] \quad \text{and} \quad
\bar g = A \otimes B \xrightarrow{g \otimes 1} [B,C] \otimes B \xrightarrow{\mathfrak e} C\rlap{ .}
\end{equation*}
We write $\omega \colon \bar{\bar f} \Rightarrow f$ to denote the invertible $2$-cell relating a morphism and its double transpose. As usual, the functors $[B, \thg]$ assemble to give a functor $[\thg, \thg] \colon \V^\op \times \V \to \V$ in such a way that the equivalences~\eqref{eq:exptransp} become pseudonatural in $B$ as well as $A$ and $C$. In the first argument, the action of this functor on $f \colon B' \to B$ is the map $[f, C] \colon [B,C] \to [B',C]$ obtained as the transpose of
\begin{equation*}
[B,C] \otimes B' \xrightarrow{1 \otimes f} [B,C] \otimes B \xrightarrow{\mathfrak e} C\rlap{ .}
\end{equation*}

The coherence data of the monoidal bicategory $\V$ may be recast in terms of the right hom $[\thg, \thg]$. First, for each $A \in \V$, we have an equivalence $\tilde{\mathfrak r}_A \defeq \mathfrak e_{IA} \circ\mathfrak r^\centerdot_A \colon [I,A] \to A$; a suitable pseudoinverse is given by the exponential transpose of $\mathfrak r \colon A \otimes I \to A$. Next, we have for each $A, B, C \in \V$ an equivalence $\tilde{\mathfrak a}_{ABC} \colon [A, [B,C]] \to [A \otimes B, C]$, obtained as the exponential transpose of
\begin{equation*}
[A, [B,C]] \otimes (A \otimes B) \xrightarrow{\mathfrak a^\centerdot} ([A,[B,C]] \otimes A) \otimes B \xrightarrow{\mathfrak e \otimes 1} [B,C] \otimes B \xrightarrow{\mathfrak e} C\rlap{ .}
\end{equation*}
A suitable pseudoinverse is obtained by transposing twice the composite\begin{equation*}
([A \otimes B, C] \otimes A) \otimes B \xrightarrow{\mathfrak a} [A \otimes B, C] \otimes (A \otimes B) \xrightarrow{\mathfrak e} C\rlap{ .}
\end{equation*}
The pseudonaturality of $\mathfrak r^\centerdot$ and $\mathfrak a^\centerdot$ immediately implies a corresponding pseudonaturality for $\tilde{\mathfrak r}$ and $\tilde{\mathfrak a}$. Finally, we have coherence $2$-cells
\begin{equation*}
\cd[@C-3em]{
 & [A \otimes I, C] \ar[dr]^{[\mathfrak r^\centerdot, 1]} \dtwocell{d}{\tilde{\rho}} \\
 [A,[I,C]] \ar[rr]_-{[1, \tilde{\mathfrak r}]} \ar[ur]^{\tilde{\mathfrak a}} & &
 [A,C]
} \quad \cd[@C-1em]{
 [A,[B,[C,D]] \ar[d]_{[1, \tilde{\mathfrak a}]} \ar[rr]^-{\tilde{\mathfrak a}}
  \rtwocell{drr}{\tilde \pi} & &
 [A \otimes B, [C,D]] \ar[d]^{\tilde{\mathfrak a}} \\
 [A,[B \otimes C, D]] \ar[r]_-{\tilde{\mathfrak a}} &
 [A \otimes (B \otimes C), D] \ar[r]_-{[\mathfrak a, 1]} &
 [(A \otimes B) \otimes C, D]\rlap{ ;}
}
\end{equation*}
to obtain these, it suffices to give $2$-cells between the adjoint transposes of their respective domains and codomains. We obtain such as the composites
\begin{equation*}
\vcenter{\hbox{\begin{tikzpicture}[y=0.80pt, x=0.8pt,yscale=-1, inner sep=0pt, outer sep=0pt, every text node part/.style={font=\tiny} ]
\path[draw=black,line join=miter,line cap=butt,line width=0.650pt]
(18.0000,1032.3622) .. controls (37.1610,1032.3622) and (45.3064,1027.8382) ..
node[above right=0.12cm,at start] {\!$\mathfrak e$}(48.3368,1024.1764);
\path[draw=black,line join=miter,line cap=butt,line width=0.650pt]
(18.0000,1012.3622) .. controls (34.9927,1012.3622) and (45.4296,1015.3926) ..
node[above right=0.12cm,at start] {\!$[\mathfrak r^\centerdot,\! 1] \t 1$}(48.4848,1019.4333);
\path[draw=black,line join=miter,line cap=butt,line width=0.650pt]
(52.5279,1022.0576) .. controls (58.5994,1027.9505) and (82.0211,1028.0244) ..
node[below=0.07cm] {$\mathfrak e$}(87.3782,1022.3101);
\path[draw=black,line join=miter,line cap=butt,line width=0.650pt]
(52.0229,1019.2797) .. controls (55.5202,1015.0504) and (62.4505,1009.0897) ..
node[below right,pos=0.65] {$1 \t \mathfrak r^\centerdot$}(71.4244,1003.8521)
(76.6096,1001.0332) .. controls (86.3936,996.0966) and (98.0015,992.3622) ..
node[above left=0.02cm,pos=0.5] {$1 \t \mathfrak r^\centerdot$}(110.0000,992.3622)
(18.0000,992.3622) .. controls (88.1675,992.3622) and (86.6180,1018.6261) ..
node[above right=0.12cm,at start] {\!$\tilde{\mathfrak a} \t 1$}(86.6180,1018.6261);
\path[draw=black,line join=miter,line cap=butt,line width=0.650pt]
(89.3985,1019.2797) .. controls (89.3985,1019.2797) and (140.4892,999.2657) ..
node[above left=0.06cm,pos=0.45] {$\mathfrak a^\centerdot$}(112.6320,993.0157);
\path[draw=black,line join=miter,line cap=butt,line width=0.650pt]
(112.6753,990.2684) .. controls (143.2111,983.1256) and (150.0077,1023.0244) ..
node[above right=0.05cm,pos=0.2] {$\mathfrak r^\centerdot$}(150.0077,1023.0244)
(90.3167,1021.0654) .. controls (112.4837,1019.7691) and (123.8685,1007.5305) ..
node[below right=0.04cm,pos=0.6] {$\mathfrak e \t 1$}(137.6038,1000.9195)
(142.3886,998.9109) .. controls (146.2013,997.5586) and (150.2818,996.7298) ..
(154.8720,996.7298) .. controls (176.8057,996.7298) and (171.0692,1022.3622) ..
node[above left=0.12cm,at end] {$\mathfrak e$\!}(200.0000,1022.3622)
(90.3571,1023.2550) .. controls (101.9643,1027.8979) and (128.6097,1027.4550) ..
node[below=0.07cm,pos=0.55] {$\mathfrak e$}
(150.0000,1023.2550) .. controls (160.1525,1021.2615) and (166.8612,1016.4538) ..
node[below=0.07cm,pos=0.5] {$\tilde{\mathfrak r}$}(173.5146,1011.9523)
(177.7103,1009.1982) .. controls (183.7963,1005.3938) and (190.4087,1002.3622) ..
node[above left=0.12cm,at end] {$[1,\! \tilde{\mathfrak r}] \t 1$\!\!}(200.0000,1002.3622);
\path[shift={(73.400512,91.955576)},draw=black,fill=black] (77.8571,931.1122)arc(0.000:180.000:1.250)arc(-180.000:0.000:1.250) -- cycle;
\path[fill=black] (49.750011,1021.5525) node[circle, draw, line width=0.65pt, minimum width=5mm, fill=white, inner sep=0.25mm] (text3199) {$\omega$};
\path[fill=black] (88.388344,1021.3) node[circle, draw, line width=0.65pt, minimum width=5mm, fill=white, inner sep=0.25mm] (text3203) {$\omega$};
\path[fill=black] (107.11678,992.00555) node[circle, draw, line width=0.65pt, minimum width=5mm, fill=white, inner sep=0.25mm] (text3207) {$\rho$};
\end{tikzpicture}}}
\quad\text{and}\quad\vcenter{\hbox{\begin{tikzpicture}[y=0.9pt, x=0.9pt,yscale=-1,xscale=-1, inner sep=0pt, outer sep=0pt, every text node part/.style={font=\tiny} ]
  \path[draw=black,line join=miter,line cap=butt,line width=0.650pt]
  (12.0000,1032.3622) .. controls (32.0000,1032.3622) and (36.4289,1026.2453) ..
  node[above left=0.12cm,at start] {$\mathfrak e$\!}(39.1075,1023.9239);
  \path[draw=black,line join=miter,line cap=butt,line width=0.650pt]
  (12.0000,1012.3622) .. controls (32.5388,1012.3622) and (36.6075,1017.6169) ..
  node[above left=0.12cm,at start] {$\tilde{\mathfrak a} \t 1$\!}(39.1075,1020.4740);
  \path[draw=black,line join=miter,line cap=butt,line width=0.650pt]
  (42.2145,1021.8051) .. controls (53.4645,1025.9122) and (80.7225,1026.4480) ..
  node[above=0.085cm,pos=0.45] {$\mathfrak e \t 1$}(89.6510,1021.8051);
  \path[draw=black,line join=miter,line cap=butt,line width=0.650pt]
  (12.0000,992.3622) .. controls (25.7076,992.3622) and (38.8505,995.2828) ..
  node[above left=0.12cm,at start] {$\tilde{\mathfrak a} \t 1$\!}(50.2588,999.2079)
  (55.6973,1001.2003) .. controls (74.5305,1008.5200) and (87.5211,1017.9397) ..
  node[above left=0.02cm,pos=0.13,rotate=25] {$(\tilde{\mathfrak a} \t 1) \t 1$}(88.4848,1019.3317)
  (41.7094,1019.2797) .. controls (43.8382,992.9095) and (95.0912,974.8997) ..
  node[above=0.1cm] {$\mathfrak a^\centerdot$}(128.9286,993.0765);
  \path[draw=black,line join=miter,line cap=butt,line width=0.650pt]
  (91.9239,1019.7848) .. controls (96.9239,1009.2491) and (116.0546,995.5411) ..
  node[above right=0.05cm,pos=0.6] {$\mathfrak a^\centerdot \t 1$}(129.8046,995.5411);
  \path[draw=black,line join=miter,line cap=butt,line width=0.650pt]
  (91.9799,1022.0882) .. controls (105.3388,1027.8007) and (152.1335,1033.7198) ..
  node[above=0.1cm] {$\mathfrak e \t 1$}(162.8871,1027.5525);
  \path[draw=black,line join=miter,line cap=butt,line width=0.650pt]
  (41.2170,1024.5217) .. controls (46.4646,1043.7508) and (145.1058,1049.3267) ..
  node[above=0.06cm] {$\mathfrak e$}(162.6039,1029.9860);
  \path[draw=black,line join=miter,line cap=butt,line width=0.650pt]
  (164.6429,1029.7908) .. controls (176.8158,1034.8769) and (192.1339,1031.4582) ..
  node[below=0.1cm] {$\mathfrak e$}(210.0000,1024.6837);
  \path[draw=black,line join=miter,line cap=butt,line width=0.650pt]
  (211.9643,1023.4337) .. controls (218.0357,1030.9337) and (242.3214,1035.5764) ..
  node[below=0.1cm] {$\mathfrak e$}(248.9286,1032.8979);
  \path[draw=black,line join=miter,line cap=butt,line width=0.650pt]
  (251.4286,1030.9335) .. controls (258.7500,1024.1478) and (263.9992,1022.3622) ..
  node[above right=0.12cm,at end] {\!\!$[\mathfrak a,\! 1] \t 1$}(284.0000,1022.3622);
  \path[draw=black,line join=miter,line cap=butt,line width=0.650pt]
  (284.0000,1042.3622) .. controls (264.6936,1042.3622) and (256.7857,1042.0051) ..
  node[above right=0.12cm,at start] {\!$\mathfrak e$}(251.2500,1033.2551);
  \path[draw=black,line join=miter,line cap=butt,line width=0.650pt]
  (164.6429,1027.2908) .. controls (169.3188,1020.2804) and (173.8557,1014.3278) ..
  node[below left=0.03cm] {$\tilde{\mathfrak a} \t 1$}(178.3267,1009.2772)
  (182.2599,1005.0472) .. controls (187.4695,999.7265) and (192.6178,995.5851) ..
  node[above right=0.02cm,pos=0.38] {$([1,\!\tilde{\mathfrak a}] \t 1) \t 1$}(197.8229,992.3712)
  (202.4037,989.7848) .. controls (224.8517,978.2646) and (248.9707,982.3622) ..
  node[above right=0.12cm,at end] {\!\!$[1,\! \tilde{\mathfrak a}] \t 1$}(284.0000,982.3622)
  (131.4286,992.3622) .. controls (192.6247,968.9983) and (211.0815,992.0040) ..
  node[above=0.1cm] {$\mathfrak a^\centerdot$}(209.8214,1019.5051)
  (96.4286,1019.5050) .. controls (113.3467,1013.0884) and (129.6258,1008.3527) ..
  node[above=0.02cm,rotate=-18,pos=0.45] {$(\mathfrak e \t 1) \t 1$}(144.8396,1006.2758)
  (151.1153,1005.5670) .. controls (173.7825,1003.5608) and (193.8505,1008.0138) ..
  node[above=0.1cm,pos=0.73] {$\mathfrak e \t 1$}(209.8214,1022.3622)
  (131.3198,994.2784) .. controls (142.9270,997.1355) and (159.0019,1017.5884) ..
  node[below right=0.03cm,pos=0.6] {$\mathfrak a^\centerdot$}(162.2161,1025.4456)
  (128.8378,990.2251) .. controls (109.7452,976.3230) and (152.2005,968.4289) ..
  node[above=0.1cm,at end] {$1 \t \mathfrak a$}
  (168.3870,968.6070) .. controls (195.5344,968.5224) and (212.2036,973.4190) .. (222.7914,980.6171)
  (228.1811,984.9350) .. controls (234.5310,990.9321) and (237.9115,997.9442) ..
  node[left=0.12cm] {$1 \t \mathfrak a$}(240.1931,1004.8309)
  (241.6368,1009.6244) .. controls (243.8409,1017.5982) and (245.1263,1025.0277) ..
  node[right=0.08cm,pos=0.4] {$1 \t \mathfrak a$}(248.5550,1030.0465)
  (211.7857,1020.9336) .. controls (224.1071,1008.7908) and (254.5162,1002.3622) ..
  node[above right=0.12cm,at end] {\!$\tilde{\mathfrak a} \t 1$}(284.0000,1002.3622);
  \path[fill=black] (38.285713,1022.005) node[circle, draw, line width=0.65pt, minimum width=5mm, fill=white, inner sep=0.25mm] (text3395) {$\omega^{-1}$     };
  \path[fill=black] (92.64286,1020.5765) node[circle, draw, line width=0.65pt, minimum width=5mm, fill=white, inner sep=0.25mm] (text3399) {$\omega^{-1} \t 1$     };
  \path[fill=black] (129.64285,993.07648) node[circle, draw, line width=0.65pt, minimum width=5mm, fill=white, inner sep=0.25mm] (text3403) {$\pi^\centerdot$     };
  \path[fill=black] (163.21428,1030.005) node[circle, draw, line width=0.65pt, minimum width=5mm, fill=white, inner sep=0.25mm] (text3407) {$\omega$     };
  \path[fill=black] (211.16881,1022.7904) node[circle, draw, line width=0.65pt, minimum width=5mm, fill=white, inner sep=0.25mm] (text3407-1) {$\omega$     };
  \path[fill=black] (249.74023,1031.3618) node[circle, draw, line width=0.65pt, minimum width=5mm, fill=white, inner sep=0.25mm] (text3407-4) {$\omega$   };
\end{tikzpicture}\rlap{ .}}}
\end{equation*}

We may similarly recast the data of a right module $W \colon \bullet \tor \B$ in terms of the internal hom. For all $x,y \in \B$, the action morphism of $W$ at $(x,y)$ transposes to yield $\bar m \colon Wy \to [\B(x,y), Wx]$. Furthermore, we obtain coherence $2$-cells
\begin{equation*}
\cd{
[\B(x,x), Wx] \ar[r]^-{[j,1]} \dtwocell{dr}{\tilde \ru} &
[I, Wx] \ar[d]^{\tilde{\mathfrak r}} \\
Wx \ar[u]^{\bar m} \ar[r]_1 & Wx
}\quad\ \text{and}\ \
\cd[@C-0.1em]{
Wz \ar[r]^-{\bar m} \ar[dd]_{\bar m}
 &
[\B(y,z), Wy] \ar[d]^{[1,\bar m]} \\
\rtwocell[0.35]{r}{\tilde \ass} &
[\B(y,z), [\B(x,y), Wx]] \ar[d]^{\tilde{\mathfrak a}} \\
[\B(x,z), Wx] \ar[r]_-{[m,1]} &
[\B(y,z) \otimes \B(x,y), Wx]
}
\end{equation*}
where $\tilde \tau$ is as on the left below, and $\tilde \alpha$ is uniquely determined by the $2$-cell between the exponential transpose of its domain and codomain displayed on the right below.
\begin{equation*}
\vcenter{\hbox{\begin{tikzpicture}[y=0.9pt, x=0.9pt,yscale=-1, inner sep=0pt, outer sep=0pt, every text node part/.style={font=\tiny} ]
  \path[draw=black,line join=miter,line cap=butt,line width=0.650pt]
  (10.0000,1042.3622) .. controls (40.3056,1042.3622) and (43.3115,1041.1784) ..
  node[above right=0.12cm,at start] {\!$\mathfrak e$}(47.8571,1038.1479);
  \path[draw=black,line join=miter,line cap=butt,line width=0.650pt]
  (51.0281,1037.8954) .. controls (57.3416,1044.2089) and (80.3856,1041.3235) ..
  node[above=0.07cm,pos=0.5] {$\mathfrak e$}(88.2143,1034.5050);
  \path[draw=black,line join=miter,line cap=butt,line width=0.650pt]
  (92.1429,1032.3622) .. controls (107.2952,1032.3622) and (115.1151,1012.9485) ..
  node[below right=0.07cm,pos=0.5] {$m$}(119.1557,1004.6147);
  \path[draw=black,line join=miter,line cap=butt,line width=0.650pt]
  (50.8368,1034.4718) .. controls (54.3045,1031.0042) and (60.8230,1024.3395) ..
  node[above left=0.07cm,pos=0.6] {$1 \t j$}(68.8833,1017.7209)
  (72.8087,1014.6036) .. controls (87.1254,1003.6286) and (105.1323,994.2330) ..
  node[above=0.07cm,pos=0.45] {$1 \t j$}(119.6429,1001.8622)
  (10.0000,1022.3622) .. controls (47.1358,1022.3622) and (70.5197,977.8073) ..
  node[above right=0.12cm,at start] {\!$\mathfrak r^\centerdot$}
  (104.0432,978.6342) .. controls (118.3571,978.9873) and (147.6615,994.1455) .. (122.1552,1003.2369)
  (10.0000,1002.3622) .. controls (21.8965,1002.3622) and (29.1430,1007.3616) ..
  node[above right=0.12cm,at start] {\!\!$[j,\!1]$}(34.0413,1013.4319)
  (36.9563,1017.5498) .. controls (41.7961,1025.2793) and (43.9074,1033.3153) ..
  node[left=0.07cm,pos=0.5] {$[j,\!1] \t 1$}(47.4387,1034.5765)
  (10.0000,982.3622) .. controls (29.6009,982.3622) and (43.6462,990.3287) ..
  node[above right=0.12cm,at start] {\!$\bar m$}(54.6778,999.8377)
  (58.4572,1003.2539) .. controls (70.6313,1014.7545) and (79.0550,1027.3404) ..
  node[above right=0.07cm,pos=0.6] {$\bar m \t 1$}(87.8571,1030.5765);
  \path[fill=black] (49.244938,1036.6642) node[circle, draw, line width=0.65pt, minimum width=5mm, fill=white, inner sep=0.25mm] (text3340) {$\omega$ };
  \path[fill=black] (89.398499,1032.9167) node[circle, draw, line width=0.65pt, minimum width=5mm, fill=white, inner sep=0.25mm] (text3344) {$\omega$  };
  \path[fill=black] (119.95561,1003.9561) node[circle, draw, line width=0.65pt, minimum width=5mm, fill=white, inner sep=0.25mm] (text3348) {$\ru$};
  \end{tikzpicture}}}\ \ \text{and}\quad\ \ \
  \vcenter{\hbox{\begin{tikzpicture}[y=0.9pt, x=0.9pt,yscale=-1, inner sep=0pt, outer sep=0pt, every text node part/.style={font=\tiny} ]
  \path[draw=black,line join=miter,line cap=butt,line width=0.650pt]
  (10.0000,1022.3622) .. controls (29.1610,1022.3622) and (37.3064,1017.8382) ..
  node[above right=0.12cm,at start] {\!$\mathfrak e$}(40.3368,1014.1764);
  \path[draw=black,line join=miter,line cap=butt,line width=0.650pt]
  (10.0000,1002.3622) .. controls (26.9927,1002.3622) and (37.4296,1005.3926) ..
  node[above right=0.12cm,at start] {\!\!$[m,\!1] \t 1$}(40.4848,1009.4333);
  \path[draw=black,line join=miter,line cap=butt,line width=0.650pt]
  (43.5178,1012.3101) .. controls (49.5892,1018.2030) and (74.0211,1018.0244) ..
  node[below=0.07cm,pos=0.5] {$\mathfrak e$}(79.3782,1012.3101);
  \path[draw=black,line join=miter,line cap=butt,line width=0.650pt]
  (44.0229,1009.2797) .. controls (48.2846,1001.2761) and (52.8467,995.5068) ..
  node[below right=0.03cm,pos=0.5] {$1 \t m$}(57.6078,991.5294)
  (62.4665,988.1127) .. controls (77.8423,979.2206) and (94.8533,986.6208) ..
  node[above=0.07cm,pos=0.2] {$1 \t m$}(110.3338,996.5043)
  (10.0000,982.3622) .. controls (80.1675,982.3622) and (79.6282,1009.8888) ..
  node[above right=0.12cm,at start] {\!$\bar m \t 1$}(79.6282,1009.8888);
  \path[draw=black,line join=miter,line cap=butt,line width=0.650pt]
  (82.0749,1010.6934) .. controls (97.2272,1010.4409) and (105.8135,1004.6325) ..
  node[below right=0.07cm,pos=0.35] {$m$}(110.8643,998.3190);
  \path[draw=black,line join=miter,line cap=butt,line width=0.650pt]
  (112.7583,996.1725) .. controls (122.2285,989.8590) and (147.3560,993.2683) ..
  node[above=0.07cm,pos=0.45] {$m \t 1$}(152.7856,997.5614);
  \path[draw=black,line join=miter,line cap=butt,line width=0.650pt]
  (112.8846,998.6979) .. controls (119.7031,1011.5773) and (130.0845,1026.7849) ..
  node[above right=0.07cm,pos=0.4] {$m$}(145.7439,1026.7849);
  \path[draw=black,line join=miter,line cap=butt,line width=0.650pt] (147.3024,1028.2875) .. controls (151.8227,1034.0600) and (214.5159,1032.3018) .. node[below=0.07cm] {$\mathfrak e$}(220.8694,1027.3609);
  \path[draw=black,line join=miter,line cap=butt,line width=0.650pt]
  (111.7212,992.7452) .. controls (87.4520,972.0066) and (167.3866,969.1566) ..
  node[above=0.06cm,at end] {$\mathfrak a^\centerdot$}(176.2562,970.1305) .. controls (206.4389,973.4448) and (211.8253,1003.0561) .. (220.8694,1022.0576)
  (154.3314,996.5821) .. controls (156.8346,990.4589) and (173.2702,983.9896) ..
  node[above=0.03cm,pos=0.64,rotate=15] {$(\bar m \t 1) \t 1$}(193.6962,979.2788)
  (199.8722,977.9331) .. controls (216.3240,974.5508) and (234.5999,972.3622) ..
  node[above left=0.12cm,at end] {$\bar m \t 1$\!}(250.0000,972.3622)
  (147.1978,1025.0481) .. controls (147.9681,1022.3949) and (155.1489,1017.7858) ..
  node[below right=0.06cm,pos=0.57,rotate=23] {$\bar m \t 1$}(165.9589,1012.8513)
  (172.0426,1010.2020) .. controls (182.3257,1005.9223) and (194.9072,1001.6163) ..
  node[above=0.04cm,pos=0.45,rotate=16] {$([1,\!\bar m] \t 1)\t 1$}(208.0415,998.3071)
  (213.0118,997.1098) .. controls (225.5347,994.2331) and (238.3559,992.3624) ..
  node[above left=0.12cm,at end] {$[1,\!\bar m] \t 1$\!\!}(250.0000,992.3624)
  (154.5533,998.8241) .. controls (157.0608,1013.4205) and (208.4433,1024.7310) ..
  node[below left=0.06cm,pos=0.6] {$\mathfrak e \t 1$}(221.1219,1025.0881);
  \path[draw=black,line join=miter,line cap=butt,line width=0.650pt] (223.1639,1026.9350) .. controls (228.4672,1031.4807) and (234.0901,1032.3622) .. node[above left=0.12cm,at end] {$\mathfrak e$\!}(250.0000,1032.3622);
  \path[draw=black,line join=miter,line cap=butt,line width=0.650pt] (222.9636,1023.6520) .. controls (227.8881,1019.3589) and (233.3312,1012.3622) .. node[above left=0.12cm,at end] {$\tilde{\mathfrak a} \t 1$\!}(250.0000,1012.3622);
  \path[fill=black] (41.750008,1011.5526) node[circle, draw, line width=0.65pt, minimum width=5mm, fill=white, inner sep=0.25mm] (text3199) {$\omega$     };
  \path[fill=black] (111.11678,996.55127) node[circle, draw, line width=0.65pt, minimum width=5mm, fill=white, inner sep=0.25mm] (text3620) {$\ass$     };
  \path[fill=black] (153.79572,997.81396) node[circle, draw, line width=0.65pt, minimum width=5mm, fill=white, inner sep=0.25mm] (text3624) {$\omega^{-1} \t 1$     };
  \path[fill=black] (142.38202,1026.6033) node[circle, draw, line width=0.65pt, minimum width=5mm, fill=white, inner sep=0.25mm] (text3628) {$\omega^{-1}$     };
  \path[fill=black] (222.13203,1026.8558) node[circle, draw, line width=0.65pt, minimum width=5mm, fill=white, inner sep=0.25mm] (text3632) {$\omega^{-1}$   };
  \path[fill=black] (80.239113,1012.8493) node[circle, draw, line width=0.65pt, minimum width=5mm, fill=white, inner sep=0.25mm] (text3199-4) {$\omega$   };
\end{tikzpicture}\rlap{ .}}}
\end{equation*}

\subsection{Construction of right homs}
\label{subsec:constr-right-homs}
Given modules $V, W \colon \bullet \tor \B$, we now consider the existence of the right hom $\h V W$. This right hom is, by definition, a birepresentation of the functor $\cat{Bimor}_{\bullet \bullet \B}(\thg, V; W) \colon \V^\op \to \cat{CAT}$; but by the remarks at the start of Section~\ref{subsec:construction-of-tensor-products}, this is equally a birepresentation of the functor ${}_\bullet \cat{Mod}_\B(\thg \otimes V, W) \colon \V^\op \to \cat{CAT}$. Thus, we seek an object $\h V W$ of $\V$ and a module morphism $\xi \colon \h V W \otimes V \to W$ such that, for every $B \in \V$, the functor
\begin{equation}\label{eq:compinduced}
\begin{aligned}
\V(B, \h V W) & \to \cat{Mod}_\B(B \otimes V, W) \\
f & \mapsto \xi \circ (f \otimes V)
\end{aligned}
\end{equation}
is an equivalence of categories.

Working still in the context of a right closed $\V$, we will describe a limit in $\V$ which, if it exists, must underlie this right hom, and whose existence, conversely, is guaranteed by the existence of the hom. To this end, let $\D$ be the category defined in Section~\ref{subsec:construction-of-tensor-products}, and consider the functor $F^{VW} \colon \D^\op \to \V$ that on objects, is given by:
\begin{itemize}
\item $F^{VW}(x) = [Vx, Wx]$;
\item $F^{VW}(x,y) = [Vy \otimes \B(x,y), Wx]$;
\item $F^{VW}(x,y,z) = [(Vz \otimes \B(y,z)) \otimes \B(x,y), Wx]$
\end{itemize}
and on generating morphisms by:
\begin{itemize}
\item $F^{VW}(i_x) =  [\mathfrak r^\centerdot,1] \circ [1 \otimes j, 1]$;
\item $F^{VW}(d_{xy}) = [m, 1]$;
\item $F^{VW}(c_{xy}) =  \bar{\mathfrak a} \circ [1, \bar m] $;
\item $F^{VW}(p_{xyz}) = [m \otimes 1,\, 1]$;
\item $F^{VW}(n_{xyz}) = [\mathfrak a, 1] \circ [1 \otimes m,\, 1]$;
\item $F^{VW}(q_{xyz}) = \bar{\mathfrak a} \circ [1, \bar m]$.
\end{itemize}
To extend these assignations to a pseudofunctor, it suffices to exhibit pseudofunctoriality cells with respect to the generating simplicial identities $di = 1$, $ci = 1$, $cp = dq$, $dp = dm$ and $cm = cq$. These are given by the respective string diagrams:
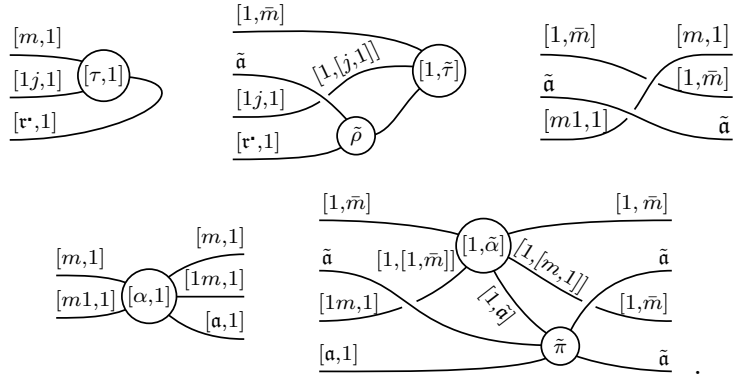
\begin{equation*}
\vcenter{\hbox{
\begin{tikzpicture}[y=0.80pt, x=0.9pt, inner sep=0pt, outer sep=0pt, every text node part/.style={font=\tiny} ]
\path[use as bounding box] (8, 960) rectangle (85,1039);
  \path[draw=black,line join=miter,line cap=butt,line width=0.650pt]
  (10.0000,1012.3622) .. controls (31.9384,1012.3622) and (44.6754,1009.3804) ..
  node[above right=0.12cm,at start] {\!$[m,\!1]$}(48.2961,1005.5467);
  \path[draw=black,line join=miter,line cap=butt,line width=0.650pt]
  (10.0000,992.3622) .. controls (30.6605,992.3622) and (45.3144,994.7050) ..
  node[above right=0.12cm,at start] {\!$[1 \t j,\! 1]$}(48.7221,999.6036);
  \path[draw=black,line join=miter,line cap=butt,line width=0.650pt]
  (10.0000,972.3622) .. controls (61.5891,972.3622) and (99.1240,1001.2238) ..
  node[above right=0.12cm,at start] {\!$[\mathfrak r^\centerdot,\!1]$}(51.2779,1002.6059);
  \path[fill=black] (49.096386,1002.6634) node[circle, draw, line width=0.65pt, minimum width=5mm, fill=white, inner sep=0.25mm] (text6241) {$[\tau,\!1]$   };
\end{tikzpicture}}} \quad
\vcenter{\hbox{
\begin{tikzpicture}[y=0.80pt, x=0.7pt,yscale=-1, inner sep=0pt, outer sep=0pt, every text node part/.style={font=\tiny} ]
  \path[draw=black,line join=miter,line cap=butt,line width=0.650pt]
  (60.0000,942.3622) .. controls (95.0000,941.6479) and (152.1429,939.5050) ..
  node[above right=0.12cm,at start] {\!$[1,\!\bar m]$}(169.2857,956.6479);
  \path[draw=black,line join=miter,line cap=butt,line width=0.650pt]
  (60.0000,982.3622) .. controls (63.3448,982.3622) and (67.2811,982.3622) ..
  node[above right=0.12cm,at start] {\!$[1 \t j,\!1]$}
  (70.0000,982.3622) .. controls (89.6712,982.3622) and (99.6685,979.1375) .. (106.6624,975.0544)
  (111.1725,972.0713) .. controls (121.9797,964.1628) and (128.0416,955.2531) ..
  node[above left=0.15cm,pos=0.75,rotate=25] {$[1,\![j,\!1]]$}(169.2857,959.5050)
  (60.0000,962.3622) .. controls (100.0000,962.3622) and (109.6429,972.7193) ..
  node[above right=0.12cm,at start] {\!$\tilde{\mathfrak a}$}(126.0714,988.4336);
  \path[draw=black,line join=miter,line cap=butt,line width=0.650pt]
  (60.0000,1002.3622) .. controls (110.0000,1002.3622) and (114.2857,999.5050) ..
  node[above right=0.12cm,at start] {\!$[\mathfrak r^\centerdot,\!1]$}(126.4286,992.0051);
  \path[draw=black,line join=miter,line cap=butt,line width=0.650pt]
  (128.2143,991.6478) .. controls (147.5000,992.3621) and (147.8572,974.5050) .. (169.2857,962.3622);
  \path[fill=black] (170.35713,960.57648) node[circle, draw, line width=0.65pt, minimum width=5mm, fill=white, inner sep=0.25mm] (text3445) {$[1,\!\tilde \ru]$     };
  \path[fill=black] (125.97502,991.23486) node[circle, draw, line width=0.65pt, minimum width=5mm, fill=white, inner sep=0.25mm] (text3449) {$\tilde \rho$   };
\end{tikzpicture}}}\ \ \ \quad
\vcenter{\hbox{
\begin{tikzpicture}[y=0.8pt, x=1pt,inner sep=0pt, yscale=-1,outer sep=0pt, every text node part/.style={font=\scriptsize} ]
  \path[draw=black,line join=miter,line cap=butt,line width=0.650pt]
  (44.0000,952.3622) .. controls (65.7595,952.3622) and (74.4984,958.2806) ..
  node[above right=0.12cm,at start] {\!$[1,\!\bar m]$}(83.7389,963.6784)
  (88.1039,966.1173) .. controls (94.7498,969.5940) and (102.6067,972.3622) ..
  node[above left=0.12cm,at end] {$[1,\!\bar m]$\!}(116.0000,972.3622)
  (44.0000,972.3622) .. controls (84.0000,972.3622) and (80.0000,992.3622) ..
  node[above right=0.12cm,at start] {\!$\tilde{\mathfrak a}$}
  node[above left=0.12cm,at end] {$\tilde{\mathfrak a}$\!}(116.0000,992.3622)
  (44.0000,992.3622) .. controls (44.9510,992.3622) and (53.1023,992.3622) ..
  node[above right=0.12cm,at start] {\!$[m \t 1,\!1]$}
  (54.0000,992.3622) .. controls (65.7332,992.3622) and (72.1551,988.1132) .. (76.5202,982.3848)
  (79.2680,978.1916) .. controls (85.8514,966.6868) and (88.0472,952.3622) ..
  (106.0000,952.3622) .. controls (107.9261,952.3622) and (115.3396,952.3622) ..
  node[above left=0.12cm,at end] {$[m,\!1]$\!}(116.0000,952.3622);
\end{tikzpicture}}}
\end{equation*}
\begin{equation*}
\vcenter{\hbox{
\begin{tikzpicture}[y=0.80pt, x=1.1pt,yscale=-1, inner sep=0pt, outer sep=0pt, every text node part/.style={font=\tiny} ]
  \path[draw=black,line join=miter,line cap=butt,line width=0.650pt]
  (122.0000,962.3622) .. controls (104.2362,962.3622) and (93.5355,974.6350) ..
  node[above left=0.12cm,at start] {$[m,\!1]$\!\!}(91.7678,979.4333);
  \path[draw=black,line join=miter,line cap=butt,line width=0.650pt]
  (122.0000,982.3622) --
  node[above left=0.12cm,at start] {$[1 \t m,\!1]$\!\!}(92.7779,982.3622);
  \path[draw=black,line join=miter,line cap=butt,line width=0.650pt]
  (122.0000,1002.3622) .. controls (103.4643,1002.3622) and (95.4296,990.3419) ..
  node[above left=0.12cm,at start] {$[\mathfrak a,\!1]$\!\!}(91.6415,985.0386);
  \path[draw=black,line join=miter,line cap=butt,line width=0.650pt]
  (58.0000,972.3622) .. controls (67.6795,972.3622) and (80.2525,972.1096) ..
  node[above right=0.12cm,at start] {\!\!$[m,\!1]$}(88.3338,980.4434);
  \path[draw=black,line join=miter,line cap=butt,line width=0.650pt]
  (58.0000,992.3622) .. controls (66.6752,992.3622) and (80.2525,992.6147) ..
  node[above right=0.12cm,at start] {\!\!$[m \t 1,\!1]$}(87.8287,983.5234);
  \path[fill=black] (89.651039,982.15656) node[circle, draw, line width=0.65pt, minimum width=5mm, fill=white, inner sep=0.25mm] (text6017) {$[\ass,\!1]$ };
\end{tikzpicture}}}\ \ \ \quad
\vcenter{\hbox{
\begin{tikzpicture}[y=0.95pt, x=1.1pt, inner sep=0pt, outer sep=0pt, every text node part/.style={font=\tiny} ]
  \path[draw=black,line join=miter,line cap=butt,line width=0.650pt]
  (70.0000,922.3622) .. controls (101.2330,922.3622) and (143.6145,921.7598) ..
  node[above right=0.12cm,at start] {\!\!$[\mathfrak a,\!1]$}(151.3870,930.2041);
  \path[draw=black,line join=miter,line cap=butt,line width=0.650pt]
  (126.7723,969.6175) .. controls (128.2875,962.0413) and (140.0945,940.6924) ..
  node[below left=0.08cm,pos=0.627,rotate=-52] {$[1,\! \tilde{\mathfrak a}]$}(151.2391,934.9696);
  \path[draw=black,line join=miter,line cap=butt,line width=0.650pt]
  (127.6786,973.3558) .. controls (136.1123,980.1329) and (156.4071,982.3622) ..
  node[above left=0.12cm,at end] {$[1, \bar m]$\!}(190.0000,982.3622);
  \path[draw=black,line join=miter,line cap=butt,line width=0.650pt]
  (152.7723,930.8750) .. controls (155.6057,928.5814) and (169.8180,922.9558) ..
  node[above left=0.12cm,at end] {$\tilde{\mathfrak a}$\!}(190.0000,922.3622);
  \path[draw=black,line join=miter,line cap=butt,line width=0.650pt]
  (153.1294,933.7322) .. controls (155.4312,936.0339) and (161.1229,962.3622) ..
  node[above left=0.12cm,at end] {$\tilde{\mathfrak a}$\!}(190.0000,962.3622)
  (127.8878,970.8146) .. controls (136.2744,968.2835) and (145.9898,957.8154) ..
  node[above right=0.05cm,pos=0.25,rotate=-32] {$[1,\![m,\!1]]$}(159.5167,950.3078)
  (164.0504,947.9860) .. controls (171.4928,944.4990) and (180.0280,942.1210) ..
  node[above left=0.12cm,at end] {$[1,\! \bar m]$\!}(190.0000,942.3622);
  \path[draw=black,line join=miter,line cap=butt,line width=0.650pt]
  (70.0000,982.3622) .. controls (94.7819,982.3622) and (118.6172,978.5843) ..
  node[above right=0.12cm,at start] {\!$[1,\!\bar m]$}(125.9408,972.7759);
  \path[draw=black,line join=miter,line cap=butt,line width=0.650pt]
  (70.0000,962.3622) .. controls (104.0413,962.3622) and (95.4124,932.6539) ..
  node[above right=0.12cm,at start] {\!$\tilde{\mathfrak a}$}(151.2919,932.6539)
  (70.0000,942.3622) .. controls (81.3043,942.3622) and (90.5727,944.4985) ..
  node[above right=0.12cm,at start] {\!\!$[1 \t m,\!1]$}(98.2107,947.7957)
  (103.1302,950.1901) .. controls (113.3893,955.7771) and (120.2698,963.4297) ..
  node[above left=0.05cm,pos=0.57] {$[1,\![1,\!\bar m]]$}(125.0912,969.9746);
  \path[fill=black] (152.0495,932.40131) node[circle, draw, line width=0.65pt, minimum width=5mm, fill=white, inner sep=0.25mm] (text3313) {$\tilde \pi$     };
  \path[fill=black] (125.89511,972.28424) node[circle, draw, line width=0.65pt, minimum width=5mm, fill=white, inner sep=0.25mm] (text3317) {$[1,\!\tilde \ass]$   };
\end{tikzpicture} \rlap{ .}}}
\end{equation*}

\begin{Prop}\label{prop:right-hom-from-bilimit}
With notation as above, the right hom $\h V W$ and the  conical bilimit $\{\Delta 1 , F^{VW}\}$ represent pseudonaturally equivalent functors $\V \to \cat{CAT}$; in particular, the one exists if and only if the other does.
\end{Prop}
\begin{proof}
It suffices to exhibit an equivalence, pseudonatural in $A$, between transformations $\Delta 1 \to \V(A, F^{VW})$ and module morphisms $A \otimes V \to W$. As before, we abbreviate $F^{VW}$ to $F$.
The arguments of Proposition~\ref{prop:tensor-product-from-bicolimit} dualise without difficulty to show that the category of transformations $\Delta 1 \to \V(A, F)$ is equivalent to the category $\E_A$ whose objects are given by a family of $1$-cells $f_x \colon A \to [Vx, Wx]$ and a family of $2$-cells
\begin{equation*}
\cd{
 A \ar[d]_{f_y} \ar[rr]^-{f_x}
  \rtwocell{drr}{\bar f_{xy}} & &
 [Vx,Wx] \ar[d]^{[m,1]} \\
 [Vy,Wy] \ar[r]_-{[1, \bar m]} &
 [Vy,[\B(x,y),Wx]] \ar[r]_-{\tilde a} &
 [Vy \otimes B(x,y), Wx]
}
\end{equation*}
satisfying two coherence axioms; and whose morphisms $f \to g$ are given by families of $2$-cells $f_x \Rightarrow g_x$ satisfying one coherence axiom. We claim that $\E_A$ is in turn equivalent to the category ${}_\bullet \cat{Mod}_\B(A \otimes V, W)$.

Indeed, from an object $f$ of $\E_A$, we obtain a module morphism $\phi \colon A \otimes V \to W$ whose $1$-cell component at $x$ is the transpose $\mathfrak e \circ (f_x \otimes 1)$ of $f_x$; and whose $2$-cell component at $(x,y)$ is the composite
\begin{equation*}
\begin{tikzpicture}[y=0.80pt, x=0.9pt,yscale=-1, inner sep=0pt, outer sep=0pt, every text node part/.style={font=\tiny} ]
  \path[draw=black,line join=miter,line cap=butt,line width=0.650pt]
  (40.0000,1032.3622) --
  node[above right=0.12cm,at start] {\!$m$}(70.0000,1032.3622);
  \path[draw=black,line join=miter,line cap=butt,line width=0.650pt]
  (71.5852,1033.9566) .. controls (78.9088,1041.7853) and (128.0694,1041.7241) ..
  node[above=0.07cm] {$\mathfrak e$}(137.9184,1033.3903);
  \path[draw=black,line join=miter,line cap=butt,line width=0.650pt]
  (40.0000,1012.3622) .. controls (74.1434,1012.3861) and (67.7142,1005.8572) ..
  node[above right=0.12cm,at start] {\!$\mathfrak e \t 1$}
  (87.2500,1005.8622) .. controls (112.4830,1005.8686) and (102.7727,1027.6019) .. (137.4567,1030.1147)
  (40.0000,992.3622) .. controls (74.5048,992.3622) and (124.0884,991.0930) ..
  node[above right=0.12cm,at start] {\!\!$(f \t 1) \t 1$}(162.4490,993.1396)
  (173.6035,993.8584) .. controls (198.1894,995.7487) and (216.2591,999.4143) ..
  node[above=0.07cm,pos=0.58] {\,\,$f \t 1$}(219.3290,1006.3341)
  (71.4657,1030.0893) .. controls (74.1681,1024.9449) and (87.9498,1020.8943) ..
  node[above=0.1cm,pos=0.61] {$\bar m \t 1$}(106.1154,1017.7783)
  (111.9550,1016.8306) .. controls (121.1256,1015.4209) and (131.1534,1014.2249) ..
  node[above=0.07cm,rotate=6] {$([1,\!\bar m] \t 1) \t 1$}(141.2867,1013.2247)
  (147.3783,1012.6475) .. controls (180.7568,1009.6123) and (213.6487,1008.6577) ..
  node[above=0.07cm,pos=0.32,rotate=2] {$[1,\!\bar m] \t 1$}(219.2857,1009.1479)
  (139.4337,1027.1530) .. controls (146.6903,982.1388) and (210.3907,982.3622) ..
  (280.0000,982.3622) .. controls (280.6531,982.3622) and (309.3412,982.3444) ..
  node[above left=0.12cm,at end] {$\mathfrak a$\!}(310.0000,982.3622);
  \path[draw=black,line join=miter,line cap=butt,line width=0.650pt]
  (140.4797,1030.2014) .. controls (152.3490,1022.3727) and (210.9340,1029.0378) ..
  node[above=0.07cm,rotate=6] {$\tilde a \t 1$}(218.9898,1011.4566);
  \path[draw=black,line join=miter,line cap=butt,line width=0.650pt]
  (221.7424,1010.6377) .. controls (229.3185,1020.9918) and (253.5714,1021.4693) ..
  node[above=0.07cm,rotate=-20,pos=0.53] {$[m,\!1]\t 1$}(259.1071,1033.7907);
  \path[draw=black,line join=miter,line cap=butt,line width=0.650pt]
  (140.7322,1033.2423) .. controls (150.5812,1039.8083) and (250.5357,1041.2907) ..
  node[above=0.07cm] {$\mathfrak e$}(258.9286,1035.9336);
  \path[draw=black,line join=miter,line cap=butt,line width=0.650pt]
  (221.0714,1007.8979) .. controls (254.9515,976.5311) and (276.1064,1022.2496) ..
  (303.8661,1022.3182) .. controls (304.5893,1022.3171) and (309.3003,1022.3622) ..
  node[above left=0.12cm,at end] {$f \t 1$\!}(310.0000,1022.3622)
  (261.2500,1033.6122) .. controls (265.6661,1030.2297) and (269.6107,1021.4001) ..
  node[below right=0.07cm,pos=0.63] {$1 \t m$}(276.7107,1013.9931)
  (280.4826,1010.4873) .. controls (286.1573,1005.8331) and (293.6134,1002.3905) ..
  (304.0378,1002.4072) .. controls (304.6784,1002.4082) and (309.3316,1002.3622) ..
  node[above left=0.12cm,at end] {$1 \t m$\!}(310.0000,1002.3622);
  \path[draw=black,line join=miter,line cap=butt,line width=0.650pt]
  (261.2500,1035.7550) .. controls (263.8409,1040.0732) and (291.4268,1042.3691) ..
  (300.0000,1042.3622) .. controls (300.4830,1042.3618) and (309.5321,1042.3622) ..
  node[above left=0.12cm,at end] {$\mathfrak e$\!}(310.0000,1042.3622);
  \path[fill=black] (67.071426,1032.7194) node[circle, draw, line width=0.65pt, minimum width=5mm, fill=white, inner sep=0.25mm] (text9821) {$\varepsilon^{-1}$ };
  \path[fill=black] (138.57143,1030.9336) node[circle, draw, line width=0.65pt, minimum width=5mm, fill=white, inner sep=0.25mm] (text9825) {$\varepsilon^{-1}$  };
  \path[fill=black] (220.35715,1009.505) node[circle, draw, line width=0.65pt, minimum width=5mm, fill=white, inner sep=0.25mm] (text9829) {$\bar f$  };
  \path[fill=black] (260,1035.5765) node[circle, draw, line width=0.65pt, minimum width=5mm, fill=white, inner sep=0.25mm] (text9833) {$\varepsilon$ };
\end{tikzpicture}
\end{equation*}
in $\V((A \otimes Vy) \otimes \B(x,y), Wx)$. This assignation provides the action on objects of a functor $G_A \colon \E_A \to {}_\bullet \cat{Mod}_\B(A \otimes V, W)$, which on morphisms acts by taking componentwise transposes.
Conversely, given a module morphism $\phi \colon A \otimes V \to W$, we obtain an object of $\E_A$ whose  $1$-cell component at $x$ is the transpose of $\phi_x$, and whose $2$-cell component at $(x,y)$ is determined by the $2$-cell
\begin{equation*}
  \vcenter{\hbox{\begin{tikzpicture}[y=0.9pt, x=0.9pt,yscale=-1, xscale=-1,inner sep=0pt, outer sep=0pt, every text node part/.style={font=\tiny} ]
  \path[draw=black,line join=miter,line cap=butt,line width=0.650pt]
  (10.0000,1022.3622) .. controls (29.1610,1022.3622) and (37.3064,1017.8382) ..
  node[above left=0.12cm,at start] {$\mathfrak e$\!}(40.3368,1014.1764);
  \path[draw=black,line join=miter,line cap=butt,line width=0.650pt]
  (10.0000,1002.3622) .. controls (26.9927,1002.3622) and (37.4296,1005.3926) ..
  node[above left=0.12cm,at start] {$[m,\!1] \t 1$\!\!}(40.4848,1009.4333);
  \path[draw=black,line join=miter,line cap=butt,line width=0.650pt]
  (43.5178,1012.3101) .. controls (49.5892,1018.2030) and (74.0211,1018.0244) ..
  node[below=0.07cm,pos=0.5] {$\mathfrak e$}(79.3782,1012.3101);
  \path[draw=black,line join=miter,line cap=butt,line width=0.650pt]
  (44.0229,1009.2797) .. controls (48.2846,1001.2761) and (52.8467,995.5068) ..
  node[below left=0.03cm,pos=0.5] {$1 \t m$}(57.6078,991.5294)
  (62.4665,988.1127) .. controls (77.8423,979.2206) and (94.8533,986.6208) ..
  node[above=0.07cm,pos=0.2] {$1 \t m$}(110.3338,996.5043)
  (10.0000,982.3622) .. controls (80.1675,982.3622) and (79.6282,1009.8888) ..
  node[above left=0.12cm,at start] {$\bar \phi \t 1$\!}(79.6282,1009.8888);
  \path[draw=black,line join=miter,line cap=butt,line width=0.650pt]
  (82.0749,1010.6934) .. controls (97.2272,1010.4409) and (105.8135,1004.6325) ..
  node[below left=0.07cm,pos=0.35] {$\phi$}(110.8643,998.3190);
  \path[draw=black,line join=miter,line cap=butt,line width=0.650pt]
  (112.7583,996.1725) .. controls (122.2285,989.8590) and (147.3560,993.2683) ..
  node[above=0.07cm,pos=0.45] {$\phi \t 1$}(152.7856,997.5614);
  \path[draw=black,line join=miter,line cap=butt,line width=0.650pt]
  (112.8846,998.6979) .. controls (119.7031,1011.5773) and (130.0845,1026.7849) ..
  node[above left=0.07cm,pos=0.4] {$m$}(145.7439,1026.7849);
  \path[draw=black,line join=miter,line cap=butt,line width=0.650pt] (147.3024,1028.2875) .. controls (151.8227,1034.0600) and (214.5159,1032.3018) .. node[below=0.07cm] {$\mathfrak e$}(220.8694,1027.3609);
  \path[draw=black,line join=miter,line cap=butt,line width=0.650pt]
  (111.7212,992.7452) .. controls (87.4520,972.0066) and (167.3866,969.1566) ..
  node[above=0.06cm,at end] {$\mathfrak a^\centerdot$}(176.2562,970.1305) .. controls (206.4389,973.4448) and (211.8253,1003.0561) .. (220.8694,1022.0576)
  (154.3314,996.5821) .. controls (156.8346,990.4589) and (173.2702,983.9896) ..
  node[above=0.03cm,pos=0.62,rotate=-17] {$(\bar \phi \t 1) \t 1$}(193.6962,979.2788)
  (199.8722,977.9331) .. controls (216.3240,974.5508) and (234.5999,972.3622) ..
  node[above right=0.12cm,at end] {\!$\bar \phi \t 1$}(250.0000,972.3622)
  (147.1978,1025.0481) .. controls (147.9681,1022.3949) and (155.1489,1017.7858) ..
  node[below left=0.06cm,pos=0.57,rotate=-23] {$\bar m \t 1$}(165.9589,1012.8513)
  (172.0426,1010.2020) .. controls (182.3257,1005.9223) and (194.9072,1001.6163) ..
  node[above=0.04cm,pos=0.45,rotate=-16] {$([1,\!\bar m] \t 1)\t 1$}(208.0415,998.3071)
  (213.0118,997.1098) .. controls (225.5347,994.2331) and (238.3559,992.3624) ..
  node[above right=0.12cm,at end] {\!\!$[1,\!\bar m] \t 1$}(250.0000,992.3624)
  (154.5533,998.8241) .. controls (157.0608,1013.4205) and (208.4433,1024.7310) ..
  node[below right=0.06cm,pos=0.6] {$\mathfrak e \t 1$}(221.1219,1025.0881);
  \path[draw=black,line join=miter,line cap=butt,line width=0.650pt] (223.1639,1026.9350) .. controls (228.4672,1031.4807) and (234.0901,1032.3622) .. node[above right=0.12cm,at end] {\!$\mathfrak e$}(250.0000,1032.3622);
  \path[draw=black,line join=miter,line cap=butt,line width=0.650pt] (222.9636,1023.6520) .. controls (227.8881,1019.3589) and (233.3312,1012.3622) .. node[above right=0.12cm,at end] {\!$\tilde{\mathfrak a} \t 1$}(250.0000,1012.3622);
  \path[fill=black] (41.750008,1011.5526) node[circle, draw, line width=0.65pt, minimum width=5mm, fill=white, inner sep=0.25mm] (text3199) {$\varepsilon^{-1}$     };
  \path[fill=black] (111.11678,996.55127) node[circle, draw, line width=0.65pt, minimum width=5mm, fill=white, inner sep=0.25mm] (text3620) {$\bar \phi^{-1}$     };
  \path[fill=black] (153.79572,997.81396) node[circle, draw, line width=0.65pt, minimum width=5mm, fill=white, inner sep=0.25mm] (text3624) {$\varepsilon \t 1$     };
  \path[fill=black] (142.38202,1026.6033) node[circle, draw, line width=0.65pt, minimum width=5mm, fill=white, inner sep=0.25mm] (text3628) {$\varepsilon$     };
  \path[fill=black] (222.13203,1026.8558) node[circle, draw, line width=0.65pt, minimum width=5mm, fill=white, inner sep=0.25mm] (text3632) {$\varepsilon$   };
  \path[fill=black] (80.239113,1012.8493) node[circle, draw, line width=0.65pt, minimum width=5mm, fill=white, inner sep=0.25mm] (text3199-4) {$\varepsilon^{-1}$   };
\end{tikzpicture}}}
\end{equation*}
between the exponential transposes of its domain and codomain. This provides the action on objects of a functor $H_A \colon {}_\bullet \cat{Mod}_\B(A \otimes V, W) \to \E_A$, which on morphisms again acts by taking componentwise exponential transposes. The units and counits of the adjoint equivalences~\eqref{eq:exptransp} now provide the components of $2$-cells witnessing $G_A$ and $H_A$ as pseudoinverse to each other; thus $\E_A$ and ${}_\bullet \cat{Mod}_\B(A \otimes V, W)$ are equivalent, as claimed. Composing with the equivalence between $\E_A$ and $\cat{Hom}(\Delta 1, \V(A, F))$, we obtain equivalences of categories
\begin{equation*}
\cat{Hom}(\Delta 1, \V(A, F)) \to {}_\bullet \cat{Mod}_\B(A \otimes V, W)\rlap{ ,}
\end{equation*}
which we may without difficulty verify are pseudonatural in $A$, as required.
\end{proof}

\begin{Cor}\label{cor:existence-of-right-homs}
If $\B$ is a small $\V$-bicategory, and $\V$ is complete and right closed, then the right hom $\h V W$ exists for all right $\B$-modules $V$ and $W$.
\end{Cor}

We could now, as we did for tensor products, go on to discuss the construction of right homs between general bimodules; as before, we would find that these can be constructed pointwise: given $N \colon \B \tor \C$ and $P \colon \A \tor \C$, the right hom $\h N P \colon \A \tor \B$ may be defined by $\h N P(b,a) = \h {N(\thg, b)} {P(\thg, a)}$, so long as each of these pointwise right homs exists in $\V$. We shall not prove this fact, since we do not need it in what follows; nonetheless, it will be useful to have a recognition principle for the more general right homs.
\begin{Prop}\label{prop:right-hom-pointwise}
Let $\phi \colon M, N \to P$ be an $(\A,\B,\C)$-module bimorphism. $\phi$ exhibits $M$ as $\h N P$ if and only if, for each $a \in \A$ and $b \in \B$, $\phi_{ab\thg} \colon M(b,a) \otimes N(\thg, b) \to P(\thg, a)$ exhibits $M(b,a)$ as $\h {N(\thg, b)} {P(\thg, a)}$.\end{Prop}
\begin{proof}
The ``only if'' direction (of which we shall not make any serious use in what follows) follows from the omitted construction above; for the ``if'' direction, we must show that, for each $L \colon A \tor B$, the functor
\begin{equation}\label{eq:right-hom-univ}
\begin{aligned}
{}_\A \cat{Mod}_\B(L, M) & \to \cat{Bimor}_{\A\B\C}(L, N; P)\\
\gamma &\mapsto \phi \circ (\gamma, 1)
\end{aligned}
\end{equation}
is an equivalence of categories.
Now, given a bimorphism $\psi \colon L, N \to P$, the universal property of each $\phi_{ab}$ induces morphisms $\theta_{ab} \colon L(b,a) \to M(b,a)$ and invertible module transformations $\bar \theta \colon \phi_{ab} \circ (\theta_{ab} \otimes 1) \Rightarrow \psi_{ab}$. We will show that the $\theta_{ab}$'s comprise the components of an $\A$-$\B$-bimodule morphism. To obtain the $2$-cells exhibiting compatibility with the left $\A$-\ and right $\B$-actions, it suffices by universality of the $\phi_{ab\thg}$'s (more precisely, by the fully faithfulness of~\eqref{eq:compinduced}) to give invertible module transformations
\begin{equation*}
\cd{
(\A(a,a') \otimes L(b,a)) \otimes N(c,b)
 \ar[r]^-{m \otimes 1} \dtwocell{drr}{}
 \ar[d]_{(1 \otimes \theta) \otimes 1} &
L(b,a') \otimes N(c,b) \ar[r]^-{\theta \otimes 1} &
M(b,a') \otimes N(c,b) \ar[d]^{\phi} \\
(\A(a,a') \otimes M(b,a)) \otimes N(c,b)
 \ar[r]^-{m \otimes 1} &
M(b,a') \otimes N(c,b)
 \ar[r]^-{\phi} &
P(c,a')
}
\end{equation*}
and
\begin{equation*}
\cd{
(L(b,a) \otimes \B(b',b)) \otimes N(c,b')
 \ar[r]^-{m \otimes 1} \rtwocell{drr}{}
 \ar[d]_{(\theta \otimes 1) \otimes 1} &
L(b',a) \otimes N(c,b') \ar[r]^-{\theta \otimes 1} &
M(b',a) \otimes N(c,b') \ar[d]^{\phi} \\
(M(b,a) \otimes \B(b',b)) \otimes N(c,b')
 \ar[r]^-{m \otimes 1} &
M(b',a) \otimes N(c,b')
 \ar[r]^-{\phi} &
P(c,a)
}
\end{equation*}
respectively.
We obtain such by taking their respective components to be
\begin{equation*}
\vcenter{\hbox{\begin{tikzpicture}[y=0.8pt, x=0.8pt,yscale=-1, inner sep=0pt, outer sep=0pt, every text node part/.style={font=\tiny} ]
  \path[draw=black,line join=miter,line cap=butt,line width=0.650pt]
  (4.0000,1012.3622) .. controls (31.5000,1012.6122) and (40.6250,1006.6122) ..
  node[above right=0.12cm,at start] {\!$\phi$}(42.6250,1003.8622);
  \path[draw=black,line join=miter,line cap=butt,line width=0.650pt]
  (4.0000,992.3622) .. controls (33.6250,992.3622) and (40.5000,997.8622) ..
  node[above right=0.12cm,at start] {\!$\theta \t 1$}(42.3750,1000.8622);
  \path[draw=black,line join=miter,line cap=butt,line width=0.650pt]
  (4.0000,972.3622) .. controls (39.8759,972.3622) and (75.2093,982.6336) ..
  node[above right=0.12cm,at start] {\!$m \t 1$}(86.7500,989.6122);
  \path[draw=black,line join=miter,line cap=butt,line width=0.650pt]
  (46.5000,1002.1122) .. controls (64.3750,1002.7372) and (84.6250,996.4872) ..
  node[above left=0.07cm,pos=0.53] {$\psi$}(87.0000,992.8622);
  \path[draw=black,line join=miter,line cap=butt,line width=0.650pt]
  (89.4090,988.4480) .. controls (95.9497,982.7911) and (128.2500,989.8622) ..
  node[above=0.07cm,pos=0.45] {$1 \t \psi$}(128.2500,989.8622);
  \path[draw=black,line join=miter,line cap=butt,line width=0.650pt]
  (130.9090,991.7764) .. controls (136.3891,996.9029) and (171.5022,1001.0086) ..
  node[above=0.07cm,pos=0.5] {$1 \t \phi$}(178.7500,1001.3622);
  \path[draw=black,line join=miter,line cap=butt,line width=0.650pt]
  (89.3964,991.7409) .. controls (89.8611,1011.4166) and (165.4110,1015.2115) ..
  node[above=0.09cm,pos=0.35] {$m$}(178.2500,1003.3622);
  \path[draw=black,line join=miter,line cap=butt,line width=0.650pt]
  (130.9090,988.9103) .. controls (136.5587,981.1558) and (158.0798,976.6216) ..
  node[above=0.03cm,rotate=14,pos=0.5] {$1 \t (\theta \t 1)$}(180.9029,974.3170)
  (189.5240,973.5530) .. controls (200.3042,972.7247) and (210.9681,972.3622) ..
  node[above left=0.12cm,at end] {$(1 \t \theta) \t 1$\!}(220.0000,972.3622)
  (86.7500,986.1122) .. controls (48.2860,961.2680) and (129.7386,960.4395) ..
  node[above=0.07cm,pos=0.7] {$\mathfrak a$}(145.7430,962.5766) .. controls (162.9479,964.8740) and (218.6948,976.3184) .. (180.1250,999.2372);
  \path[draw=black,line join=miter,line cap=butt,line width=0.650pt]
  (220.0000,992.3622) .. controls (200.5514,992.3622) and (190.4712,998.1850) ..
  node[above left=0.12cm,at start] {$m \t 1$\!}(181.1465,1001.3873);
  \path[draw=black,line join=miter,line cap=butt,line width=0.650pt]
  (181.4419,1003.9328) .. controls (185.7730,1010.1200) and (199.4939,1012.3622) ..
  node[above left=0.12cm,at end] {$\phi$\!}(220.0000,1012.3622);
  \path[fill=black] (43.497475,1002.8647) node[circle, draw, line width=0.65pt, minimum width=5mm, fill=white, inner sep=0.25mm] (text3131) {$\bar\theta^{-1}$     };
  \path[fill=black] (129.75409,990.13678) node[circle, draw, line width=0.65pt, minimum width=5mm, fill=white, inner sep=0.25mm] (text3135) {$1 \t \bar\theta$     };
  \path[fill=black] (88.034798,990.13678) node[circle, draw, line width=0.65pt, minimum width=5mm, fill=white, inner sep=0.25mm] (text3139) {$\bar \psi$     };
  \path[fill=black] (179.35207,1003.2183) node[circle, draw, line width=0.65pt, minimum width=5mm, fill=white, inner sep=0.25mm] (text3143) {$\bar \phi^{-1}$   };
\end{tikzpicture}}}\qquad \text{and}\qquad
\vcenter{\hbox{\begin{tikzpicture}[y=0.8pt, x=0.8pt,yscale=-1, inner sep=0pt, outer sep=0pt, every text node part/.style={font=\tiny} ]
  \path[draw=black,line join=miter,line cap=butt,line width=0.650pt]
  (4.0000,1012.3622) .. controls (27.1731,1012.3622) and (35.5935,1007.7769) ..
  node[above right=0.12cm,at start] {\!$\phi$}(38.5858,1004.2304);
  \path[draw=black,line join=miter,line cap=butt,line width=0.650pt]
  (4.0000,992.3622) .. controls (24.9809,992.3622) and (36.1109,997.0468) ..
  node[above right=0.12cm,at start] {\!$m \t 1$}(38.5858,1000.6707);
  \path[draw=black,line join=miter,line cap=butt,line width=0.650pt]
  (41.4420,1003.5718) .. controls (50.5649,1015.0696) and (106.0244,1014.4661) ..
  node[above=0.07cm,pos=0.45] {$\phi$}(112.2116,1007.7486);
  \path[draw=black,line join=miter,line cap=butt,line width=0.650pt]
  (113.9793,1006.8647) .. controls (129.0054,1007.5718) and (151.0262,995.4401) ..
  node[above=0.07cm,pos=0.35] {$\psi$}(154.3849,989.4297);
  \path[draw=black,line join=miter,line cap=butt,line width=0.650pt]
  (4.0000,972.3622) .. controls (19.6411,972.3622) and (36.4594,973.5749) ..
  node[above right=0.12cm,at start] {\!\!$(\theta \t 1) \t 1$}(52.0410,975.9283)
  (59.8046,977.2155) .. controls (89.0898,982.5231) and (112.1889,992.0778) ..
  node[above right=0.07cm,pos=0.25] {$\theta \t 1$}(111.4767,1005.3536)
  (41.2652,1000.0363) .. controls (36.2891,964.0822) and (131.9065,952.6088) ..
  node[above=0.07cm,pos=0.5] {$\mathfrak a$}(154.0314,984.4799)
  (42.5027,1001.4505) .. controls (47.0795,996.7712) and (69.2468,992.0843) ..
  node[above left=0.07cm,pos=0.65] {$1 \t m$}(92.8413,989.0302)
  (101.7278,987.9643) .. controls (125.0478,985.3860) and (147.7926,984.5874) ..
  node[above=0.07cm,pos=0.35] {$1 \t m$}(154.5617,987.1316);
  \path[draw=black,line join=miter,line cap=butt,line width=0.650pt]
  (156.3295,985.7174) .. controls (162.2500,977.2372) and (205.1172,972.3622) ..
  node[above left=0.12cm,at end] {$m \t 1$\!}(230.0000,972.3622);
  \path[draw=black,line join=miter,line cap=butt,line width=0.650pt]
  (156.5062,987.8387) .. controls (159.1312,994.8387) and (185.9272,1002.3344) ..
  node[above right=0.07cm,pos=0.55] {$\psi$}(196.3052,1002.3344);
  \path[draw=black,line join=miter,line cap=butt,line width=0.650pt]
  (198.6033,1001.9808) .. controls (201.2283,995.9808) and (219.6220,992.3622) ..
  node[above left=0.12cm,at end] {$\theta \t 1$\!}(230.0000,992.3622);
  \path[draw=black,line join=miter,line cap=butt,line width=0.650pt]
  (198.6033,1004.2789) .. controls (200.9783,1010.2789) and (217.9896,1012.3622) ..
  node[above left=0.12cm,at end] {$\phi$\!}(230.0000,1012.3622);
  \path[fill=black] (40,1002.3622) node[circle, draw, line width=0.65pt, minimum width=5mm, fill=white, inner sep=0.25mm] (text3446) {$\bar \phi$     };
  \path[fill=black] (197.25,1003.3622) node[circle, draw, line width=0.65pt, minimum width=5mm, fill=white, inner sep=0.25mm] (text3450) {$\bar \theta$     };
  \path[fill=black] (155.5,987.36218) node[circle, draw, line width=0.65pt, minimum width=5mm, fill=white, inner sep=0.25mm] (text3454) {$\bar\psi^{-1}$     };
  \path[fill=black] (112.25,1007.1122) node[circle, draw, line width=0.65pt, minimum width=5mm, fill=white, inner sep=0.25mm] (text3458) {$\bar \theta^{-1}$   };
\end{tikzpicture}\rlap{ .}}}
\end{equation*}

The assignation $\psi \mapsto \theta$ just described provides the action on objects of a pseudoinverse to the functor~\eqref{eq:right-hom-univ}. The remaining data of this functor, together with the $2$-cells exhibiting it as pseudoinverse are obtained directly from the universality of each $M(b,a)$.
\end{proof}

\section{The Yoneda Lemma}
\label{sec:yoneda}

\subsection{Representables and the Yoneda lemma}
Let $\B$ be a $\V$-bicategory and let $W$ be a right $\B$-module. As in Example~\ref{ex:yoneda-bimorphisms}, the morphisms $m_{bb'} \colon Wb \otimes \B(b',b) \to Wb'$ comprise the $1$-cell components of a module bimorphism $W, \B \to W$, whose corresponding $2$-cell components are given by the associativity constraints of the action of $\B$ on $W$.

\begin{Prop}[Yoneda lemma] For each $b \in \B$, the right module morphism $m_{b\thg} \colon Wb \otimes \B(\thg, b) \to W$ exhibits $Wb$ as $\h{\B(\thg, b)} W$.\label{prop:yoneda}
\end{Prop}
\begin{proof}
We must show that for each $A \in \V$, the functor $F \defeq m \circ (\thg \otimes 1) \colon \V(A, Wb) \to \cat{Mod}_\B(A \otimes \B(\thg, b), W)$ is an equivalence of categories. We will do so by exhibiting an explicit pseudoinverse functor $G \colon {}_\bullet \cat{Mod}_\B(A \otimes \B(\thg, b), W) \to \V(A,Wb)$. On objects, this functor sends a module morphism $\phi \colon A \otimes \B(\thg,b) \to W$ to the $1$-cell
\begin{equation*}
A \xrightarrow{\mathfrak r^\centerdot} A \otimes I \xrightarrow {1 \otimes j} A \otimes \B(b,b) \xrightarrow{\phi_b} Wb\rlap{ ;}
\end{equation*}
on morphisms, it sends a module transformation $\Gamma \colon \phi \Rightarrow \psi$ to the $2$-cell obtained by whiskering $\Gamma_b$ with $(1 \otimes j) \circ \mathfrak r^\centerdot$. We now show that this $G$ is pseudoinverse to $F$. On the one hand, if given $f \in \V(A,Wb)$ then $GFf$ is the composite morphism
\begin{equation*}
A \xrightarrow{\mathfrak r^\centerdot} A \otimes I \xrightarrow {1 \otimes j} A \otimes \B(b,b) \xrightarrow{f \otimes 1} Wb \otimes \B(b,b) \xrightarrow{m} Wb\rlap{ ,}
\end{equation*}
which is isomorphic to $f$ via the $2$-cell
\begin{equation*}
\vcenter{\hbox{
\begin{tikzpicture}[y=0.65pt, x=0.7pt,yscale=-1, inner sep=0pt, outer sep=0pt, every text node part/.style={font=\tiny} ]
  \path[draw=black,line join=miter,line cap=butt,line width=0.650pt]
  (0.7143,1032.7193) .. controls (33.4402,1032.9653) and (50.7053,1013.5463) ..
  node[above right=0.15cm,at start] {$\!f \t 1$}(73.9303,1001.6762)
  (80.6246,998.5644) .. controls (89.4024,994.9073) and (99.2003,992.5432) ..
  node[above left=0.15cm,at end] {$f\!$}(111.0000,992.7193)
  (0.7143,1011.6479) .. controls (18.2407,1011.6479) and (35.1220,1012.0913) ..
  node[above right=0.15cm,at start] {$\!1 \t j$}(49.5738,1013.5828)
  (57.3368,1014.5222) .. controls (74.8965,1016.9929) and (87.9017,1021.4249) ..
  node[below left=0.08cm,pos=0.45] {$1 \t j$}(93.2143,1029.1479)
  (0.7143,990.5765) .. controls (102.7520,989.9005) and (134.2353,1027.0446) ..
  node[above right=0.15cm,at start] {$\!\mathfrak r^\centerdot$}(98.7857,1030.2194);
  \path[draw=black,line join=miter,line cap=butt,line width=0.650pt]
  (0.7143,1052.0051) .. controls (56.4286,1052.3622) and (84.4286,1044.8621) ..
  node[above right=0.15cm,at start] {$\!m$}(93.7143,1031.6479);
  \path[fill=black] (94.071426,1029.8622) node[circle, draw, line width=0.65pt, minimum width=5mm, fill=white, inner sep=0.25mm] (text3445) {$\ru$};
\end{tikzpicture}
}}\rlap{ .}
\end{equation*}
The naturality of these isomorphisms in $f$ is immediate, and so we have $GF \cong 1$. On the other hand, given $\phi \colon A \otimes \B(\thg,b) \to W$ a module morphism, $FG\phi$ is the morphism whose $1$-cell component at $x$ is (after using pseudofunctoriality of $\thg\otimes 1$):
{\small
\begin{equation*}
 A \otimes \B(x,b) \xrightarrow{\mathfrak r^\centerdot \otimes 1}
(A \otimes I) \otimes \B(x,b) \xrightarrow {(1 \otimes j) \otimes 1}
(A \otimes \B(b,b)) \otimes \B(x,b) \xrightarrow{\phi_b \otimes 1}
Wb \otimes \B(x,b) \xrightarrow{m} Wx\rlap{.}
\end{equation*}
}
There is now an invertible module transformation $FG\phi \Rightarrow \phi$ with components
\begin{equation*}
\vcenter{\hbox{
\begin{tikzpicture}[y=0.80pt, x=0.8pt,yscale=-1, inner sep=0pt, outer sep=0pt, every text node part/.style={font=\tiny} ]
  \path[draw=black,line join=miter,line cap=butt,line width=0.650pt]
    (8.2143,947.7193) .. controls (43.2143,947.0050) and (140.3571,944.8622) ..
    node[above right=0.15cm,at start] {$\mathfrak r^\centerdot \t 1$}(157.5000,962.0050);
  \path[draw=black,line join=miter,line cap=butt,line width=0.650pt]
    (8.2143,968.7908) .. controls (23.3018,968.7908) and (45.8868,971.1394) ..
    node[above right=0.15cm,at start] {$\!\!(1 \t j) \t 1$}(66.7018,975.4731)
    (74.3039,977.1662) .. controls (91.8110,981.3341) and (107.1382,986.9540) ..
    node[above right=0.06cm,pos=0.45]{$1 \t (j \t 1)$} (114.2857,993.7908)
    (38.9286,997.7193) .. controls (55.7143,987.3622) and (80.7143,955.2193) ..
    node[above left=0.12cm,pos=0.72] {$\mathfrak a$}(157.5000,964.8622);
  \path[draw=black,line join=miter,line cap=butt,line width=0.650pt]
    (38.5714,1000.9337) .. controls (58.2143,1006.2907) and (102.5000,1004.8622)
    .. node[above=0.13cm,pos=0.45]{$1 \t m$} (114.6429,997.3622);
  \path[draw=black,line join=miter,line cap=butt,line width=0.650pt]
    (39.2857,1003.7908) .. controls (45.2089,1021.2793) and (140.6406,1018.4336)
    .. node[above left=0.12cm,at end]{$\phi$} (166.2143,1018.4336);
  \path[draw=black,line join=miter,line cap=butt,line width=0.650pt]
    (116.4286,997.0050) .. controls (135.7143,997.7193) and (136.0714,979.8622) ..
    node[below right=0.07cm,pos=0.55]{$1 \t \mathfrak l$}(158.2143,967.7193);
  \path[draw=black,line join=miter,line cap=butt,line width=0.650pt]
    (8.2143,989.1479) .. controls (33.6008,989.0186) and (35.7143,998.0765) ..
    node[above right=0.15cm,at start] {$\phi \t 1$}(35.7143,998.0765);
  \path[draw=black,line join=miter,line cap=butt,line width=0.650pt]
    (8.2143,1010.9336) .. controls (33.2143,1011.2907) and (36.0714,1002.3622) ..
    node[above right=0.15cm,at start] {$m$}(36.0714,1002.3622);
  \path[fill=black] (158.57143,965.93359) node[circle, draw, line width=0.65pt,
    minimum width=5mm, fill=white, inner sep=0.25mm] (text3445) {$\nu$    };
  \path[fill=black] (114.1893,996.59198) node[circle, draw, line width=0.65pt,
    minimum width=5mm, fill=white, inner sep=0.25mm] (text3449) {$1 \t \lu$    };
  \path[fill=black] (36.42857,1001.2907) node[circle, draw, line width=0.65pt,
    minimum width=5mm, fill=white, inner sep=0.25mm] (text4000) {$\bar \phi$    };
\end{tikzpicture}
}}
\end{equation*}
and the naturality of these isomorphisms in $\phi$ is immediate using the module transformation axiom.
We thus have an isomorphism $FG \cong 1$ as required.
\end{proof}

\begin{Cor}\label{cor:yoneda-with-naturality}
The module bimorphism $m \colon W, \B \to W$ exhibits $W$ as $\h \B W$.
\end{Cor}
\begin{proof}
By the preceding result and Proposition~\ref{prop:right-hom-pointwise}.
\end{proof}

We now use the Yoneda lemma to give a number of different descriptions of $\V$-transformations in terms of module transformations.
\begin{Prop}\label{prop:vnats-as-bimods}
For all $\V$-functors $F, G \colon \B \to \C$, the functor
\begin{equation}\label{eq:vnats-as-bimods}
\begin{aligned}
\V\text-\cat{Bicat}(\B,\C)(F,G) &\to {}_\B \cat{Mod}_{\B}(\B,\, \C(F,G)) \\
\gamma & \mapsto \B \xrightarrow{\hat F} \C(F,F) \xrightarrow{\C(F, \gamma)} \C(F,G)
\end{aligned}
\end{equation}
is an equivalence of categories.
\end{Prop}
\begin{proof}
By the Yoneda lemma, for each $b, b' \in \B$, the composite
\begin{equation}\label{eq:universality-vnat-bimod}
\C(Fb,Gb') \otimes \B(\thg, b) \xrightarrow{1 \otimes F} \C(Fb',Gb') \otimes \C(F,Fb') \xrightarrow{m}
\C(F,Gb')
\end{equation}
exhibits $\C(Fb,Gb')$ as $\h{\B(\thg, b)}{\C(F,Gb')}$. We use this universal property repeatedly in what follows. Let now $\phi \colon \B \to \C(F,G)$ be a bimodule morphism. For each $b \in \B$, we have a right $\B$-module morphism
\begin{equation*}
I \otimes \B(\thg, b) \xrightarrow{\mathfrak l} \B(\thg, b) \xrightarrow{\phi_{b\thg}} \C(F, Gb)\rlap{ ;}
\end{equation*}
applying universality of~\eqref{eq:universality-vnat-bimod}, we obtain morphisms $\gamma_b \colon I \to \C(Fb,Gb)$, together with invertible $2$-cells $\theta \colon m \circ (1 \otimes F) \circ (\gamma \otimes 1) \Rightarrow \phi \circ \mathfrak l$. We claim that the morphisms $\gamma_b$ are the $1$-cell components of a $\V$-transformation $\gamma \colon F \Rightarrow G$. To obtain the corresponding $2$-cell components $\C(\gamma, 1) \circ G \Rightarrow \C(1, \gamma) \circ F \colon \B(b,b') \to \C(Fb,Gb')$, we will show that both the domain and codomain $1$-cells are isomorphic to $\phi_{bb'} \colon \B(b,b') \to \C(Fb,Gb')$. On the one hand, we have a $2$-cell $\C(1, \gamma) \circ F \Rightarrow \phi$ given as on the left below; on the other, we obtain a $2$-cell $\C(\gamma, 1) \circ G \Rightarrow \phi$ by universality of~\eqref{eq:universality-vnat-bimod} applied to the invertible module modification with components as on the right below.
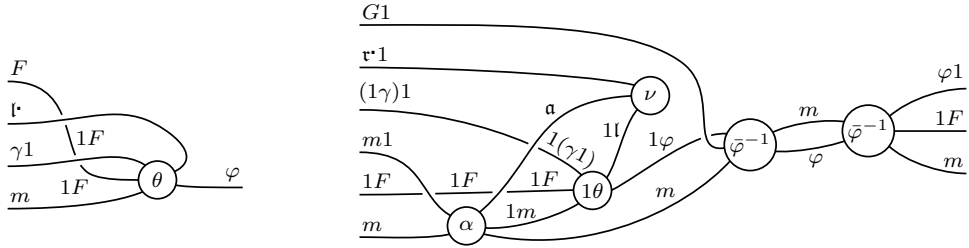
\begin{equation*}
\vcenter{\hbox{
\begin{tikzpicture}[y=0.8pt, x=0.8pt,yscale=-1, inner sep=0pt, outer sep=0pt, every text node part/.style={font=\tiny} ]
  \path[draw=black,line join=miter,line cap=butt,line width=0.650pt] 
  (70.0000,1002.3622) .. controls (129.4049,1002.3622) and (139.6988,990.7056) .. 
  node[above right=0.12cm,at start] {\!$m$}(139.6988,990.7056);
  \path[draw=black,line join=miter,line cap=butt,line width=0.650pt] 
  (70.0000,982.3622) .. controls (104.4486,982.3622) and (125.1417,970.4665) .. 
  node[above right=0.12cm,at start] {\!$\gamma \t 1$}(139.3089,985.9929)
  (70.0000,962.3622) .. controls (101.0319,962.5887) and (117.0120,953.7201) .. 
  node[above right=0.12cm,at start] {\!$\mathfrak l^\centerdot$}
  (134.6386,960.3866) .. controls (144.1646,963.9894) and (168.1766,980.2405) .. (141.6001,987.2706)
  (70.0000,942.3622) .. controls (85.4535,942.3622) and (91.7725,949.4899) .. 
  node[above right=0.12cm,at start] {\!$F$}(95.4263,958.1095)
  (97.1911,962.9167) .. controls (98.7742,967.8175) and (99.9127,972.8580) .. 
  node[ right=0.08cm,pos=0.42] {$1 \t F$}(101.6599,977.1209)
  (104.1343,981.8311) .. controls (105.0706,983.1875) and (106.1630,984.3705) .. 
  (107.4736,985.3261) .. controls (114.3409,990.3335) and (138.1999,988.4863) .. 
  node[below left=0.08cm,pos=0.1] {$1 \t F$}(138.1999,988.4863);
  \path[draw=black,line join=miter,line cap=butt,line width=0.650pt] 
  (141.4111,989.7449) .. controls (150.8946,992.2860) and (156.9820,992.3622) .. 
  node[above left=0.12cm,at end] {$\phi$\!}(180.0000,992.3622);
  \path[fill=black] (140.12048,988.80798) node[circle, draw, line width=0.65pt, minimum width=5mm, fill=white, inner sep=0.25mm] (text3056) {$\theta$   };
\end{tikzpicture}
}} \qquad \quad
\vcenter{\hbox{
\begin{tikzpicture}[y=0.8pt, x=0.8pt,yscale=-1, inner sep=0pt, outer sep=0pt, every text node part/.style={font=\tiny} ]
  \path[draw=black,line join=miter,line cap=butt,line width=0.650pt] 
  (10.0000,922.3622) .. controls (53.0233,922.3622) and (135.6002,926.4432) .. 
  node[above right=0.12cm,at start] {\!$\mathfrak r^\centerdot \t 1$}(144.7585,932.8328);
  \path[draw=black,line join=miter,line cap=butt,line width=0.650pt] 
  (10.0000,1002.3622) .. controls (24.4650,1002.3622) and (54.0663,1004.1032) .. 
  node[above right=0.12cm,at start] {\!$m$}(58.8855,997.9285);
  \path[draw=black,line join=miter,line cap=butt,line width=0.650pt] 
  (61.0843,994.9164) .. controls (97.3219,967.1572) and (86.5944,935.5427) .. 
  node[above left=0.07cm,pos=0.64] {$\mathfrak a$}(145.1656,935.5427)
  (10.0000,962.3622) .. controls (43.8137,962.3622) and (36.4037,988.5007) .. 
  node[above right=0.12cm,at start] {\!$m \t 1$}(58.4337,994.6754)
  (10.0000,982.3622) .. controls (17.7130,982.3622) and (28.1520,982.2148) .. 
  node[above right=0.12cm,at start] {\!$1 \t F$}(39.4966,981.9929)
  (45.4532,981.8714) .. controls (54.3068,981.6841) and (63.5044,981.4608) .. 
  node[above=0.08cm] {$1 \t F$}(72.2374,981.2338)
  (78.9663,981.0560) .. controls (100.0594,980.4892) and (116.9277,979.9525) .. 
  node[above=0.07cm,pos=0.3] {$1 \t F$}(116.9277,979.9525)
  (10.0000,942.3622) .. controls (41.6733,942.3622) and (69.6366,949.9342) .. 
  node[above right=0.12cm,at start] {\!\!$(1 \t \gamma) \t 1$}(89.0103,958.3033)
  (93.6559,960.4007) .. controls (106.7826,966.5962) and (115.0794,972.9177) .. 
  node[above right=0.04cm,rotate=-29,pos=0.025] {$1 \t (\gamma \t 1)$}(116.6265,976.6995);
  \path[draw=black,line join=miter,line cap=butt,line width=0.650pt] 
  (61.8997,997.1647) .. controls (80.0510,1001.4868) and (111.0438,990.4878) .. 
  node[above left=0.05cm] {$1 \t m$}(118.0723,982.1815);
  \path[draw=black,line join=miter,line cap=butt,line width=0.650pt] 
  (61.2973,999.2603) .. controls (119.8335,1015.5447) and (170.9694,983.9529) .. 
  node[above left=0.07cm,pos=0.66] {$m$}(187.5331,961.5530);
  \path[draw=black,line join=miter,line cap=butt,line width=0.650pt] 
  (120.7831,978.8682) .. controls (131.6453,974.6085) and (132.0122,943.8428) .. 
  node[above left=0.05cm,pos=0.55] {$1 \t \mathfrak l$}(144.7913,938.5182);
  \path[draw=black,line join=miter,line cap=butt,line width=0.650pt] 
  (121.0843,981.5791) .. controls (125.9583,985.5096) and (150.1606,963.2547) .. 
  node[above left=0.02cm,pos=0.8] {$1 \t \phi$}(169.7816,955.8733)
  (174.5539,954.3674) .. controls (179.6347,953.1084) and (184.2308,953.1992) .. (187.8624,955.4636)
  (10.0000,902.3622) .. controls (67.7187,902.3622) and (134.6510,900.0402) .. 
  node[above right=0.12cm,at start] {\!$G \t 1$}
  (161.0756,920.2173) .. controls (180.9449,935.3890) and (159.6548,968.9639) .. (187.4664,958.5277);
  \path[draw=black,line join=miter,line cap=butt,line width=0.650pt] 
  (197.9771,956.6423) .. controls (203.3988,950.0158) and (229.6944,944.5923) .. 
  node[above=0.08cm,pos=0.55] {$m$}(243.0954,949.3385);
  \path[draw=black,line join=miter,line cap=butt,line width=0.650pt] 
  (198.4913,960.0123) .. controls (203.0094,965.2834) and (238.2762,959.1069) .. 
  node[below=0.08cm,pos=0.55] {$\phi$}(242.9448,953.9864);
  \path[draw=black,line join=miter,line cap=butt,line width=0.650pt] 
  (254.4559,949.1694) .. controls (261.0584,940.2241) and (279.0896,932.3622) .. 
  node[above left=0.12cm,at end] {$\phi \t 1$\!}(294.0000,932.3622);
  \path[draw=black,line join=miter,line cap=butt,line width=0.650pt] 
  (254.0000,952.3622) -- 
  node[above left=0.12cm,at end] {$1 \t F$\!}(294.0000,952.3622);
  \path[draw=black,line join=miter,line cap=butt,line width=0.650pt] 
  (254.6841,955.5730) .. controls (258.9437,963.8794) and (278.6592,972.3622) .. 
  node[above left=0.12cm,at end] {$m$\!}(294.0000,972.3622);
  \path[fill=black] (60.071724,997.06952) node[circle, draw, line width=0.65pt, minimum width=5mm, fill=white, inner sep=0.25mm] (text3358) {$\ass$     };
  \path[fill=black] (146.10699,935.64697) node[circle, draw, line width=0.65pt, minimum width=5mm, fill=white, inner sep=0.25mm] (text3362) {$\nu$     };
  \path[fill=black] (119.05803,980.3736) node[circle, draw, line width=0.65pt, minimum width=5mm, fill=white, inner sep=0.25mm] (text3366) {$1 \t \theta$     };
  \path[fill=black] (192.83363,958.86224) node[circle, draw, line width=0.65pt, minimum width=5mm, fill=white, inner sep=0.25mm] (text3370) {$\bar \phi^{-1}$     };
  \path[fill=black] (248.11607,952.25977) node[circle, draw, line width=0.65pt, minimum width=5mm, fill=white, inner sep=0.25mm] (text3374) {$\bar \phi^{-1}$   };
\end{tikzpicture}
}}
\end{equation*}
Composing one with the inverse of the other, we obtain the $2$-cell components of the $\V$-transformation $\gamma$ as desired. 
The assignation $\phi \mapsto \gamma$ just described is the action on objects of a functor $H \colon {}_\B \cat{Mod}_{\B}(\B,\, \C(F,G)) \to \V\text-\cat{Bicat}(\B,\C)(F,G)$ that on morphisms, sends a bimodule transformation $\Gamma \colon \phi \Rightarrow \phi'$  to the $\V$-modification whose $2$-cell components are induced by universality of~\eqref{eq:universality-vnat-bimod} applied to the module transformations $\Gamma_{b\thg} \circ \mathfrak l$. 

Finally, we show that this functor $H$ is pseudoinverse to~\eqref{eq:vnats-as-bimods}. On the one hand, starting from $\phi \colon \B \to \C(F,G)$, forming $\gamma = H \phi$ and then applying~\eqref{eq:vnats-as-bimods}, the resultant $\C(F, \gamma) \circ \hat F$ admits an invertible modification to $\phi$ whose $2$-cell components are as on the left above. On the other, given $\gamma \colon F \Rightarrow G$, we may apply~\eqref{eq:vnats-as-bimods} to yield $\C(F, \gamma) \circ \hat F$ and then apply $H$ to this to obtain $\delta \colon F \Rightarrow G$; we now have an invertible $\V$-modification $\delta \Rightarrow \gamma$ with components
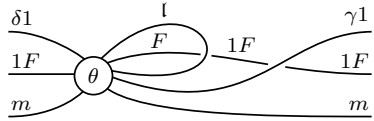
\begin{equation*}
\begin{tikzpicture}[y=0.8pt, x=0.8pt,yscale=-1, inner sep=0pt, outer sep=0pt, every text node part/.style={font=\tiny} ]
  \path[draw=black,line join=miter,line cap=butt,line width=0.650pt] 
  (10.0000,992.3622) .. controls (32.5763,992.5752) and (46.7950,981.0161) .. 
  node[above right=0.12cm,at start] {\!$m$}(48.2859,975.9624);
  \path[draw=black,line join=miter,line cap=butt,line width=0.650pt] 
  (10.0000,972.3622) -- 
  node[above right=0.12cm,at start] {\!$1 \t F$}(47.8599,972.3622);
  \path[draw=black,line join=miter,line cap=butt,line width=0.650pt] 
  (10.0000,952.3622) .. controls (26.6127,951.9362) and (46.5820,964.2893) .. 
  node[above right=0.12cm,at start] {\!$\delta \t 1$}(48.4989,968.3360);
  \path[draw=black,line join=miter,line cap=butt,line width=0.650pt] 
  (51.5421,973.9841) .. controls (60.0000,992.3622) and (132.7077,992.3622) .. 
  node[above left=0.12cm,at end] {$m$\!}(180.0000,992.3622);
  \path[draw=black,line join=miter,line cap=butt,line width=0.650pt] 
  (51.7551,971.2153) .. controls (124.4829,1000.4556) and (136.1518,952.8662) .. 
  node[above left=0.12cm,at end] {$\gamma \t 1$\!}(180.0000,952.3622)
  (51.1162,967.3816) .. controls (99.7584,916.3920) and (137.2462,988.6006) .. 
  node[above=0.1cm,pos=0.24] {$\mathfrak l$}(53.2460,969.7244)
  (52.1811,968.6595) .. controls (67.8278,961.9169) and (83.3621,961.3872) .. 
  node[above=0.07cm,pos=0.6] {$ F$}(99.9407,963.0838)
  (105.0501,963.6685) .. controls (113.5230,964.7321) and (122.2913,966.2582) .. 
  node[above=0.1cm,pos=0.54] {$1 \t F$}(131.5034,967.7358)
  (139.5938,968.9908) .. controls (152.1021,970.8471) and (165.4574,972.3622) .. 
  node[above left=0.12cm,at end] {$1 \t F$\!}(180.0000,972.3622);
  \path[fill=black] (49.838245,973.13214) node[circle, draw, line width=0.65pt, minimum width=5mm, fill=white, inner sep=0.25mm] (text4277) {$\theta$   };
\end{tikzpicture}\rlap{ .}\qedhere
\end{equation*}
\end{proof}

\begin{Prop}
For all $\V$-functors $F, G \colon \B \to \C$, the functors
\begin{gather*}
\phantom{and\ }\C(1, \thg) \colon \V\text-\cat{Bicat}(\B,\C)(F,G) \to {}_\B \cat{Mod}_{\C}(\C(1,F),\, \C(1,G))\\
\text{and}\ \C(\thg, 1) \colon \V\text-\cat{Bicat}(\B,\C)(F,G) \to {}_\C \cat{Mod}_{\B}(\C(G,1),\, \C(F,1))
\end{gather*}
are equivalences of categories.
\end{Prop}
\begin{proof}
We prove only the case of $\C(1, \thg)$; the other is dual.
By Corollary~\ref{cor:yoneda-with-naturality} and Proposition~\ref{prop:right-hom-pointwise}, the bimorphism
$m \colon \C(F,G), \C(1,F) \to \C(1,G)$ exhibits $\C(F,G)$ as $\h{\C(1,F)}{\C(1,G)}$; by the dual of Corollary~\ref{cor:yoneda-with-naturality}, the bimorphism
\begin{equation*}
\phi = \B, \C(1,G) \xrightarrow{\hat G, 1} \C(G,G), \C(1,G) \xrightarrow{m} \C(1,G)
\end{equation*}
exhibits $\C(1,G)$ as $\h{\B}{\C(1,G)}_\ell$.
Thus, in the diagram of categories and functors
\begin{equation*}
\cd[@C+1.5em]{
\V\text-\cat{Bicat}(\B,\C)(F,G) \ar[d]_{\C(1, \thg)}
\ar[r]^-{\text{\eqref{eq:vnats-as-bimods}}} &
{}_\B \cat{Mod}_{\B}(\B, \C(F,G))
\ar[d]^{m \circ (\thg, 1)} \\
{}_\B \cat{Mod}_{\C}(\C(1,F), \C(1,G)) \ar[r]_-{\phi \circ (1, \thg)} &
\cat{Bimor}_{\B}(\B, \C(1,F); \C(1,G))
}
\end{equation*}
the top, bottom and right sides are equivalences; it will follow that the left side is too, so long as we can show that the square commutes to within natural isomorphism. To obtain the component of this at some $\V$-transformation $\gamma \colon F \Rightarrow G$, observe that the two sides of the square
send $\gamma$ to the respective bimorphisms
\begin{gather*}
\B, \C(1,F) \xrightarrow{\hat F, 1} \C(F,F), \C(1,F) \xrightarrow{\C(F,\gamma), 1} \C(F, G), \C(1,F) \xrightarrow m \C(1,G)\\
\text{and} \quad \B, \C(1,F) \xrightarrow{1, \C(1, \gamma)} \B, \C(1,G) \xrightarrow{\hat G, 1} \C(G,G), \C(1,G) \xrightarrow m \C(1,G)\rlap{ .}
\end{gather*}
We obtain the required invertible transformation between these bimorphisms by taking as components the $2$-cells\begin{equation*}
\begin{tikzpicture}[y=0.8pt, x=0.9pt,yscale=-1, inner sep=0pt, outer sep=0pt, every text node part/.style={font=\tiny} ]
  \path[draw=black,line join=miter,line cap=butt,line width=0.650pt] 
  (180.0000,952.3622) .. controls (145.8263,952.3622) and (134.1532,957.4738) .. 
  node[above left=0.12cm,at start] {$F\t 1$\!}(131.2779,960.8815);
  \path[draw=black,line join=miter,line cap=butt,line width=0.650pt] 
  (180.0000,972.3622) .. controls (149.3230,972.3622) and (135.1116,968.4220) .. 
  node[above left=0.12cm,at start] {$\C(1,\!\gamma)\t 1$\!}(131.3844,963.8428);
  \path[draw=black,line join=miter,line cap=butt,line width=0.650pt] 
  (180.0000,992.3622) .. controls (138.4604,992.3622) and (82.6973,987.9386) .. 
  node[above left=0.12cm,at start] {$m$\!}(77.6922,983.7854);
  \path[draw=black,line join=miter,line cap=butt,line width=0.650pt] 
  (128.3064,964.1623) .. controls (126.9220,967.7830) and (84.6324,973.3597) .. 
  node[above left=0.05cm,pos=0.32,rotate=11] {$\C(\gamma,\!1)\t 1$}(77.7104,980.8142);
  \path[draw=black,line join=miter,line cap=butt,line width=0.650pt] 
  (74.3035,980.6894) .. controls (64.5633,977.5586) and (50.8545,972.2616) .. 
  node[below left=0.03cm,pos=0.38,rotate=-18] {$1 \t \C(1,\!\gamma)$}(36.4904,967.0894)
  (27.7302,963.9883) .. controls (9.6569,957.7183) and (-8.3945,952.3622) .. 
  node[above right=0.12cm,at end] {\!\!$1 \t \C(1,\!\gamma)$}(-20.0000,952.3622)
  (128.3064,961.0945) .. controls (87.9459,937.9858) and (42.6070,972.3622) .. 
  node[above right=0.12cm,at end] {\!$G \t 1$}(-20.0000,972.3622);
  \path[draw=black,line join=miter,line cap=butt,line width=0.650pt] 
  (74.4100,983.9984) .. controls (70.5763,988.6840) and (3.8580,992.3622) .. 
  node[above right=0.12cm,at end] {\!$m$}(-20.0000,992.3622);
  \path[fill=black] (75.006805,982.50342) node[circle, draw, line width=0.65pt, minimum width=5mm, fill=white, inner sep=0.25mm] (text4505) {$\alpha^{-1}$     };
  \path[fill=black] (129.04779,963.3349) node[circle, draw, line width=0.65pt, minimum width=5mm, fill=white, inner sep=0.25mm] (text4509) {$\bar \gamma$   };
\end{tikzpicture}\rlap{ .}\qedhere
\end{equation*}
\end{proof}

In discussing Kan extensions, we will need a mild generalisation of this result; the proof is identical in form to the one just given and hence omitted.
\begin{Prop}\label{prop:vnat-for-kan}
For all $\V$-functors $F \colon \A \to \C$, $G \colon \A \to \B$ and $H \colon \B \to \C$, the functors
\begin{align*}
\V\text-\cat{Bicat}(\B,\C)(HG,F) &\to {}_\A \cat{Mod}_{\B}(\B(1,G),\, \C(H,F))\\
\gamma & \mapsto \B(1,G) \xrightarrow{\hat H} \C(H, HG) \xrightarrow{\C(1, \gamma)} \C(H,F) \\[0.3\baselineskip]
\text{and} \quad \quad \V\text-\cat{Bicat}(\B,\C)(F,HG) &\to {}_\B \cat{Mod}_{\A}(\B(G,1),\, \C(F,H))\\
\gamma & \mapsto \B(G,1) \xrightarrow{\hat H} \C(HG, H) \xrightarrow{\C(\gamma, 1)} \C(F,H)
\end{align*}
are equivalences of categories. \end{Prop}

\section{$\V$-categories of right modules}
\label{sec:vcatsofmodules}

In this section, we describe how the right modules over some $\V$-bicategory $\B$ may themselves be formed into a $\V$-bicategory. If we were identifying right modules with $\V$-functors $\B^\op \to \V$, then the construction we give would be an instance of the more general construction of functor $\V$-bicategories; but since we are treating right modules as basic, rather than derived, structures, the construction becomes significantly simpler, and may be carried out, again, with fewer assumptions on our base bicategory $\V$; in particular, no symmetry is needed.

\subsection{Moderate right modules}
\label{subsec:vcats-of-modules}
Let $\B$ be a $\V$-bicategory. A right $\B$-module $V$ is said to be \emph{moderate} if, for every right $\B$-module $W$, the right hom $\h V W$ exists, with universal morphism
$\xi_{VW} \colon \h V W \otimes V \to W$, say. By the Yoneda lemma~\ref{prop:yoneda}, every representable right module $\B(\thg, b)$ is moderate; moreover, by Corollary~\ref{cor:existence-of-right-homs}, we have:
\begin{Prop}\label{prop:every-module-over-small-moderate}
If $\V$ is complete and right closed, and $\B$ is a small $\V$-bicategory, then every right $\B$-module is moderate.
\end{Prop}

\subsection{The $\V$-category of moderate modules}
For any $\V$-bicategory $\B$, we now define a $\V$-bicategory $\M \B$ whose objects are the moderate right $\B$-modules, and whose hom-objects are the $\h V W$'s. For each $W \in \M\B$, the unit morphism $j_W \colon I \to \h W W$ is obtained by applying the universality of $\xi_{WW}$ to the module morphism $\mathfrak l \colon I \otimes W \to W$. This universality also yields an invertible module transformation
\begin{equation}\label{eq:jbar}
\cd[@!C@C-4.5em]{
 & \h W W \otimes W \ar[dr]^-{\xi} \dtwocell[0.6]{d}{\bar \jmath_W} \\
 I \otimes W \ar[rr]_{\mathfrak l} \ar[ur]^-{j \otimes 1} & &
 W\rlap{ .}
}
\end{equation}
Given $U,V,W \in \M\B$, the composition morphism $m_{UVW} \colon \h V W \otimes \h U V \to \h U W$ is obtained by applying universality of $\xi_{UW}$ to the module morphism \begin{equation*}
(\h V W \otimes \h U V) \otimes U \xrightarrow{\mathfrak a}
\h V W \otimes (\h U V \otimes U) \xrightarrow{1 \otimes \xi_{UV}} \
\h V W \otimes V \xrightarrow{\xi_{VW}}
W\rlap{ .}
\end{equation*}
This universality also yields an invertible module transformation
\begin{equation}\label{eq:mbar}
\cd[@C-1.5em]{
 (\h V W \otimes \h U V) \otimes U \ar[rr]^{\mathfrak a} \ar[d]_{m \otimes 1}  & {} \rtwocell[0.55]{d}{\bar m_{UVW}} &
 \h V W \otimes (\h U V \otimes U) \ar[d]^{1 \otimes \xi} \\
 \h U W \otimes U \ar[r]_-{\xi} & W & \h V W \otimes V \rlap{ .}
 \ar[l]^-{\xi}
}
\end{equation}
To obtain the left unit constraint $\lu_{VW} \colon m \circ (j \otimes 1) \Rightarrow \mathfrak l \colon I \otimes \h V W \to \h V W$ of $\M\B$,
it suffices by universality of $\xi_{VW}$ to give an invertible module transformation between the composites of $\xi_{VW}$ with $(m \circ (j \otimes 1)) \otimes 1$ and $\mathfrak l \otimes 1$. We obtain such with $2$-cell components as on the left in:
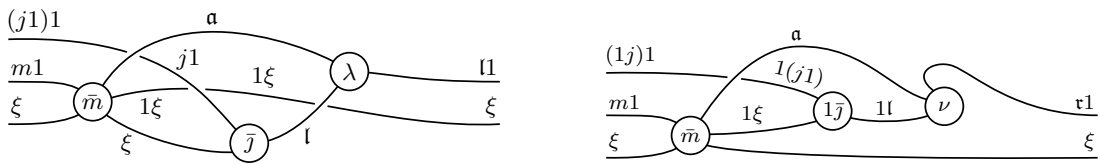
\begin{equation*}
\begin{tikzpicture}[y=0.80pt, x=0.8pt,yscale=-1, inner sep=0pt, outer sep=0pt, every text node part/.style={font=\scriptsize} ]
  \path[draw=black,line join=miter,line cap=butt,line width=0.650pt]
  (10.0000,1032.3622) .. controls (28.6481,1032.6606) and (42.6853,1030.8704) ..
  node[above right=0.12cm,at start] {\!$\xi$}(48.6527,1023.5604);
  \path[draw=black,line join=miter,line cap=butt,line width=0.650pt]
  (10.0000,1012.3622) .. controls (29.8416,1012.6606) and (39.8082,1011.2753) ..
  node[above right=0.12cm,at start] {\!$m \t 1$}(48.6527,1020.1198);
  \path[draw=black,line join=miter,line cap=butt,line width=0.650pt]
  (50.7125,1023.4890) .. controls (67.7767,1041.7087) and (118.0760,1050.4154) ..
  node[left=0.52cm,pos=0.45] {$\xi$}(122.8223,1043.0423);
  \path[draw=black,line join=miter,line cap=butt,line width=0.650pt]
  (170.9870,1007.7713) .. controls (188.1724,1007.7713) and (194.7164,1012.3622) ..
  node[above left=0.12cm,at end] {$\mathfrak l \t 1$\!}(240.0000,1012.3622);
  \path[draw=black,line join=miter,line cap=butt,line width=0.650pt]
  (52.1452,1021.9337) .. controls (65.2033,1017.1378) and (79.3246,1015.4846) ..
  node[below=0.1cm,pos=0.6] {$1 \t \xi$}(94.1198,1015.6823)
  (102.1856,1015.9621) .. controls (118.7052,1016.8669) and (135.9547,1019.7272) ..
  node[above=0.1cm,pos=0.54] {$1 \t \xi$}(153.4166,1022.8260)
  (159.3539,1023.8850) .. controls (186.5094,1028.7421) and (214.0353,1033.7052) ..
  node[above left=0.12cm,at end] {$\xi$\!}(240.0000,1032.3622)
  (10.0000,992.3622) .. controls (33.6375,991.6289) and (51.5811,994.1408) ..
  node[above right=0.12cm,at start] {\!\!$(j \t 1) \t 1$}(65.5095,998.3851)
  (71.4182,1000.3722) .. controls (103.1494,1012.0977) and (111.6557,1033.1261) ..
  node[above=0.17cm,pos=0.3] {$j \t 1$}(122.6313,1040.3004)
  (123.9221,1041.3205) .. controls (144.1543,1041.2970) and (161.5971,1017.3972) ..
  node[below right=0.07cm,pos=0.4] {$\mathfrak l$}(168.9655,1009.0195)
  (51.0090,1019.4542) .. controls (80.1930,964.5804) and (158.8751,998.5958) ..
  node[above=0.1cm,pos=0.42] {$\mathfrak a$}(168.1893,1006.5667);
  \path[fill=black] (122.74157,1041.3734) node[circle, draw, line width=0.65pt, minimum width=5mm, fill=white, inner sep=0.25mm] (text3353) {$\bar{\jmath}$};
  \path[fill=black] (49.5,1021.8622) node[circle, draw, line width=0.65pt, minimum width=5mm, fill=white, inner sep=0.25mm] (text3357) {$\bar m$};
  \path[fill=black] (169.59831,1007.5729) node[circle, draw, line width=0.65pt, minimum width=5mm, fill=white, inner sep=0.25mm] (text3361) {$\lambda$};
\end{tikzpicture} \qquad \ \quad
\begin{tikzpicture}[y=0.80pt, x=0.8pt,yscale=-1, inner sep=0pt, outer sep=0pt, every text node part/.style={font=\tiny} ]
  \path[draw=black,line join=miter,line cap=butt,line width=0.650pt]
  (10.0000,1032.3622) .. controls (28.6481,1032.6606) and (42.6853,1030.8704) ..
  node[above right=0.12cm,at start] {\!$\xi$}(48.6527,1023.5604);
  \path[draw=black,line join=miter,line cap=butt,line width=0.650pt]
  (10.0000,1012.3622) .. controls (29.8416,1012.6606) and (39.8082,1011.2753) ..
  node[above right=0.12cm,at start] {\!$m \t 1$}(48.6527,1020.1198);
  \path[draw=black,line join=miter,line cap=butt,line width=0.650pt]
  (52.0608,1021.3318) .. controls (75.0576,1021.7538) and (108.3681,1017.5166) ..
  node[above=0.1cm,pos=0.35] {$1 \t \xi$}(115.5414,1012.0311);
  \path[draw=black,line join=miter,line cap=butt,line width=0.650pt]
  (118.2592,1010.8486) .. controls (123.6600,1016.4881) and (159.4398,1017.3972) ..
  node[above=0.085cm,pos=0.47] {$1 \t \mathfrak l$}(166.8082,1009.0195);
  \path[draw=black,line join=miter,line cap=butt,line width=0.650pt]
  (10.0000,992.3622) .. controls (31.6003,991.9643) and (50.8012,992.7258) ..
  node[above right=0.12cm,at start] {\!\!$(1 \t j) \t 1$}(66.7503,994.4162)
  (73.1361,995.1647) .. controls (95.5024,998.0515) and (110.4183,1002.9226) ..
  node[above=0.07cm,rotate=-14,pos=0.46] {$1 \t (j \t 1)$}(115.0807,1009.0196)
  (51.0090,1019.4542) .. controls (64.7186,996.7679) and (81.6324,980.5234) ..
  (98.9213,979.9352) .. controls (133.9430,978.7437) and (145.9460,1005.9545) ..
  node[above=0.085cm,at start] {$\mathfrak a$}(166.5713,1005.7577);
  \path[draw=black,line join=miter,line cap=butt,line width=0.650pt]
  (167.2117,1003.4567) .. controls (149.6598,993.4283) and (163.2963,985.6240) ..
  (172.3462,989.1172) .. controls (182.3874,992.9932) and (201.8066,1012.3622) ..
  node[above left=0.12cm,at end] {$\mathfrak r \t 1$\!}(240.0000,1012.3622);
  \path[draw=black,line join=miter,line cap=butt,line width=0.650pt]
  (51.6059,1024.0910) .. controls (58.3519,1034.5666) and (240.0000,1032.3622) ..
  node[above left=0.12cm,at end] {$\xi$\!}(240.0000,1032.3622);
  \path[fill=black] (116,1010.3622) node[circle, draw, line width=0.65pt, minimum width=5mm, fill=white, inner sep=0.25mm] (text3353) {$1 \bar{\jmath}$     };
  \path[fill=black] (49.5,1021.8622) node[circle, draw, line width=0.65pt, minimum width=5mm, fill=white, inner sep=0.25mm] (text3357) {$\bar m$     };
  \path[fill=black] (168.25,1008.1122) node[circle, draw, line width=0.65pt, minimum width=5mm, fill=white, inner sep=0.25mm] (text3361) {$\nu$   };
\end{tikzpicture}
\end{equation*}
To give the right unit constraint $\ru_{VW} \colon m \circ (1 \otimes j) \Rightarrow \mathfrak r$, it likewise suffices to give an invertible module transformation between the composites of $\xi_{VW}$ with $(m \circ (1 \otimes j)) \otimes 1$ and $\mathfrak r \otimes 1$. We obtain such by taking its $2$-cell component to be as on the right above.
Finally, to give the associativity constraint $\ass_{UVWX}$, it suffices to give an invertible module transformation between the composites of $\xi_{UX}$ with $(m \circ (1 \otimes m) \circ \mathfrak a) \otimes 1$ and $(m \circ (m \otimes 1)) \otimes 1$. The $2$-cell components of this transformation are given by
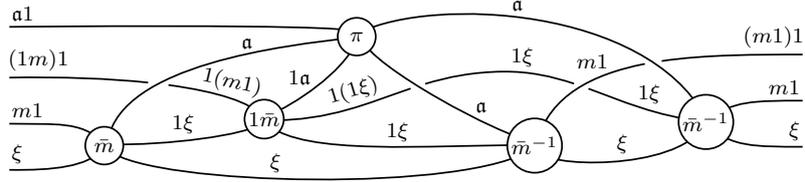
\begin{equation*}
\begin{tikzpicture}[y=0.87pt, x=0.9pt,yscale=-1, inner sep=0pt, outer sep=0pt, every text node part/.style={font=\tiny} ]
  \path[draw=black,line join=miter,line cap=butt,line width=0.650pt]
  (10.0000,1032.3622) .. controls (28.6481,1032.6606) and (42.6853,1030.8704) ..
  node[above right=0.12cm,at start] {\!$\xi$}(48.6527,1023.5604);
  \path[draw=black,line join=miter,line cap=butt,line width=0.650pt]
  (10.0000,1012.3622) .. controls (29.8416,1012.6606) and (39.8082,1011.2753) ..
  node[above right=0.12cm,at start] {\!$m \t 1$}(48.6527,1020.1198);
  \path[draw=black,line join=miter,line cap=butt,line width=0.650pt]
  (10.0000,970.3622) .. controls (10.0000,970.3622) and (146.3527,968.8266) ..
  node[above right=0.12cm,at start] {\!$\mathfrak a \t 1$}(154.1309,972.3622);
  \path[draw=black,line join=miter,line cap=butt,line width=0.650pt]
  (52.0608,1021.3318) .. controls (75.0576,1021.7538) and (108.3681,1017.5166) ..
  node[above=0.1cm,pos=0.4] {$1 \t \xi$}(115.5414,1012.0311);
  \path[draw=black,line join=miter,line cap=butt,line width=0.650pt]
  (117.2592,1009.0986) .. controls (129.2850,1003.6131) and (148.1898,990.5222) ..
  node[above left=0.07cm,pos=0.47] {$1 \t \mathfrak a$}(154.3082,977.0195);
  \path[draw=black,line join=miter,line cap=butt,line width=0.650pt]
  (51.0090,1019.4542) .. controls (56.8485,986.8815) and (123.4933,978.4403) ..
  node[above left=0.07cm,pos=0.65] {$\mathfrak a$}(153.3213,975.0077)
  (10.0000,992.3622) .. controls (32.0443,991.9561) and (51.5895,992.7576) ..
  node[above right=0.12cm,at start] {\!\!$(1 \t m) \t 1$}(67.7297,994.5216)
  (76.1783,995.5761) .. controls (96.8544,998.5027) and (110.6330,1003.2034) ..
  node[above=0.07cm,rotate=-13] {$1 \t (m \t 1)$}(115.0807,1009.0196);
  \path[draw=black,line join=miter,line cap=butt,line width=0.650pt]
  (117.4339,1012.6940) .. controls (132.8354,1026.4077) and (183.4823,1021.6122) ..
  node[above=0.08cm,pos=0.58] {$1 \t \xi$}(226.7500,1021.6122);
  \path[draw=black,line join=miter,line cap=butt,line width=0.650pt]
  (51.6059,1024.0910) .. controls (70.0000,1039.3622) and (201.5000,1042.3622) ..
  node[above=0.08cm,pos=0.43] {$\xi$}(227.7500,1023.1122);
  \path[draw=black,line join=miter,line cap=butt,line width=0.650pt]
  (229.5000,1019.8622) .. controls (234.9965,1000.2319) and (252.4130,990.8776) ..
  node[above left=0.08cm,pos=0.785] {$m \t 1$}(274.3623,986.4200)
  (281.3726,985.1785) .. controls (299.7606,982.3622) and (320.5923,982.3622) ..
  node[above left=0.12cm,at end] {$(m \t 1) \t 1$\!\!}(340.0000,982.3622)
  (156.4196,975.6489) .. controls (161.4435,988.0195) and (213.2044,1017.1122) ..
  node[above right=0.08cm,pos=0.65] {$\mathfrak a$}(227.2500,1019.8622)
  (118.4527,1010.8888) .. controls (143.5901,1010.9989) and (159.3833,1003.2413) ..
  node[above=0.08cm,pos=0.57,rotate=17] {$1 \t (1 \t \xi)$}(176.9367,997.0963)
  (183.0085,995.0768) .. controls (192.4298,992.1325) and (202.6702,989.9173) ..
  (215.2966,989.7688) .. controls (226.4532,989.6377) and (235.8159,991.8750) ..
  node[above=0.08cm,pos=0.25] {$1 \t \xi$}(244.8824,994.9818)
  (250.8838,997.1595) .. controls (264.8803,1002.4709) and (278.9634,1009.0791) ..
  node[above right=0.05cm,pos=0.45] {$1 \t \xi$}(298.7500,1011.3622)
  (155.8980,972.6975) .. controls (187.1970,957.1265) and (274.5018,962.1122) ..
  node[above=0.08cm,pos=0.45] {$\mathfrak a$}(299.2500,1009.1122);
  \path[draw=black,line join=miter,line cap=butt,line width=0.650pt]
  (230.5000,1023.8622) .. controls (230.5000,1031.3622) and (287.0000,1033.3622) ..
  node[above=0.1cm,pos=0.55] {$\xi$}(299.0000,1013.3622);
  \path[draw=black,line join=miter,line cap=butt,line width=0.650pt]
  (340.0000,1002.3622) .. controls (320.0000,1002.3622) and (310.0000,1002.3622) ..
  node[above left=0.12cm,at start] {$m \t 1$\!}(301.2500,1010.8622);
  \path[draw=black,line join=miter,line cap=butt,line width=0.650pt]
  (340.0000,1022.3622) .. controls (320.0000,1022.3622) and (310.0000,1022.3622) ..
  node[above left=0.12cm,at start] {$\xi$\!}(301.2500,1013.3622);
  \path[fill=black] (154.5,974.61218) node[circle, draw, line width=0.65pt, minimum width=5mm, fill=white, inner sep=0.25mm] (text3361) {$\pi$     };
  \path[fill=black] (49.5,1021.8622) node[circle, draw, line width=0.65pt, minimum width=5mm, fill=white, inner sep=0.25mm] (text3357) {$\bar m$     };
  \path[fill=black] (228.52539,1021.8378) node[circle, draw, line width=0.65pt, minimum width=5mm, fill=white, inner sep=0.25mm] (text3357-1) {${\bar m}^{-1}$     };
  \path[fill=black] (299.77539,1011.5878) node[circle, draw, line width=0.65pt, minimum width=5mm, fill=white, inner sep=0.25mm] (text3357-5) {${\bar m}^{-1}$     };
  \path[fill=black] (116,1010.3622) node[circle, draw, line width=0.65pt, minimum width=5mm, fill=white, inner sep=0.25mm] (text3353) {$1 \t \bar m$   };
\end{tikzpicture}\rlap{ .}
\end{equation*}

Observe that for all $V,W \in \M\B$, we have by~\eqref{eq:compinduced} an equivalence of categories
$\V(I, \h V W) \to {}_\bullet\cat{Mod}_\B(I \otimes V, W)$; and on composing this with the equivalence of categories ${}_\bullet\cat{Mod}_\B(\mathfrak l^\centerdot, W) \colon {}_\bullet\cat{Mod}_\B(I \otimes V, W) \to {}_\bullet\cat{Mod}_\B(V, W)$, we conclude that the
underlying bicategory of $\M\B$ is, up to an identity-on-objects biequivalence, the full sub-bicategory of ${}_\bullet\cat{Mod}_\B$ on the moderate right modules.

\subsection{The $\M\B$-module induced by a $\B$-module}\label{subsec:modules-over-module-categories}
Given a right $\B$-module $W$, the right hom $\h V W$ exists by definition for every moderate $\B$-module $V$. This allows us to define from $W$ a right $\M\B$-module whose component at $V \in \M \B$ is $\h V W$, and whose action morphisms and coherence $2$-cells defined just as in the preceding section. By abuse of notation, we shall denote this right module by $\M\B(\thg, W)$, although in general $W$ need not be an object of $\M \B$.

\subsection{Yoneda embedding}
As observed above, every representable right $\B$-module is moderate. We may thus define the \emph{Yoneda embedding} $Y \colon \B \to \M \B$ to be the $\V$-functor with action on objects $Yb = \B(\thg, b)$, and action on homs $Y_{bb'} \colon \B(b,b') \to \M\B(Yb,Yb')$ given by applying the universal property of $\M \B(Yb', Yb) = \h{Yb}{Yb'}$ to the right module morphism
\begin{equation}\label{eq:yon-embedding-module-maps}
\B(b,b') \otimes \B(\thg, b) \xrightarrow m \B(\thg, b')\rlap{ .}
\end{equation}
This universality also provides invertible module modifications $\Gamma \colon m \Rightarrow \xi \circ ( Y_{bb'} \otimes 1)$. To obtain the unit isomorphism $\fu_{b} \colon j_{Yb} \Rightarrow Y_{bb} \circ j_b \colon I \to \h{Yb}{Yb}$, it suffices by universality of $\h{Yb}{Yb}$ to construct an invertible module transformation $\xi_{Yb, Yb} \circ (j_{Yb} \otimes 1)  \Rightarrow \xi_{Yb, Yb} \circ ((Y_{bb} \circ j_b) \otimes 1)$. We obtain such by taking its $2$-cell components to be 
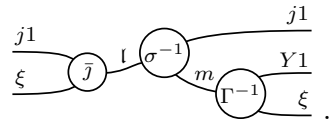
\begin{equation*}
\begin{tikzpicture}[y=0.80pt, x=0.7pt,yscale=-1, inner sep=0pt, outer sep=0pt, every text node part/.style={font=\tiny} ]
  \path[draw=black,line join=miter,line cap=butt,line width=0.650pt]
  (170.0000,1012.3622) .. controls (150.0000,1012.3622) and (140.0000,1012.3622) ..
  node[above left=0.12cm,at start] {$\xi$\!}(131.5327,1005.0719);
  \path[draw=black,line join=miter,line cap=butt,line width=0.650pt]
  (170.0000,992.3622) .. controls (150.0000,992.3622) and (140.0000,992.3622) ..
  node[above left=0.12cm,at start] {$Y \t 1$\!}(131.5854,999.4988);
  \path[draw=black,line join=miter,line cap=butt,line width=0.650pt]
  (170.0000,972.3622) .. controls (128.8033,972.3622) and (95.5339,974.3721) ..
  node[above left=0.12cm,at start] {$j \t 1$\!}(91.5372,985.1324);
  \path[draw=black,line join=miter,line cap=butt,line width=0.650pt]
  (128.2181,1001.9424) .. controls (113.4611,1001.9424) and (97.1669,995.7970) ..
  node[above right=0.07cm,pos=0.53] {$m$}(91.9404,988.1110);
  \path[draw=black,line join=miter,line cap=butt,line width=0.650pt]
  (10.0000,1002.3622) .. controls (37.6941,1002.3622) and (44.5136,1000.1886) ..
  node[above right=0.12cm,at start] {\!$\xi$}(48.5917,994.1027);
  \path[draw=black,line join=miter,line cap=butt,line width=0.650pt]
  (48.2843,990.4135) .. controls (41.8281,983.6499) and (36.7487,982.3622) ..
  node[above right=0.12cm,at end] {\!$j \t 1$}(10.0000,982.3622);
  \path[draw=black,line join=miter,line cap=butt,line width=0.650pt]
  (89.7718,986.5738) .. controls (70.9954,986.5738) and (75.8267,992.6696) ..
  node[above left=0.07cm] {$\mathfrak l$}(51.8446,992.6696);
  \path[fill=black] (132.07925,1002.1044) node[circle, draw, line width=0.65pt, minimum width=5mm, fill=white, inner sep=0.25mm] (text3463) {$\Gamma^{-1}$     };
  \path[fill=black] (50.082329,991.94238) node[circle, draw, line width=0.65pt, minimum width=5mm, fill=white, inner sep=0.25mm] (text3850) {$\bar \jmath$     };
  \path[fill=black] (91.309006,982.88129) node[circle, draw, line width=0.65pt, minimum width=5mm, fill=white, inner sep=0.25mm] (text3854) {$\lu^{-1}$   };
\end{tikzpicture}\rlap{ .}
\end{equation*}
Similarly, to define the binary coherence $2$-cell $\fm_{bb'b''}$, it suffices to give an invertible module transformation between the composites of $\xi_{Yb, Yb''}$ with $(m \circ (Y \otimes Y)) \otimes 1$ and $(Y \circ m) \otimes 1$. We obtain such by taking its $2$-cell components
to be 
\begin{equation*}
\begin{tikzpicture}[y=0.8pt, x=0.8pt,yscale=-1, inner sep=0pt, outer sep=0pt, every text node part/.style={font=\tiny} ]
  \path[draw=black,line join=miter,line cap=butt,line width=0.650pt]
  (210.0000,1022.3622) .. controls (190.1148,1022.3622) and (180.0710,1020.7287) ..
  node[above left=0.12cm,at start] {$\xi$\!}(175.3804,1014.3135);
  \path[draw=black,line join=miter,line cap=butt,line width=0.650pt]
  (210.0000,1002.3622) .. controls (190.2564,1002.3622) and (181.4818,1002.6917) ..
  node[above left=0.12cm,at start] {$Y \t 1$\!}(175.7206,1011.0399);
  \path[draw=black,line join=miter,line cap=butt,line width=0.650pt]
  (210.0000,982.3622) .. controls (188.6368,982.3622) and (147.4904,981.0509) ..
  node[above left=0.12cm,at start] {$m \t 1$\!}(137.7449,990.9489);
  \path[draw=black,line join=miter,line cap=butt,line width=0.650pt]
  (173.7193,1012.3338) .. controls (158.9623,1012.3338) and (143.0872,1001.6134) ..
  node[above right=0.07cm,pos=0.55] {$m$}(137.8607,993.9275);
  \path[draw=black,line join=miter,line cap=butt,line width=0.650pt]
  (-63.8940,1011.8479) .. controls (-48.6766,1017.2827) and (5.2359,1013.7260) ..
  node[above=0.07cm,pos=0.6] {$1 \t \xi$}(25.8880,1008.5086);
  \path[draw=black,line join=miter,line cap=butt,line width=0.650pt]
  (-63.7585,1013.7885) .. controls (-58.4669,1027.1301) and (78.7117,1027.1104) ..
  node[above=0.07cm,pos=0.7] {$\xi$}(87.2609,1011.0578);
  \path[draw=black,line join=miter,line cap=butt,line width=0.650pt]
  (-65.9288,1013.6412) .. controls (-72.1775,1020.6122) and (-79.4913,1022.3622) ..
  node[above right=0.12cm,at end] {\!$\xi$}(-100.0000,1022.3622);
  \path[draw=black,line join=miter,line cap=butt,line width=0.650pt]
  (29.5924,1006.1428) .. controls (38.9961,1001.4682) and (51.7367,996.1376) ..
  node[above left=0.07cm,pos=0.65] {$1 \t m$}(67.0724,992.5178)
  (74.6202,990.9232) .. controls (92.1366,987.6550) and (112.5680,986.8456) ..
  node[above=0.07cm] {$1 \t m$}(134.9672,991.5170)
  (86.9957,1007.4072) .. controls (84.3219,990.2175) and (17.4013,964.6290) ..
  node[above=0.1cm,pos=0.5] {$Y \t 1$}(-12.8508,989.2007);
  \path[draw=black,line join=miter,line cap=butt,line width=0.650pt]
  (135.3974,994.2526) .. controls (122.7080,1007.7308) and (100.8214,1001.0956) ..
  node[below right=0.07cm,pos=0.4] {$\xi$}(89.8828,1009.3566);
  \path[draw=black,line join=miter,line cap=butt,line width=0.650pt]
  (-66.0797,1010.7451) .. controls (-71.6848,1004.2777) and (-79.3419,1002.3622) ..
  node[above right=0.12cm,at end] {\!$m \t 1$}(-100.0000,1002.3622);
  \path[draw=black,line join=miter,line cap=butt,line width=0.650pt]
  (-100.0000,982.3622) .. controls (-65.5323,982.3622) and (-88.3974,982.3702) ..
  node[above right=0.12cm,at start] {\!\!$(Y \t Y) \t 1$}
  (-59.0146,982.6198) .. controls (-58.1893,982.6268) and (-57.3727,982.6386) .. (-56.5645,982.6551)
  (-49.5370,982.9251) .. controls (-15.0388,984.8937) and (3.4438,995.7568) ..
  node[below=0.07cm,pos=0.15,rotate=-9] {$Y \t (Y \t 1)$}
  node[above=0.02cm,pos=0.73,rotate=-17] {$1 \t (Y \t 1)$} (26.1000,1004.2487)
  (-64.1352,1009.3284) .. controls (-62.9885,966.5154) and (-11.4257,965.5628) ..
  (35.8893,964.8030) .. controls (71.6867,964.2281) and (128.1646,976.5675) ..
  node[above=0.07cm,at start] {$\mathfrak a$}(135.6458,989.7161);
  \path[fill=black] (-64.886513,1012.3622) node[circle, draw, line width=0.65pt, minimum width=5mm, fill=white, inner sep=0.25mm] (text3838) {$\ass$     };
  \path[fill=black] (27.869759,1006.2278) node[circle, draw, line width=0.65pt, minimum width=5mm, fill=white, inner sep=0.25mm] (text3842) {$1 \Gamma$     };
  \path[fill=black] (87.669739,1008.9585) node[circle, draw, line width=0.65pt, minimum width=5mm, fill=white, inner sep=0.25mm] (text3846) {$\Gamma$     };
  \path[fill=black] (175.18152,1012.7709) node[circle, draw, line width=0.65pt, minimum width=5mm, fill=white, inner sep=0.25mm] (text3850) {$\Gamma^{-1}$     };
  \path[fill=black] (136.59808,992.43054) node[circle, draw, line width=0.65pt, minimum width=5mm, fill=white, inner sep=0.25mm] (text3854) {$\ass$     };
  \path[shift={(-89.194619,58.292418)},draw=black,fill=black] (77.8600,931.1000)arc(0.000:180.000:1.250)arc(-180.000:0.000:1.250) -- cycle;
\end{tikzpicture}\rlap{ .}
\end{equation*}

\begin{Prop}
The Yoneda embedding is fully faithful.
\end{Prop}
\begin{proof}
By the Yoneda lemma, each morphism~\eqref{eq:yon-embedding-module-maps} exhibits $\B(b,b')$ as $\h{Yb}{Yb'}$; whence each $Y_{bb'}$ is an equivalence in $\V$ as required.
\end{proof}
The following further reformulation of the Yoneda lemma will prove useful in what follows. Given any right $\B$-module $W$, we may as in Section~\ref{subsec:modules-over-module-categories} form the right $\B$-module $\M \B(Y,W)$. For each $b \in \B$, the right module morphism $m_{\thg b} \colon Wb \otimes \B(\thg, b) \to W$ induces by the universal property of $\h{Yb}{W}$ a morphism
\begin{equation*}
\upsilon_b \colon Wb \to \h{Yb}{W}\rlap{ ,}
\end{equation*}
together with an invertible module transformation $\Delta \colon m_{\thg b} \Rightarrow \xi_{Yb,W} \circ (\upsilon_b \otimes 1)$. We claim that the morphisms $\upsilon_b$ constitute the components of a module morphism $W \to \M\B(Y,W)$. To give the $2$-cells $m \circ (1 \otimes Y) \circ (\upsilon \otimes 1) \Rightarrow \upsilon \circ m$ verifying compatibility with the right $\B$-actions, it suffices to give an invertible module transformation between the composites of $(m \circ (1 \otimes Y) \circ (\upsilon \otimes 1)) \otimes 1$ and $(\upsilon \circ m) \otimes 1$ with $\xi$. We obtain such by taking its $2$-cell components to be

\begin{equation*}
\begin{tikzpicture}[y=0.8pt, x=0.8pt,yscale=-1, inner sep=0pt, outer sep=0pt, every text node part/.style={font=\tiny} ]
  \path[draw=black,line join=miter,line cap=butt,line width=0.650pt] 
  (220.0000,1012.3622) .. controls (200.1148,1012.3622) and (190.0710,1010.7287) .. 
  node[above left=0.12cm,at start] {$\xi$\!}(185.3804,1004.3135);
  \path[draw=black,line join=miter,line cap=butt,line width=0.650pt] 
  (220.0000,992.3622) .. controls (200.2564,992.3622) and (191.4818,992.6917) .. 
  node[above left=0.12cm,at start] {$\upsilon \t 1$\!}(185.7206,1001.0399);
  \path[draw=black,line join=miter,line cap=butt,line width=0.650pt] 
  (220.0000,972.3622) .. controls (198.6368,972.3622) and (150.0000,972.3622) .. 
  node[above left=0.12cm,at start] {$m \t 1$\!}(137.7449,990.9489);
  \path[draw=black,line join=miter,line cap=butt,line width=0.650pt] 
  (183.7193,1002.3338) .. controls (168.9623,1002.3338) and (143.0872,1001.6134) .. 
  node[above=0.07cm,pos=0.44] {$m$}(137.8607,993.9275);
  \path[draw=black,line join=miter,line cap=butt,line width=0.650pt] 
  (-63.8940,1011.8479) .. controls (-41.3596,1013.2177) and (-8.1788,1010.4739) .. 
  node[above=0.07cm,pos=0.45] {$1 \t \xi$}(2.7173,1000.7850);
  \path[draw=black,line join=miter,line cap=butt,line width=0.650pt] 
  (-63.7585,1013.7885) .. controls (-58.4669,1027.1301) and (78.7117,1027.1104) .. 
  node[above=0.07cm,pos=0.5] {$m$}(87.2609,1011.0578);
  \path[draw=black,line join=miter,line cap=butt,line width=0.650pt] 
  (-65.9288,1013.6412) .. controls (-72.1775,1020.6122) and (-79.4913,1022.3622) .. 
  node[above right=0.12cm,at end] {\!$\xi$}(-100.0000,1022.3622);
  \path[draw=black,line join=miter,line cap=butt,line width=0.650pt] 
  (135.3974,994.2526) .. controls (122.7080,1007.7308) and (100.8214,1001.0956) .. 
  node[above=0.07cm,pos=0.55] {$m$}(89.8828,1009.3566);
  \path[draw=black,line join=miter,line cap=butt,line width=0.650pt] 
  (-66.0797,1010.7451) .. controls (-71.6848,1004.2777) and (-79.3419,1002.3622) .. 
  node[above right=0.12cm,at end] {\!$m \t 1$}(-100.0000,1002.3622);
  \path[draw=black,line join=miter,line cap=butt,line width=0.650pt] 
  (4.7957,998.4192) .. controls (14.5930,993.5490) and (34.9309,991.5572) .. 
  node[below=0.07cm,pos=0.66] {$1 \t m$}(57.8164,990.8932)
  (68.8728,990.6712) .. controls (92.2664,990.3908) and (116.9286,991.1831) .. 
  node[above=0.07cm,pos=0.45] {$1 \t m$}(134.9672,991.5170)
  (86.9957,1007.4072) .. controls (77.5214,994.9986) and (44.7536,978.8523) .. 
  node[above=0.1cm,pos=0.6] {$\upsilon \t 1$}(-8.5578,969.6447)
  (-23.3698,967.3460) .. controls (-45.8177,964.2412) and (-71.4172,962.3622) .. 
  node[above right=0.12cm,at end] {\!\!$(\upsilon \t 1) \t 1$}(-100.0000,962.3622)
  (-100.0000,982.3622) .. controls (-84.6299,982.3622) and (-69.6827,982.8087) .. 
  node[above right=0.12cm,at start] {\!\!$(1 \t Y) \t 1$}(-56.0746,983.7823)
  (-49.2204,984.3280) .. controls (-25.1185,986.4513) and (-6.0612,990.3991) .. 
  node[above=0.07cm,pos=0.52,rotate=-8] {$1 \t (Y \t 1)$}(2.2764,996.6711)
  (-64.1352,1009.3284) .. controls (-62.9885,966.5155) and (1.0189,964.8030) .. 
  node[above=0.07cm,at end] {$\mathfrak a$}
  (35.8893,964.8030) .. controls (71.6913,964.8030) and (128.1646,976.5676) .. (135.6458,989.7161);
  \path[fill=black] (-64.886513,1012.3622) node[circle, draw, line width=0.65pt, minimum width=5mm, fill=white, inner sep=0.25mm] (text3838) {$\bar m$     };
  \path[fill=black] (3.0894313,998.29712) node[circle, draw, line width=0.65pt, minimum width=5mm, fill=white, inner sep=0.25mm] (text3842) {$1 \Gamma$     };
  \path[fill=black] (87.669739,1008.9585) node[circle, draw, line width=0.65pt, minimum width=5mm, fill=white, inner sep=0.25mm] (text3846) {$\Delta$     };
  \path[fill=black] (185.18152,1002.7709) node[circle, draw, line width=0.65pt, minimum width=5mm, fill=white, inner sep=0.25mm] (text3850) {$\Delta^{-1}$     };
  \path[fill=black] (136.59808,992.43054) node[circle, draw, line width=0.65pt, minimum width=5mm, fill=white, inner sep=0.25mm] (text3854) {$\ass^{-1}$   };
\end{tikzpicture}\rlap{ .}
\end{equation*}

Since, by the Yoneda lemma, each $m_{\thg b}$ exhibits $Wb$ as $\h {Yb}{W}$, we deduce that the morphisms $\upsilon_b$ are equivalences, and thus conclude that:
\begin{Prop}\label{prop:reformulate-yoneda}
For every right $\B$-module $W$, there is an equivalence of right $\B$-modules $\upsilon \colon W \to \M\B(Y,W)$, defined as above.
\end{Prop}

\section{Colimits and left Kan extensions}\label{sec:colimits}
We have now developed enough theory to describe and discuss colimits and left Kan extensions for enriched bicategories. By colimit, we of course mean \emph{weighted} colimit; $\V$-bicategories of right modules will play the important role of indexing the weights for such colimits. We have chosen to discuss colimits rather than limits since it is these that we will require in our further development; note that defining limits and right Kan extensions would instead require us to make use of $\V$-bicategories of \emph{left} modules.

\subsection{Weighted colimits}\label{subsec:wclim}
Suppose that we are given $F \colon \B \to \C$ and a right $\B$-module $W$. A \emph{$W$-weighted cylinder} under $F$ is given by an object $v \in \C$ (the \emph{vertex} of the cylinder) together with a morphism of right modules $\phi \colon W \to \C(F, v)$.
Given such a cylinder, we obtain for each $c \in \C$ a right module morphism
\begin{equation}
\label{eq:inducedphi}
\phi^c = \C(v,c) \otimes W \xrightarrow{1 \otimes \phi} \C(v,c) \otimes \C(F,v) \xrightarrow{m} \C(F,c)\end{equation}
whose second constituent  is the right module morphism obtained from the bimodule structure of $\C(F,1)$ as in Example~\ref{ex:yoneda-bimorphisms}.
The cylinder $\phi$ is said to be \emph{colimiting} if each induced $\phi^c$ exhibits $\C(v,c)$ as $\h{W}{\C(F,c)}$.
We write $W \star F$ for the vertex of a colimiting cylinder, and call it the \emph{colimit of $F$ weighted by $W$}.
Note that by Proposition~\ref{prop:right-hom-pointwise}, $\phi$ is colimiting just when the bimodule morphism
\begin{equation*}
\C(v,\thg) \otimes W \xrightarrow{1, \phi} \C(v,\thg) \otimes \C(F,v) \xrightarrow{m} \C(F,1)
\end{equation*}
exhibits $\C(v, \thg)$ as $\h{W}{\C(F,1)}$.

\begin{Ex}\label{ex:yonedacolimit}
Let $F \colon \B \to \C$ be any $\V$-functor, and $b \in \B$. Partially evaluating the bimodule morphism $\hat F \colon \B \to \C(F,F)$ of Example~\ref{ex:modulemorphisminducedbyfunctor} yields a morphism of right $\B$-modules $\zeta = \hat F_{b\thg} \colon \B(\thg, b) \to \C(F, Fb)$; now for each $c \in \C$, the induced morphism $\zeta^c \colon \C(Fb, c) \otimes \B(\thg, b) \to \C(F,c)$ is precisely that which is asserted to be universal by the Yoneda lemma, so that $\zeta$ exhibits $Fb$ as $\B(\thg, b) \star F$.
\end{Ex}

\subsection{Uniqueness of colimits}\label{subsec:uniqueness-of-colimits}
Given $F \colon \B \to \C$ a $\V$-functor, we have the bimodule $\C(F, 1)$ which, as in Example~\ref{ex:homs-induced-by-modules}, induces a functor $\C_0 \to {}_\bullet\cat{Mod}_\B$ sending $c$ to $\C(F,c)$ and sending $f \colon c \to d$ to the module morphism
\begin{equation*}
\C(F, c) \xrightarrow{\mathfrak l^\centerdot} I \otimes \C(F, c) \xrightarrow{f \otimes 1} \C(c,d) \otimes \C(F, c) \xrightarrow{m} \C(F, d)\rlap{ .}
\end{equation*}
Now let $W$ be a right $\B$-module, and let the cylinder $\phi \colon W \to \C(F,v)$ exhibit $v$ as $W \star F$. For each $c \in \C$, the universality of $\phi^c$ yields an equivalence of categories
\begin{equation*}
\C_0(v, c) = \V(I, \C(v,c)) \xrightarrow{\phi^c \circ (\thg \otimes 1)} {}_\bullet\cat{Mod}_\B(I \otimes W, \C(F,c)) \xrightarrow{(\thg) \circ \mathfrak l^\centerdot}
{}_\bullet\cat{Mod}_\B(W, \C(F,c))\rlap{ .}
\end{equation*}
To within isomorphism, this functor sends $f \colon v \to c$ to the composite morphism $\C(F,f) \circ \phi \colon W \to \C(F,d)$; and we therefore conclude that:
\begin{Prop}\label{prop:uniquenessofcolimits}
  If the cylinder $\phi \colon W \to \C(F,v)$ exhibits $v$ as $W \star F$, then it also exhibits $v$ as birepresenting object for the functor ${}_\bullet\cat{Mod}_\B(W, \C(F, \thg)) \colon \C_0 \to \cat{CAT}$.
  Consequently, any two $W$-weighted colimits of $F$ are related in an essentially-unique way by an equivalence in $\C_0$ commuting with the universal cylinders.
\end{Prop}

\subsection{Preservation of colimits}
Given $F \colon \B \to \C$, a right $\B$-module $W$, and a colimiting cylinder $\phi \colon W \to \C(F,W \star F)$ in $\C$, a functor $G \colon \C \to \D$ is said to \emph{preserve} the colimit $W \star F$ if the composite cylinder
\begin{equation*}
W \xrightarrow{\phi} \C(F,W \star F) \xrightarrow{\hat G} \D(GF,G(W \star F))
\end{equation*}
is colimiting. 

\subsection{Functoriality of colimits}
The simplest way of discussing the functoriality of taking colimits is to generalise the basic notions to depend on an indexing $\V$-bicategory $\A$. Suppose that we are given $F \colon \B \to \C$ and an $\A$-$\B$-module $M$. An \emph{$M$-weighted cylinder} under $F$ is given by a $\V$-functor $V \colon \A \to \C$ (the \emph{vertex} of the cylinder) together with a bimodule morphism $\phi \colon M \to \C(F, V)$.
Given such a cylinder, we obtain for each $c \in \C$ a module bimorphism
\begin{equation}
\label{eq:inducedphi2}
\phi^c = \C(V,c) , M \xrightarrow{1,\phi} \C(V,c), \C(F,V) \xrightarrow{m} \C(F,c)
\end{equation}
whose second component is obtained as in Examples~\ref{ex:yoneda-bimorphisms} and \ref{ex:pullback-of-bimorphisms}.
The cylinder $\phi$ is said to be \emph{colimiting} if each induced $\phi^c$ exhibits $\C(V,c)$ as $\h{M}{\C(F,c)}$.
We write $M \star F \colon \A \to \C$ for the vertex of a colimiting cylinder, and call it the \emph{colimit of $F$ weighted by $M$}.
By Proposition~\ref{prop:right-hom-pointwise}, to ask that $\phi$ be colimiting is equally to ask that
\begin{equation*}
\C(V,1), M \xrightarrow{1, \phi} \C(V,1), \C(F,V) \xrightarrow{m} \C(F,1)
\end{equation*}
exhibit $\C(V, 1)$ as $\h{M}{\C(F,1)}$.

By Proposition~\ref{prop:right-hom-pointwise} again, if $\phi \colon M \to \C(F,V)$ is a colimiting cylinder, then so is each $\phi_{a \thg} \colon M(\thg, a) \to \C(F,Va)$, so that the existence of $M \star F$ implies that of each $M(\thg, a) \star F$; conversely, we have:

\begin{Prop}\label{prop:functcolimits}
Let $M$ be an $\A$-$\B$-bimodule and let $F \colon \B \to \C$. If for all $a \in \A$ the weighted colimit $M(\thg, a) \star F$ exists in $\C$, then the weighted colimit $M \star F \colon \A \to \C$ does too, and may be computed pointwise in the sense that $(M
\star F)(a) = M(\thg, a) \star F$ with the colimiting cylinder for $M \star F$ being in each component the colimiting cylinder for $M(\thg, a) \star F$.
\end{Prop}
\begin{proof}
Suppose that each $M(\thg, a) \star F$ exists, with universal cylinder $\eta_a \colon M(\thg, a) \to \C(F, M(\thg, a) \star F)$, say.
We define a $\V$-functor $L \colon \A \to \C$ that will be the desired $M$-weighted colimit of $F$. On objects, we take $La = M(\thg, a) \star F$, as anticipated. To define the action on homs $L_{aa'} \colon \A(a,a') \to \C(La, La')$, consider the composite module morphism
\begin{equation*}
\A(a,a') \otimes M(\thg, a) \xrightarrow m M(\thg, a') \xrightarrow{\eta_{a'}} \C(Fb, L a')\rlap{ .}
\end{equation*}
The morphism $\eta_{a}^{L a'} \colon \C(L a, L a') \otimes M(\thg, a) \to \C(F, L a')$ induced by $\eta_{a}$ as in~\eqref{eq:inducedphi} exhibits $\C(La, La')$ as $\h{M(\thg, a)}{\C(F, La')}$, and thus we 
induce the desired morphism $L_{aa'} \colon \A(a,a') \to \C(L a, L a')$, together with an invertible module transformation
 \begin{equation}\label{eq:the-theta}
\cd[@C+1em]{
 \A(a,a') \otimes M(\thg,a) \ar[rr]^-{m} \ar[d]_{L_{aa'} \otimes 1} \dtwocell{drr}{\Theta_{aa'}} & &
 M(\thg,a') \ar[d]^{\eta_{a'}} \\
 \C(L a, L a') \otimes M(\thg,a) \ar[r]_-{1 \otimes \eta_{a}}  &
  \C(L a, L a') \otimes \C(F, L a') \ar[r]_-{m} &
 \C(F,L a')\rlap{ .}
}
\end{equation}

We now construct the functoriality coherence cells for $L$. To obtain the unit coherence isomorphism $\fu_{a} \colon j_{L a} \Rightarrow L_{aa} \circ j_a \colon I \to \C(La, La)$, it suffices by universality to construct an invertible module transformation $(j_{L a} \otimes 1) \circ \eta_a^{L a} \Rightarrow ((L_{aa} \circ j_a) \otimes 1) \circ \eta_a^{La}$. We obtain such by taking its $2$-cell components to be
\begin{equation*}
\begin{tikzpicture}[y=0.80pt, x=0.8pt,yscale=-1, inner sep=0pt, outer sep=0pt, every text node part/.style={font=\tiny} ]
  \path[draw=black,line join=miter,line cap=butt,line width=0.650pt]
  (170.0000,1022.3622) .. controls (146.5217,1022.3622) and (138.8102,1012.3494) ..
  node[above left=0.12cm,at start]{$m$\!}(131.5327,1005.0719);
  \path[draw=black,line join=miter,line cap=butt,line width=0.650pt]
  (170.0000,1002.3622) --
  node[above left=0.12cm,at start]{$1 \t \eta_{a}$\!}(131.7391,1002.3622);
  \path[draw=black,line join=miter,line cap=butt,line width=0.650pt]
  (170.0000,982.3622) .. controls (146.9508,982.3622) and (140.6522,990.4320) ..
  node[above left=0.12cm,at start]{$L \t 1$\!}(131.5854,999.4988);
  \path[draw=black,line join=miter,line cap=butt,line width=0.650pt]
  (170.0000,962.3622) .. controls (128.8033,962.3622) and (95.5339,970.3721) ..
  node[above left=0.12cm,at start]{$j \t 1$\!}(91.5372,981.1324);
  \path[draw=black,line join=miter,line cap=butt,line width=0.650pt]
  (128.2181,1001.9424) .. controls (113.4611,1001.9424) and (97.1669,991.7970) ..
  node[above right=0.07cm]{$m$}(91.9404,984.1110);
  \path[draw=black,line join=miter,line cap=butt,line width=0.650pt]
  (-30.0000,1014.3622) .. controls (-2.6134,1014.3622) and (22.6690,1016.4828) ..
  node[above right=0.12cm,at start]{\!$m$}(29.8215,1004.0944);
  \path[draw=black,line join=miter,line cap=butt,line width=0.650pt]
  (-30.0000,992.3622) .. controls (-14.3347,992.3622) and (-5.3471,989.2214) ..
  node[above right=0.12cm,at start]{\!$1 \t \eta_a$}(1.7530,985.6088)
  (6.5497,982.9947) .. controls (13.7510,978.8946) and (20.0223,974.9941) ..
  (31.9557,974.9662) .. controls (44.4292,974.9370) and (53.0039,982.1314) ..
  node[above=0.07cm,pos=0.1]{$1 \t \eta_a$}(60.7347,990.6127)
  (64.1378,994.4625) .. controls (79.1386,1011.8194) and (92.6002,1030.7726) ..
  node[above=0.07cm,pos=0.63]{$\eta_a$}(128.5090,1004.7094)
  (28.8991,1000.4052) .. controls (10.1454,996.1011) and (2.8974,972.3622) ..
  node[above right=0.12cm,at end]{\!$j \t 1$}(-30.0000,972.3622)
  (89.7718,982.5738) .. controls (70.9954,982.5738) and (55.6482,1002.2498) ..
  node[below right=0.07cm,pos=0.75]{$\mathfrak l$}(31.6661,1002.2498);
  \path[fill=black] (130.07925,1002.1044) node[circle, draw, line width=0.65pt, minimum width=5mm, fill=white, inner sep=0.25mm] (text3463) {$\Theta$};
  \path[fill=black] (29.821461,1001.9424) node[circle, draw, line width=0.65pt, minimum width=5mm, fill=white, inner sep=0.25mm] (text3850) {$\lu$};
  \path[fill=black] (91.309006,982.88129) node[circle, draw, line width=0.65pt, minimum width=5mm, fill=white, inner sep=0.25mm] (text3854) {$\lu^{-1}$};
\end{tikzpicture}\rlap{ .}
\end{equation*}
Similarly, to define the composition coherence constraint $\fm_{aa'a''}$, it suffices to give an invertible module transformation between the composites of $\eta_a^{La''}$ with $(m \circ (L \otimes L)) \otimes 1$ and $(L \circ m) \otimes 1$. We obtain such by taking its $2$-cell components to be
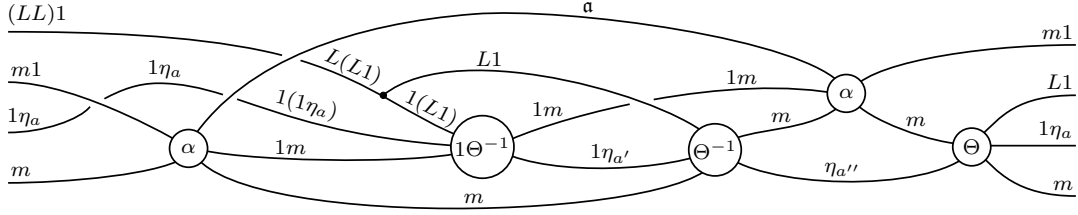
\begin{equation*}
\begin{tikzpicture}[y=0.95pt, x=0.95pt,yscale=-1, inner sep=0pt, outer sep=0pt, every text node part/.style={font=\tiny} ]
  \path[draw=black,line join=miter,line cap=butt,line width=0.650pt]
  (284.7826,1027.5796) .. controls (261.3043,1027.5796) and (253.5928,1017.5668) ..
  node[above left=0.12cm,at start]{$m$\!}(246.3153,1010.2893);
  \path[draw=black,line join=miter,line cap=butt,line width=0.650pt]
  (284.7826,1007.5796) --
  node[above left=0.12cm,at start]{$1 \t \eta_a$\!}(246.5217,1007.5796);
  \path[draw=black,line join=miter,line cap=butt,line width=0.650pt]
  (284.7826,987.5796) .. controls (261.7334,987.5796) and (255.4348,995.6494) ..
  node[above left=0.12cm,at start]{$L \t 1$\!}(246.3680,1004.7162);
  \path[draw=black,line join=miter,line cap=butt,line width=0.650pt]
  (284.7826,967.5796) .. controls (243.5859,967.5796) and (200.3165,975.5895) ..
  node[above left=0.12cm,at start]{$m \t 1$\!}(196.3198,986.3498);
  \path[draw=black,line join=miter,line cap=butt,line width=0.650pt]
  (243.0007,1007.1598) .. controls (228.2437,1007.1598) and (201.9495,997.0143) ..
  node[above right=0.07cm]{$m$}(196.7230,989.3284);
  \path[draw=black,line join=miter,line cap=butt,line width=0.650pt]
  (58.1710,1009.8937) .. controls (78.4866,1020.0391) and (125.4615,1017.7035) ..
  node[above=0.07cm]{$1 \t \eta_{a'}$}(143.0669,1009.0044);
  \path[draw=black,line join=miter,line cap=butt,line width=0.650pt]
  (193.6849,989.9651) .. controls (180.9955,1003.4434) and (155.7228,999.3710) ..
  node[above left=0.07cm,pos=0.4]{$m$}(144.7842,1007.6319);
  \path[draw=black,line join=miter,line cap=butt,line width=0.650pt]
  (243.4783,1010.1883) .. controls (226.5217,1027.1448) and (158.3728,1024.5207) ..
  node[above=0.07cm]{$\eta_{a''}$}(144.8945,1009.7381);
  \path[draw=black,line join=miter,line cap=butt,line width=0.650pt]
  (142.7594,1006.5449) .. controls (130.3295,990.1682) and (42.0840,962.9596) ..
  node[above right=0.09cm,pos=0.73]{$L \t 1$}(11.8318,987.5314)
  (58.2412,1006.8614) .. controls (67.7385,1002.1403) and (86.7266,995.7900) ..
  node[above left=0.07cm,pos=0.6]{$1 \t m$}(108.9637,991.1125)
  (118.4142,989.2668) .. controls (143.3560,984.7724) and (171.0928,982.7810) ..
  node[above=0.07cm,pos=0.45]{$1 \t m$}(193.6858,987.4928);
  \path[draw=black,line join=miter,line cap=butt,line width=0.650pt]
  (-63.0317,1008.3986) .. controls (-47.8143,1013.8334) and (23.2494,1015.3069) ..
  node[above=0.07cm,pos=0.4]{$1\t m$}(43.9015,1010.0895);
  \path[draw=black,line join=miter,line cap=butt,line width=0.650pt]
  (-62.6087,1011.0578) .. controls (-42.3135,1042.4784) and (142.3552,1036.7664) ..
  node[above=0.07cm,pos=0.5]{$m$}(143.2609,1011.0578);
  \path[draw=black,line join=miter,line cap=butt,line width=0.650pt]
  (-64.3478,1010.6231) .. controls (-72.6087,1019.3188) and (-115.4913,1022.3622) ..
  node[above right=0.12cm,at end]{\!$m$}(-136.0000,1022.3622);
  \path[draw=black,line join=miter,line cap=butt,line width=0.650pt]
  (-63.5604,1005.0168) .. controls (-42.2916,973.3827) and (0.5004,959.0455) ..
  (47.6744,955.3175) .. controls (84.0467,952.4430) and (178.6911,965.2136) ..
  node[above=0.07cm,pos=0.3]{$\mathfrak a$}(193.0709,984.1110)
  (-65.2174,1008.0144) .. controls (-65.2174,1008.0144) and (-107.2912,982.3622) ..
  node[above right=0.12cm,at end]{\!$m \t 1$}(-136.0000,982.3622)
  (-136.0000,1002.3622) .. controls (-122.0625,1002.3622) and (-107.6294,997.1434) ..
  node[above right=0.12cm,at start]{\!\!$1 \t \eta_a$}(-103.7218,992.0973)
  (-98.0477,988.9402) .. controls (-91.6444,985.5261) and (-85.0719,982.8030) ..
  (-76.8306,982.8030) .. controls (-67.8171,982.8030) and (-59.9275,984.6325) ..
  node[above=0.07cm,pos=0.1]{$1 \t \eta_a$}(-51.3302,987.3211)
  (-45.7679,989.1316) .. controls (-26.7271,995.5119) and (-2.7630,1004.8340) ..
  node[above=0.07cm,pos=0.4,rotate=-15]{$1 \t (1 \t \eta_a)$}(43.5941,1007.8410)
  (-136.0000,962.3622) .. controls (-101.5323,962.3622) and (-64.6834,963.6861) ..
  node[above right=0.12cm,at start]{\!\!$(L \t L) \t 1$}
  (-35.3006,969.9724) .. controls (-32.7815,970.5113) and (-30.3633,971.0833) .. (-28.0369,971.6843)
  (-20.7604,973.7434) .. controls (8.5916,982.8012) and (22.4311,996.5171) ..
  node[above=0.03cm,pos=0.245,rotate=-23]{$L \t (L \t 1)$}
  node[above=0.03cm,pos=0.75,rotate=-27]{$1 \t (L \t 1)$}(43.8261,1004.5361);
  \path[fill=black] (-64.88929,1008.671) node[circle, draw, line width=0.65pt, minimum width=5mm, fill=white, inner sep=0.25mm] (text3838) {$\ass$     };
  \path[fill=black] (50.913471,1008.0961) node[circle, draw, line width=0.65pt, minimum width=5mm, fill=white, inner sep=0.25mm] (text3842) {$1 \Theta^{-1}$     };
  \path[fill=black] (142.57112,1009.2459) node[circle, draw, line width=0.65pt, minimum width=5mm, fill=white, inner sep=0.25mm] (text3846) {$\Theta^{-1}$     };
  \path[fill=black] (243.75063,1008.0961) node[circle, draw, line width=0.65pt, minimum width=5mm, fill=white, inner sep=0.25mm] (text3850) {$\Theta$     };
  \path[fill=black] (194.31064,986.82544) node[circle, draw, line width=0.65pt, minimum width=5mm, fill=white, inner sep=0.25mm] (text3854) {$\ass$     };
  \path[shift={(-64.732665,56.496783)},draw=black,fill=black] (77.8600,931.1000)arc(0.000:180.000:1.250)arc(-180.000:0.000:1.250) -- cycle;
\end{tikzpicture}\rlap{ .}
\end{equation*}
This completes the definition of the $\V$-functor $L \colon \A \to \C$; we shall now make the morphisms $(\eta_a)_b \colon M(b,a) \to \C(Fb, La)$ into the components of a bimodule morphism $M \to \C(F,L)$. Each $\eta_a$ is already a right $\B$-module morphism; while the $2$-cells making each $(\eta_{\thg})_b$ into a left $\A$-module morphism, are obtained by pasting the components of $\Theta_{aa'}$ with pseudofunctoriality of $\otimes$.
Finally, we check that $\eta \colon M \to \C(F,L)$ exhibits $L$ as $M \star F$, which is equally to check that
\begin{equation*}
\C(L, 1), M \xrightarrow{1, \eta} \C(L, 1), \C(F, L) \xrightarrow \C(F,1)
\end{equation*}
exhibits $\C(L,1)$ as $\h{M}{\C(F,1)}$. But the morphism obtained by evaluating at any $a \in \A$ and $c \in \C$ exhibits $\C(La, c)$ as $\h{M(\thg,a)}{\C(F, c)}$, since each $\eta_a$ exhibits $La$ as $M(\thg, a) \star F$; and so the result follows by Proposition~\ref{prop:right-hom-pointwise}.
\end{proof}

We now discuss the uniqueness of these more general kinds of colimits. 
\begin{Prop}\label{prop:uniqueness-indexed-colimits}
Let $M$ be an $\A$-$\B$-bimodule. If the cylinder $\eta \colon M \to \C(F,V)$ exhibits $V$ as $M \star F$, then it also exhibits $V$ as birepresenting object for the functor ${}_\A \cat{Mod}_\B(M, \C(F, \thg)) \colon \V\text-\cat{Bicat}(\A,\C) \to \cat{CAT}$.
  Consequently, any two $M$-weighted colimits of $F$ are related  essentially-uniquely by an equivalence in $\V\text-\cat{Bicat}(\A, \C)$ commuting with the universal cylinders.
\end{Prop}
\begin{proof}
We must show that, for each $H \colon \A \to \C$, the functor 
\begin{equation}\label{eq:colimiting-cylinder-ordinary-universal}
\begin{aligned}
\V\text-\cat{Bicat}(\A,\C)(V, H) & \to {}_\A \cat{Mod}_{\B}(M, \C(F,H)) \\
\gamma & \mapsto \C(1, \gamma) \circ \eta
\end{aligned}
\end{equation}
is an equivalence of categories.
By definition of colimit and  Proposition~\ref{prop:right-hom-pointwise}, the bimorphism
\begin{equation*}
\psi = \C(V,H), M \xrightarrow{1, \eta} \C(V,H), \C(F,V) \xrightarrow{m} \C(F,H)
\end{equation*}
exhibits $\C(V, H)$ as $\h{M}{\C(F,H)}$; while by the dual of Corollary~\ref{cor:yoneda-with-naturality}, the bimorphism
\begin{equation*}
\phi = \A, \C(F,H) \xrightarrow{\hat H, 1} \C(H,H), \C(F,H) \xrightarrow{m} \C(F,H)
\end{equation*}
exhibits $\C(F,H)$ as $\h{\A}{\C(F,H)}_\ell$.
Thus, in the diagram of categories and functors
\begin{equation*}
\cd[@C+1.5em]{
\V\text-\cat{Bicat}(\A,\C)(V, H) \ar[d]_{\eqref{eq:colimiting-cylinder-ordinary-universal}}
\ar[r]^-{\text{\eqref{eq:vnats-as-bimods}}} &
{}_\A \cat{Mod}_{\A}(\A, \C(V,H))
\ar[d]^{\psi \circ (\thg, 1)} \\
{}_\A \cat{Mod}_{\B}(M, \C(F,H)) \ar[r]_-{\phi \circ (1, \thg)} &
\cat{Bimor}_{\A\A\B}(\A, M; \C(F,H))
}
\end{equation*}
the top, bottom and right sides are equivalences; it will follow that the left side is too, as required, so long as we can show that the square commutes to within natural isomorphism. Evaluating at some $\V$-transformation $\gamma \colon V \Rightarrow H$, the two sides of the square
yield the respective bimorphisms
\begin{gather*}
\A, M \xrightarrow{\hat V, 1} \C(V,V), M \xrightarrow{\C(1,\gamma), 1} \C(V, H), M \xrightarrow{1, \eta} \C(V,H), \C(F,V) \xrightarrow{m} \C(F,H) \\
\text{and} \ \
\A, M \xrightarrow{1, \eta} \A, \C(F,V) \xrightarrow{1, \C(1, \gamma)} \A, \C(F, H) \xrightarrow{\hat H, 1} \C(H, H), \C(F, H) \xrightarrow{m} \C(F,H)\rlap{ ,}
\end{gather*}
between which we obtain an invertible modification with components:\begin{equation*}
\begin{tikzpicture}[y=0.8pt, x=0.8pt,yscale=-1, inner sep=0pt, outer sep=0pt, every text node part/.style={font=\tiny} ]
  \path[draw=black,line join=miter,line cap=butt,line width=0.650pt] 
  (-10.0000,932.3622) .. controls (16.5129,932.3622) and (43.9759,935.3742) .. 
  node[above right=0.12cm,at start] {\!$V \t 1$}(47.8916,940.2537);
  \path[draw=black,line join=miter,line cap=butt,line width=0.650pt] 
  (-10.0000,952.3622) .. controls (13.2006,952.3622) and (43.9759,949.3501) .. 
  node[above right=0.12cm,at start] {\!$\C(1,\!\gamma) \t 1$}(47.5904,944.5309);
  \path[draw=black,line join=miter,line cap=butt,line width=0.650pt] 
  (-10.0000,992.3622) .. controls (23.7470,992.3622) and (93.3735,990.5549) .. 
  node[above right=0.12cm,at start] {\!$m$}(98.1928,984.8321);
  \path[draw=black,line join=miter,line cap=butt,line width=0.650pt] 
  (102.4096,984.1694) .. controls (112.0482,991.6995) and (161.3156,992.3622) .. 
  node[above left=0.12cm,at end] {$m$\!}(180.0000,992.3622);
  \path[draw=black,line join=miter,line cap=butt,line width=0.650pt] 
  (102.4096,980.3742) .. controls (109.4513,972.8026) and (117.1751,967.1803) .. 
  node[below right=0.04cm,rotate=0,pos=0.45] {$1 \t \C(1,\!\gamma)$\!}(125.3980,963.0445)
  (131.2682,960.3720) .. controls (146.3991,954.1680) and (162.9937,952.3848) .. 
  node[above left=0.12cm,at end] {$1 \t \C(1,\!\gamma)$\!}(180.0000,952.3622)
  (52.7108,944.8321) .. controls (55.1205,954.4706) and (91.5663,969.2297) .. 
  node[below left=0.04cm,pos=0.73] {$\C(\gamma,\!1)\t 1$\!}(97.2892,980.3742)
  (-10.0000,972.3622) .. controls (26.8989,972.3622) and (50.3715,967.8320) .. 
  node[above right=0.12cm,at start] {\!$1 \t \eta$}(68.6778,961.8207)
  (74.7470,959.7125) .. controls (84.5785,956.1137) and (93.0354,952.1460) .. 
  node[below right=0.07cm,pos=0.5] {$1 \t \eta$}(101.6505,948.3754)
  (107.2561,945.9715) .. controls (125.1243,938.5024) and (145.0276,932.3622) .. 
  node[above left=0.12cm,at end] {$1 \t \eta$\!}(180.0000,932.3622)
  (52.4096,940.3140) .. controls (111.2967,925.2388) and (104.5418,972.4979) .. 
  node[above left=0.12cm,at end] {$H \t 1$\!}(180.0000,972.3622);
  \path[fill=black] (49.397587,942.72363) node[circle, draw, line width=0.65pt, minimum width=5mm, fill=white, inner sep=0.25mm] (text3071) {$\bar \gamma^{-1} \t 1$     };
  \path[fill=black] (99.397583,982.78387) node[circle, draw, line width=0.65pt, minimum width=5mm, fill=white, inner sep=0.25mm] (text3075) {$\alpha$   };
\end{tikzpicture}\rlap{ .}\qedhere
\end{equation*}
\end{proof}

\subsection{Left Kan extensions}
Suppose we are given a diagram
\begin{equation*}
\cd[@+0.5em]{
\A \ar[r]^-F \ar[d]_G \dtwocell[0.3]{dr}{\gamma} & \C \\
\B \ar[ur]_H & {}
}
\end{equation*}
of $\V$-bicategories, $\V$-functors and a $\V$-transformation. We say that \emph{$\gamma$ exhibits $H$ as the left Kan extension of $F$ along $G$}, or that \emph{$\gamma$ exhibits $H$ as $\Lan_G F$} if the morphism of $\B$-$\A$-bimodules
\begin{equation}\label{eq:inducedbykan}
\B(G,1) \xrightarrow{\hat H} \C(HG, H) \xrightarrow{\C(\gamma, 1)} \C(F, H)
\end{equation}
exhibits $H$ as $\B(G,1) \star F$. So by Proposition~\ref{prop:right-hom-pointwise}, the existence of the left Kan extension implies the existence of the colimit $\B(G,b) \star F$ for each $b \in \B$. 
(Our terminology follows that of~\cite{kelly:enriched}; according to~\cite{maclane} these would be ``pointwise'' Kan extensions.)

Conversely, if $\B(G,b) \star F$ exists for each $b \in \B$, then on applying Proposition~\ref{prop:functcolimits} to the $\B$-$\A$-bimodule $\B(G,1)$, we obtain a $\V$-functor $H \colon \B \to \C$ and a morphism $\eta \colon \B(G, 1) \to \C(F, H)$ of $\B$-$\A$-bimodules exhibiting $H$ as $\B(G,1) \star F$. Applying Proposition~\ref{prop:vnat-for-kan} to $\eta$, we obtain a $\V$-transformation $\gamma \colon F \Rightarrow HG$ such that the composite~\eqref{eq:inducedbykan} is isomorphic to $\eta$. Since $\eta$ is colimiting, so too is this~\eqref{eq:inducedbykan}, whence this $\gamma$ exhibits $H$ as $\Lan_G F$.

\begin{Prop}\label{prop:left-kan-ext-left-biadj}
Suppose that $\gamma \colon F \Rightarrow HG$ exhibits $H$ as $\Lan_G
F$. Then $\gamma$ exhibits $H$ as the value at $F$ of a (partial) left biadjoint to the functor $(\thg) G \colon \V\text-\cat{Bicat}(\B,\C) \to
\V\text-\cat{Bicat}(\A,\C)$.
\end{Prop}
\begin{proof}
Given any $K \colon \B \to \C$, we have for each $\delta \colon H
\Rightarrow K$  the module morphisms
\begin{align*}
\B(G,1) \xrightarrow{\hat K} \C(KG, K) \xrightarrow{\C(\delta G, 1)} \C(HG, K) \xrightarrow{\C(\gamma, 1)} \C(F,K)\\
\text{and} \ \
\B(G,1) \xrightarrow{\hat H} \C(HG, H) \xrightarrow{\C(\gamma,1)} \C(F, H) \xrightarrow{\C(1,\delta)} \C(F,K)
\end{align*}
and between these, we may construct an invertible modification with components:
\begin{equation*}
\begin{tikzpicture}[y=0.8pt, x=0.8pt,yscale=-1, inner sep=0pt, outer sep=0pt, every text node part/.style={font=\tiny} ]
  \path[draw=black,line join=miter,line cap=butt,line width=0.650pt] 
  (-5.0000,1002.3622) .. controls (11.3877,1002.3622) and (44.8795,1005.0730) .. 
  node[above right=0.12cm,at start] {\!$K$}(47.5904,1010.1935);
  \path[draw=black,line join=miter,line cap=butt,line width=0.650pt] 
  (-5.0000,1022.3622) .. controls (13.2006,1022.3622) and (44.2771,1018.7478) .. 
  node[above right=0.12cm,at start] {\!$\C(\delta, 1)$}(48.1928,1015.4345);
  \path[draw=black,line join=miter,line cap=butt,line width=0.650pt] 
  (52.1687,1009.8923) .. controls (57.8916,1005.0730) and (104.0106,1002.3622) .. 
  node[above left=0.12cm,at end] {$H$\!}(120.0000,1002.3622);
  \path[draw=black,line join=miter,line cap=butt,line width=0.650pt] 
  (-5.0000,1042.3622) .. controls (69.5193,1042.3622) and (85.9622,1022.3622) .. 
  node[above right=0.12cm,at start] {\!$\C(\gamma, 1)$}
  node[above left=0.12cm,at end] {$\C(\gamma, 1)$\!}(120.0000,1022.3622)
  (53.0723,1014.8321) .. controls (54.2816,1018.6112) and (56.3967,1025.5311) .. (65.6000,1031.4874)
  (72.1112,1034.9519) .. controls (81.8507,1039.2172) and (96.7970,1042.3622) .. 
  node[above left=0.12cm,at end] {$\C(1, \delta)$\!}(120.0000,1042.3622);
  \path[fill=black] (49.698792,1012.6031) node[circle, draw, line width=0.65pt, minimum width=5mm, fill=white, inner sep=0.25mm] (text3490) {$\bar \delta$   };
\end{tikzpicture}\rlap{ .}
\end{equation*}
As $\delta$ varies, these invertible modifications provide the
components of a natural isomorphism filling the triangle:
\begin{equation*}
\cd[@C-3em]{
\V\text-\cat{Bicat}(\B,\C)(H, K) \ar[rr]^{(\thg)G \circ \gamma} \ar[dr]_(0.35){\C(1, \thg) \circ \C(\gamma, 1) \circ \hat H\ \ \ } & \twocong{d} & 
\V\text-\cat{Bicat}(\A,\C)(F, KG)\rlap{ .} \ar[dl]^(0.4){\ \ \ \C(\thg, 1) \circ \hat K} \\ &
{}_\B \cat{Mod}_{\B}(\B(G,1),\, \C(F,K))
}
\end{equation*}
The left arrow therein is invertible by
Proposition~\ref{prop:uniqueness-indexed-colimits}, while the right
one is so by Proposition~\ref{prop:vnat-for-kan}; whence the top arrow
is also an equivalence, as required.
\end{proof}

\begin{Prop}\label{prop:kan-ff-ff}
If $\gamma \colon F \Rightarrow HG$ exhibits $H$ as $\Lan_G F$, and $G$ is fully faithful, then $\gamma$ is an equivalence in $\V\text-\cat{Bicat}(\A,\C)$.
\end{Prop}
\begin{proof}
By Proposition~\ref{prop:equivifcomponentsequiv}, it suffices to show that each component $\gamma_a \colon Fa \to HGa$ is an equivalence in $\C_0$. As in Example~\ref{ex:yonedacolimit}, the functors $F$ and $G$ induce module morphisms $\zeta_G \colon \A(\thg, a) \to \B(G, Ga)$ and $\zeta_F \colon \A(\thg, a) \to \C(F,Fa)$. Because $G$ is fully faithful, the former of these has equivalences as its $1$-cell components, and so by Example~\ref{ex:pseudoinversemodulemorphism} admits a pseudoinverse $\zeta^\centerdot_G \colon \B(G,Ga) \to \A(\thg, a)$. Now consider the composite
\begin{equation*}
\phi \defeq \B(G,Ga) \xrightarrow{\zeta^\centerdot_G} \A(\thg, a) \xrightarrow{\zeta_F} \C(F,Fa)\rlap{ .}
\end{equation*}
Since $\zeta_F$ is colimiting, and $\zeta^\centerdot_G$ is an equivalence, it follows easily that their composite is also colimiting, and so exhibits $Fa$ as $\B(G,Ga) \star F$. On the other hand, by the definition of Kan extension, the module morphism
\begin{equation*}
\phi_a = \B(G,Ga) \xrightarrow{\hat H} \C(HG, HGa) \xrightarrow{\C(\gamma, Hb)} \C(F, HGa)
\end{equation*}
exhibits $HGa$ as $\B(G,Ga) \star F$. Thus by Proposition~\ref{prop:uniquenessofcolimits}, there exists an equivalence $\delta_a \colon Fa \to HGa$ in $\C_0$ and an invertible module transformation
\begin{equation*}
\cd[@C+1.5em]{
 \B(G,Ga) \ar[r]^-{\zeta_G^\centerdot} \ar[d]_{\hat H} \rtwocell{drr}{} &
 \A(\thg, a) \ar[r]^-{\zeta_F} &
 \C(F,Fa) \ar[d]^{\C(F,\delta_a)} \\
 \C(HG, HGa) \ar[rr]_-{\C(\gamma, Hb)} && \C(F, HGa)\rlap{ .}
}
\end{equation*}
To conclude that $\gamma_a$ is an equivalence, it now suffices to show that it is isomorphic to $\delta_a$ in $\C_0$. By Proposition~\ref{prop:uniquenessofcolimits}, it is enough to show that, on replacing $\delta_a$ with $\gamma_a$ in the preceding display, there is still an invertible module transformation mediating the centre. We obtain such by taking its components to be the $2$-cells
\begin{equation*}
\begin{tikzpicture}[y=0.9pt, x=0.9pt,yscale=-1, inner sep=0pt, outer sep=0pt, every text node part/.style={font=\tiny} ]
  \path[draw=black,line join=miter,line cap=butt,line width=0.650pt]
  (10.0000,982.3622) --
  node[above right=0.12cm,at start]{\!$1 \t \gamma$}
  node[above left=0.09cm,pos=0.93]{$(H \! \cdot \! G) \t \gamma$}(129.8178,983.2141);
  \path[draw=black,line join=miter,line cap=butt,line width=0.650pt]
  (10.0000,1002.3622) .. controls (50.0000,1002.3622) and (119.4507,1006.5630) ..
  node[above right=0.12cm,at start]{\!$m$}(128.9658,985.9727);
  \path[draw=black,line join=miter,line cap=butt,line width=0.650pt]
  (10.0000,962.3622) .. controls (60.0000,962.3622) and (116.6335,958.3362) ..
  node[above right=0.12cm,at start]{\!$\mathfrak r^\centerdot$}(128.9658,979.4009)
  (10.0000,942.3622) .. controls (23.9734,942.3622) and (31.2479,950.4224) ..
  node[above right=0.12cm,at start]{\!$H$}(39.2605,959.6210)
  (43.6356,964.5849) .. controls (50.1243,971.7417) and (57.8095,978.8465) ..
  node[left=0.12cm,pos=0.45]{$H \t 1$}(69.9898,982.8291)
  (70.0717,982.7681) .. controls (66.7974,976.9470) and (63.8342,970.7229) ..
  node[right=0.12cm,pos=0.55]{$G \t 1$}(62.0354,964.6202)
  (60.8298,959.6041) .. controls (59.5128,952.3637) and (60.2374,945.5115) ..
  (64.5065,939.9709) .. controls (79.5118,920.4960) and (164.5489,922.3622) ..
  node[above left=0.12cm,at end]{$G^\centerdot$\!}(220.0000,922.3622);
  \path[draw=black,line join=miter,line cap=butt,line width=0.650pt]
  (130.7517,985.4852) .. controls (142.0004,1007.3832) and (180.0000,1002.3622) ..
  node[above left=0.12cm,at end]{$m$\!}(220.0000,1002.3622);
  \path[draw=black,line join=miter,line cap=butt,line width=0.650pt]
  (131.6036,982.5034) --
  node[above right=0.12cm,pos=0.12]{$\gamma \t F$}
  node[above left=0.12cm,at end]{$\gamma \t 1$\!}(220.0000,982.3622);
  \path[draw=black,line join=miter,line cap=butt,line width=0.650pt]
  (131.6036,979.0957) .. controls (139.1864,959.1598) and (180.0000,962.3622) ..
  node[above left=0.12cm,at end]{$\mathfrak l^\centerdot$\!}(220.0000,962.3622)
  (165.2551,982.0775) .. controls (165.2551,982.0775) and (166.3450,973.9768) ..
  node[below right=0.06cm,pos=0.85]{$1 \t F$}(172.0911,965.0925)
  (175.9700,959.9066) .. controls (183.8099,950.8059) and (197.2215,942.3622) ..
  node[above left=0.12cm,at end]{$F$\!}(220.0000,942.3622);
  \path[shift={(-6.8577508,51.616449)},draw=black,fill=black] (77.8600,931.1000)arc(0.000:180.000:1.250)arc(-180.000:0.000:1.250) -- cycle;
  \path[shift={(88.645067,50.977493)},draw=black,fill=black] (77.8600,931.1000)arc(0.000:180.000:1.250)arc(-180.000:0.000:1.250) -- cycle;
  \path[fill=black] (129.89972,982.50342) node[circle, draw, line width=0.65pt, minimum width=5mm, fill=white, inner sep=0.25mm] (text6012) {$\bar \gamma$};
\end{tikzpicture}\rlap{ .} \qedhere
\end{equation*}
\end{proof}

\section{Colimits in categories of right modules}
\label{sec:moderate}
In this section we shall prove that---over a suitably well-behaved base $\V$---every category of moderate right modules $\M\B$ is cocomplete; and moreover, that for any $\V$-functor $F \colon \B \to \C$, the weighted colimit functor $(\thg) \star F \colon \M \B \to \C$ is cocontinuous insofar as it is defined.

\subsection{Cocompleteness of moderate right module categories}
We begin with a construction that we will use extensively in this section. Given a $\V$-functor $F \colon \A \to \M \B$, we have the $\A$-$\B$-bimodule $\M\B(Y,F)$, where $Y \colon \B \to \M \B$ is, as before, the Yoneda embedding. Evaluating at any $a \in \A$ yields the right $\B$-module $\M\B(Y, Fa)$, which by Proposition~\ref{prop:reformulate-yoneda} is equivalent to $Fa$. We may transport the bimodule structure of $\M\B(Y,F)$ along these equivalences to  yield an $\A$-$\B$-bimodule $\tilde F$ with $\tilde F(\thg,a) = Fa$ and with left $\A$-actions induced by the action of $F$ on homs.

Now for any right $\B$-module $W$, the hom-object $\h V W$ exists by assumption for any $V \in \M \B$; so by the construction of Section~\ref{subsec:modules-over-module-categories}, we may form the right $\A$-module $\M \B(F, W)$. There is now a module bimorphism $\varepsilon \colon \M\B(F,W), \tilde F \to W$ with $\A$-components $\varepsilon_{a \thg} = \xi_{Fa,W} \colon \h {Fa} W \otimes Fa \to W$, and with bimodule compatibility $2$-cells
given by
\begin{equation*}
\begin{tikzpicture}[y=0.75pt, x=0.75pt,yscale=-1, inner sep=0pt, outer sep=0pt, every text node part/.style={font=\tiny} ]
  \path[draw=black,line join=miter,line cap=butt,line width=0.650pt]
  (192.2929,1020.4333) .. controls (200.0000,1012.3622) and (230.0000,1012.3622) ..
  node[above left=0.12cm,at end] {$1 \t \xi$\!}(240.0000,1012.3622);
  \path[draw=black,line join=miter,line cap=butt,line width=0.650pt]
  (240.0000,1032.3622) .. controls (230.0000,1032.3622) and (200.0000,1032.3622) ..
  node[above left=0.12cm,at start] {$\xi$\!}(192.3107,1023.6330);
  \path[draw=black,line join=miter,line cap=butt,line width=0.650pt]
  (189.0792,1023.6520) .. controls (180.0000,1032.3622) and (160.7590,1032.3622) ..
  node[above right=0.12cm,at end] {\!$\xi$}(150.0000,1032.3622);
  \path[draw=black,line join=miter,line cap=butt,line width=0.650pt]
  (188.7267,1021.2262) .. controls (180.0000,1012.3622) and (160.2273,1012.3622) ..
  node[above right=0.12cm,at end] {\!$m \t 1$}(150.0000,1012.3622);
  \path[draw=black,line join=miter,line cap=butt,line width=0.650pt]
  (189.6201,1019.2846) .. controls (200.0000,1002.3622) and (188.5587,972.3622) ..
  node[above left=0.12cm,at end] {$\mathfrak a$\!}(240.0000,972.3622)
  (150.0000,992.3622) --
  node[above right=0.12cm,at start] {\!\!$(1 \t F) \t 1$}(195.5661,992.3622)
  (200.8457,992.3622) --
  node[above left=0.12cm,at end] {$1 \t (F \t 1)$\!\!}(240.0000,992.3622);
  \path[fill=black] (191.06065,1023.8874) node[circle, draw, line width=0.65pt, minimum width=5mm, fill=white, inner sep=0.25mm] (text3153) {$\bar m$
     };
\end{tikzpicture}\rlap{ .}
\end{equation*}
Note that by the universality of each component $\varepsilon_{a \thg}$ and Proposition~\ref{prop:right-hom-pointwise}, the bimorphism $\varepsilon$ expresses $\M\B(F,W)$ as $\h{\tilde F} W$.

\begin{Prop}\label{prop:mb-is-cocomplete}
If $\V$ is complete, cocomplete and left and right closed, then every $\V$-category $\M \B$ is cocomplete.
\end{Prop}
\begin{proof}
Let $\A$ be a small $\V$-bicategory, $V$ be a right $\A$-module and $F \colon \A \to \M \B$. We must show that the colimit $V \star F$ exists in $\M \B$. Let $\tilde F$ be the $\A$-$\B$-bimodule constructed from $F$ as above. Since $\A$ is small, and $\V$ biclosed and cocomplete, the tensor product $V \otimes_\A \tilde F$ exists; for brevity, we shall denote it by $C$, and write $\phi \colon V, \tilde F \to C$ for its universal bimorphism. We first show that $C$ is a moderate $\B$-module; which is to say that, for every right $\B$-module $W$, the functor ${}_\bullet \cat{Mod}_\B(\thg \otimes C, W) \colon \B_0 \to \cat{Cat}$ is birepresentable.

Given a right $\B$-module $W$, let $\varepsilon_W \colon \M\B(F, W),\,\tilde F \to W$ be as defined in the preceding section. Observe also that by Proposition~\ref{prop:copowers-of-universal}, for any $A \in \V$, the bimorphism $A \otimes \phi \colon A \otimes V, \tilde F \to A \otimes C$ is universal because $\phi$ is. Thus, for each $A \in \V$, we have equivalences of categories
\begin{gather}
\label{eq:first-equivalence-colims}
{}_\bullet \cat{Mod}_\B(A \otimes C, W) \xrightarrow{(\thg) \circ (A \otimes \phi)}
\cat{Bimor}(A \otimes V, \tilde F; W) \\
\label{eq:second-equivalence-colims}
{}_\bullet \cat{Mod}_\A(A \otimes V, \M \B(F, W)) \xrightarrow{\varepsilon_W \circ (\thg, 1)}
\cat{Bimor}(A \otimes V, \tilde F; W)
\end{gather}
pseudonatural in $A \in \V$. Composing the first with the pseudoinverse of the second, we conclude that $
{}_\bullet \cat{Mod}_\A(\thg \otimes V, \M \B(F, W)) \simeq {}_\bullet \cat{Mod}_\B(\thg \otimes C, W)$. Since $\A$ is small, every right $\A$-module is moderate by Proposition~\ref{prop:every-module-over-small-moderate}, and so the former functor is birepresented by $\h{V}{\M\B(F,W)}$; whence also the latter. This proves that $C$ is a moderate $\B$-module.

We now show that $C = V \otimes_\A \tilde F$ is in fact the colimit $V \star F$. Since $\varepsilon_C$ exhibits $\M\B(F,C)$ as $\h {\tilde F} C$, applying this universal property to the bimorphism $\phi \colon V, \tilde F \to C$ yields a module morphism $\psi \colon V \to \M\B(F,C)$ together with an invertible transformation $\Gamma \colon \varepsilon_C \circ (\psi, 1) \Rightarrow
\phi$. We claim that $\psi$ exhibits $C$ as $V \star F$. Thus, for each $W \in \M \B$, we must show that the module morphism
\begin{equation}\label{eq:morphism-to-be-shown-universal}
\M\B(C, W) \otimes V \xrightarrow{1 \otimes \psi}
\M\B(C, W) \otimes \M\B(F, C) \xrightarrow{m} \M\B(F, W)
\end{equation}
exhibits $\M\B(C,W)$ as $\h V{\M\B(F,W)}$. For this, it suffices to show that its image under~\eqref{eq:second-equivalence-colims} is isomorphic
to the image under~\eqref{eq:first-equivalence-colims} of the universal module morphism $\xi_{CW} \colon \M\B(C,W) \otimes C \to W$. Thus we must show that the bimorphisms
\begin{gather*}
\M\B(C,W) \otimes V,\, \tilde F \xrightarrow{\M\B(C,W) \otimes \phi} \M\B(C,W) \otimes C \xrightarrow{\xi} W \qquad\text{and}
\\
\M\B(C,W) \otimes V,\, \tilde F \xrightarrow{(1 \otimes \psi),\, 1}
\M\B(C,W) \otimes \M\B(F,C),\, \tilde F \xrightarrow{m, 1}
\M\B(F,W),\, \tilde F \xrightarrow{\varepsilon} W
\end{gather*}
are isomorphic. The $2$-cells witnessing this to be so are given by the composites
\begin{equation*}
\begin{tikzpicture}[y=0.9pt, x=0.9pt,yscale=-1, inner sep=0pt, outer sep=0pt, every text node part/.style={font=\tiny} ]
  \path[draw=black,line join=miter,line cap=butt,line width=0.650pt]
  (192.2929,1020.4333) .. controls (200.0000,1012.3622) and (233.0052,1019.7868) ..
  node[below=0.04cm] {$1 \t \xi$}(238.5858,1011.6551);
  \path[draw=black,line join=miter,line cap=butt,line width=0.650pt]
  (270.0000,1032.3622) .. controls (260.0000,1032.3622) and (200.0000,1032.3622) ..
  node[above left=0.12cm,at start] {$\xi$\!}(192.3107,1023.6330);
  \path[draw=black,line join=miter,line cap=butt,line width=0.650pt]
  (189.0792,1023.6520) .. controls (180.0000,1032.3622) and (160.7590,1032.3622) ..
  node[above right=0.12cm,at end] {\!$\xi$}(150.0000,1032.3622);
  \path[draw=black,line join=miter,line cap=butt,line width=0.650pt]
  (188.7267,1021.2262) .. controls (180.0000,1012.3622) and (160.2273,1012.3622) ..
  node[above right=0.12cm,at end] {\!$m \t 1$}(150.0000,1012.3622);
  \path[draw=black,line join=miter,line cap=butt,line width=0.650pt]
  (150.0000,992.3622) .. controls (150.0000,992.3622) and (177.0413,993.8235) ..
  node[above right=0.12cm,at start] {\!$(1 \t \psi) \t 1$}(201.4808,997.6226)
  (211.0916,999.2773) .. controls (223.9625,1001.7400) and (234.8121,1004.9577) ..
  node[below left=0.04cm,rotate=-10,pos=0.57] {$1 \t (\psi \t 1)$}(238.5858,1009.0797)
  (189.6201,1019.2846) .. controls (190.0000,992.3622) and (231.5587,992.3622) ..
  node[above left=0.12cm,at end] {$\mathfrak a$\!}(270.0000,992.3622);
  \path[draw=black,line join=miter,line cap=butt,line width=0.650pt]
  (240.1633,1009.9358) .. controls (250.0000,1012.3622) and (259.3934,1012.3622) ..
  node[above left=0.12cm,at end] {$1 \t \phi$\!}(270.0000,1012.3622);
  \path[fill=black] (239.1026,1010.2893) node[circle, draw, line width=0.65pt, minimum width=5mm, fill=white, inner sep=0.25mm] (text3572) {$1\t \Gamma$   };
  \path[fill=black] (191.06065,1023.8874) node[circle, draw, line width=0.65pt, minimum width=5mm, fill=white, inner sep=0.25mm] (text3153) {$\bar m$     };
\end{tikzpicture}\rlap{ .} \qedhere
\end{equation*}
\end{proof}

\subsection{Cocontinuity of colimits in the weight}
We now give results proving the cocontinuity of the colimit operation, insofar as it is defined, in the weight. To this end, let $\A$, $\B$ and $\C$ be $\V$-bicategories with $\A$ and $\B$ small, let $V$ be a right $\A$-module, let $F \colon \A \to \M \B$, and let $G \colon \B \to \C$. Suppose that 
the colimit $Fa \star G$ exists in $\C$ for every $a \in \A$; then by Proposition~\ref{prop:functcolimits}, the assignation $a \mapsto Fa \star G$ is the object assignation of a functor $\tilde F \star G \colon \A \to \C$.

\begin{Prop}\label{prop:cocontinuity-of-colimits}
In the circumstances just described, if the colimit $V \star (\tilde F \star G)$ exists in $\C$, then so too does the colimit $(V \star F) \star G$, and they are equivalent.
\end{Prop}
\newcommand{\ts}{.}
\begin{proof}
For the sake of brevity, we shall in this proof write a weighted colimit $W \star D$ simply as $W \ts D$; since we do not require functor composition, no confusion should arise. Let $\gamma \colon \tilde F \to \C(G, \tilde F \ts G)$ and $\delta \colon V \to \C(\tilde F \ts G, V \ts (\tilde F \ts G))$ be colimiting cylinders for $\tilde F \ts G$ and $V \ts (\tilde F \ts G)$ respectively. As in the preceding proof, we construct the colimit $V \ts F$ in $\M\B$ as the tensor product $V \otimes_\A \tilde F$. Let $\phi \colon V, \tilde F \to V \ts F$ be the universal bimorphism of this tensor product, and $\psi \colon V \to \M\B(F,V \ts F)$ the corresponding colimiting cylinder; as before, it comes equipped with an invertible transformation $\Gamma \colon \varepsilon \circ (\psi, 1) \Rightarrow
\phi$.

Forming the bimorphism
\begin{equation*}
V, \tilde F \xrightarrow{\delta, \gamma} \C(\tilde F \ts G, V \ts (\tilde F \ts G)),\, \C(G, \tilde F \ts G)
\xrightarrow{m} \C(G, V \ts (\tilde F \ts G))
\end{equation*}
and applying the universality of $\phi$ yields a module morphism $\eta \colon V \ts F \to \C(G, V \ts (\tilde F \ts G))$ and an invertible module transformation $\Delta 
\colon m \circ (\delta, \gamma) \Rightarrow \eta \circ \phi$. We claim that $\eta$ exhibits $V \ts (\tilde F \ts G)$ as $(V \ts F) \ts G$. Thus for each $X \in \C$, we must show that the 
module morphism
\begin{equation}\label{eq:cocont-colim-universal}
\C(V \ts (\tilde F \ts G), X) \otimes (V \ts F) \xrightarrow{1 \otimes \eta} \C(V \ts (\tilde F \ts G), X) \otimes \C(G, V \ts (\tilde F \ts G)) \xrightarrow m \C(G, X)
\end{equation}
exhibits $\C(V \ts (\tilde F \ts G), X)$ as $\h{V \ts F}{\,\C(G,X)}$. To do so, note first that, under the equivalences~\eqref{eq:first-equivalence-colims} and~\eqref{eq:second-equivalence-colims}, the morphism \eqref{eq:cocont-colim-universal} corresponds to a morphism $\C(V \ts (\tilde F \ts G), X) \otimes V \to \M\B(F, \C(G,X))$. But since the two morphisms
\begin{gather*}
\C(\tilde F \ts G, X), \tilde F \xrightarrow{1, \gamma} \C(\tilde F \ts G, X), \C(G, \tilde F \ts G) \xrightarrow{m} \C(G, X) \\
\text{and}\qquad \M\B(F, \C(G,X)), \tilde F \xrightarrow{\varepsilon} \C(G,X)
\end{gather*}
both exhibit their first argument as the hom $\h{\tilde F}{\C(G,X)}$, we conclude that there is an equivalence of modules $\zeta \colon \C(\tilde F \ts G, X) \to \M\B(F, \C(G,X))$ and an invertible module transformation $\Upsilon \colon \varepsilon \circ (\zeta, 1) \Rightarrow m \circ (1,\gamma)$. Thus, on applying~\eqref{eq:first-equivalence-colims} and~\eqref{eq:second-equivalence-colims} and composing with $\zeta$,~\eqref{eq:cocont-colim-universal} corresponds to a morphism $\C(V \ts (\tilde F \ts G), X) \otimes V \to \C(\tilde F \ts G, X)$, and it suffices to show that this morphism is universal. By the universal property of $V \ts (\tilde F \ts G)$, we have the the universal morphism
\begin{equation*}
\C(V \ts (\tilde F \ts G),X) \otimes V \xrightarrow{1 \otimes \delta} \C(V \ts (\tilde F \ts G), X) \otimes \C(\tilde F \ts G, V \ts (\tilde F \ts G)) \xrightarrow m \C(\tilde F \ts G, X) \ \rlap;
\end{equation*}
composing with $\zeta$ yields a universal morphism $\C(V \ts (\tilde F \ts G), X) \rightarrow \M\B(F, \C(G,X))$, and we will be done if we can show that the image of this under~\eqref{eq:second-equivalence-colims} is isomorphic to the image under~\eqref{eq:first-equivalence-colims} of~\eqref{eq:cocont-colim-universal}. The invertible module modification witnessing this to be the case has components
\begin{equation*}
\begin{tikzpicture}[y=0.8pt, x=0.8pt,yscale=-1, inner sep=0pt, outer sep=0pt, every text node part/.style={font=\tiny} ]
  \path[draw=black,line join=miter,line cap=butt,line width=0.650pt] 
  (24.0000,962.3622) .. controls (49.6183,962.3622) and (57.0120,966.8803) .. 
  node[above right=0.12cm,at start] {\!$\zeta \t 1$}(58.5181,969.2899);
  \path[draw=black,line join=miter,line cap=butt,line width=0.650pt] 
  (24.0000,982.3622) .. controls (47.4940,982.3622) and (56.4096,976.9405) .. 
  node[above right=0.12cm,at start] {\!$\xi$}(58.5181,974.5309);
  \path[draw=black,line join=miter,line cap=butt,line width=0.650pt] 
  (24.0000,942.3622) .. controls (70.7838,942.3622) and (120.1722,955.3166) .. 
  node[above right=0.12cm,at start] {\!$m \t 1$}(126.4975,961.9431)
  (62.7952,968.9887) .. controls (63.1262,963.7860) and (68.1114,956.5989) .. 
  node[below right=0.06cm,pos=0.8] {$1 \t \gamma$}(73.4083,949.0290)
  (76.5699,944.4820) .. controls (80.9385,938.1081) and (84.9253,931.6894) .. 
  node[below right=0.06cm,pos=0.6] {$1 \t \gamma$}(86.0273,926.1491);
  \path[draw=black,line join=miter,line cap=butt,line width=0.650pt] 
  (62.4940,972.8321) .. controls (70.1186,978.8776) and (121.5581,973.4597) .. 
  node[above=0.09cm,pos=0.5] {$m$}(126.2479,964.9831);
  \path[draw=black,line join=miter,line cap=butt,line width=0.650pt] 
  (129.7248,965.8351) .. controls (138.3564,977.0955) and (185.3562,982.3622) .. 
  node[above left=0.12cm,at end] {$m$\!}(220.0000,982.3622);
  \path[draw=black,line join=miter,line cap=butt,line width=0.650pt] 
  (129.6797,963.6664) .. controls (140.8243,956.7387) and (170.3614,961.7598) .. 
  node[above=0.07cm,pos=0.5] {$1 \t m$}(178.1928,953.9284);
  \path[draw=black,line join=miter,line cap=butt,line width=0.650pt] 
  (129.3719,961.1642) .. controls (133.6105,937.0273) and (141.9880,922.3622) .. 
  node[above left=0.12cm,at end] {$\mathfrak a$\!}(220.0000,922.3622)
  (24.0000,922.3622) .. controls (69.9157,922.3622) and (110.2579,928.8927) .. 
  node[above right=0.12cm,at start] {\!\!$(1 \t \delta) \t 1$}
  node[above right=0.07cm,pos=0.63] {$(1 \t \delta) \t \gamma$}(138.0800,935.9561)
  (143.9935,937.5155) .. controls (163.1004,942.7501) and (175.1434,948.0353) .. 
  node[above right=0.01cm,pos=0.24] {$1 \t (\delta \t \gamma)$}(177.3494,950.9766);
  \path[draw=black,line join=miter,line cap=butt,line width=0.650pt] 
  (182.1687,950.4947) .. controls (188.1928,945.3742) and (200.1913,942.3622) .. 
  node[above left=0.12cm,at end] {$1 \t \phi$\!}(220.0000,942.3622);
  \path[draw=black,line join=miter,line cap=butt,line width=0.650pt] 
  (182.1687,954.2296) .. controls (185.4819,959.3501) and (200.1913,962.3622) .. 
  node[above left=0.12cm,at end] {$1 \t \eta$\!}(220.0000,962.3622);
  \path[cm={{-1.0,0.0,0.0,1.0,(162.43967,-4.8365028)}},draw=black,fill=black] (77.8571,931.1122)arc(0.000:180.000:1.250)arc(-180.000:0.000:1.250) -- cycle;
  \path[fill=black] (60.430504,972.28021) node[circle, draw, line width=0.65pt, minimum width=5mm, fill=white, inner sep=0.25mm] (text4251) {$\Upsilon$     };
  \path[fill=black] (127.79037,964.61279) node[circle, draw, line width=0.65pt, minimum width=5mm, fill=white, inner sep=0.25mm] (text4255) {$\alpha^{-1}$     };
  \path[fill=black] (180.18443,952.2597) node[circle, draw, line width=0.65pt, minimum width=5mm, fill=white, inner sep=0.25mm] (text4259) {$1 \t \Delta$   };
\end{tikzpicture}\rlap{ .}\qedhere
\end{equation*}
\end{proof}

We now give a general result expressing the cocontinuity of taking colimits in the weight. Let $G \colon \B \to \C$ be a $\V$-functor,  let $\D$ be the full sub-$\V$-category of $ \M \B$ on those weights $W$ for which $W \star G$ exists in $\C$, and let $J \colon \D \hookrightarrow \M \B$ be the inclusion $\V$-functor. Since for every $W \in \D$, the weighted colimit $\tilde J(
\thg, W) \star G = W \star G$ exists in $\C$, so too does the $\tilde J$-weighted colimit of $G$; we write the underlying functor of this colimit as $(\thg) \star G \colon \D \to \C$, since on objects it sends $W$ to $W \star G$.

\begin{Prop}\label{prop:colims-cocont-in-weight}
In the situation just described, the functor $(\thg) \star G \colon \D \to \C$ preserves every colimit in $\D$ that is preserved by the inclusion $\D \hookrightarrow \M \B$.
\end{Prop}
\begin{proof}
As in the preceding proof, we write weighted colimits as $W \ts D$ rather than $W \star D$. To say that the $V$-weighted colimit of a diagram $F \colon \A \to \D$ (with $\A$ small) exists in $\D$ and is preserved by the inclusion $\D \hookrightarrow \M \B$ is to say that the colimit $V \ts F$, computed in $\M \B$, again lies in $\D$. We show that any such colimit is preserved by $(\thg) \ts G$. For such $V$ and $F$, we have as in the preceding proof $\phi \colon V, \tilde F \to V \ts F$ exhibiting $V \ts F$ as $V \star \tilde F$, inducing $\psi \colon V \to \D(F, V \ts F)$ the colimiting cylinder, seen now as landing in $\D \subset \M \B$. We must show that the composite
\begin{equation}\label{eq:cocont-weight-tobeshownuniversal}
V \xrightarrow{\psi} \D(F, V \ts F) \xrightarrow{(\thg) \ts G} \C(F(\thg) \ts G, (V \ts F) \ts G)
\end{equation}
is again colimiting. Observe first that the functor $F(\thg) \ts G$ has the same action on objects as $\tilde F \ts G$, and has its action on morphisms determined by the same universal property; we may therefore without loss of generality, assume that $F(\thg) \ts G = \tilde F \ts G$. We will also assume without loss of generality that $(V \ts F) \ts G$ has been constructed as $V \ts (\tilde F \ts G)$, with the colimiting cylinder $\eta \colon V \ts F \to \C(G, V \ts (\tilde F \ts G))$ as in the previous proof. Under these assumptions,~\eqref{eq:cocont-weight-tobeshownuniversal} becomes a morphism $V \to \C(\tilde F \ts G, V \ts (\tilde F \ts G))$; and to show this is colimiting, it suffices to show that it is isomorphic to the colimiting cylinder $\delta$ of the previous proof. Let $\gamma \colon \tilde F \to \C(G, \tilde F \ts G)$ be, as before, a colimiting cylinder for $\tilde F \ts G$. By definition of colimit, the module bimorphism
\begin{equation*}
\C(\tilde F \ts G, V \ts (\tilde F \ts G)), \tilde F \xrightarrow{1, \gamma} \C(\tilde F \ts G, V \ts (\tilde F \ts G)), \C(G, \tilde F \ts G) \xrightarrow{m} \C(G, V \ts (\tilde F \ts G)
\end{equation*}
exhibits its first argument as $\h{\tilde F}{\C(G, V \ts (\tilde F \ts G))}$. Thus to show that~\eqref{eq:cocont-weight-tobeshownuniversal} and $\delta$ are isomorphic, it suffices to construct an isomorphism between their composites with the above-displayed module bimorphism. We obtain such by taking its $2$-cell components to be
\begin{equation*}
\begin{tikzpicture}[y=0.8pt, x=0.96pt,yscale=-1, inner sep=0pt, outer sep=0pt, every text node part/.style={font=\tiny} ]
  \path[draw=black,line join=miter,line cap=butt,line width=0.650pt] 
  (10.0000,922.3622) .. controls (24.4679,922.3622) and (36.0843,934.7718) .. 
  node[above right=0.12cm,at start] {\!$\delta 1$}(37.8916,939.2899);
  \path[draw=black,line join=miter,line cap=butt,line width=0.650pt] 
  (10.0000,942.3622) -- 
  node[above right=0.12cm,at start] {\!$1 \gamma$}(37.5904,942.3622);
  \path[draw=black,line join=miter,line cap=butt,line width=0.650pt] 
  (10.0000,962.3622) .. controls (25.0633,962.3622) and (36.3855,949.7116) .. 
  node[above right=0.12cm,at start] {\!$m$}(37.8916,945.4947);
  \path[draw=black,line join=miter,line cap=butt,line width=0.650pt] 
  (74.3976,928.0248) .. controls (82.2289,922.0007) and (105.9519,912.3622) .. 
  node[above left=0.12cm,at end] {$\psi \t 1$\!}(145.0000,912.3622);
  \path[draw=black,line join=miter,line cap=butt,line width=0.650pt] 
  (101.5060,948.7477) .. controls (104.5181,944.2297) and (121.3229,932.3622) .. 
  node[above left=0.12cm,at end] {$(\thg \star G)1$\!}(145.0000,932.3622);
  \path[draw=black,line join=miter,line cap=butt,line width=0.650pt] 
  (102.4096,952.3622) -- 
  node[above left=0.12cm,at end] {$1 \t \gamma$\!}(145.0000,952.3622);
  \path[draw=black,line join=miter,line cap=butt,line width=0.650pt] 
  (102.1084,955.1935) .. controls (103.9157,960.3140) and (120.1000,972.3622) .. 
  node[above left=0.12cm,at end] {$m$\!}(145.0000,972.3622);
  \path[draw=black,line join=miter,line cap=butt,line width=0.650pt] 
  (99.0964,950.2538) .. controls (92.1687,949.6513) and (76.2048,940.6152) .. 
  node[above=0.1cm,pos=0.35] {$\varepsilon$}(73.7952,930.9766);
  \path[draw=black,line join=miter,line cap=butt,line width=0.650pt] 
  (39.4578,944.2297) .. controls (43.6747,955.6754) and (75.0000,966.5188) .. 
  node[above=0.08cm,pos=0.5] {$\eta$}(99.0964,953.2658);
  \path[draw=black,line join=miter,line cap=butt,line width=0.650pt] 
  (39.7590,941.2176) .. controls (40.9639,935.7959) and (55.4217,922.5429) .. 
  node[above=0.08cm,pos=0.5] {$\phi$}(72.5904,928.8682);
  \path[fill=black] (38.25301,942.72363) node[circle, draw, line width=0.65pt, minimum width=5mm, fill=white, inner sep=0.25mm] (text4678) {$\Delta$     };
  \path[fill=black] (73.192764,929.47064) node[circle, draw, line width=0.65pt, minimum width=5mm, fill=white, inner sep=0.25mm] (text4682) {$\Gamma^{-1}$     };
  \path[fill=black] (100.6024,952.66339) node[circle, draw, line width=0.65pt, minimum width=5mm, fill=white, inner sep=0.25mm] (text4686) {$\Theta$   };
\end{tikzpicture}
\end{equation*}
where $\Gamma$ and $\Delta$ are as in the preceding proof, and $\Theta$ is the invertible $2$-cell of~\eqref{eq:the-theta}.
\end{proof}

\section{Free cocompletions}
\label{sec:free-cocompletions}
We are now ready to define the free cocompletion of a $\V$-bicategory under a class of colimits. Through this section, we assume that our base bicategory $\V$ is left and right closed, complete and cocomplete.

By a \emph{class of weights} $\Phi$, we mean a collection of pairs $(\A, W)$ with $\A$ a small $\V$-bicategory and $W$ a right $\A$-module. A $\V$-bicategory $\B$ is said to be \emph{$\Phi$-cocomplete} if, for every $(\A, W) \in \Phi$ and every $\V$-functor $F \colon \A \to \B$, the weighted colimit $W \star F$ exists in $\B$; a $\V$-functor $G \colon \B \to \C$ between $\Phi$-cocomplete categories is said to be \emph{$\Phi$-cocontinuous} if it preserves every such colimit. Given $\Phi$-cocomplete $\V$-bicategories $\B$ and $\C$, we write $\Phi\text-\cat{Cocts}(\B, \C)$ for the bicategory of $\Phi$-cocontinuous $\V$-functors together with all transformations and modifications between them.

Given a class of weights $\Phi$ and a $\V$-bicategory $\C$, let us write $\Phi(\C)$ for the closure of the representables in $\M \B$ under $\Phi$-weighted colimits; that is, the smallest full, equivalence-replete, sub-$\V$-bicategory of $\M \B$ that contains the representables and that, for any $(\A, W) \in \Phi$, contains $W \star F$ whenever it contains each $Fb$. By Proposition~\ref{prop:mb-is-cocomplete} and our standing hypotheses, $\M \B$ is cocomplete, and it follows immediately that $\Phi(\B)$ is $\Phi$-cocomplete, and that the inclusion $J \colon \Phi(\B) \to \M \B$ is $\Phi$-cocontinuous. Moreover, the Yoneda embedding $Y \colon \B \to \M \B$ factors through $\Phi(\B)$, yielding a fully faithful functor $Z \colon \B \to \Phi(\B)$.

\begin{Thm}
$Z \colon \B \to \Phi(\B)$ exhibits $\Phi(\B)$ as the free completion of $\B$ under $\Phi$-colimits; which is to say that for every $\Phi$-cocomplete $\C$, the functor
\begin{equation}\label{eq:free-cocomp-functor}
\Phi\text-\cat{Cocts}(\Phi(\B), \C) \to \V\text-\cat{Bicat}(\B, \C)
\end{equation}
induced by composition with $Z$ is a biequivalence, with a suitable weak inverse being given by left Kan extension along $Z$.
\end{Thm} 
\begin{proof}
Let $\C$ be a $\Phi$-cocomplete $\V$-bicategory. First we show that, for any $F \colon \B \to \C$, the left Kan extension $\Lan_Z F \colon \Phi(\B) \to \C$ exists; for this we must show that the colimit $\Phi(\B)(Z, W) \star F$ exists in $\C$ for any $W \in \Phi(\B)$. Since $\Phi(\B)(Z, W) = \M \B(Y,W) \simeq W$ by Proposition~\ref{prop:reformulate-yoneda}, it suffices to show that $W \star F$ exists in $\C$ for every $W \in \Phi(\B)$. So consider the class of all $W \in \Phi(\B)$ for which $W \star F$ exists in $\C$; this class is clearly equivalence-closed, contains the representables by Example~\ref{ex:yonedacolimit}, and is closed under $\Phi$-weighted colimits by Proposition~\ref{prop:cocontinuity-of-colimits}; and so is all of $\Phi(\B)$ as required.
Thus $\Lan_Z F$ exists for all $F \colon \B \to \C$, and it follows from Proposition~\ref{prop:colims-cocont-in-weight} that it is $\Phi$-cocontinuous. Given this, it follows by Proposition~\ref{prop:left-kan-ext-left-biadj}, that the functor~\eqref{eq:free-cocomp-functor} has a left biadjoint given by left Kan extension along $Z$. Since $Z$ is fully faithful, the unit of this biadjunction is, by Proposition~\ref{prop:kan-ff-ff}, a pseudonatural equivalence. It remains to show that the counit is likewise a pseudonatural equivalence.

So let $H \colon \Phi(\B) \to \C$ be a $\Phi$-cocontinuous $\V$-functor, and write $H' = \Lan_Z(HZ)$. We must show that the counit component $\varepsilon_H \colon H' \Rightarrow H \colon \Phi(\B) \to \C$ is an equivalence. This is equally to show that each $(\varepsilon_H)_W \colon H'W \to HW$ is an equivalence in $\C_0$. To this end, let $\U$ be the collection of $W \in \Phi(\B)$ such that $(\varepsilon_H)_W$ is an equivalence. Now, because the unit of $\Lan_Z \dashv Z^\ast$ is an equivalence, it follows from the triangle identities that the composite $(\varepsilon_H) \circ Z$ is also an equivalence. Thus $\U$ contains the representables, and it is clearly equivalence-closed; it will thus suffice to show that $\U$ is also closed under $\Phi$-colimits. To this end, let $F \colon \A \to \Phi(\B)$ take values in $\U$, and let $(\A, V) \in \Phi$; we must show that $V \star F$ also lies in $\U$. Consider the diagram
\begin{equation*}
\cd{
 V \ar[r]^-\eta \ar@{=}[d] &
 \Phi(\B)(F, V \star F) \ar[rr]^-{H'} & & \C(H'F, H'(V \star F)) \ar[d]^{\C(1, (\varepsilon_H)_{V \star F})} \\
 V \ar[r]_-{\eta} & \Phi(\B)(F, V \star F) \ar[r]_-{H} & \C(HF, H(V \star F)) \ar[r]_-{\C(\varepsilon_H F, 1)} & \C(H'F, H(V \star F))
}
\end{equation*}
where $\eta$ is a colimiting cylinder for $V \star F$. This rectangle commutes up to isomorphism by the $\V$-naturality of $\varepsilon_H$. Since $H'$ preserves $\Phi$-colimits, the top row is a colimiting cylinder. Since $H$ preserves $\Phi$-colimits, and $\varepsilon_H F \colon H'F \Rightarrow HF$ is by assumption an equivalence, the bottom row is also colimiting. It follows from Proposition~\ref{prop:uniquenessofcolimits} that $(\varepsilon_H)_{V \star F}$ is an equivalence, and so that $V \star F \in \U$ as required.
\end{proof}

\section{Change of base}
\label{sec:monoidal-adjunctions}
In this section, we consider \emph{change of base} for enriched bicategories. From a lax monoidal functor $\V \to \W$ between monoidal bicategories, we induce a trifunctor $\V\text-\cat{Bicat} \to \W\text-\cat{Bicat}$; while from monoidal transformations and modifications between the former, we induce tritransformations and trimodifications between the latter. We shall be particularly interested in the situation where we have a \emph{monoidal biadjunction} of monoidal bicategories, and in particular, the effect that this has on weighted colimits over the two bases for enrichment.

\subsection{The tricategory of monoidal bicategories}
Given $\V$ and $\W$ monoidal bicategories, a \emph{lax monoidal functor} $(L,\laxu,\laxc,\omega,\delta,\gamma)\colon\V\to\W$ is a lax functor in the sense of~\cite{Gordon1995Coherence} between the corresponding one-object tricategories, all of whose $3$-cell coherence data is invertible; thus the underlying lax functor of bicategories is strong, the lax-natural functoriality constraints $\laxu$ and $\laxc$ (called $\iota$ and $\chi$ in the terminology of~\cite{Gordon1995Coherence}) are pseudonatural, and the modifications $\omega$, $\delta$, and $\gamma$ are invertible.
In the terminology of~\cite[Definition~1]{gg:ldstr-tricat}, this is a ``lax homomorphism'' between one-object tricategories.
If $\laxu$ and $\laxc$ are additionally equivalences, so that $L$ becomes a strong functor of tricategories, we will call it a \emph{strong monoidal functor}.

A \emph{monoidal transformation} $L \Rightarrow K$ between lax monoidal functors is a ``pseudo-icon'' in the sense of~\cite[Definition~5]{gg:ldstr-tricat}; it comprises a transformation $\alpha \colon L \Rightarrow K$ between the underlying functors together with suitably coherent invertible modifications $M$ and $\Pi$ asserting the compatibility of $\alpha$ with $\laxu$ and $\laxc$ for $L$ and $K$. A \emph{monoidal modification} $\alpha \Rrightarrow \beta$ is, in the terminology of~\cite[Definition 6]{gg:ldstr-tricat}, a ``pseudo-icon modification''; it comprises a modification $\Gamma$ between underlying transformations satisfying compatibility axioms with $M$ and $\Pi$ for $\alpha$ and $\beta$.
Monoidal bicategories, lax monoidal functors, monoidal transformations and monoidal modifications comprise a tricategory $\cat{MonBicat}$; it is constructed, for example, in~\cite[Corollary~27]{gg:ldstr-tricat}.

\subsection{Change of base along a lax monoidal functor}\label{subsec:change-of-base-trihom}
Given $L \colon \V \to \W$ a lax monoidal functor, we now describe the induced change of base trifunctor $L 
\colon \V\text-\cat{Bicat} \to \W\text-\cat{Bicat}$.
First, given a $\V$-bicategory $\B$, the $\W$-bicategory $L\B$ has the same objects as $\B$, hom-objects $L\B(x,y) = L(\B(x,y))$, and remaining data obtained as follows:
\begin{itemize}
\item The unit map $j_x\colon I \to L\B(x,x)$ is the composite
  \[ I \xrightarrow{\laxu} L I \xrightarrow{L(j_x)} L(\B(x,x))\rlap{ ;}\]
\item The composition map $m_{xyz} \colon L \B(y,z) \otimes L\B(x,y) \to L\B(x,z)$ is the composite
  \[ L(\B(y,z)) \otimes L(\B(x,y)) \xrightarrow{\laxc} L(\B(y,z) \otimes \B(x,y)) \xrightarrow{L(m_{xyz})} L(\B(x,z))\rlap{ ;}\]
\item The associativity constraint $\alpha_{wxyz}$ is the composite
  \begin{equation*}
\vcenter{\hbox{
\begin{tikzpicture}[y=0.8pt, x=0.8pt,yscale=-1, inner sep=0pt, outer sep=0pt, every text node part/.style={font=\tiny} ]
  \path[draw=black,line join=miter,line cap=butt,line width=0.650pt] 
  (15.0000,992.3622) .. controls (70.8214,992.3622) and (66.7500,959.3622) .. 
  node[above right=0.12cm,at start] {\!$\laxc$}(140.7143,964.5765)
  (15.0000,972.3622) .. controls (35.1208,972.3622) and (47.9256,975.2715) .. 
  node[above right=0.12cm,at start] {\!\!$(Lm)1$}(59.8629,979.3192)
  (66.1342,981.5765) .. controls (83.6854,988.2585) and (98.1848,996.8425) .. 
  node[above right=0.04cm,pos=0.27] {$L(m \t 1)$}(100.6747,1001.1108);
  \path[draw=black,line join=miter,line cap=butt,line width=0.650pt] (15.0000,1012.3622) .. controls (46.2976,1012.3622) and (91.3909,1013.5024) .. 
  node[above right=0.12cm,at start] {\!$Lm$}(100.2436,1004.5306);
  \path[draw=black,line join=miter,line cap=butt,line width=0.650pt] (15.0000,952.3622) .. controls (77.7500,952.3622) and (124.9398,948.5405) .. 
  node[above right=0.12cm,at start] {\!$\laxc \t 1$}(140.0000,961.0050);
  \path[draw=black,line join=miter,line cap=butt,line width=0.650pt] (103.1493,1000.0969) .. controls (109.0839,994.5411) and (132.4643,977.6122) .. 
  node[below right=0.04cm] {$L\mathfrak a$}(140.5000,967.0765);
  \path[draw=black,line join=miter,line cap=butt,line width=0.650pt] (143.2143,961.3622) .. controls (143.2143,961.3622) and (170.0000,942.3622) .. 
  node[above left=0.12cm,at end] {$\mathfrak a$\!}(205.0000,942.3622);
  \path[draw=black,line join=miter,line cap=butt,line width=0.650pt] (144.2857,964.5765) .. controls (144.2857,964.5765) and (170.0000,962.3622) .. 
  node[above left=0.12cm,at end] {$1 \t \laxc$\!}(205.0000,962.3622);
  \path[draw=black,line join=miter,line cap=butt,line width=0.650pt] (102.8571,1005.5408) .. controls (109.0325,1016.4964) and (150.0000,1022.3622) .. 
  node[above left=0.12cm,at end] {$Lm$\!}(205.0000,1022.3622);
  \path[draw=black,line join=miter,line cap=butt,line width=0.650pt] 
  (143.5714,968.1479) .. controls (143.5714,968.1479) and (150.2500,1002.3622) .. 
  node[above left=0.12cm,at end] {$\laxc$\!}(205.0000,1002.3622)
  (103.6163,1003.2193) .. controls (126.4617,1008.5114) and (141.9441,1000.3538) .. 
  node[below right=0.04cm,pos=0.55] {$L(1 \t m)$}(158.9681,992.7424)
  (164.3364,990.4025) .. controls (174.7225,986.0275) and (185.9835,982.3622) .. 
  node[above left=0.12cm,at end] {$1 \t (Lm)$\!\!}(205.0000,982.3622);
  \path[fill=black] (142.25,964.86218) node[circle, draw, line width=0.65pt, minimum width=5mm, fill=white, inner sep=0.25mm] (text4677) {$\omega$
     };
  \path[fill=black] (101.5,1002.6122) node[circle, draw, line width=0.65pt, minimum width=5mm, fill=white, inner sep=0.25mm] (text4681) {$L\ass$
   };
  \end{tikzpicture}}}\rlap{ ,}
  \end{equation*}
 where here, and throughout what follows, we silently suppress applications of the binary and nullary functoriality constraints for $L$;
\vskip0.5\baselineskip
\item The left and right unit constraints $\sigma_{xy}$ and $\tau_{xy}$ are the respective composites
  \begin{equation*}
\vcenter{\hbox{
\begin{tikzpicture}[y=0.8pt, x=0.8pt,yscale=-1, inner sep=0pt, outer sep=0pt, every text node part/.style={font=\tiny} ]
  \path[draw=black,line join=miter,line cap=butt,line width=0.650pt] 
  (30.0000,992.3622) .. controls (79.0000,992.3622) and (65.2500,963.1122) .. 
  node[above right=0.12cm,at start] {\!$\laxc$}(147.4643,981.8265)
  (30.0000,972.3622) .. controls (43.6066,972.3622) and (57.0476,975.5313) .. 
  node[above right=0.12cm,at start] {\!\!$(Lj) \t 1$}(69.4222,979.8561)
  (76.8432,982.6323) .. controls (93.1593,989.1386) and (106.3141,997.0642) .. 
  node[below left=0.02cm,pos=0.38] {$L(j \t 1)$}(108.6747,1001.1108);
  \path[draw=black,line join=miter,line cap=butt,line width=0.650pt] (30.0000,1012.3622) .. controls (60.0000,1012.3622) and (99.3909,1013.5024) .. 
  node[above right=0.12cm,at start] {\!$Lm$}(108.2436,1004.5306);
  \path[draw=black,line join=miter,line cap=butt,line width=0.650pt] (30.0000,952.3622) .. controls (80.0000,952.3622) and (134.2500,958.1122) .. 
  node[above right=0.12cm,at start] {\!$\laxu \t 1$}(148.2500,979.0050);
  \path[draw=black,line join=miter,line cap=butt,line width=0.650pt] (111.1493,1002.0969) .. controls (122.8339,1002.5411) and (140.4643,995.1122) .. 
  node[below right=0.09cm,pos=0.5] {$L \mathfrak l$}(148.5000,984.5765);
  \path[draw=black,line join=miter,line cap=butt,line width=0.650pt] (150.0000,982.3622) .. controls (150.0000,982.3622) and (170.0000,982.3622) .. 
  node[above left=0.12cm,at end] {$\mathfrak l$\!}(180.0000,982.3622);
  \path[fill=black] (149,981.86218) node[circle, draw, line width=0.65pt, minimum width=5mm, fill=white, inner sep=0.25mm] (text4677) {$\gamma$
     };
  \path[fill=black] (109.5,1002.6122) node[circle, draw, line width=0.65pt, minimum width=5mm, fill=white, inner sep=0.25mm] (text4681) {$L\lu$
     };
  \end{tikzpicture}}}\qquad \text{and} \qquad
\vcenter{\hbox{
\begin{tikzpicture}[y=0.8pt, x=0.8pt,yscale=-1, inner sep=0pt, outer sep=0pt, every text node part/.style={font=\tiny} ]
  \path[draw=black,line join=miter,line cap=butt,line width=0.650pt] 
  (30.0000,972.3622) .. controls (44.5544,972.3622) and (61.0063,976.0700) .. 
  node[above right=0.07cm,pos=0.2] {$1 \t (Lj)$}(76.2189,980.9376)
  (86.8712,984.6017) .. controls (103.3429,990.6702) and (116.5479,997.4649) .. 
  node[below left=0.05cm,pos=0.41] {$L(1 \t j)$}(118.6747,1001.1108)
  (30.0000,992.3622) .. controls (69.2166,992.3622) and (73.5528,973.7400) .. 
  node[above right=0.12cm,at start] {\!$\laxc$}(147.9643,984.3265);
  \path[draw=black,line join=miter,line cap=butt,line width=0.650pt] (30.0000,1012.3622) .. controls (50.0000,1012.3622) and (109.3909,1013.5024) .. 
  node[above right=0.12cm,at start] {\!$Lm$}(118.2436,1004.5306);
  \path[draw=black,line join=miter,line cap=butt,line width=0.650pt] (30.0000,952.3622) .. controls (120.0000,952.3622) and (120.0000,981.8622) .. 
  node[above right=0.12cm,at start] {\!$1 \t \laxu$}(146.7500,980.7550);
  \path[draw=black,line join=miter,line cap=butt,line width=0.650pt] (121.1493,1002.0969) .. controls (131.5735,1002.4932) and (144.8262,1000.9348) .. 
  node[below right=0.06cm,pos=0.66] {$L\mathfrak r$}(157.5389,996.4729) .. controls (168.6280,992.5808) and (166.6562,984.8590) .. 
  (152.8505,981.2476);
  \path[draw=black,line join=miter,line cap=butt,line width=0.650pt] (146.3761,977.5008) .. controls (133.5541,974.3625) and (128.7974,963.4580) .. (141.0638,963.1074) .. controls (153.3301,962.7569) and (159.5670,981.3899) .. 
  node[above left=0.12cm,at end] {$\mathfrak r$}(190.1993,981.3899);
  \path[fill=black] (149,981.86218) node[circle, draw, line width=0.65pt, minimum width=5mm, fill=white, inner sep=0.25mm] (text4677) {$\delta^{-1}$
     };
  \path[fill=black] (119.5,1002.6122) node[circle, draw, line width=0.65pt, minimum width=5mm, fill=white, inner sep=0.25mm] (text4681) {$L\ru$
     };
  \end{tikzpicture}}}\rlap{ .}
  \end{equation*}
 \end{itemize}

Given a $\V$-functor $F\colon \B\to\C$, the $\W$-functor $L F \colon L \B \to L \C$ has the same action on objects as $F$,  
action on homs $(LF)_{xy} = L(F_{xy}) \colon L \B(x,y) \to L \C(Fx,Fy)$, and coherence $2$-cells $\fu_{x}$ and $\fm_{xyz}$ given by the respective composites
  \begin{equation*}
\vcenter{\hbox{
\begin{tikzpicture}[y=0.8pt, x=0.85pt,yscale=-1, inner sep=0pt, outer sep=0pt, every text node part/.style={font=\tiny} ]
  \path[draw=black,line join=miter,line cap=butt,line width=0.650pt] (126.0000,972.3622) .. controls (110.0000,972.3622) and (99.9765,972.3169) .. 
  node[above left=0.12cm,at start] {$Lj$\!}(92.5625,979.7309);
  \path[draw=black,line join=miter,line cap=butt,line width=0.650pt] (126.0000,992.3622) .. controls (110.0000,992.3622) and (99.9298,992.1317) .. 
  node[above left=0.12cm,at start] {$LF$\!}(92.6156,984.8176);
  \path[draw=black,line join=miter,line cap=butt,line width=0.650pt] (90.0000,982.3622) -- 
  node[above right=0.12cm,at end] {\!$Lj$}(60.0000,982.3622);
  \path[draw=black,line join=miter,line cap=butt,line width=0.650pt] (60.0000,962.3622) .. controls (100.0000,962.3622) and (78.7500,952.3622) .. 
  node[above left=0.12cm,at end] {$\laxu$\!} node[above right=0.12cm,at start] {\!$\laxu$}(126.0000,952.3622);
  \path[fill=black] (89.651039,982.15656) node[circle, draw, line width=0.65pt, minimum width=5mm, fill=white, inner sep=0.25mm] (text6017) {$L\iota$
     };
  \end{tikzpicture}}}\qquad \text{and} \qquad
\vcenter{\hbox{
\begin{tikzpicture}[y=0.8pt, x=0.85pt,yscale=-1, inner sep=0pt, outer sep=0pt, every text node part/.style={font=\tiny} ]
  \path[draw=black,line join=miter,line cap=butt,line width=0.650pt] 
  (60.0000,912.3622) .. controls (72.1177,912.3622) and (85.4282,915.1636) .. 
   node[above right=0.12cm,at start] {\!$(LF)\t(LF)$}(97.5234,919.0912)
  (105.3673,921.8531) .. controls (121.3611,927.9286) and (133.9336,935.4795) .. 
  node[below left=0.045cm,pos=0.45] {$L(F\t F)$}(136.5000,939.9247)
  (60.0000,932.3622) .. controls (89.8608,932.3622) and (89.2443,912.3622) .. 
  node[above left=0.12cm,at end] {$\laxc$\!} node[above right=0.12cm,at start] {\!$\laxc$}(175.0000,912.3622);
  \path[draw=black,line join=miter,line cap=butt,line width=0.650pt] (139.5000,940.8622) .. controls (146.2500,932.8622) and (160.2372,932.3622) .. 
  node[above left=0.12cm,at end] {$Lm$\!}(175.0000,932.3622);
  \path[draw=black,line join=miter,line cap=butt,line width=0.650pt] (175.0000,952.3622) .. controls (159.9969,952.3622) and (148.0321,952.1443) .. 
   node[above left=0.12cm,at start] {$LF$\!}(139.2500,943.3622);
  \path[draw=black,line join=miter,line cap=butt,line width=0.650pt] (60.0000,952.3622) .. controls (90.0000,952.3622) and (124.3655,950.7623) .. 
  node[above right=0.12cm,at start] {\!$Lm$}(136.7500,943.6122);
  \path[fill=black] (137.90105,941.15656) node[circle, draw, line width=0.65pt, minimum width=5mm, fill=white, inner sep=0.25mm] (text6017) {$L\fm$
     };
  \end{tikzpicture}}}\ \rlap{ .}
  \end{equation*}

  Next, for a $\V$-transformation $\gamma \colon F \Rightarrow G$, the $\W$-transformation $L\gamma$ has $1$-cell components
\[(L\gamma)_x = I \xrightarrow{\laxu} L I \xrightarrow{L(\gamma_x)} L(\C(Fx,Gx))\rlap{ ,}\]
and $2$-cell components given by
\begin{equation*}
\begin{tikzpicture}[y=0.8pt, x=0.9pt,yscale=-1, inner sep=0pt, outer sep=0pt, every text node part/.style={font=\tiny} ]
  \path[draw=black,line join=miter,line cap=butt,line width=0.650pt] 
  (42.0000,932.3622) .. controls (90.6380,932.3622) and (138.6005,949.9627) .. 
  node[above right=0.12cm,at start] {\!$LG$}(146.7322,976.9835);
  \path[draw=black,line join=miter,line cap=butt,line width=0.650pt] 
  (42.0000,952.3622) .. controls (45.0028,952.3622) and (86.0000,958.3622) .. 
  node[above right=0.12cm,at start] {\!$\mathfrak r^\centerdot$}(96.0000,967.3622);
  \path[draw=black,line join=miter,line cap=butt,line width=0.650pt] 
  (42.0000,972.3622) -- 
  node[above right=0.12cm,at start] {\!$1 \laxu$}(94.0000,972.3622);
  \path[draw=black,line join=miter,line cap=butt,line width=0.650pt] 
  (42.0000,1032.3622) .. controls (97.0088,1032.3622) and (136.0399,1023.3285) .. 
  node[above right=0.12cm,at start] {\!$Lm$}(147.0000,991.8622);
  \path[draw=black,line join=miter,line cap=butt,line width=0.650pt] 
  (42.0000,992.3622) .. controls (52.4852,992.3622) and (66.2025,992.2505) .. 
  node[above right=0.12cm,at start] {\!$1\t(L\gamma)$}(80.4234,991.8885)
  (87.9501,991.6738) .. controls (112.0986,990.9073) and (136.3162,989.3589) .. 
  node[below=0.07cm,pos=0.4] {$L(1 \t \gamma)$}(147.5000,986.3622)
  (42.0000,1012.3622) .. controls (70.0000,1012.3622) and (88.0000,988.8622) .. 
  node[above right=0.12cm,at start] {\!$\laxc$}(95.5000,976.3622);
  \path[draw=black,line join=miter,line cap=butt,line width=0.650pt] 
  (102.0000,972.3622) .. controls (122.0000,972.3622) and (139.9645,977.1300) .. 
  node[above=0.12cm,pos=0.45] {$L \mathfrak r^\centerdot$}(146.4645,982.1300);
  \path[draw=black,line join=miter,line cap=butt,line width=0.650pt] 
  (280.0000,932.3622) .. controls (229.3620,932.3622) and (161.3995,954.9627) .. 
  node[above left=0.12cm,at start] {$LF$\!}(153.2678,976.9835);
  \path[draw=black,line join=miter,line cap=butt,line width=0.650pt] 
  (280.0000,952.3622) .. controls (246.7452,952.3622) and (238.3036,950.6657) .. 
  node[above left=0.12cm,at start] {$\mathfrak l^\centerdot$\!}(214.7644,954.3685) .. controls (192.7505,957.8314) and (184.2579,966.5484) .. (199.7574,970.8977);
  \path[draw=black,line join=miter,line cap=butt,line width=0.650pt] 
  (280.0000,972.3622) .. controls (245.3518,972.7157) and (218.7868,949.2807) .. 
  node[above left=0.12cm,at start] {$\laxu 1$\!}(203.5251,968.8266);
  \path[draw=black,line join=miter,line cap=butt,line width=0.650pt] 
  (280.0000,1032.3622) .. controls (222.9912,1032.3622) and (166.7886,1022.2678) .. 
  node[above left=0.12cm,at start] {$Lm$\!}(153.0000,991.8622);
  \path[draw=black,line join=miter,line cap=butt,line width=0.650pt] 
  (280.0000,1012.3622) .. controls (206.9167,1012.3622) and (265.0330,972.9523) .. 
  node[above left=0.12cm,at start] {$\laxc$\!}(205.2071,972.1195)
  (280,992.3622) .. controls (272.8182,992.3622) and (257.2795,992.2620) .. 
  node[above left=0.12cm,at start] {$(L\gamma)1$}(240.2431,991.9439)
  (234.5645,991.8306) .. controls (201.4469,991.1276) and (164.5075,989.5796) .. 
  node[below=0.09cm,pos=0.4] {$L(\gamma 1)$}(152.5000,986.3622);
  \path[draw=black,line join=miter,line cap=butt,line width=0.650pt] 
  (205.0711,975.1906) .. controls (222.1819,977.4572) and (214.8962,986.7288) .. (206.8472,986.1992) .. controls (194.0653,985.3581) and (165.9847,977.1094) .. 
  node[above=0.12cm,pos=0.65] {$L \mathfrak l^\centerdot$}(153.5355,982.1299);
  \path[fill=black] (99.102539,972.51746) node[circle, draw, line width=0.65pt, minimum width=5mm, fill=white, inner sep=0.25mm] (text5074-2) {$\delta^{-1}$     };
  \path[fill=black] (200.27373,972.51746) node[circle, draw, line width=0.65pt, minimum width=5mm, fill=white, inner sep=0.25mm] (text5074-2-1) {$\gamma$     };
  \path[fill=black] (149.10254,985.01746) node[circle, draw, line width=0.65pt, minimum width=5mm, fill=white, inner sep=0.25mm] (text5074-8) {$\L\bar \gamma$   };
  \end{tikzpicture}\ \rlap{ .}
\end{equation*}

Finally, for a $\V$-modification $\Gamma \colon \gamma \Rrightarrow \delta$, the $\W$-modification $L\Gamma \colon L\gamma \Rrightarrow L\delta$ has components given by $(L\Gamma)_x = L(\Gamma_x) \circ \laxu$.

Given $\V$-functors $F \colon \B \to \C$ and $G \colon \C \to \D$, the composites $(LG)(LF)$ and $L(GF) \colon L\B \to L\D$ agree on objects, and differ on hom-objects only up to binary functoriality constraints for $L$; these constraints assemble into an invertible $\W$-icon $(LG)(LF) \Rightarrow L(GF)$. Similarly, for every $\V$-bicategory $\B$, we have an invertible $\W$-icon $1_{L\B} \Rightarrow L(1_\B) \colon L\B \to L \B$. Using Proposition~\ref{thm:icons} and arguing as in Section~\ref{sec:tricat-vbicat}, we may derive from these $\W$-icons the coherence data making change of base into a trifunctor $L \colon \V\text-\cat{Bicat} \to \W\text-\cat{Bicat}$.
We may make this more precise, as in Remark~\ref{rk:tricat-vbicat}, by constructing change of base first as a morphism of locally cubical bicategories $\V\text-\underline{\cat{Bicat}} \to \W\text-\underline{\cat{Bicat}}$, and recovering its instantiation as a trifunctor from the ``locally horizontal'' structure.
\begin{Ex}\label{ex:underlying-change-of-base}
For a monoidal bicategory $\V$, the functor {$V = \V(I, \thg) \colon \V \to \cat{Cat}$} becomes lax monoidal in a canonical way; the binary and nullary monoidality constraints have $1$-cell components given by
\begin{gather*}
\V(I, A) \times \V(I,B) \xrightarrow{\otimes} \V(I \otimes I, A \otimes B) \xrightarrow{\V(\mathfrak l^\centerdot, 1)} \V(I, A \otimes B) \quad
\text{and} \quad 1 \xrightarrow{\id_I} \V(I,I)
\end{gather*}
respectively. The change of base trifunctor $\V\text-\cat{Bicat} \to \cat{Cat}\text-\cat{Bicat} = \cat{Bicat}$ is in this case the underlying ordinary bicategory trifunctor $(\thg)_0$ of Section~\ref{subsec:ordinary-trifunctor}.
\end{Ex}

\subsection{Change of base along higher monoidal cells}

Given a monoidal transformation $\alpha \colon L \Rightarrow L' \colon \V \to \W$, we induce a tritransformation $\alpha \colon L \Rightarrow L' \colon \V\text-\cat{Bicat} \to \W\text-\cat{Bicat}$ as follows. Its $1$-cell component at a $\V$-bicategory $\B$ is the identity-on-objects $\W$-functor $\alpha_\B \colon L\B \to L'\B$ with action on homs $(\alpha_\B)_{xy} = \alpha_{\B(x,y)} \colon L(\B(x,y)) \rightarrow L'(\B(x,y))$, and with functoriality constraint cells built from pseudonaturality of $\alpha$ together with $M$ and $\Pi$. Now for any $\V$-functor $F \colon \B \to \C$, there is an invertible $\W$-icon $\alpha_\C \circ LF \Rightarrow L'F \circ \alpha_\B$ whose components are built from pseudonaturality $2$-cells for $\alpha$; and, arguing as before, we may construct the remaining data of the tritransformation $\alpha$ from these invertible icons, either directly using Proposition~\ref{thm:icons}, or by deriving the structure from a transformation of locally cubical bicategories.

Similarly, given a monoidal modification $\Gamma \colon \alpha \Rightarrow \beta$, we have for each $\V$-bicategory $\B$ a $\W$-icon $\Gamma_\B$ whose $2$-cell components are given by $(\Gamma_\B)_{xy} = \Gamma_{\B(x,y)}$. The $\W$-transformations corresponding to these $\W$-icons constitute the components of a trimodification, whose 
remaining coherence data may be obtained, as before, either by direct construction, or more precisely via locally cubical bicategories.

\subsection{Functoriality of change of base}
The operations described in the preceding two sections are, in a suitable sense, functorial; the precise nature of this functoriality is a little delicate. They ought to comprise a morphism of tetracategories (weak $4$-categories) $(\thg)\text-\cat{Bicat} \colon \cat{MonBicat} \to \cat{Tricat}$, but the sheer quantity of coherence that would be involved in making this precise leads us to consider a more refined approach. We begin from the assignation $\V \mapsto \V\text-\underline{\cat{Bicat}}$, and view the change of base operations as landing in locally cubical bicategories; whereupon we obtain a trifunctor $(\thg)\text-\underline{\cat{Bicat}} \colon \cat{MonBicat} \to \cat{DblCat}\text-\cat{Bicat}$. We leave a detailed description of the coherence constraints of this trifunctor to the reader.

\subsection{Monoidal adjunctions}
We will be particularly concerned with change of base in the situation of an adjunction of the following sort.

\begin{Lemma}[Doctrinal adjunction]\label{lem:docadjn}
  Suppose $L\colon \V\to\W$ is a strong monoidal functor that has a right adjoint $R\colon \W\to\V$.
  Then $R$ is canonically a lax monoidal functor.
\end{Lemma}
\begin{proof}
  Since $L$ is strong, $\laxu$ and $\laxc$ have inverse adjoint equivalences $\laxu^\centerdot$ and $\laxc^\centerdot$.
  By the usual mates correspondence, the constraints $\omega$, $\gamma$, and $\delta$ induce analogous constraints for $\laxu^\centerdot$ and $\laxc^\centerdot$, making $L$ into an ``oplax monoidal functor'' in an obvious sense.

  Now we can transfer this structure across the biadjunction $L\dashv R$ using a categorified mates correspondence.
  That is, if $\eta\colon 1_{\V} \to R L$ and $\epsilon \colon  L R \to 1_{\W}$ are the pseudonatural unit and counit of the biadjunction, then $\laxc$ and $\laxu$ for $R$ are the respective composites
\begin{gather*}
R x \otimes R y \xrightarrow{\eta_{R\otimes R}}
  R L (R x \otimes R y) \xrightarrow{R \laxc^\centerdot_{R,R}}
  R (L R x \otimes L R y) \xrightarrow{R(\epsilon \otimes \epsilon)}
  R (x \otimes y) \\
\text{and} \quad \ \ 
I \xrightarrow{\eta_I} R L I \xrightarrow{R \laxu^\centerdot} R I\rlap{ .}
\end{gather*}
  The categorified mates correspondence is stated precisely in~\cite[\S3]{lauda:faaa}, as an equivalence of categories of 2-cells in a Gray-category, but we may apply it in any tricategory by the coherence theorem for tricategories.
  Thus we can obtain the constraint isomorphisms $\omega$, $\gamma$, and $\delta$ for $R$, and their axioms, by the functoriality of the mates correspondence on 3-cells.
\end{proof}

We call this situation a \emph{monoidal biadjunction}; by a categorification of the argument of~\cite{kelly:doc-adjn}, the unit and counit of such a biadjunction are monoidal transformations, and the modifications witnessing the coherent satisfaction of the triangle identities are monoidal modifications. Thus the entire biadjunction $L \dashv R$ lifts to the tricategory $\cat{MonBicat}$; applying the trifunctor $(\thg)\text-\underline{\cat{Bicat}}$, we obtain a biadjunction $\V\text-\underline{\cat{Bicat}} \rightleftarrows \W\text-\underline{\cat{Bicat}}$ which, since both domain and codomain are locally fibrant, induces in turn a triadjunction between underlying tricategories $\V\text-\cat{Bicat} \rightleftarrows \W\text-\cat{Bicat}$. We shall write $\eta_\B \colon \B \to RL\B$ and $\varepsilon_\C \colon LR\C \to \C$ for the identity-on-objects $\V$- and $\W$-functors forming the unit and counit $1$-cells of this triadjunction.

\begin{Ex}\label{ex:underlying-monoidal-biadj}
If $\V$ is cocomplete and left and right closed, then the lax monoidal functor $V = \V(I,A) \colon \V \to \cat{Cat}$ is the right adjoint of a monoidal biadjunction; the 
strong monoidal left adjoint $L$ is defined by taking $L C$ to be the copower of $I\in\V$ by $C\in\cat{Cat}$. It is evident that $L1 \simeq I$, while the binary monoidality constraints
$LA \otimes LB \simeq L(A \times B)$ are obtained via the chain of equivalences
\[
LA \otimes LB = (I \cdot A) \otimes (I \cdot B) \simeq (I \otimes (I \cdot B)) \cdot A \simeq ((I \otimes I) \cdot A) \cdot B \simeq I \cdot (A \times B) = L(A\times B)
\] 
using the fact that each functor $X \otimes (\thg)$ and $(\thg) \otimes Y$ preserves colimits, in particular copowers. \end{Ex}

\subsection{Monoidal adjunctions and weighted colimits}
For the rest of this section, we investigate the effect of a monoidal biadjunction $\V \rightleftarrows \W$ on weighted colimits in $\V$- and $\W$-bicategories.
Observe first that, given a lax monoidal functor $L \colon \V \to \W$, the constructions of Section~\ref{subsec:change-of-base-trihom} carry over \emph{mutatis mutandis} to bimodules; so we induce from any $\V$-bimodule $M \colon \A \tor \B$ a $
\W$-bimodule $LM \colon L\A \tor L\B$ with components $(LM)(b, a) = L(M(b,a))$, and so on, and correspondingly for bimodule morphisms and transformations. In this way, we obtain a functor $L \colon {}_\A \cat{Mod}_\B \to {}_{L\A} \cat{Mod}_{L \B}$; and the same is evidently true for left and right modules, which we can include in this notation with the convention $L \bullet = \bullet$. 

The functors just described are also ``lax'' with respect to copowers of modules: given $A\in\V$ and a right $\B$-module $W$, we have a morphism of right $L\B$-modules $\laxc \colon L A \otimes L W \to L(A\otimes W)$, with $1$-cell components $\laxc_{A,Wx} \colon LA \otimes L(Wx) \to L(A \otimes Wx)$ and $2$-cell components obtained from $\omega$ together with functoriality constraints for $L$.
There are now invertible $L\B$-module modifications $\omega$ and $\gamma$ whose components are those of the corresponding coherence constraints for the lax monoidal functor $L$; it follows that these satisfy axioms corresponding to those for a lax monoidal functor.
(In fact, this is part of the structure of a lax functor of tricategories $\V\text-\cat{Mod} \to \W\text-\cat{Mod}$.)

Suppose now that we have a monoidal biadjunction $L \dashv R \colon \W \to \V$. For any $\V$-bicategory $\B$ and any right $L \B$-module $W$, we can restrict the $R L \B$-module $R W$ along $\eta_\B \colon \B \to R L \B$ to obtain the right-$\B$-module  $\hat{R} W \defeq RW(\eta_\B)$.

\begin{Lemma}\label{lem:rhat}
  $\hat{R} \colon {}_\bullet \cat{Mod}_{L \B} \to {}_\bullet \cat{Mod}_\B$ is right adjoint to $L \colon {}_\bullet \cat{Mod}_\B \to {}_\bullet \cat{Mod}_{L \B}$.
\end{Lemma}
\begin{proof}
  Let $V$ be a $\B$-module and $W$ an $L\B$-module.
  A morphism $LV \to W$ consists of morphisms $\phi_x\colon L(Vx) \to Wx$ for each $x\in\B$, together with invertible 2-cells
  \begin{equation*}
    \cd[@C+2.5em]{
      L(Vy) \otimes L\B(x,y) \ar[r]^-{\laxc}
      \ar[d]_{\phi \otimes 1} \rtwocell{drr}{\bar \phi_{xy}} &
      L(Vy \otimes \B(x,y))
      \ar[r]^-{Lm}  &
      L(Vx) \ar[d]^{\phi} \\
      W y \otimes L\B(x,y) \ar[rr]_-{m} && W x
    }
  \end{equation*}
  satisfying appropriate axioms.
  On the other hand, a morphism $V \to \hat{R}W$ consists of morphisms $\psi_x\colon Vx \to R(Wx)$ for each $x\in \B$, together with invertible 2-cells
  \begin{equation*}
    \cd[@C+0em]{
      Vy \otimes \B(x,y) \ar[rrr]^-{m} \ar[d]_{\psi \otimes 1} \rtwocell{drrr}{\bar \psi_{xy}} &&&
      Vx \ar[d]^{\psi} \\
      R(W y) \otimes \B(x,y) \ar[r]_-{1\otimes \eta} &
      R(Wy) \otimes RL\B(x,y) \ar[r]_-{\laxc} &
      R(Wy) \otimes L\B(x,y)) \ar[r]_-{Rm} & R(W x)
    }
  \end{equation*}
  satisfying appropriate axioms.
  Of course the biadjunction $L\dashv R$ gives us equivalences of categories $\W(LVx,Wx) \simeq \V(Vx,RWx)$, so it remains to show that these lift compatibly to the additional 2-cell data.
  Thus, given on the one hand $\phi$ and $\bar\phi$ as above, we take $\psi_x = R\phi_x\circ \eta_{Vx}$ and $\bar\psi_{x,y}$ to be the composite
\begin{equation*}
\begin{tikzpicture}[y=0.8pt, x=1pt,yscale=-1, inner sep=0pt, outer sep=0pt, every text node part/.style={font=\tiny} ]
  \path[draw=black,line join=miter,line cap=butt,line width=0.650pt] 
  (40.0000,952.3622) .. controls (82.7279,952.3622) and (107.9069,965.6269) .. 
  node[above right=0.12cm,at start] {\!$\eta \t 1$}(107.9069,965.6269);
  \path[draw=black,line join=miter,line cap=butt,line width=0.650pt] 
  (40.0000,972.3622) .. controls (60.9164,972.3622) and (77.6588,984.4538) .. 
  node[above right=0.12cm,at start] {\!\!$(R \phi) \t 1$}(90.8895,996.5881)
  (94.8667,1000.3257) .. controls (104.8607,1009.9212) and (112.6587,1018.7526) .. 
  node[left=0.19cm,pos=0.55] {$R( \phi \t 1)$}(118.5979,1020.6828)
  (40.0000,992.3622) .. controls (47.1114,992.3622) and (56.7036,987.8205) .. 
  node[above right=0.12cm,at start] {\!$1 \t \eta$}(67.6039,982.3563)
  (72.9581,979.6642) .. controls (96.8589,967.6670) and (125.7766,954.1429) .. 
  node[below=0.09cm,pos=0.25] {$1 \t \eta$} node[above=0.09cm,pos=0.7] {$\eta \eta$}(148.5858,973.4228)
  (40.0000,1012.3622) .. controls (80.3068,1012.1305) and (131.4142,975.7452) .. 
  node[above right=0.12cm,at start] {\!$\laxc$}(148.5616,975.9220);
  \path[draw=black,line join=miter,line cap=butt,line width=0.650pt] 
  (40.0000,1032.3622) .. controls (69.7032,1032.3622) and (113.0936,1027.7660) .. 
  node[above right=0.12cm,at start] {\!$Rm$}(118.9272,1023.4350);
  \path[draw=black,line join=miter,line cap=butt,line width=0.650pt] 
  (121.4456,1019.6585) .. controls (137.0019,1004.6325) and (135.1508,987.3119) .. 
  node[above left=0.07cm,pos=0.25] {$R \laxc$}(148.4090,977.9428);
  \path[draw=black,line join=miter,line cap=butt,line width=0.650pt] 
  (121.6794,1021.6551) .. controls (139.9715,1018.1289) and (162.5983,999.8455) .. 
  node[below right=0.04cm,pos=0.67] {$RLm$\!}(183.9450,986.5727)
  (189.4718,983.2423) .. controls (200.3804,976.9017) and (210.8158,972.3622) .. 
  node[above left=0.12cm,at end] {$m$\!}(220.0000,972.3622)
  (151.8512,975.6411) .. controls (169.8824,974.9340) and (195.4274,992.3622) .. 
  node[above left=0.12cm,at end] {$\eta$\!}(220.0000,992.3622);
  \path[draw=black,line join=miter,line cap=butt,line width=0.650pt] 
  (121.2688,1024.0779) .. controls (159.4223,1034.6134) and (187.2844,1012.3622) .. 
  node[above left=0.12cm,at end] {$R\phi$\!}(220.0000,1012.3622);
  \path[fill=black] (149.55309,975.64111) node[circle, draw, line width=0.65pt, minimum width=5mm, fill=white, inner sep=0.25mm] (text3214) {$\Pi^{-1}$     };
  \path[fill=black] (119.8546,1022.6637) node[circle, draw, line width=0.65pt, minimum width=5mm, fill=white, inner sep=0.25mm] (text3218) {$R\bar \phi$     };
  \path[shift={(30.703863,34.624305)},draw=black,fill=black] (77.8600,931.1000)arc(0.000:180.000:1.250)arc(-180.000:0.000:1.250) -- cycle;
  \end{tikzpicture}\ \rlap{ ,}
\end{equation*}
wherein $\Pi$ is part of the $2$-cell data exhibiting $\eta \colon 1 \Rightarrow RL$ as a monoidal transformation.
  On the other hand, given $\psi$ and $\bar\psi$ as above, we take $\phi_x = \epsilon_{Wx} \circ L\psi_x$ and $\bar\phi_{x,y}$ to be the composite
\begin{equation*}
\begin{tikzpicture}[y=0.9pt, x=1pt,yscale=-1, inner sep=0pt, outer sep=0pt, every text node part/.style={font=\tiny} ]
  \path[draw=black,line join=miter,line cap=butt,line width=0.650pt] 
  (20.0000,992.3622) .. controls (51.0091,992.3622) and (68.2500,973.8622) .. 
  node[above right=0.12cm,at start] {\!$\varepsilon 1$}node[below right=0.07cm,pos=0.45] {$\varepsilon \varepsilon$} (68.2500,973.8622);
  \path[draw=black,line join=miter,line cap=butt,line width=0.650pt] (140.0763,970.9480) .. controls (151.0364,961.0485) and (161.7913,957.3622) .. 
  node[above left=0.12cm,at end] {$Lm$\!}(180.0000,957.3622);
  \path[draw=black,line join=miter,line cap=butt,line width=0.650pt] (139.6464,973.7001) .. controls (149.1989,983.2526) and (162.2910,987.3622) .. 
  node[above left=0.12cm,at end] {$L\psi$\!}(180.0000,987.3622);
  \path[draw=black,line join=miter,line cap=butt,line width=0.650pt] 
  (71.2936,974.3457) .. controls (77.5029,1005.8466) and (132.2441,1012.3622) .. 
  node[above left=0.12cm,at end] {$\varepsilon$\!}(180.0000,1012.3622)
  (20.0000,1012.3622) .. controls (50.7709,1012.3622) and (74.5419,1007.6894) .. node[above right=0.12cm,at start] {\!$m$}(92.4909,1001.4817)
  (99.2219,998.9804) .. controls (120.8538,990.3605) and (132.6832,979.7941) .. 
  node[below right=0.06cm,pos=0.37] {$LRm$}(137.3884,974.4167);
  \path[draw=black,line join=miter,line cap=butt,line width=0.650pt] (29.3878,963.0167) .. controls (40.2696,966.7134) and (45.0685,979.7551) .. 
  node[above right=0.04cm,pos=0.5] {$1 \varepsilon$}(46.5685,987.0051);
  \path[draw=black,line join=miter,line cap=butt,line width=0.650pt] 
  (31.9234,958.9812) .. controls (41.5638,950.4454) and (64.7797,948.1121) .. 
  node[above left=0.05cm,pos=0.6] {$1(L \eta)$}(87.2831,950.8926)
  (95.2485,952.0981) .. controls (112.9732,955.2855) and (129.1524,961.6697) .. 
  node[below=0.14cm,pos=0.28] {$L(1 \eta)$}(136.4474,970.6913)
  (71.0865,969.9822) .. controls (88.2816,941.2373) and (144.9964,932.3622) .. 
  node[above left=0.12cm,at end] {$\laxc$\!}(180.0000,932.3622)
  (20.0000,932.3622) .. controls (59.2825,932.3622) and (89.1218,933.9996) .. node[above right=0.12cm,at start] {\!$(L\psi)1$} (109.0963,940.7915)
  (116.3626,943.6988) .. controls (128.2758,949.3003) and (135.5448,957.5284) .. 
  node[above right=0.02cm,pos=0.5] {$L(\psi 1)$}(138.0348,969.5095);
  \path[draw=black,line join=miter,line cap=butt,line width=0.650pt] (73.5391,972.4591) .. controls (91.7623,976.3187) and (111.8280,977.4145) .. 
  node[below=0.07cm,pos=0.5] {$L\laxc$}(135.6922,972.8127);
  \path[fill=black] (139.33125,971.75201) node[circle, draw, line width=0.65pt, minimum width=5mm, fill=white, inner sep=0.25mm] (text4788) {$L\bar \psi$
     };
  \path[fill=black] (70.003571,972.10553) node[circle, draw, line width=0.65pt, minimum width=5mm, fill=white, inner sep=0.25mm] (text4792) {$\Pi^{-1}$
     };
  \path[fill=black] (26.343145,962.36218) node[circle, draw, line width=0.65pt, minimum width=5mm, fill=white, inner sep=0.25mm] (text4796) {$1 \Delta$
   };
     \path[shift={(-29.940952,56.208382)},draw=black,fill=black] (77.8600,931.1000)arc(0.000:180.000:1.250)arc(-180.000:0.000:1.250) -- cycle;
  \end{tikzpicture}\ \rlap{ ,}
\end{equation*}
where now $\Pi$ is part of the $2$-cell data exhibiting $\varepsilon \colon LR \Rightarrow 1$ as a monoidal transformation, and $\Delta$ is a triangle identity $2$-cell for the biadjunction $L \dashv R$.
\end{proof}

Since $L$ is strong monoidal, each of the comparison morphisms $\laxc \colon LA \otimes LW \to L(A \otimes W)$ described above is an equivalence of modules; combining these equivalences with the adjunction just described, we obtain equivalences of categories
\begin{equation}
  {}_\bullet \cat{Mod}_\B(A \otimes V, \hat{R} W) \simeq
  {}_\bullet \cat{Mod}_{L \B} (L A \otimes L V, W)\rlap{ ,}\label{eq:indrhat}
\end{equation}
pseudonatural in $A$. We write the left-to-right version of such an equivalence as $(-)^\sharp$ and the right-to-left version as $(-)^\flat$.

\begin{Lemma}\label{lem:musicaluniv}
  Suppose that the morphism $\phi\colon B \otimes L V \to W$ exhibits $B$ as $\h {L V} W$; then $(\phi \circ (\varepsilon_B \otimes 1))^\flat \colon R B \otimes V \to \hat{R} W$ exhibits $R B$ as $\h V {\hat{R} W}$.
\end{Lemma}
\begin{proof}
  Using the adjunction $L\dashv R$, the assumption, and~\eqref{eq:indrhat}, we have equivalences
  \[ \V(A,R B)
  \simeq \W(L A, B)
  \simeq {} \cat{Mod}_{L\B}(L A \otimes L V, W)
  \simeq {} \cat{Mod}_\B(A\otimes V,\hat{R} W)
  \]
  which are pseudonatural in $A\in\V$.
  Thus, by the Yoneda lemma, the image of $1_{R B}\in\V(R B, R B)$ under these equivalences is a universal morphism. It is easy to see that this image is isomorphic to $(\phi \circ (\varepsilon_B \otimes 1))^\flat$, which is thus is itself universal, as required.
\end{proof}

Now suppose that $\V$ and $\W$ are right closed.
The equivalences $L(A\otimes B) \simeq L A \otimes L B$ coming from the strong monoidal structure of $L$ can be regarded as a pseudonatural equivalence between the composite functors
\begin{gather*}
  \V \xrightarrow{\thg \otimes B} \V \xrightarrow{L} \W \qquad \quad \text{and} \qquad \quad
  \V \xrightarrow{L} \W \xrightarrow{\thg \otimes LB} \W\rlap{ .}
\end{gather*}
Therefore, the composite right biadjoints
\begin{gather*}
  \W \xrightarrow{R} \V \xrightarrow{[B,\thg]} \V \qquad \quad \text{and} \qquad \quad  \W \xrightarrow{[L B, \thg]} \W \xrightarrow{R} \V\rlap{ .}
\end{gather*}
are also equivalent; i.e.\ we have pseudonatural equivalences
\begin{equation}
  [B,R C] \simeq R[L B, C].\label{eq:strong-closed-adjointness}
\end{equation}
This yields a converse to Lemma~\ref{lem:musicaluniv}.

\begin{Lemma}\label{lem:muconv}
  Suppose $\V$ and $\W$ are complete and right closed, and $R$ reflects bilimits.
  Then if $\psi\colon A\otimes V \to \hat{R} W$ exhibits $A$ as $\h V {\hat{R} W}$, and $\phi\colon B \otimes L V \to W$ is an $L \B$-module morphism with $R B = A$ and $(\phi \circ (\varepsilon_B \otimes 1)^\flat\cong \psi$, then $\phi$ exhibits $B$ as $\h {L V} W$.
\end{Lemma}
\begin{proof}
  Since $\B$ and $L\B$ have the same objects, the category $\D$ from Section~\ref{subsec:construction-of-tensor-products} is the same whether we define it for $\B$ or for $L\B$.
  Let $F^{VW} \colon \D^\op \to \V$ be the functor defined as in Section~\ref{subsec:constr-right-homs} for $V$ and $\hat{R} W$; so we have:
  \begin{align*}
  F^{VW}(x) &= [Vx, RWx];\\
  F^{VW}(x,y) &= [Vy \otimes \B(x,y), RWx];\\
  F^{VW}(x,y,z) &= [(Vz \otimes \B(y,z)) \otimes \B(x,y), RWx];
  \end{align*}
  and so on.
 Let $G^{VW} \colon \D^\op \to \W$ be the analogous functor defined for $L V$ and $W$; so we have:
  \begin{align*}
    G^{VW}(x) &=[L V(x),Wx]\\
    &= [L(Vx), Wx];\\
    G^{VW}(x,y) &= [L V(y) \otimes L\B(x,y), Wx]\\
    &\simeq [L(Vy \otimes \B(x,y)), Wx];\\
    G^{VW}(x,y,z) &= [(L V(z) \otimes L \B(y,z)) \otimes L \B(x,y),Wx]\\
    &\simeq [L((Vz \otimes \B(y,z)) \otimes \B(x,y)), Wx];
  \end{align*}
  and so on.
  Now using~\eqref{eq:strong-closed-adjointness}, we can construct an equivalence $R \circ G^{VW} \simeq F^{VW}$.
  Thus the claim follows from Proposition~\ref{prop:right-hom-from-bilimit} and the assumption that $R$ reflects bilimits.
\end{proof}

Now let $\B$ be a $\V$-bicategory, $\C$ a $\W$-bicategory, $W$ a right $\B$-module, and $F\colon L\B \to \C$ a $\W$-functor.
 $F$ has the adjunct $\bar{F} \defeq RF \circ \eta_\B \colon \B \to R\C$ under the triadjunction $\V\text-\cat{Bicat} \rightleftarrows \W\text-\cat{Bicat}$;
 now for any object $v \in \C$, the $L \B$-module $\C(F,v)$ satisfies $\hat{R}(\C(F,v)) = (R\C)(\bar{F},v)$, whence
by Lemma~\ref{lem:rhat}, any morphism $\phi\colon L(W)\to \C(F,v)$ has an adjunct $\bar \phi \colon W \to R\C(\bar{F},v)$.

\begin{Thm}\label{thm:adjoints}
  If $\phi\colon L(W)\to \C(F,v)$ is an $L W$-weighted colimit of $F$, then $\bar \phi \colon W \to R\C(\bar{F},v)$ is a $W$-weighted colimit of $\bar{F}$.
  The converse is true if $\V$ and $\W$ are complete and right closed and $R$ reflects limits.
\end{Thm}
\begin{proof}
  By Lemmas~\ref{lem:musicaluniv} and~\ref{lem:muconv}, it suffices to verify that the operation $(\thg \circ (\varepsilon_\B \otimes 1))^\flat$ takes the transformation of~\eqref{eq:inducedphi} for $\phi$ to the analogous transformation for $\bar \phi$.
\end{proof}

\begin{Ex}
Let $\V$ be cocomplete and left and right closed, so that as in Example~\ref{ex:underlying-monoidal-biadj}, we have a monoidal biadjunction $L \dashv V \colon \V \to \cat{Cat}$;
thus for any ordinary bicategory $\B$ and any $\cat{Cat}$-weight $W\colon \B\to \cat{Cat}$, we have a weight $L W \colon L \B \to \V$ such that $L W$-weighted colimits in a $\V$-bicategory $\C$ are, in particular, $W$-weighted colimits in its underlying ordinary bicategory $\C_0$.
If $R$ happens to reflect limits (such as if $\V$ is complete and $\V(I, \thg)$ is conservative), then there is no difference between $L W$-weighted colimits in $\C$ and $W$-weighted colimits in $R \C$.
This is a categorification of~\cite[Section~3.8]{kelly:enriched}.
\end{Ex}

\section{Constructions on monoidal bicategories}
\label{sec:constructions}

In this section we describe two methods of constructing new monoidal bicategories from old that will be useful later on.

\subsection{Comma bicategories}
\label{sec:comma}

If $\V$ and $\W$ are bicategories and $R\colon \V\to\W$ is a functor, we denote by $(\W\downarrow R)$ the \emph{comma bicategory} defined as follows:
\begin{itemize}
\item Its objects are triples $(V,W,f)$ where $V\in \V$, $W\in\W$, and $f\colon W \to R V$.
\item Its morphisms are triples $(p,q,\gamma)$ where $p\in\V(V,V')$, $q\in \W(W,W')$, and $\gamma\colon R p \circ f \cong f' \circ q$.
\item Its 2-cells are pairs $(\alpha,\beta)$ where $\alpha\colon p\Rightarrow p'$ in $\V$, $\beta\colon q\Rightarrow q'$ in $\W$, and the evident cylinder commutes.
\end{itemize}
There are forgetful functors $U_\V\colon (\W\downarrow R) \to \V$ and $U_\W\colon (\W\downarrow R) \to \W$, with a canonical transformation $U_R\colon U_\W \to R U_\V$.
It is straightforward to verify that $U_\V$ and $U_\W$ jointly create bicolimits; if $R$ preserves bilimits, then they jointly create bilimits as well.
(Here, a functor $F \colon \A \to \B$ is said to \emph{create bicolimits} if it preserves and reflects them, and whenever given $D \colon \I \to \A$ and $W$ a right $\I$-module such that $W \star FD$ exists in $\B$, also $W \star D$ exists in $\A$; likewise for creation of bilimits).

\begin{Thm}\label{thm:comma}
  If $\V$ and $\W$ are monoidal bicategories and $R\colon \V\to\W$ is a lax monoidal functor, then $(\W\downarrow R)$ is a monoidal bicategory, $U_\V$ and $U_\W$ are strong monoidal, and $U_R$ is a monoidal transformation.
  Moreover, a $(\W\downarrow R)$-bicategory $\A$ is determined exactly by the data:
  \begin{itemize}
  \item A $\V$-bicategory $\A_\V = U_\V(\A)$;
  \item A $\W$-bicategory $\A_\W = U_\W(\A)$ with the same objects as $\A_\V$; and
  \item An identity-on-objects functor $J_\A = U_R(\A) \colon \A_\W \to R(\A_\V)$.
  \end{itemize}
  A $(\W\downarrow R)$-functor $F\colon \A\to \B$ is determined exactly by the data:
  \begin{itemize}
  \item A $\V$-functor $F_\V = U_\V(F) \colon \A_\V \to \B_\V$;
  \item A $\W$-functor $F_\W = U_\W(F) \colon \A_\W \to \B_\W$; and
  \item An invertible $\W$-icon
    \[\vcenter{\xymatrix{
        \A_\W \ar[rr]^{J_\A}\ar[d]_{F_\W} \dtwocell{drr}{U_R(F)} &&
        R(\A_\V)\ar[d]^{R(F_\V)}\\
        \B_\W \ar[rr]_{J_\B} &&
        R(\B_\V)\rlap{ .}
      }}
    \]
  \end{itemize}
  A $(\W\downarrow R)$-transformation $\gamma \colon F \Rightarrow G$ is determined exactly by the data:
  \begin{itemize}
  \item A $\V$-transformation $\gamma_\V = U_\V(\gamma) \colon  F_\V \to G_\V$;
  \item A $\W$-transformation $\gamma_\W = U_\W(\gamma) \colon  F_\W \to G_\W$; and
  \item An invertible $\W$-modification
    \[\vcenter{\xymatrix{
        R (F_\V) \circ J_\A \ar[rr]^{U_R(F)} \ar[d]_{R (\gamma_\V) \circ J_\A}
        \dtwocell{drr}{U_R(\gamma)} &&
        J_\B \circ F_\W \ar[d]^{J_\B \circ \gamma_\W}\\
        R(G_\V) \circ J_\A \ar[rr]_{U_R (G)} &&
        J_\B \circ G_\W\rlap{ .}
      }}
    \]
  \end{itemize}
  And a $(\W\downarrow R)$-modification $\Gamma \colon \gamma \Rrightarrow \delta$ is determined exactly by the data:
  \begin{itemize}
  \item A $\V$-modification $\Gamma_\V = U_\V (\Gamma) \colon \gamma_\V \Rrightarrow \delta_\V$; and
  \item A $\W$-modification $\Gamma_\W = U_\W( \Gamma) \colon \gamma_\W \Rrightarrow \delta_\W$,
  \end{itemize}
  such that the obvious diagram of $2$-cells commutes.
\end{Thm}
\begin{proof}
  The tensor product of $W_1 \to R V_1$ and $W_2 \to R V_2$ is defined to be the composite
  \[ W_1 \otimes W_2 \to R V_1 \otimes R V_2 \xrightarrow{\laxc} R(V_1 \otimes V_2) \]
  This extends to a functor $(\W\downarrow R) \times (\W\downarrow R)\to (\W\downarrow R)$, whose action on morphisms involves pseudofunctoriality of $\otimes$ and the pseudonaturality constraint of $\laxc$.
  The unit object of $(\W\downarrow R)$ is defined to be $\laxu\colon I_\W \to R(I_\V)$.
  Using interchange isomorphisms, the left-associated triple tensor product of $W_1 \to R V_1$ and $W_2 \to R V_2$ and $W_3 \to R V_3$ is isomorphic to the composite
  \[ (W_1 \otimes W_2) \otimes W_3 \to (R V_1 \otimes R V_2)\otimes RV_3
  \xrightarrow{\laxc\otimes 1} R(V_1 \otimes V_2) \otimes R V_3
  \xrightarrow{\laxc}  R((V_1 \otimes V_2) \otimes V_3)
  \]
  and similarly for other multiple tensor products.
  Now the tensor associativity constraint $\mathfrak a$ for $(\W\downarrow R)$ has 1-cell components taken from $\mathfrak a$ in $\V$ and $\W$, and 2-cell components built out of these together with the modification $\omega$.
  The unitality constraints $\mathfrak l$ and $\mathfrak r$ for $(\W\downarrow R)$ are likewise built out of $\mathfrak l$ and $\mathfrak r$ in $\V$ and $\W$ together with the modifications $\gamma$ and $\delta$.
  The pseudonaturality constraint 2-isomorphisms for $\mathfrak a$, $\mathfrak l$, and $\mathfrak r$ have their components taken from those of $\mathfrak a$, $\mathfrak l$, and $\mathfrak r$ in $\V$ and $\W$, and their axioms follow from those for $\V$ and $\W$ together with the modification axiom for $\omega$, $\gamma$, and $\delta$.

  The modifications $\pi$, $\nu$, $\lambda$, and $\rho$ for $(\W\downarrow R)$ have their (2-cell) components taken from $\pi$, $\nu$, $\lambda$, and $\rho$ in $\V$ and $\W$.
  The modification axioms for $\pi$ and $\nu$ are precisely the axioms demanded of a lax functor (see~\cite[pp.~17--18]{Gordon1995Coherence}), while those for $\lambda$ and $\rho$ can be verified by hand (they should also follow from a coherence theorem for lax functors).
  The tricategory axioms for $\pi$, $\nu$, $\lambda$, and $\rho$ for $(\W\downarrow R)$ follow immediately from their counterparts in $\V$ and $\W$.

  This completes the proof that $(\W\downarrow R)$ is a monoidal bicategory, and it is obvious on inspection that the functors $U_\V$ and $U_\W$ are in fact \emph{strict} monoidal. The modifications $M$ and $\Pi$ making $U_R$ into a monoidal transformation are simply coherence isomorphisms in $\W$, and their axioms follow from the coherence theorem for tricategories.

  Now, the hom-objects of a $(\W\downarrow R)$-bicategory $\A$ are precisely determined by the hom-objects of $\A_\V$ and $\A_\W$ and the action of $J_\A$ on homs.
  The composition in $\A$ is similarly precisely determined by the compositions in $\A_\V$ and $\A_\W$ and the functoriality isomorphisms of $J_\A$, and correspondingly for the unit.
  The associativity and unitality constraints of $\A$ are determined by those in $\A_\V$ and $\A_\W$, with the requisite axioms (making them 2-cells in $(\W\downarrow R)$) reducing to the functoriality axioms of $J_\A$.
  And the pentagon and triangle axioms of $\A$ are exactly those of $\A_\V$ and $\A_\W$.
  Functors, transformations, and modifications behave similarly; we leave them to the reader.
\end{proof}

\begin{Rk}
  Note that we used essentially all the data and axioms of a lax functor in proving Theorem~\ref{thm:comma}.
  This is not a coincidence.
  It is shown in~\cite{Bourke2012TwoDMonadicity} that in the 2-dimensional case, the correctness of a notion of ``lax morphism'' (in the sense that it coincides with the general definition for algebras over a 2-monad) is essentially determined by doctrinal adjunction (Lemma~\ref{lem:docadjn}) and the existence of colax limits of arrows (comma objects under the identity, as in Theorem~\ref{thm:comma}).
  Thus, one may expect these two theorems to eventually imply that lax monoidal functors of monoidal bicategories coincide with a general 3-monadic definition of lax morphism.
\end{Rk}

\subsection{Reflective sub-bicategories}
\label{sec:reflective}

Suppose $\V$ is a monoidal bicategory and $\W\subseteq\V$ a full, replete, reflective sub-bicategory with inclusion $R \colon \W \to \V$ and reflector $L \colon \V \to \W$. Write $^{\perp}\W$ for the class of morphisms of $\V$ that are sent to equivalences by $L$.

\begin{Thm}\label{thm:reflect}
  If each functor $X \otimes (\thg)$ and $(\thg) \otimes Y$ preserves $^{\perp}\W$, then:
  \begin{itemize}
  \item $\W$ inherits a monoidal structure.
  \item The biadjunction $L\dashv R$ is monoidal.
  \item For a $\W$-bicategory $\B$, the induced counit $LR \B \to \B$ is an identity-on-objects biequivalence.
    In particular, $\W\text-\cat{Bicat}$ is a reflective sub-tricategory of $\V\text-\cat{Bicat}$.
  \item A $\V$-bicategory lies in the triessential image of $R$ if and only if its hom-objects all lie in $\W$.
  \end{itemize}
\end{Thm}
\begin{proof}
  By the Yoneda lemma, an equivalent definition of $^{\perp}\W$ is as the class of morphisms $V_1\to V_2$ such that the induced functor $\V(V_2,W) \to \V(V_1,W)$ is an equivalence for all $W\in\W$.
  We can then characterize $L V$, for $V\in \V$, as an object of $\W$ equipped with a morphism $\eta_V \colon V\to L V$ that lies in $^{\perp}\W$.

  We define the tensor product of $\W$ by $W_1 \botimes W_2 = L(W_1 \otimes W_2)$, where $\otimes$ denotes the tensor product in $\V$, and its unit object by $I_\W = L (I_\V)$.
  The assumption implies that the induced morphism
  \[ (W_1 \otimes W_2) \otimes W_3 \xrightarrow{\eta \otimes 1} (W_1 \botimes W_2) \otimes W_3 \xrightarrow{\eta} (W_1 \botimes W_2) \botimes W_3  \]
  lies in $^{\perp}\W$; thus, by the alternate characterisation of $^{\perp}\W$, we obtain a morphism $\bar{\mathfrak a}_{W_1,W_2,W_3}$ and invertible $2$-cell as in
  \begin{equation}
    \vcenter{\xymatrix{
        (W_1 \otimes W_2) \otimes W_3\ar[r]\ar[d]_{\mathfrak a_{W_1,W_2,W_3}} \dtwocell{dr}{\cong} &
        (W_1 \botimes W_2) \botimes W_3\ar[d]^{\bar{\mathfrak a}_{W_1,W_2,W_3}}\\
        W_1 \otimes (W_2 \otimes W_3) \ar[r] &
        W_1 \botimes (W_2 \botimes W_3)\rlap{ .}
      }}\label{eq:reflect-omega}
  \end{equation}
  Applying the same argument to $\mathfrak a^\centerdot_{W_1,W_2,W_3}$ yields a pseudoinverse for $\bar{\mathfrak a}_{W_1,W_2,W_3}$; now  corresponding arguments allow us to choose the components of $\bar{\mathfrak l}$ and $\bar{\mathfrak r}$ and their adjoint inverses; the $2$-cells making $\bar{\mathfrak a}$, $\bar{\mathfrak l}$ and $\bar{\mathfrak r}$ pseudonatural; and the 2-cell constraints $\bar\pi$, $\bar\nu$, $\bar\lambda$, and $\bar\rho$.
  The coherence axioms for $\V$ directly imply the coherence axioms for $\W$ (using again the definition of $^{\perp}\W$).
  We make $L$ a strong monoidal functor with $\laxc$ and $\laxu$ identities; the constraint $\omega$ is induced by the isomorphism~\eqref{eq:reflect-omega}, while $\gamma$ and $\delta$ are similarly induced by the analogous isomorphisms defining $\bar{\mathfrak l}$ and $\bar{\mathfrak r}$ respectively.
  It follows by Theorem~\ref{thm:adjoints} that $R$ is lax monoidal; its constraints $\laxc$ and $\laxu$ are (up to isomorphism) the defining morphisms $\eta_{W_1 \otimes W_2}\colon  W_1 \otimes W_2 \to  L(W_1 \otimes W_2)$ and $\eta_{I_\V}\colon  I_\V \to L (I_\V)$.

  The counit $\varepsilon_\B \colon LR\B\to \B$ is, by definition, the identity on objects, and given by the counit of $L\dashv R$ on hom-objects.
  Thus, since the latter counit consists of equivalences, each $\varepsilon_\B$ is an identity-on-objects biequivalence. Finally, if a $\V$-bicategory $\B$ admits a biequivalence $\theta \colon R\C \to \B$, then, since such biequivalences are easily shown to be fully faithful, each hom-object $\B(x,y)$ is equivalent to an object of $\W$, and hence by repleteness must actually lie in $\W$. Conversely,  if a $\V$-bicategory $\A$ has hom-objects in $\W$, then the unit $\A \to R L \A$ is the identity on objects and an equivalence on hom-objects, hence a biequivalence.
\end{proof}

\section{Enriched categories and modules as a free cocompletion}\label{sec:collages}
We have now developed enough enriched bicategory theory to embark upon the second main objective of this paper, that of showing that the basic structures into which enriched (one-dimensional) categories form themselves can be obtained as free cocompletions of certain kinds of enriched bicategory.

In this section, we consider the construction which to a bicategory $\C$ assigns the bicategory of $\C$-enriched categories and modules between them.\footnote{As observed in the introduction, we must be careful to distinguish between a category enriched over a bicategory $\C$, and a bicategory enriched over a monoidal bicategory $\V$; the former notion has homs given by \emph{morphisms} of $\C$, while the latter has homs given by \emph{objects} of $\V$.} Some care is needed in setting up this construction, since composition of $\C$-modules is by tensor product---decategorifying Section~\ref{sec:tensorproductofmodules}---and such tensor products need not always exist. To make our construction well-defined, therefore, we impose a cardinality bound $\kappa$, of a kind to be discussed in the following section, which we use in two ways:
\begin{enumerate}[(i)]
\item We restrict the $\C$-categories over which we consider $\C$-modules to those with a $\kappa$-small set of objects; and
\vskip0.5\baselineskip
\item We assume that the bicategory $\C$ admits reflexive coequalisers and $\kappa$-small coproducts in each hom, preserved by whiskering on each side; we call such a $\C$ \emph{locally $\kappa$-cocomplete}.
\end{enumerate}

Under the assumption (ii) on $\C$, we have enough tensor products to form the bicategory of $\kappa$-small $\C$-categories and $\C$-modules as in (i), yielding an endo-operation $\C \mapsto \cat{Mod}_\kappa(\C)$ on locally $\kappa$-cocomplete bicategories. We will show that this operation is one of free cocompletion; more precisely, we exhibit locally $\kappa$-cocomplete bicategories as bicategories enriched over a suitable monoidal bicategory $\cat{Colim}_\kappa$, and the operation $\cat{Mod}_\kappa(\thg)$ as given by free cocompletion under a class of $\cat{Colim}_\kappa$-enriched colimits.

\subsection{Locally $\kappa$-cocomplete bicategories}\label{subsec:loc-kappa-cocomp}
Throughout the rest of the paper,  $\kappa$ will be an \emph{arity class} in the sense of~\cite{shulman:exsite}.
This means that it is either:
\begin{enumerate}
\item The set $\{1\}$, or
\item The set of cardinal numbers less than some regular cardinal.
\end{enumerate}
The most interesting cases are when $\kappa=\{1\}$, or when it is the set of cardinalities below some inaccessible cardinal.
A set will be called \emph{$\kappa$-small} if its cardinality belongs to $\kappa$. A category will be called \emph{$\kappa$-cocomplete} when it admits reflexive coequalisers and $\kappa$-small coproducts; a \emph{$\kappa$-cocontinuous} functor is one preserving such colimits. As above, a bicategory $\C$ will be called \emph{locally $\kappa$-cocomplete} when its hom-categories are $\kappa$-cocomplete, and its composition functors are $\kappa$-cocontinuous in each variable.

\subsection{$\C$-categories}
\label{subsec:Ccats}

Let $\C$ be a bicategory.
A \emph{$\kappa$-ary $\C$-category} $\f A$, consists of:
\begin{itemize}
\item A $\kappa$-small set $\ob{\f{A}}$;
\item For each $x\in\ob{\f{A}}$, an object $\epsilon x$ of $\C$ (its \emph{extent});
\item For each $x,y\in\ob{\f{A}}$, a morphism ${\f{A}}(x,y) \in \C(\epsilon y, \epsilon x)$;
\item For each $x,y,z\in\ob{\f{A}}$, a 2-cell $m_{xyz}\colon {\f{A}}(x,y) \circ {\f{A}}(y,z) \Rightarrow {\f{A}}(x,z)$;
\item For each $x\in\ob{\f{A}}$, a 2-cell $j_x\colon 1_{\epsilon x} \Rightarrow {\f{A}}(x,x)$;
\end{itemize}
such that the following diagrams commute for all $w,x,y,z \in \ob \f A$:
\begin{equation*}
\cd{
 {\f{A}}(w, x) \circ ({\f{A}}(x, y) \circ {\f{A}}(y, z)) \ar[rr]^{\cong} \ar[d]_{1 \circ m}  & &
 ({\f{A}}(w, x) \circ {\f{A}}(x, y)) \circ {\f{A}}(y, z) \ar[d]^{m \circ 1} \\
 {\f{A}}(w, x) \circ {\f{A}}(x, z) \ar[r]_-m & {\f{A}}(w, z) & {\f{A}}(w, y) \circ {\f{A}}(y, z) \ar[l]^-{\ m}
}
\end{equation*}
\begin{equation*}
\cd[@!C@C-4.5em]{
 & {\f{A}}(x, y) \circ {\f{A}}(y, y) \ar[dr]^{m} \\
 {\f{A}}(x, y) \circ 1_{\epsilon y}  \ar[rr]_{\cong} \ar[ur]^{1 \circ j} & & 
 {\f{A}}(x, y)
}\quad \text{and} \quad
\cd[@!C@C-4.5em]{
 &  {\f{A}}(x, x) \circ {\f{A}}(x, y)  \ar[dr]^{m}  \\
 1_{\epsilon x} \circ {\f{A}}(x, y) \ar[rr]_{\cong} \ar[ur]^{j \circ 1} & & 
 {\f{A}}(x, y)\rlap{ .}
}
\end{equation*}

\begin{Ex}\label{ex:monad}
  If $\ob{\f{A}}$ is a singleton, then ${\f{A}}$ is just a \emph{monad} in $\C$.
  (For this reason,~\cite{benabou:bicats} referred to $\C$-categories as \emph{polyads} in $\C$.)
  In particular, if $\C$ is the bicategory $\cat{Span}({\mathbf{C}})$ of spans in a category $\mathbf{C}$ with pullbacks, then a $\{1\}$-ary $\C$-category is an internal category in $\mathbf{C}$.
\end{Ex}

\begin{Ex}\label{ex:monoidally-enriched}
  If $\C$ has only one object, hence is the delooping of a monoidal category $\mathbf{V}$, then a $\kappa$-ary $\C$-category is a $\mathbf{V}$-enriched category in the usual sense (with a $\kappa$-small set of objects).
\end{Ex}

\begin{Ex}\label{ex:represented-Catcat}
  Suppose ${\f{A}}$ is a $\C$-category and $c\in\C$.
  Then we have a $\cat{Cat}$-category $\C(c,{\f{A}})$ whose objects are those of ${\f{A}}$, with extents $\C(c,\epsilon x)\in\cat{Cat}$.
  The functor $\C(c,{\f{A}})(x,y) \colon \C(c,\epsilon y) \to \C(c,\epsilon x)$ is just postcomposition with ${\f{A}}(x,y)$, and similarly for $m$ and $j$. We emphasise that here, and subsequently, we use $\cat{Cat}$-category in accordance with the preceding definition to mean ``category enriched over the bicategory $\cat{Cat}$'', rather than ``category enriched over the monoidal category $\cat{Cat}$''.
\end{Ex}

\subsection{$\C$-modules}
\label{sec:Cmodules}

If ${\f{A}}$ and ${\f{B}}$ are $\kappa$-ary $\C$-categories, then an \emph{${\f{A}}$-${\f{B}}$-module} $T$ is given by the following data:
\begin{itemize}
\item For each $x\in \ob{\f{A}}$ and $u\in\ob{\f{B}}$, a morphism $T(u,x)\in\C(\epsilon x, \epsilon u)$;
\item For each $x,y\in\ob{\f{A}}$ and $u\in\ob{\f{B}}$, a 2-cell $m_{xyu}\colon T(u,x) \circ {\f{A}}(x,y) \Rightarrow T(u,y)$; and
\item For each $x\in\ob{\f{A}}$ and $u,v\in\ob{\f{B}}$, a 2-cell $m_{xuv}\colon {\f{B}}(u,v) \circ T(v,x) \Rightarrow T(u,x)$,
\end{itemize}
such that diagrams analogous to those in Section~\ref{subsec:Ccats} commute.

If $T$ and $S$ are two ${\f{A}}$-${\f{B}}$-modules, then an \emph{${\f{A}}$-${\f{B}}$-module morphism} $f\colon T\to S$ consists of $2$-cells $f_{ux}\colon T(u,x) \Rightarrow S(u,x)$ such that the following diagrams commute:
\begin{equation*}
\cd[@C+2em]{
 T(u,x) \circ {\f{A}}(x,y) \ar[r]^-{m} \ar[d]_{f_{ux} \circ 1} &
 T(u,y) \ar[d]^{f_{uy}} \\
 S(u,x) \circ {\f{A}}(x,y) \ar[r]_-{m} & S(u,y)
}
\qquad\text{and}\qquad
\cd[@C+2em]{
 {\f{B}}(u,v) \circ T(v,x) \ar[r]^-{m} \ar[d]_{1 \circ f_{vx}} &
 T(u,x) \ar[d]^{f_{ux}} \\
 {\f{B}}(u,v) \circ S(v,x) \ar[r]_-{m} & S(u,x)\rlap{ .}
}\end{equation*}

Thus the $\f A$-$\f B$-modules and their morphisms form a category $\cat{Mod}_{\kappa}(\C)(\f A, \f B)$. If $\C$ is locally $\kappa$-cocomplete, then such categories  comprise the homs of a bicategory $\cat{Mod}_\kappa(\C)$; the composite of an ${\f{A}}$-${\f{B}}$-module $S$ and a ${\f{B}}$-${\f{C}}$-module $T$ is the ${\f{A}}$-${\f{C}}$-module $T\circ S$ defined by the coequaliser:
\begin{equation}
  \xymatrix@C=1pc{ \displaystyle\sum_{x,y\in\ob{\f{B}}}^{\phantom{{\f{B}}}} \scriptstyle T(u,x) \circ {\f{B}}(x,y) \circ S(y,w)
    \ar@<1mm>[r] \ar@<-1mm>[r] &
    \displaystyle\sum_{x\in\ob{\f{B}}}^{\phantom{{\f{B}}}} \scriptstyle T(u,x) \circ S(x,w)
    \ar[r] &
    (T\circ S)(u,w)\rlap{ ,}
  }\label{eq:mod-coeq}
\end{equation}
with the obvious extension of this action to module morphisms. The identity module of ${\f{A}}$ is defined by $1_{\f{A}}(x,y) = {\f{A}}(x,y)$, with actions given by the composition of ${\f{A}}$. The composite $T \circ S$ ``classifies bilinear $\C$-transformations'' out of $S$ and $T$, in the sense that $\f A$-$\f C$-module morphisms $T \circ S \to U$ are in natural bijection with families of maps $T(u,x) \circ S(x,w) \to U(u,w)$ that respect the outer $\f A$- and $\f C$-actions and are ``dinatural'' with respect to the $\f B$ action (a decategorification of Section~\ref{subsec:modulebimorphisms}). It follows easily from this that both triple composites $(T \circ S) \circ R$ and $T \circ (S \circ R)$ ``classify trilinear maps'' out of $R$, $S$ and $T$, whence are canonically isomorphic; these isomorphisms provide the associativity constraints of the bicategory $\cat{Mod}_\kappa(\C)$, and the corresponding classifying properties of four-fold composites now imply the pentagon axiom. The unit constraints and unit coherence axiom follow similarly.

Every $v\in\ob\C$ induces a $\kappa$-ary $\C$-category $\hat{v}$ that has one object $\star$ with $\epsilon(\star) = v$, $\hat{v}(\star,\star) = 1_v$, $m_{\star \star \star} = \mathfrak l = \mathfrak r$, and $j_\star = 1_{1_v}$. Similarly, every $f \colon v \to w$ in $\C$ induces a $\hat v$-$\hat w$-module $\hat f$ with $\hat f(\star,\star) = f$; and every $\alpha \colon f \Rightarrow g$ induces a module transformation $\hat \alpha$ with $\hat \alpha_{\star \star} = \alpha$. Clearly $\widehat{1_v} = 1_{\hat v}$, while for a composite $\hat g \circ \hat f$, the coequaliser in~\eqref{eq:mod-coeq} is trivial, so that $\hat g \circ \hat f \cong \widehat{g \circ f}$. Thus we obtain a fully faithful functor $\hat{(\thg)}\colon \C \to \cat{Mod}_\kappa(\C)$.

\begin{Rk}\label{rk:order}
  Observe that the components of an $\f A$-$\f B$-module $T$ appear to exhibit $\f A$ as acting ``on the right'' and $\f B$ as acting ``on the left''; 
  despite this, we call $T$ an $\f A$-$\f B$-module rather than a $\f B$-$\f A$-module, because we consider a left action to be one that is \emph{covariant} in the variable, and a right action as one that is \emph{contravariant}, regardless of their typographical disposition on the page.
  This discrepancy would be resolved if we were to write composition in $\C$ in ``diagrammatic'' order rather than ``applicative'' order, since then the contravariant $\f B$-action would be displayed  to act on the right and $\f A$ on the left. What would not change is that a morphism $T$ from $\f A$ to $\f B$ in $\cat{Mod}_\kappa(\C)$ involves 1-morphisms for which the following strings are composable in $\C$:
  \[ \xymatrix@C=3pc{ \epsilon v \ar[r]^{\f A(u,v)} & \epsilon u \ar[r]^{T(y,u)} & \epsilon y \ar[r]^{\f B(x,y)} & \epsilon x. } \]
  This is mandated, among other things, by the need for the functor $\hat{(\thg)}$ to be covariant rather than contravariant, since we shall later want to identify it with the embedding of $\C$ into its free cocompletion under $\kappa$-ary collages.

  On the other hand, we appear to be writing composition \emph{in} a $\C$-category in diagrammatic order, since the domain of $m_{xyz}$ is $\f A(x,y)\circ \f A(y,z)$ rather than $\f A(y,z) \circ \f A(x,y)$.
  Of course, if we wrote composition in $\C$ itself in diagrammatic order, then this appearance would be reversed;
  what would not change is that the morphism in $\C$ denoted by $\f A(x,y)$ has domain $\epsilon y$ and codomain $\epsilon x$, rather than vice versa.
  As we will see in Section~\ref{sec:tight-collages}, this is necessary in order to identify $\C$-transformations with 2-cells in the free cocompletion under tight collages.
\end{Rk}

\subsection{Locally $\kappa$-cocomplete bicategories via enrichment}
\label{subsec:loc-cocplt}
We shall now show that the embedding $\hat{(\thg)} \colon \C \to \cat{Mod}_\kappa(\C)$ exhibits $\cat{Mod}_\kappa(\C)$ as the free completion of $\C$ under a suitable class of enriched bicolimits.
The appropriate base for enrichment will be the 2-category $\cat{Colim}_\kappa$ of $\kappa$-cocomplete categories and $\kappa$-cocontinuous functors, equipped with a tensor product which we now construct.

There is a 2-monad $P_\kappa$ on $\cat{Cat}$ whose 2-category $P_\kappa\text{-}\cat{Alg}$ of strict algebras and algebra pseudomorphisms is biequivalent to $\cat{Colim}_\kappa$.
The 2-monad $P_\kappa$ is \emph{lax-idempotent} in the sense of~\cite{kl:property-like}, and thus by~\cite{ilf:pscom-kz} it is pseudo-commutative. Since $\cat{Cat}$ is cocomplete and $P_\kappa$ has a rank, it follows by~\cite{hp:pscom-mnds,ilf:pscom-kz} that $P_\kappa\text{-}\cat{Alg}$, hence also $\cat{Colim}_\kappa$, is a left and right closed monoidal bicategory.
By~\cite{bkp:2dmonads}, it is also complete and cocomplete as a bicategory, and the forgetful functor $\cat{Colim}_\kappa \to \cat{Cat}$ preserves and reflects bilimits.
Moreover, the biadjunction $L \colon \cat{Cat} \rightleftarrows \cat{Colim}_\kappa\colon R$ is monoidal in the sense of Section~\ref{sec:monoidal-adjunctions}.
By construction, the tensor product in $\cat{Colim}_\kappa$ represents functors of multiple variables that are $\kappa$-cocontinuous in each variable separately, and it follows that $\cat{Colim}_\kappa$-enriched bicategories are precisely locally $\kappa$-cocomplete bicategories; likewise, $\cat{Colim}_\kappa$-functors are precisely locally $\kappa$-cocontinuous functors, and so on.

\subsection{Collages}
\label{subsec:collages}
We will now  describe the kinds of $\cat{Colim}_\kappa$-enriched colimits under which the locally $\kappa$-cocomplete $\C$ is to be completed to obtain $\cat{Mod}_\kappa(\C)$. We first give an elementary presentation of these colimits. Given a locally $\kappa$-cocomplete $\C$ and a $\kappa$-ary $\C$-category $\f A$, we may form the functor
\[ \cat{Mod}_\kappa(\C)({\f{A}},\hat{(\thg)}) \colon \C \to \cat{Cat}\rlap{ ;} \]
now a \emph{collage}~\cite{street:cauchy-enr} of ${\f{A}}$ is a birepresentation of this functor.
It comprises an object $v \in \C$ together with a universal ${\f{A}}$-$\hat v$-module $T$. We call the morphisms $T(\star,x) \colon \epsilon x \to v$ the \emph{coprojections} of the collage, and may write the underlying object $v$ as $|{\f{A}}|$. 

Note that a $\C$-category with object set $O$ is precisely a lax functor $\f A \colon \nabla_O \to \C$, where $\nabla_O$ is the chaotic (=~indiscrete) category with object set $O$.
Similarly, if $v\in\ob\C$, then a $\f A$-$\hat v$-module is precisely an oplax natural transformation from the lax functor representing ${\f{A}}$ to the constant functor with value $v$; it follows that a collage of $\f A$ is equally a \emph{lax colimit} of ${\f{A}}\colon \nabla_O \to \C$, that is, a birepresentation of the functor $\cat{Oplax}(\f A, \Delta(\thg)) \colon \C \to \cat{Cat}$. Thus the notion of collage, which was defined above only for locally $\kappa$-cocomplete bicategories, in fact makes sense for any bicategory.

Two special cases of collages are worth singling out.

\begin{Ex}\label{ex:kleisli}
  If ${\f{A}}$ has one object, hence is just a monad in $\C$, then a collage of ${\f{A}}$ is a \emph{Kleisli object} for this monad.
\end{Ex}

\begin{Ex}\label{ex:coproducts}
  Given a $\kappa$-small set $O$ and a function $\epsilon\colon O\to \ob\C$, define
  \[ {\f{A}_\epsilon}(x,y) =
  \begin{cases}
    1_{\epsilon x} &\quad x=y\\
    \emptyset &\quad x\neq y,
  \end{cases}\]
  where $\emptyset$ denotes the initial object of $\C(\epsilon y,\epsilon x)$.
  Then ${\f{A}_\epsilon}$ is a $\kappa$-ary $\C$-category, and a collage of ${\f{A}_\epsilon}$ is just a coproduct $\sum_{x\in O} \epsilon x$.
\end{Ex}

In the locally $\kappa$-cocomplete case, we can construct all collages from these two types.

\begin{Thm}\label{thm:collages}
  If $\C$ is locally $\kappa$-cocomplete, then the following are equivalent:
  \begin{enumerate}
  \item $\C$ admits all lax colimits of lax functors with $\kappa$-small domain.\label{item:c1}
  \item $\C$ admits all $\kappa$-ary collages.\label{item:c2}
  \item $\C$ admits Kleisli objects and $\kappa$-small coproducts.\label{item:c3}
  \end{enumerate}
\end{Thm}
\begin{proof}
  Clearly~(\ref{item:c1})$\Rightarrow$(\ref{item:c2}), while we have just observed that~(\ref{item:c2})$\Rightarrow$(\ref{item:c3}).
  And~(\ref{item:c3})$\Rightarrow$(\ref{item:c1}) is~\cite[Proposition~2.2(a)]{street:cauchy-enr}.
\end{proof}

\subsection{Collages as $\cat{Colim}_\kappa$-colimits}
Let $O$ be a $\kappa$-small set; by the general theory of~\cite{bkp:2dmonads}, there is a bicategory $O^\diamond$ such that lax functors $\nabla_O \to \C$ correspond bijectively with strict functors $O^\diamond \to \C$; and this $O^\diamond$ is ``flexible'', so that any functor $O^\diamond \to \C$ is equivalent to a strict one. Furthermore, there is as in~\cite[\S4]{street:catval-2lim}, a right $O^\diamond$-module $V_O$ such that the lax colimit of $\f A \colon \nabla_O \to \C$ can be identified with the $V_O$-weighted colimit of the corresponding $O^\diamond \to \C$.
Letting $W_O = L V_O$, where $L$ is the free-$\kappa$-cocomplete-category functor from Section~\ref{subsec:loc-cocplt}, it follows from Theorems~\ref{thm:collages} and~\ref{thm:adjoints}
that a locally $\kappa$-cocomplete bicategory admits collages of $\kappa$-ary categories if and only if, when regarded as a $\cat{Colim}_\kappa$-bicategory, it is $\Phi_\kappa$-cocomplete for the class of weights $\Phi_\kappa = \{W_O \mid \text{$O$ a $\kappa$-small set}\}$. We aim to show that $\cat{Mod}_\kappa(\C) \simeq \Phi_\kappa(\C)$; the following few sections gather some preparatory material we shall need for this result.

\subsection{Collages as absolute colimits}
Recall that a \emph{map} in a bicategory is defined to be a left adjoint morphism.
It is thus natural to call a right adjoint morphism a \emph{comap}.
Note also that a $\C$-category $\f A$ can equivalently be regarded as a $\C^\op$-category $\f A^\op$ with $\f A^\op(y,x) = \f A(x,y)$; we refer to collages in $\C^\op$ as \emph{cocollages} in $\C$, and their coprojections in $\C^\op$ as \emph{projections} in $\C$.

\begin{Thm}\label{thm:collages2}
  If $\C$ is locally $\kappa$-cocomplete with $\kappa$-ary collages, then:
  \begin{enumerate}
  \item The coprojections into any collage are maps, and the projections out of any cocollage are comaps.\label{item:c4}
  \item A morphism out of a collage is a map if and only if its composite with each coprojection is a map, and a morphism into a cocollage is a comap if and only if its composite with each projection is a comap.\label{item:c5}
  \item For any $\C$-category ${\f{A}}$ and ${\f{A}}$-$\hat{v}$-module $T$, the following are equivalent:\label{item:c6}
    \begin{enumerate}
    \item $T$ exhibits $v$ as a collage of ${\f{A}}$.\label{item:c6a} 
    \item Each morphism in $T$ is a map, and their right adjoints form a $\hat{v}$-${\f{A}}$-module (with action defined by mates) exhibiting $v$ as a cocollage of ${\f{A}}$.\label{item:c6b} 
    \item $T$ is an equivalence in $\cat{Mod}_\kappa(\C)$.\label{item:c6c}
    \end{enumerate}
  \end{enumerate}
\end{Thm}
\begin{proof}
  The first halves of~(\ref{item:c4}) and~(\ref{item:c5}) are~\cite[Prop.~2.2(b) and (e)]{street:cauchy-enr}.
  The equivalence~(\ref{item:c6a})$\Leftrightarrow$(\ref{item:c6c}) is~\cite[Prop.~2.5]{ckw:axiom-mod}, and this implies the equivalence with~(\ref{item:c6b}) as remarked in~\cite[Remark~2.6(ii)]{ckw:axiom-mod}.
  (Note that in~\cite{ckw:axiom-mod}, our~(\ref{item:c6c}) is taken as the definition of ``collage''.)
  The second halves of~(\ref{item:c4}) and~(\ref{item:c5}) then follow by duality.
\end{proof}

\begin{Cor}
  $\kappa$-ary collages are absolute colimits in locally $\kappa$-cocomplete bicategories; that is, they are preserved by any locally $\kappa$-cocontinuous functor.
\end{Cor}
\begin{proof}
  The construction of $\cat{Mod}_\kappa$ is preserved by any locally $\kappa$-cocontinuous functor, so this follows from~(\ref{item:c6a})$\Leftrightarrow$(\ref{item:c6c}) of Theorem~\ref{thm:collages2}.
\end{proof}

\subsection{Collages in $\cat{Colim}_\kappa$}
\label{subsec:coll-cocplt}

Since $\cat{Colim}_\kappa$ is \emph{closed} monoidal, it is enriched over itself.
(In elementary terms, this is just the fact that cocontinuous functors are closed under pointwise colimits, since colimits commute with colimits.)
Thus, Theorem~\ref{thm:collages2}(\ref{item:c6b}) tells us that we can construct collages in $\cat{Colim}_\kappa$ by constructing cocollages, which (being limits) are created by the forgetful functor $R \colon \cat{Colim}_\kappa \to \cat{Cat}$.

More specifically, suppose ${\f{A}}$ is a $\cat{Colim}_\kappa$-category; thus each $\epsilon x$ is a $\kappa$-cocomplete category and each ${\f{A}}(x,y)\colon \epsilon y \to \epsilon x$ is a $\kappa$-cocontinuous functor.
Then its collage $|{\f{A}}|$ in $\cat{Colim}_\kappa$ (which is different from its collage in $\cat{Cat}$, of course) is the same as its cocollage, which is easy to describe explicitly: it is the category whose objects are given by a family of objects $(\xi_x \in \epsilon x)_{x\in\ob{\f{A}}}$ together with a family of morphisms ${\f{A}}(x,y)(\xi_y) \to \xi_x$ such that for all $x,y,z \in \ob \f A$, the diagram
\[\vcenter{\xymatrix@-.5pc{
    {\f{A}}(x,y)({\f{A}}(y,z)(\xi_z)) \ar[r]\ar[d] &
    {\f{A}}(x,y)(\xi_y)\ar[d]\\
    {\f{A}}(x,z)(\xi_z)\ar[r] &
    \xi_x
  }} \qquad\text{and}\qquad
\vcenter{\xymatrix{
    \xi_x \ar[r]\ar[dr] &
    {\f{A}}(x,x)(\xi_x)\ar[d]\\
    & \xi_x      
  }}
\]
commutes.
Its morphisms are, of course, families $(\psi_x\colon \xi_x \to \zeta_x)_{x\in\ob{\f{A}}}$ such that for all $x,y \in \ob \f A$, the diagram
\[\vcenter{\xymatrix@-.5pc{
    {\f{A}}(x,y)(\xi_y) \ar[r]\ar[d] &
    {\f{A}}(x,y)(\zeta_y)\ar[d]\\
    \xi_x\ar[r] &
    \zeta_x
  }}\]
commutes.
The category $|\f A|$ is $\kappa$-cocomplete, with colimits computed pointwise on the $\xi$'s.
The coprojection $T_x\colon \epsilon x \to |{\f{A}}|$ sends $\omega\in\epsilon x$ to the tuple $(\xi_y \coloneqq {\f{A}}(y,x)(\omega))_y$.

In particular, if ${\f{A}}$ has only one object, so that it is just a $\kappa$-cocontinuous monad on a $\kappa$-cocomplete category, then this describes its Eilenberg--Moore category, which therefore coincides with its Kleisli object in $\cat{Colim}_\kappa$.
Similarly, if each ${\f{A}}(x,y)$ is initial (i.e.\ constant at the initial object of $\epsilon y$), then the morphisms ${\f{A}}(x,y)(\xi_y) \to \xi_x$ and axioms are trivial; thus $|{\f{A}}|$ is just the product $\prod_{x\in\ob{\f{A}}} \epsilon x$, which therefore coincides with the coproduct $\sum_{x\in\ob{\f{A}}} \epsilon x$ in $\cat{Colim}_\kappa$.

Finally, the proof of Theorem~\ref{thm:collages2} supplies the following recipe for the universal property of such collages.
Suppose $B$ is $\kappa$-cocomplete and $T_x\colon \epsilon x \to B$ are $\kappa$-cocontinuous functors forming a lax cocone; thus we have transformations $T_x \circ {\f{A}}(x,y) \to T_y$ satisfying the evident axioms.
Then the induced functor $|T| \colon |{\f{A}}| \to B$ is defined by the following reflexive coequaliser
\begin{equation}
  \xymatrix{
    \displaystyle \sum_{x,y\in\ob{\f{A}}}^{\phantom{\f{A}}} T_x({\f{A}}(x,y)(\xi_y))
    \ar@<1mm>[r] \ar@<-1mm>[r] &
    \displaystyle \sum_{x\in\ob{\f{A}}}^{\phantom{\f{A}}} T_x(\xi_x) \ar[r] &
    |T|(\xi)\rlap{ .}
  }\label{eq:coll-coeq}
\end{equation}

We are finally ready to prove:

\begin{Thm}\label{thm:free-collages}
  For any locally $\kappa$-cocomplete bicategory $\C$, the free cocompletion $\Phi_\kappa(\C)$ of $\C$ under $\kappa$-ary collages is equivalent to $\cat{Mod}_\kappa(\C)$.
\end{Thm}
\begin{proof}
  By its construction in Section~\ref{sec:free-cocompletions}, $\Phi_\kappa(\C)$ is the closure of $\C$ in $\M\C$ under $\kappa$-ary collages.
  This closure certainly contains the full sub-bicategory of $\M\C$ whose objects are collages of the image of some $\kappa$-ary $\C$-category under the Yoneda embedding $Y\colon\C\to\M\C$.
  We will show that this full sub-bicategory is equivalent to $\cat{Mod}_\kappa(\C)$, and that it is closed in $\M\C$ under $\kappa$-ary collages; thus it coincides with the desired closure.

  Up to equivalence, we can certainly consider the bicategory whose objects are literally the $\kappa$-ary $\C$-categories (i.e.\ the objects of $\cat{Mod}_\kappa(\C)$), and whose hom-categories are the hom-categories between their collages in $\M\C$.
  If ${\f{A}}$ is a $\kappa$-ary $\C$-category, then by the construction of colimits in $\M\C$ in Proposition~\ref{prop:mb-is-cocomplete}, the collage $|Y{\f{A}}|$ of its image in $\M\C$ can be defined by
  \[ |Y{\f{A}}|(c) \coloneqq |\C(c,{\f{A}})|. \]
  Here $\C(c,{\f{A}})$ is the $\cat{Cat}$-category from Example~\ref{ex:represented-Catcat}---which because $\C$ is locally $\kappa$-cocomplete, is easily seen to in fact be a $\cat{Colim}_\kappa$-category---and $|\C(c,{\f{A}})|$ denotes its collage in $\cat{Colim}_\kappa$.

  Now suppose ${\f{A}}$ and ${\f{B}}$ are two $\kappa$-ary $\C$-categories.
  By the universal property of collages, a morphism $|Y{\f{A}}| \to |Y{\f{B}}|$ in $\M\C$ is determined by a lax cocone under $Y{\f{A}}$ with vertex $|Y{\f{B}}|$, i.e.\ by a collection of right $\C$-module morphisms $Y(\epsilon x) \to |Y{\f{B}}|$, for each $x\in\ob{\f{A}}$, together with associative module transformations.

  However, by the Yoneda lemma, a right $\C$-module morphism $Y(\epsilon x) \to |Y{\f{B}}|$ is uniquely determined by an object of $|Y{\f{B}}|(\epsilon x) = |\C(\epsilon x,{\f{B}})|$.
  And using the explicit description of collages in $\cat{Colim}_\kappa$ from Section~\ref{subsec:coll-cocplt}, $|\C(\epsilon x,{\f{B}})|$ is equivalent to the category of tuples
  \[\Big(T(z,x)\Big)_{z\in\ob{\f{B}}} \qquad\text{with }T(z,x) \in\C(\epsilon x, \epsilon z)\]
  equipped with associative morphisms ${\f{B}}(w,z)\circ T(z,x) \to T(w,x)$.
  Thus, to give a morphism $|Y{\f{A}}| \to |Y{\f{B}}|$ in $\M\C$ is equivalent to giving such a tuple for each $x\in\ob{\f{A}}$, together with associative morphisms $T(z,x) \circ {\f{A}}(x,y) \to T(z,y)$ that assemble into morphisms in $|\C(\epsilon x,{\f{B}})|$; which is precisely to say that they commute with the action morphisms of $\f B$.
  In sum, what we have is exactly an ${\f{A}}$-${\f{B}}$-module.

  Now the components of the actual morphism $|Y{\f{A}}| \to |Y{\f{B}}|$ can be computed using~\eqref{eq:coll-coeq}.
  It follows straightforwardly that composition of these morphisms in $\M\C$ is computed as in~\eqref{eq:mod-coeq}.
  Similarly, the identity morphism $|Y{\f{A}}| \to |Y{\f{A}}|$ corresponds to the coprojections as described in Section~\ref{subsec:coll-cocplt}, and so to the identity module $1_{\f{A}}$.

  We have shown that the ``one-step'' closure of $\C$ in $\M\C$ under $\kappa$-ary collages is equivalent to $\cat{Mod}_\kappa(\C)$.
  But by~\cite[Proposition~2.2]{ckw:axiom-mod} (plus a little bit of attention to check that $\kappa$-smallness is preserved), $\cat{Mod}_\kappa(\C)$ has $\kappa$-ary collages, which by absoluteness of collages are preserved by the inclusion into $\M\C$. Thus $\cat{Mod}_\kappa(\C)$ is closed in $\M\C$ under $\kappa$-ary collages, and so is the free cocompletion of $\C$ under $\kappa$-ary collages.
\end{proof}

Under the equivalence of the preceding theorem, the embedding of $\C$ into its free cocompletion is identified with $\hat{(\thg)}\colon \C \to \cat{Mod}_\kappa(\C)$.
If $\C$ already has $\kappa$-ary collages, then by Theorem~\ref{thm:collages2}(\ref{item:c6c}), $\hat{(\thg)}$ is an equivalence (and conversely).
This is as we expect for a cocompletion with respect to a class of absolute weights.
In particular, we have~\cite[Corollary~2.4]{ckw:axiom-mod}: a locally $\kappa$-cocomplete $\C$ is equivalent to $\cat{Mod}_\kappa(\B)$, for some locally $\kappa$-cocomplete $\B$, if and only if $\C$ has $\kappa$-ary collages (and in this case we can take $\B=\C$).

\begin{Rk}
  As is often the case for absolute colimits, if a bicategory has $\kappa$-ary collages that behave suitably nicely, then it is automatically locally $\kappa$-cocomplete.
  In this case, the requisite niceness properties are Theorem~\ref{thm:collages2}(\ref{item:c4}), (\ref{item:c5}), and (\ref{item:c6a})$\Rightarrow$(\ref{item:c6b}); see~\cite[Proposition~3.3]{ckw:axiom-mod}.
  Note the strong analogy with the fact that a category with finite coproducts is enriched over commutative monoids if and only if those coproducts are biproducts.
\end{Rk}

\begin{Rk}\label{remark:laxcolim-completion}
  By Theorem~\ref{thm:collages}, $\cat{Mod}_\kappa(\C)$ also has all lax colimits of lax functors with $\kappa$-small domain.
  Thus, it is also the free cocompletion of $\C$ under such lax colimits.
\end{Rk}

\begin{Rk}\label{remark:matmod}
  We saw in Theorem~\ref{thm:collages} that $\Phi_\kappa$-cocompleteness of a locally $\kappa$-cocomplete bicategory is equivalent to cocompleteness for $\Phi_{\{1\}}$ (consisting of Kleisli objects) and also for the class $\Phi_{\kappa\sqcup}$ of of $\kappa$-small coproducts.
  In general, if $\Phi$-cocompleteness is equivalent to $\Phi_1$-cocompleteness together with $\Phi_2$-cocompleteness, it does not follow that $\Phi(\C)$ may be computed as $\Phi_1(\Phi_2(\C))$ or $\Phi_2(\Phi_1(\C))$; a transfinite iteration is often required.
  However, in this case it is true that $\Phi_\kappa(\C) \simeq \Phi_{\{1\}}(\Phi_{\kappa\sqcup}(\C))$ (though not in the reverse order). Let us explain briefly why.

  Recall that a $\kappa$-small coproduct $\sum_{x\in O} \epsilon x$ in a locally $\kappa$-cocomplete bicategory $\C$ is equally a collage of the $\C$-category ${\f{A}_\epsilon}$ from Example~\ref{ex:coproducts}.
  Thus, the full sub-bicategory of $\M\C$ consisting of $\kappa$-small coproducts of representables is equivalent to the full sub-bicategory of $\cat{Mod}_\kappa(\C)$ determined by the $\C$-categories of this form.
  But an ${\f{A}_{\epsilon_1}}$-${\f{A}_{\epsilon_2}}$-module is nothing but an $(O_1\times O_2)$-matrix of morphisms $\epsilon_1 x \to \epsilon_2 y$ in $\C$, and the  composition of such modules is given simply by ``matrix multiplication'', so that the bicategory in question is that which~\cite[\S4.2]{ckw:axiom-mod} denotes by $\cat{Matr}_\kappa(\C)$; it has $\kappa$-small coproducts, and hence is $\Phi_{\kappa\sqcup}(\C)$.
  (In fact, this analysis can be performed with bicategories having merely local $\kappa$-small \emph{coproducts}, regarded as enriched over a monoidal bicategory $\cat{Coprod}_\kappa$.)

  Given these explicit descriptions, the results of~\cite[Sections~1--3]{bcsw:variation-enr} now show that $\cat{Mod}_\kappa(\C) \simeq \cat{Mod}_{\{1\}}(\cat{Matr}_\kappa(\C))$, whence $\Phi_\kappa(\C) \simeq \Phi_{\{1\}}(\Phi_{\kappa\sqcup}(\C))$. This decomposition can be seen as expressing a \emph{distributive law} (at a suitable level of weakness) between the cocompletion operations $\Phi_{\{1\}}$ and $\Phi_{\kappa\sqcup}$.
\end{Rk}

\section{Enriched categories, functors, and modules as a free cocompletion}
\label{sec:tight-collages}
We now refine the results of the preceding section to describe a free cocompletion that from a given enrichment base $\C$ generates not only the bicategory $\cat{Mod}_\kappa(\C)$ of $\kappa$-ary $\C$-categories and $\C$-modules, but also the bicategory $\cat{Cat}_\kappa(\C)$ of $\kappa$-ary $\C$-categories and $\C$-functors, together with the embedding $J_\kappa \colon \cat{Cat}_\kappa(\C) \to \cat{Mod}_\kappa(\C)$ of the former into the latter. Since a free cocompletion process must produce the same kind of structure as output as it consumes as input, this means that the appropriate notion of ``base for enrichment'' now changes: it will involve a pair of bicategories related by an identity on objects and locally fully faithful functor. We call instances of this notion \emph{equipments}; they are a generalisation of the proarrow equipments of~\cite{wood:proarrows-i}. In this section, we first describe the construction assigning to any (locally $\kappa$-cocomplete) equipment $\C$, the equipment $\fcat{Mod}_\kappa(\C)$ of $\kappa$-ary $\C$-enriched categories, functors and modules; we then explain how locally $\kappa$-cocomplete equipments may be viewed as bicategories enriched in a certain monoidal bicategory $\F_\kappa$, and finally, exhibit the assignation $\C \mapsto \fcat{Mod}_\kappa(\C)$ as a free cocompletion process on $\F_\kappa$-bicategories.

\subsection{Equipments}
\label{subsec:equipments}
As anticipated above, we define an \emph{equipment} $\C$ to be given by a pair of bicategories $\C\ti$ and $\C\lo$ with the same objects and a functor $J_\C \colon \C\ti \to \C\lo$ that is the identity on objects and locally fully faithful. We refer to morphisms in $\C\ti$ as \emph{tight} and morphisms in $\C\lo$ as \emph{loose}.
Every tight morphism $f$ has an underlying loose morphism $J_\C(f)$; conversely, by a \emph{tightening} of a loose morphism $g$, we mean a tight morphism $f$ and invertible $2$-cell $J_\C(f) \cong g$. Since $J_\C$ is locally fully faithful, any two tightenings of a loose morphism are uniquely isomorphic.

Note that an equipment, in our sense, is a structure satisfying the first two axioms of a proarrow equipment as defined in~\cite{wood:proarrows-i}. The third axiom given there is that each morphism in the image of $J_\C$ is a map; in this paper, we shall refer to equipments satisfying this extra axiom as \emph{map equipments}.

A \emph{morphism of equipments} $F \colon \C \to \D$ is given by a pair of functors $F\ti, F\lo$ between the respective tight and loose parts, together with an invertible icon
\begin{equation*}
\cd{
\C\ti \ar[d]_{J_\C} \ar[r]^-{F\ti} \dtwocell{dr}{J_F} &
\D\ti \ar[d]^{J_\D} \\
\C\lo \ar[r]_-{F\lo} &
\D\lo\rlap{ .}
}
\end{equation*}
Note that such an $F$ has the property that, whenever $f \colon x \to y$ in $\C\lo$ admits a tightening, so too does $F\lo(f)$; in fact, given just $F\lo$ with this property, we may recover $F\ti$ from it to within an invertible icon. There are corresponding notions of transformation, modification and icon between equipment morphisms; we do not give the details since we will not need them.

\begin{Ex}\label{ex:chordate-inchordate}
A bicategory $\C$ can be made into an equipment in two canonical ways. For the first, we take $\C\ti = \C\lo = \C$ and $J_\C = 1$, so that every morphism in $\C$ is tight in a unique way; following~\cite{ls:limlax}, we call such an equipment \emph{chordate}. For the second, we factorise the inclusion-of-objects functor $\ob \C \to \C$ as bijective on objects and $1$-cells followed by locally fully faithful; the second part of this factorisation gives an equipment $J_\C \colon \C\ti \to \C\lo$ in which only identity morphisms admit tightenings, which are unique. As in~\cite{ls:limlax}, we call such an equipment \emph{inchordate}. Both the chordate and the inchordate equippings can be made functorial with respect to $1$-, $2$- and $3$-cells between bicategories.
\end{Ex}

\subsection{Categories and functors enriched in an equipment}

Let $\C$ be an equipment. A $\kappa$-ary \emph{$\C$-category} is simply a $\kappa$-ary $\C\lo$-category in the sense of Section~\ref{subsec:Ccats}. Given $\kappa$-ary $\C$-categories $\f A$ and $\f B$, a \emph{$\C$-functor} $D\colon{\f{A}}\to {\f{B}}$ is given by the following data:
\begin{itemize}
\item For each $x\in\ob {\f{A}}$, an object $Dx\in\ob{\f{B}}$ and a \emph{tight} morphism $D_x \colon \epsilon x \to \epsilon (Dx)$;
\item For each $x,y\in\ob{\f{A}}$, a 2-cell $D_{xy}\colon D_x \circ {\f{A}}(x,y) \Rightarrow {\f{B}}(Dx,Dy) \circ D_y$,
\end{itemize}
subject to the two axioms that, for all $x,y,z \in \ob \f A$, we have:
\begin{equation*}
\vcenter{\hbox{\begin{tikzpicture}[y=0.85pt, x=1.14pt,yscale=-1, inner sep=0pt, outer sep=0pt, every text node part/.style={font=\tiny} ]
  \path[draw=black,line join=miter,line cap=butt,line width=0.650pt] (14.0000,952.3622) .. controls (45.6430,952.3622) and (73.2383,950.2136) .. 
  node[above right=0.12cm,at start] {\!$\f A(y,\!z)$}(78.9810,965.0901);
  \path[draw=black,line join=miter,line cap=butt,line width=0.650pt] (14.0000,972.3622) .. controls (36.8049,972.3622) and (38.8735,976.9584) .. 
  node[above right=0.12cm,at start] {\!$\f A(x,\!y)$}(42.2322,980.4939);
  \path[draw=black,line join=miter,line cap=butt,line width=0.650pt] (14.0000,992.3622) .. controls (36.4534,992.3622) and (39.9341,985.9982) .. 
  node[above right=0.12cm,at start] {\!$D_x$}(42.0555,983.8769);
  \path[draw=black,line join=miter,line cap=butt,line width=0.650pt] (46.2218,983.7001) .. controls (59.7943,988.9377) and (99.7285,989.7962) .. 
  node[below=0.09cm] {$\f B(Dx,\!Dy)$}(124.2426,983.7001);
  \path[draw=black,line join=miter,line cap=butt,line width=0.650pt] (45.6915,980.4939) .. controls (48.5199,974.4835) and (68.0325,967.9428) .. 
  node[above=0.08cm,pos=0.5] {$D_y$}(78.6759,967.9428);
  \path[draw=black,line join=miter,line cap=butt,line width=0.650pt] (81.4280,965.1906) .. controls (94.5975,941.3532) and (130.7507,961.4289) .. 
  node[above left=0.12cm,at end] {$D_z$\!}(168.0000,962.3622);
  \path[draw=black,line join=miter,line cap=butt,line width=0.650pt] (81.7816,968.0433) .. controls (93.9843,968.0433) and (111.6395,976.4280) .. 
  node[above right=0.04cm,pos=0.17,rotate=-12] {$\f B(Dy,\!Dz)$}(124.1143,980.3171);
  \path[draw=black,line join=miter,line cap=butt,line width=0.650pt] (127.4038,982.1091) -- 
  node[above left=0.12cm,at end] {$\f B(Dx,\!Dz)$\!}(168.0000,982.3622);
  \path[fill=black] (43.674252,982.35858) node[circle, draw, line width=0.65pt, minimum width=5mm, fill=white, inner sep=0.25mm] (text5087) {$D$
     };
  \path[fill=black] (80.069435,966.80225) node[circle, draw, line width=0.65pt, minimum width=5mm, fill=white, inner sep=0.25mm] (text5091) {$D$
     };
  \path[fill=black] (126.22198,982.18182) node[circle, draw, line width=0.65pt, minimum width=5mm, fill=white, inner sep=0.25mm] (text5095) {$m$
   };
\end{tikzpicture}}}
\qquad=\qquad
 \vcenter{\hbox{\raisebox{0.14\height}{\begin{tikzpicture}[y=0.85pt, x=0.8pt,yscale=-1, inner sep=0pt, outer sep=0pt, every text node part/.style={font=\tiny} ]
  \path[draw=black,line join=miter,line cap=butt,line width=0.650pt] (5.0000,952.3622) .. controls (31.2500,952.3622) and (43.5000,952.8622) .. 
  node[above right=0.12cm,at start] {\!$\f A(y,\!z)$}(52.0000,960.3622);
  \path[draw=black,line join=miter,line cap=butt,line width=0.650pt] (5.0000,972.3622) .. controls (31.0000,972.3622) and (42.7500,971.8622) .. 
  node[above right=0.12cm,at start] {\!$\f A(x,\!y)$}(52.0000,964.6122);
  \path[draw=black,line join=miter,line cap=butt,line width=0.650pt] (5.0000,992.3622) .. controls (47.0104,992.3622) and (84.5966,982.7655) .. 
  node[above right=0.12cm,at start] {\!$D_x$}(91.5000,974.3622);
  \path[draw=black,line join=miter,line cap=butt,line width=0.650pt] (54.0000,962.3622) .. controls (66.0154,962.3622) and (86.5000,963.8622) .. 
  node[above right=0.05cm,pos=0.22] {$\f A(x,z)$}(91.5000,971.1122);
  \path[draw=black,line join=miter,line cap=butt,line width=0.650pt] (94.7500,970.8622) .. controls (104.5000,962.3622) and (120.7095,962.3622) .. 
  node[above left=0.12cm,at end] {$D_z$\!}(160.0000,962.3622);
  \path[draw=black,line join=miter,line cap=butt,line width=0.650pt] (94.9482,974.5086) .. controls (100.0584,981.7526) and (121.2350,982.3622) .. 
  node[above left=0.12cm,at end] {$\f B(Dx,\!Dz)$\!}(160.0000,982.3622);
  \path[fill=black] (52.75,962.36218) node[circle, draw, line width=0.65pt, minimum width=5mm, fill=white, inner sep=0.25mm] (text5442) {$m$
     };
  \path[fill=black] (93,973.11218) node[circle, draw, line width=0.65pt, minimum width=5mm, fill=white, inner sep=0.25mm] (text5446) {$D$
   };
\end{tikzpicture}}}}
\end{equation*}
in $\C(\epsilon z, \epsilon(Dx))$; and that for all $x \in \ob \f A$, we have
\begin{equation*}
\vcenter{\hbox{\begin{tikzpicture}[y=0.85pt, x=0.8pt,yscale=-1, inner sep=0pt, outer sep=0pt, every text node part/.style={font=\tiny} ]
  \path[draw=black,line join=miter,line cap=butt,line width=0.650pt] (20.0000,982.3622) .. controls (70.0000,982.3622) and (83.1755,982.6866) .. 
  node[above right=0.12cm,at start] {\!$D_x$}(91.5000,974.3622);
  \path[draw=black,line join=miter,line cap=butt,line width=0.650pt] (50.0000,962.3622) .. controls (70.2500,962.3622) and (82.6241,962.2363) .. 
  node[above=0.1cm,pos=0.48] {\!$\f A(x,\!x)$}(91.5000,971.1122);
  \path[draw=black,line join=miter,line cap=butt,line width=0.650pt] (94.7500,970.8622) .. controls (104.5000,962.3622) and (120.7095,962.3622) .. 
  node[above left=0.12cm,at end] {$D_x$\!}(160.0000,962.3622);
  \path[draw=black,line join=miter,line cap=butt,line width=0.650pt] (94.9482,974.5086) .. controls (102.6256,982.1860) and (121.2350,982.3622) .. 
  node[above left=0.12cm,at end] {$\f B(Dx,\!Dx)$\!}(160.0000,982.3622);
  \path[fill=black] (46.75,962.36218) node[circle, draw, line width=0.65pt, minimum width=5mm, fill=white, inner sep=0.25mm] (text5442) {$j$
     };
  \path[fill=black] (93,973.11218) node[circle, draw, line width=0.65pt, minimum width=5mm, fill=white, inner sep=0.25mm] (text5446) {$D$
   };
\end{tikzpicture}}} \qquad = \qquad
\vcenter{\hbox{\begin{tikzpicture}[y=0.85pt, x=0.8pt,yscale=-1, inner sep=0pt, outer sep=0pt, every text node part/.style={font=\tiny} ]
  \path[draw=black,line join=miter,line cap=butt,line width=0.650pt] (50.0000,992.3622) .. controls (103.5091,992.3622) and (80.6919,972.3622) .. 
  node[above right=0.12cm,at start] {\!$D_x$}node[above left=0.12cm,at end] {$D_x$\!}(175.0000,972.3622);
  \path[draw=black,line join=miter,line cap=butt,line width=0.650pt] (112.0000,992.3622) -- 
  node[above left=0.12cm,at end] {$\f B(Dx,\!Dx)$\!}(175.0000,992.3622);
  \path[fill=black] (112,992.36218) node[circle, draw, line width=0.65pt, minimum width=5mm, fill=white, inner sep=0.25mm] (text5442) {$j$
     };
\end{tikzpicture}}}
\end{equation*}
in $\C(\epsilon x, \epsilon(Dx))$.
It is crucial in the above definition that the morphisms $D_x$ are tight. Note that in the domain and codomain of $D_{xy}$, and in the axioms, we have omitted to notate the functor $J_\C$ that ought to be applied to occurrences of $D_x$ and $D_y$. We will continue in such abuses without further comment.

If $E\colon{\f{A}}\to{\f{B}}$ is another $\C$-functor, then a \emph{$\C$-transformation} $\vartheta\colon D\Rightarrow E$ is given by 
$2$-cells $\vartheta_x \colon D_x \Rightarrow {\f{B}}(Dx,Ex) \circ E_x$ for all $x\in \ob {\f{A}}$, 
such that for all $x,y \in \ob  \f A$, we have:
\begin{equation*}
\vcenter{\hbox{\begin{tikzpicture}[y=0.85pt, x=1.14pt,yscale=-1, inner sep=0pt, outer sep=0pt, every text node part/.style={font=\tiny} ]
  \path[draw=black,line join=miter,line cap=butt,line width=0.650pt] (14.0000,972.3622) .. controls (36.8049,972.3622) and (38.8735,976.9584) .. 
  node[above right=0.12cm,at start] {\!$\f A(x,\!y)$}(42.2322,980.4939);
  \path[draw=black,line join=miter,line cap=butt,line width=0.650pt] (14.0000,992.3622) .. controls (36.4534,992.3622) and (39.9341,985.9982) .. 
  node[above right=0.12cm,at start] {\!$D_x$}(42.0555,983.8769);
  \path[draw=black,line join=miter,line cap=butt,line width=0.650pt] (46.2218,983.7001) .. controls (59.7943,988.9377) and (99.7285,989.7962) .. 
  node[below=0.09cm] {$\f B(Dx,\!Dy)$}(124.2426,983.7001);
  \path[draw=black,line join=miter,line cap=butt,line width=0.650pt] (45.6915,980.4939) .. controls (48.5199,974.4835) and (68.0325,967.9428) .. 
  node[above=0.08cm,pos=0.5] {$D_y$}(78.6759,967.9428);
  \path[draw=black,line join=miter,line cap=butt,line width=0.650pt] (81.4280,965.1906) .. controls (94.5975,941.3532) and (130.7507,961.4289) .. 
  node[above left=0.12cm,at end] {$E_y$\!}(168.0000,962.3622);
  \path[draw=black,line join=miter,line cap=butt,line width=0.650pt] (81.7816,968.0433) .. controls (93.9843,968.0433) and (111.6395,976.4280) .. 
  node[above right=0.04cm,pos=0.15,rotate=-12] {$\f B(Dy,\!Ey)$}(124.1143,980.3171);
  \path[draw=black,line join=miter,line cap=butt,line width=0.650pt] (127.4038,982.1091) -- 
  node[above left=0.12cm,at end] {$\f B(Dx,\!Ey)$\!}(168.0000,982.3622);
  \path[fill=black] (43.674252,982.35858) node[circle, draw, line width=0.65pt, minimum width=5mm, fill=white, inner sep=0.25mm] (text5087) {$D$
     };
  \path[fill=black] (80.069435,966.80225) node[circle, draw, line width=0.65pt, minimum width=5mm, fill=white, inner sep=0.25mm] (text5091) {$\vartheta$
     };
  \path[fill=black] (126.22198,982.18182) node[circle, draw, line width=0.65pt, minimum width=5mm, fill=white, inner sep=0.25mm] (text5095) {$m$
   };
\end{tikzpicture}}} \qquad = \qquad
\vcenter{\hbox{\begin{tikzpicture}[y=0.85pt, x=1.14pt,yscale=-1, inner sep=0pt, outer sep=0pt, every text node part/.style={font=\tiny} ]
  \path[draw=black,line join=miter,line cap=butt,line width=0.650pt] (14.0000,952.3622) .. controls (45.6430,952.3622) and (73.2383,950.2136) .. 
  node[above right=0.12cm,at start] {\!$\f A(x,\!y)$}(78.9810,965.0901);
  \path[draw=black,line join=miter,line cap=butt,line width=0.650pt] (14.0000,992.3622) .. controls (36.4534,992.3622) and (39.9341,985.9982) .. 
  node[above right=0.12cm,at start] {\!$D_x$}(42.0555,983.8769);
  \path[draw=black,line join=miter,line cap=butt,line width=0.650pt] (46.2218,983.7001) .. controls (59.7943,988.9377) and (99.7285,989.7962) .. 
  node[below=0.09cm] {$\f B(Dx,\!Ex)$}(124.2426,983.7001);
  \path[draw=black,line join=miter,line cap=butt,line width=0.650pt] (45.6915,980.4939) .. controls (48.5199,974.4835) and (68.0325,967.9428) .. 
  node[above=0.08cm,pos=0.5] {$E_x$}(78.6759,967.9428);
  \path[draw=black,line join=miter,line cap=butt,line width=0.650pt] (81.4280,965.1906) .. controls (94.5975,941.3532) and (130.7507,961.4289) .. 
  node[above left=0.12cm,at end] {$E_y$\!}(168.0000,962.3622);
  \path[draw=black,line join=miter,line cap=butt,line width=0.650pt] (81.7816,968.0433) .. controls (93.9843,968.0433) and (111.6395,976.4280) .. 
  node[above right=0.04cm,pos=0.15,rotate=-12] {$\f B(Ex,\!Ey)$}(124.1143,980.3171);
  \path[draw=black,line join=miter,line cap=butt,line width=0.650pt] (127.4038,982.1091) -- 
  node[above left=0.12cm,at end] {$\f B(Dx,\!Ey)$\!}(168.0000,982.3622);
  \path[fill=black] (43.674252,982.35858) node[circle, draw, line width=0.65pt, minimum width=5mm, fill=white, inner sep=0.25mm] (text5087) {$\vartheta$
     };
  \path[fill=black] (80.069435,966.80225) node[circle, draw, line width=0.65pt, minimum width=5mm, fill=white, inner sep=0.25mm] (text5091) {$E$
     };
  \path[fill=black] (126.22198,982.18182) node[circle, draw, line width=0.65pt, minimum width=5mm, fill=white, inner sep=0.25mm] (text5095) {$m$
   };
\end{tikzpicture}}}
\end{equation*}
in $\C(\epsilon y, \epsilon(Dx))$.
\begin{Ex}\label{ex:c-icon}
Let $D, E \colon \f A \to \f B$ be $\C$-functors that agree on objects. A \emph{$\C$-icon} from $D$ to $E$ is given by a family of $2$-cells $\vartheta_x \colon D_x \Rightarrow E_x \colon \epsilon x \to \epsilon(Dx) = \epsilon(Ex)$ in $\C\ti$ satisfying the axiom
\begin{equation*}
\vcenter{\hbox{\begin{tikzpicture}[y=0.95pt, x=0.8pt,yscale=-1, inner sep=0pt, outer sep=0pt, every text node part/.style={font=\tiny} ]
  \path[draw=black,line join=miter,line cap=butt,line width=0.650pt] (20.0000,982.3622) .. controls (70.0000,982.3622) and (83.1755,982.6866) .. 
  node[above right=0.12cm,at start] {\!$D_x$}node[above=0.09cm,pos=0.45] {$E_x$}(91.5000,974.3622);
  \path[draw=black,line join=miter,line cap=butt,line width=0.650pt] (20.0000,962.3622) .. controls (70.2500,962.3622) and (82.6241,962.2363) .. 
  node[above right=0.12cm,at start] {\!$\f A(x,\!y)$}(91.5000,971.1122);
  \path[draw=black,line join=miter,line cap=butt,line width=0.650pt] (94.7500,970.8622) .. controls (104.5000,962.3622) and (120.7095,962.3622) .. 
  node[above left=0.12cm,at end] {$E_y$\!}(160.0000,962.3622);
  \path[draw=black,line join=miter,line cap=butt,line width=0.650pt] (94.9482,974.5086) .. controls (102.6256,982.1860) and (121.2350,982.3622) .. 
  node[above left=0.12cm,at end] {$\f B(Ex,\!Ey)$\!}(160.0000,982.3622);
  \path[fill=black] (48,982.36218) node[circle, draw, line width=0.65pt, minimum width=5mm, fill=white, inner sep=0.25mm] (text5442) {$\vartheta$
     };
  \path[fill=black] (93,973.11218) node[circle, draw, line width=0.65pt, minimum width=5mm, fill=white, inner sep=0.25mm] (text5446) {$E$
   };
\end{tikzpicture}}} \qquad = \qquad
\vcenter{\hbox{\begin{tikzpicture}[y=0.95pt, x=0.8pt,yscale=-1, inner sep=0pt, outer sep=0pt, every text node part/.style={font=\tiny} ]
  \path[draw=black,line join=miter,line cap=butt,line width=0.650pt] (50.0000,982.3622) .. controls (70.0000,982.3622) and (83.1755,982.6866) .. 
  node[above right=0.12cm,at start] {\!$D_x$}(91.5000,974.3622);
  \path[draw=black,line join=miter,line cap=butt,line width=0.650pt] (50.0000,962.3622) .. controls (70.2500,962.3622) and (82.6241,962.2363) .. 
  node[above right=0.12cm,at start] {\!$\f A(x,\!y)$}(91.5000,971.1122);
  \path[draw=black,line join=miter,line cap=butt,line width=0.650pt] (94.7500,970.8622) .. controls (104.5000,962.3622) and (120.7095,962.3622) .. 
  node[above=0.09cm,pos=0.45] {$D_y$}node[above left=0.12cm,at end] {$E_y$\!}(170.0000,962.3622);
  \path[draw=black,line join=miter,line cap=butt,line width=0.650pt] (94.9482,974.5086) .. controls (102.6256,982.1860) and (121.2350,982.3622) .. 
  node[above left=0.12cm,at end] {$\f B(Ex,\!Ey)$\!}(170.0000,982.3622);
  \path (78,986.36218) node[circle] (text5442) {$\phantom{\vartheta}$
     };
  \path[fill=black] (138,962.36218) node[circle, draw, line width=0.65pt, minimum width=5mm, fill=white, inner sep=0.25mm] (text5442) {$\vartheta$
     };
  \path[fill=black] (93,973.11218) node[circle, draw, line width=0.65pt, minimum width=5mm, fill=white, inner sep=0.25mm] (text5446) {$E$
   };
\end{tikzpicture}\ \rlap{ .}}}
\end{equation*}
Every $\C$-icon gives rise to a $\C$-transformation $D \Rightarrow E$ with components
\[
D_x \xrightarrow{\ \ \vartheta\ \ } E_x \xrightarrow{\ \ \cong\ \ } 1_{\epsilon(Dx)} \circ E_x \xrightarrow{\ \ j \circ 1\ \ } \f B(Dx, Ex) \circ E_x\rlap{ .}
\]
\end{Ex}

\subsection{The bicategory of $\C$-categories}
\label{sec:bicat-Ccats}
For any equipment $\C$, the $\kappa$-ary $\C$-categories, $\C$-functors and $\C$-transformations form a bicategory $\cat{Cat}_\kappa(\C)$ as follows.

 Given a $\C$-functor $D \colon \f A \to \f B$, the identity $\C$-transformation $D\Rightarrow D$ is induced by the $\C$-icon whose components are identity $2$-cells 
$1_{D_x} \colon D_x \Rightarrow D_x$.
Given $\C$-transformations $\vartheta\colon D\Rightarrow E$ and $\varsigma\colon E\Rightarrow F$, their composite $\varsigma\vartheta\colon D\Rightarrow F$ has components
\begin{equation*}
\vcenter{\hbox{\begin{tikzpicture}[y=0.85pt, x=1.14pt,yscale=-1, inner sep=0pt, outer sep=0pt, every text node part/.style={font=\tiny} ]
  \path[draw=black,line join=miter,line cap=butt,line width=0.650pt] (14,982.1091) -- 
  node[above right=0.12cm,at start] {\!$D_x$}(42.0555,982.1091);
  \path[draw=black,line join=miter,line cap=butt,line width=0.650pt] (46.2218,983.7001) .. controls (59.7943,988.9377) and (99.7285,989.7962) .. 
  node[below=0.09cm] {$\f B(Dx,\!Ex)$}(124.2426,983.7001);
  \path[draw=black,line join=miter,line cap=butt,line width=0.650pt] (45.6915,980.4939) .. controls (48.5199,974.4835) and (68.0325,967.9428) .. 
  node[above=0.08cm,pos=0.5] {$E_x$}(78.6759,967.9428);
  \path[draw=black,line join=miter,line cap=butt,line width=0.650pt] (81.4280,965.1906) .. controls (94.5975,941.3532) and (130.7507,961.4289) .. 
  node[above left=0.12cm,at end] {$F_x$\!}(168.0000,962.3622);
  \path[draw=black,line join=miter,line cap=butt,line width=0.650pt] (81.7816,968.0433) .. controls (93.9843,968.0433) and (111.6395,976.4280) .. 
  node[above right=0.04cm,pos=0.15,rotate=-12] {$\f B(Ex,\!Fx)$}(124.1143,980.3171);
  \path[draw=black,line join=miter,line cap=butt,line width=0.650pt] (127.4038,982.1091) -- 
  node[above left=0.12cm,at end] {$\f B(Dx,\!Fx)$\!}(168.0000,982.3622);
  \path[fill=black] (43.674252,982.35858) node[circle, draw, line width=0.65pt, minimum width=5mm, fill=white, inner sep=0.25mm] (text5087) {$\vartheta$
     };
  \path[fill=black] (80.069435,966.80225) node[circle, draw, line width=0.65pt, minimum width=5mm, fill=white, inner sep=0.25mm] (text5091) {$\varsigma$
     };
  \path[fill=black] (126.22198,982.18182) node[circle, draw, line width=0.65pt, minimum width=5mm, fill=white, inner sep=0.25mm] (text5095) {$m$
   };
\end{tikzpicture}}} 
\end{equation*}
(Observe that the $\C$-transformations induced by $\C$-icons are stable under this composition).
This defines the hom-category $\cat{Cat}_\kappa(\C)(\f A, \f B)$.

Suppose now that $D\colon \f A \to \f B$ and $E\colon \f B \to \f C$ are $\C$-functors; we take their composite $E D \colon \f A \to \f C$ to have action on objects $(ED)(x) = E(Dx)$, $1$-cell components given by $(ED)_x = E_{D x} \circ D_x \colon \epsilon x \to \epsilon (E D x)$, and $2$-cell components $(ED)_{xy}$ given by
\begin{equation*}
\vcenter{\hbox{\begin{tikzpicture}[y=0.8pt, x=0.9pt,yscale=-1, inner sep=0pt, outer sep=0pt, every text node part/.style={font=\tiny} ]
  \path[draw=black,line join=miter,line cap=butt,line width=0.650pt] (60.0000,932.3622) .. controls (80.0000,932.3622) and (90.0505,932.6866) .. 
  node[above right=0.12cm,at start] {\!$D_x$}(98.3750,924.3622);
  \path[draw=black,line join=miter,line cap=butt,line width=0.650pt] (60.0000,912.3622) .. controls (74.2500,912.3622) and (89.4991,911.6113) .. 
  node[above right=0.12cm,at start] {\!$\f A(x,\!y)$}(98.3750,920.4872);
  \path[draw=black,line join=miter,line cap=butt,line width=0.650pt] (102.0000,920.3622) .. controls (111.7500,911.8622) and (160.7095,912.3622) .. 
  node[above left=0.12cm,at end] {$D_y$\!}(220.0000,912.3622);
  \path[draw=black,line join=miter,line cap=butt,line width=0.650pt] (102.0732,924.3836) .. controls (111.5184,933.8288) and (127.5236,930.1389) .. 
  node[above right=0.04cm,pos=0.42] {$\f B(Dx,\!Dy)$}(137.7019,940.3171);
  \path[fill=black] (100,922.36218) node[circle, draw, line width=0.65pt, minimum width=5mm, fill=white, inner sep=0.25mm] (text5446-2) {$D$
     };
  \path[draw=black,line join=miter,line cap=butt,line width=0.650pt] (60.0000,952.3622) .. controls (80.0000,952.3622) and (130.0505,952.6866) .. 
  node[above right=0.12cm,at start] {\!$E_{Dx}$}(138.3750,944.3622);
  \path[draw=black,line join=miter,line cap=butt,line width=0.650pt] (142.0000,940.3622) .. controls (151.7500,931.8622) and (160.7095,932.3622) .. 
  node[above left=0.12cm,at end] {$E_{Dy}$\!}(220.0000,932.3622);
  \path[draw=black,line join=miter,line cap=butt,line width=0.650pt] (142.0732,944.3836) .. controls (149.7506,952.0610) and (161.2350,952.3622) .. 
  node[above left=0.12cm,at end] {$\f C(EDx,\!EDy)$\!}(220.0000,952.3622);
  \path[fill=black] (140.00002,942.36224) node[circle, draw, line width=0.65pt, minimum width=5mm, fill=white, inner sep=0.25mm] (text5446-2-1) {$E$
   };
\end{tikzpicture} \ \rlap{ .}}} \end{equation*}
Given $\C$-transformations $\varsigma\colon E \Rightarrow E'$ and $\vartheta\colon D \Rightarrow D'$, we define the whiskerings $\varsigma D\colon E D \Rightarrow E' D$ and  $E \vartheta \colon E D \Rightarrow E D'$ to have respective components
\begin{equation*}
\vcenter{\hbox{\begin{tikzpicture}[y=0.8pt, x=0.9pt,yscale=-1, inner sep=0pt, outer sep=0pt, every text node part/.style={font=\tiny} ]
  \path[draw=black,line join=miter,line cap=butt,line width=0.650pt] (142.0000,940.3622) .. controls (151.7500,931.8622) and (160.7095,932.3622) .. 
  node[above left=0.12cm,at end] {$E'_{Dx}$\!}(220.0000,932.3622);
  \path[draw=black,line join=miter,line cap=butt,line width=0.650pt] (142.0732,944.3836) .. controls (149.7506,952.0610) and (161.2350,952.3622) .. 
  node[above left=0.12cm,at end] {$\f C(EDx,\!E'Dx)$\!}(220.0000,952.3622);
  \path[draw=black,line join=miter,line cap=butt,line width=0.650pt] (100.0000,942.3622) -- 
  node[above right=0.12cm,at start] {\!$E_{Dx}$}(138.1317,942.3622);
  \path[draw=black,line join=miter,line cap=butt,line width=0.650pt] (100.0000,922.3622) .. controls (149.8510,922.3622) and (140.4005,912.3622) .. 
  node[above right=0.12cm,at start] {\!$D_x$}node[above left=0.12cm,at end] {$D_x$\!}(220.0000,912.3622);
  \path[fill=black] (140.00002,942.36224) node[circle, draw, line width=0.65pt, minimum width=5mm, fill=white, inner sep=0.25mm] (text5446-2-1) {$\varsigma$
     };
\end{tikzpicture}}}\qquad \text{and}\qquad 
\vcenter{\hbox{\begin{tikzpicture}[y=0.8pt, x=0.9pt,yscale=-1, inner sep=0pt, outer sep=0pt, every text node part/.style={font=\tiny} ]
  \path[draw=black,line join=miter,line cap=butt,line width=0.650pt] (102.0000,920.3622) .. controls (111.7500,911.8622) and (160.7095,912.3622) .. 
  node[above left=0.12cm,at end] {$D'_{x}$\!}(220.0000,912.3622);
  \path[draw=black,line join=miter,line cap=butt,line width=0.650pt] (102.0732,924.3836) .. controls (111.5184,933.8288) and (127.5236,930.1389) ..
  node[above right=0.04cm,pos=0.42] {$\f B(Dx,\!D'x)$}(137.7019,940.3171);
  \path[draw=black,line join=miter,line cap=butt,line width=0.650pt] (60.0000,942.3622) .. controls (90.0000,942.3622) and (130.0505,952.6866) .. 
  node[above right=0.12cm,at start] {\!$E_{Dx}$}(138.3750,944.3622);
  \path[draw=black,line join=miter,line cap=butt,line width=0.650pt] (142.0000,940.3622) .. controls (151.7500,931.8622) and (160.7095,932.3622) .. 
  node[above left=0.12cm,at end] {$E_{D'x}$\!}(220.0000,932.3622);
  \path[draw=black,line join=miter,line cap=butt,line width=0.650pt] (142.0732,944.3836) .. controls (149.7506,952.0610) and (161.2350,952.3622) .. 
  node[above left=0.12cm,at end] {$\f C(EDx,\!ED'x)$\!}(220.0000,952.3622);
  \path[draw=black,line join=miter,line cap=butt,line width=0.650pt] (60.0000,922.3622) -- 
  node[above right=0.12cm,at start] {\!$D_x$}(97.2478,922.3622);
  \path[fill=black] (100,922.36218) node[circle, draw, line width=0.65pt, minimum width=5mm, fill=white, inner sep=0.25mm] (text5446-2) {$\vartheta$
     };
  \path[fill=black] (140.00002,942.36224) node[circle, draw, line width=0.65pt, minimum width=5mm, fill=white, inner sep=0.25mm] (text5446-2-1) {$E$
     };
\end{tikzpicture} \ \rlap{ .}}} 
\end{equation*}

It is easy to verify that these whiskering operations satisfy the middle-four interchange law, so yielding a composition functor $\cat{Cat}_\kappa(\C)(\f B, \f C) \times \cat{Cat}_\kappa(\C)(\f A, \f B) \to \cat{Cat}_\kappa(\C)(\f A, \f C)$. The identity $\C$-functor on a $\C$-category $\f A$ is the identity on objects, has $1$-cell components $1_{\epsilon x} \colon \epsilon x \to \epsilon x$, and $2$-cell components obtained from left and right unit constraints in $\C\lo$.

It remains to give the associativity and unitality constraints making  $\cat{Cat}_\kappa(\C)$ into a bicategory. Given a triple $(C,D,E)$ of composable $\C$-functors, we observe that the composites $(ED)C$ and $E(DC)$ agree on objects, and differ in their $1$-cell components only by associativity constraints in $\C\ti$; these constraints constitute the components of an invertible $\C$-icon $(ED)C \Rightarrow E(DC)$, from which we induce the $\C$-transformation giving the associativity constraint for $\cat{Cat}_\kappa(\C)$ at $(C,D,E)$. The left and right unit constraints are obtained similarly from left and right unit constraints in $\C\ti$; the pentagon and unit axioms are now easily verified using the corresponding axioms in $\C\ti$ together with the closure of $\C$-icons under composition. This completes the definition of the bicategory $\cat{Cat}_\kappa(\C)$. Note that if $\C\ti$ is a $2$-category---which it will be in many of our examples---then so too is $\cat{Cat}_\kappa(\C)$.

As in Section~\ref{sec:Cmodules}, we have for each object $v \in \C\ti$ the $\kappa$-ary $\C$-category $\hat v$ with one object $\star$ such that $\epsilon(\star) = v$ and $\hat v(v,v) = 1_v$. Each tight morphism $f \colon v \to w$ induces a $\C$-functor $\hat f \colon \hat v \to \hat w$ that is the identity on objects, with $1$-cell component  $\hat f_\star = f$ and $2$-cell component $\hat f_{\star \star}$ built from unit constraints in $\C\lo$. Similarly, each $2$-cell $\alpha \colon f \Rightarrow g$ of $\C\ti$ induces a $\C$-icon $\hat f \Rightarrow \hat g$ with unique component $\alpha$, whence a $\C$-transformation $\hat \alpha \colon \hat f \Rightarrow \hat g$. These data assemble to give a functor $\hat{(\thg)} \colon \C\ti \to \cat{Cat}_\kappa(\C)$.

\begin{Ex}\label{ex:enriched}
  Let $\C$ be the inchordate equipment on the delooping of a monoidal category $\mathbf{V}$; then $\cat{Cat}_\kappa(\C)$ is the $2$-category of ($\kappa$-small) $\mathbf{V}$-enriched categories, $\mathbf{V}$-functors and $\mathbf{V}$-transformations. More generally, if $\C$ is the inchordate equipment on any bicategory, then the notions of  $\C$-category, $\C$-functor and $\C$-transformation as defined above coincide with those defined for bicategories in, for example,~\cite{Street2005Enriched}. 
  \end{Ex}

\begin{Ex}\label{ex:internal}
  Let $\mathbf{C}$ be a category with pullbacks, and consider the equipment $\fcat{Span}(\mathbf{C})$ given by the functor $\cat C \to \cat{Span}(\mathbf{C})$ that sends $f \colon x \to y$ to the span $1_x \colon x \leftarrow x \rightarrow y \colon f$.
 Then $\{1\}$-ary $\fcat{Span}(\mathbf{C})$-categories are internal categories in $\mathbf C$, as in Example~\ref{ex:monad}, and $\cat{Cat}_{\{1\}}(\fcat{Span}(\mathbf{C}))$ is the $2$-category of internal categories, functors and transformations in $\mathbf C$. The embedding $\hat{(\thg)}$ is the ``discrete internal category'' functor.
\end{Ex}

\begin{Ex}\label{ex:mealy}
  On the other hand, if $\C$ is the \emph{chordate} equipment on the delooping of a monoidal category $\mathbf{V}$, then $\cat{Cat}_\kappa(\C)$ is the bicategory of $\mathbf{V}$-categories, \emph{Mealy morphisms}, and Mealy cells as described in~\cite{Pare2012Mealy}.
  We may extend this terminology more generally to any chordate equipment.
  For instance, if $\C$ is the chordate equipment on the bicategory of spans in a category $\mathbf{C}$ with pullbacks, then the objects of $\cat{Cat}_{\{1\}}(\C)$ are internal categories as in Example~\ref{ex:internal}, while its morphisms can be identified with spans of internal categories whose source leg is a discrete opfibration.
\end{Ex}

\subsection{The equipment of $\C$-categories}
\label{sec:equip-Ccats}

We call an equipment $\C$ \emph{locally $\kappa$-cocomplete} if $\C\lo$ is so. In this case, we have as in Section~\ref{sec:Cmodules}, a locally $\kappa$-cocomplete bicategory $\cat{Mod}_\kappa(\C\lo)$ whose objects are also $\kappa$-ary $\C$-categories, but whose morphisms are modules.
We will combine $\cat{Cat}_\kappa(\C)$ and $\cat{Mod}_\kappa(\C\lo)$ into a (locally $\kappa$-cocomplete) equipment.

Firstly, if $D\colon \f A \to \f B$ is a $\C$-functor, we have an $\f A$-$\f B$-module with $1$-cell components $\f B(u,D x) \circ D_x \in \C\lo(\epsilon x,\epsilon u)$, and  with $\f A$- and $\f B$-actions given by
\begin{equation*}
\vcenter{\hbox{\begin{tikzpicture}[y=0.8pt, x=0.9pt,xscale=-1, inner sep=0pt, outer sep=0pt, every text node part/.style={font=\tiny} ]
  \path[draw=black,line join=miter,line cap=butt,line width=0.650pt] (102.0000,920.3622) .. controls (111.7500,911.8622) and (160.7095,912.3622) .. 
  node[above right=0.12cm,at end] {\!$\f B(u,\!Dx)$}(180.0000,912.3622);
  \path[draw=black,line join=miter,line cap=butt,line width=0.650pt] (102.0732,924.3836) .. controls (111.5184,933.8288) and (127.5236,930.1389) ..
  node[above right=0.04cm,pos=0.55] {$\f B(Dx,\!Dy)$}(137.7019,940.3171);
  \path[draw=black,line join=miter,line cap=butt,line width=0.650pt] (35.0000,942.3622) .. controls (90.0000,942.3622) and (130.0505,952.6866) .. 
  node[above left=0.12cm,at start] {$D_y$\!}(138.3750,944.3622);
  \path[draw=black,line join=miter,line cap=butt,line width=0.650pt] (142.0000,940.3622) .. controls (151.7500,931.8622) and (160.7095,932.3622) .. 
  node[above right=0.12cm,at end] {\!$D_{x}$}(180.0000,932.3622);
  \path[draw=black,line join=miter,line cap=butt,line width=0.650pt] (142.0732,944.3836) .. controls (149.7506,952.0610) and (161.2350,952.3622) .. 
  node[above right=0.12cm,at end] {\!$\f A(x,\!y)$}(180.0000,952.3622);
  \path[draw=black,line join=miter,line cap=butt,line width=0.650pt] (35.0000,922.3622) -- 
  node[above left=0.12cm,at start] {$\f B(u,\!Dy)$\!}(97.2478,922.3622);
  \path[fill=black] (100,922.36218) node[circle, draw, line width=0.65pt, minimum width=5mm, fill=white, inner sep=0.25mm] (text5446-2) {$m$
     };
  \path[fill=black] (140.00002,942.36224) node[circle, draw, line width=0.65pt, minimum width=5mm, fill=white, inner sep=0.25mm] (text5446-2-1) {$D$
     };
\end{tikzpicture}}}\qquad \text{and}\qquad 
\vcenter{\hbox{\begin{tikzpicture}[y=0.8pt, x=0.9pt,xscale=-1,yscale=-1, inner sep=0pt, outer sep=0pt, every text node part/.style={font=\tiny} ]
  \path[draw=black,line join=miter,line cap=butt,line width=0.650pt] (142.0000,940.3622) .. controls (151.7500,931.8622) and (160.7095,932.3622) .. 
  node[above right=0.12cm,at end] {\!$\f B(v,\!Dx)$}(190.0000,932.3622);
  \path[draw=black,line join=miter,line cap=butt,line width=0.650pt] (142.0000,944.3836) .. controls (149.7506,952.0610) and (161.2350,952.3622) .. 
  node[above right=0.12cm,at end] {\!$\f B(u,\!v)$}(190.0000,952.3622);
  \path[draw=black,line join=miter,line cap=butt,line width=0.650pt] (87.0000,942.3622) -- 
  node[above left=0.12cm,at start] {$\f B(u,\!Dx)$\!}(138.1317,942.3622);
  \path[draw=black,line join=miter,line cap=butt,line width=0.650pt] (87.0000,922.3622) .. controls (149.8510,922.3622) and (140.4005,912.3622) .. 
  node[above left=0.12cm,at start] {$D_x$\!}node[above right=0.12cm,at end] {\!$D_x$}(190.0000,912.3622);
  \path[fill=black] (140.00002,942.36224) node[circle, draw, line width=0.65pt, minimum width=5mm, fill=white, inner sep=0.25mm] (text5446-2-1) {$\varsigma$
     };
\end{tikzpicture}}} 
\end{equation*}
respectively.
We denote this $\f A$-$\f B$-module by $\f B(1,D)$, and call a module \emph{representable} when it is isomorphic to one of this form.
Now if $E\colon \f A \to \f B$ is another $\C$-functor and $\vartheta\colon D\Rightarrow E$ is a $\C$-transformation, then we have an $\f A$-$\f B$-module morphism $\f B(1,\vartheta)\colon \f B(1,D) \to \f B(1,E)$ with components
\begin{equation*}
\vcenter{\hbox{\begin{tikzpicture}[y=1pt, x=1.4pt,yscale=-1, inner sep=0pt, outer sep=0pt, every text node part/.style={font=\tiny} ]
  \path[draw=black,line join=miter,line cap=butt,line width=0.650pt] (104.7730,922.3622) -- 
  node[above left=0.12cm,at end] {$E_x$\!}(180.0000,922.3622);
  \path[draw=black,line join=miter,line cap=butt,line width=0.650pt] (102.0732,924.3836) .. controls (111.5184,933.8288) and (127.5236,930.1389) .. 
  node[above right=0.04cm,pos=0.65] {$\f B(Dx,\!Ex)$}(137.7019,940.3171);
  \path[draw=black,line join=miter,line cap=butt,line width=0.650pt] (65.0000,922.3622) -- 
  node[above right=0.12cm,at start] {\!$D_x$}(97.2478,922.3622);
  \path[draw=black,line join=miter,line cap=butt,line width=0.650pt] (65.0000,942.3622) -- 
  node[above right=0.12cm,at start] {\!$\f B(u,\!Dx)$}(136.0000,942.3622);
  \path[draw=black,line join=miter,line cap=butt,line width=0.650pt] (144.0000,942.3622) -- 
  node[above left=0.12cm,at end] {$\f B(u,\!Ex)$\!}(180.0000,942.3622);
  \path[fill=black] (100,922.36218) node[circle, draw, line width=0.65pt, minimum width=5mm, fill=white, inner sep=0.25mm] (text5446-2) {$\vartheta$
     };
  \path[fill=black] (140.00002,942.36224) node[circle, draw, line width=0.65pt, minimum width=5mm, fill=white, inner sep=0.25mm] (text5446-2-1) {$m$
     };
\end{tikzpicture}\ \rlap{ .}}} 
\end{equation*}
On the other hand, any $\f A$-$\f B$-module morphism $f\colon \f B(1,D) \to \f B(1,E)$ induces a $\C$-transformation $\vartheta\colon D\Rightarrow E$ with components
\begin{equation*}
\vcenter{\hbox{\begin{tikzpicture}[y=0.85pt, x=1pt, inner sep=0pt, outer sep=0pt, every text node part/.style={font=\tiny} ]
  \path[draw=black,line join=miter,line cap=butt,line width=0.650pt] (0.0000,982.3622) .. controls (70.0000,982.3622) and (83.1755,982.6866) .. 
  node[above right=0.12cm,at start] {\!$D_x$}(91.5000,974.3622);
  \path[draw=black,line join=miter,line cap=butt,line width=0.650pt] (30.0000,962.3622) .. controls (70.2500,962.3622) and (82.6241,962.2363) .. 
  node[above=0.08cm,pos=0.28] {$\f B(Dx,\!Dx)$}(91.5000,971.1122);
  \path[draw=black,line join=miter,line cap=butt,line width=0.650pt] (94.7500,970.8622) .. controls (104.5000,962.3622) and (120.7095,962.3622) .. 
  node[above left=0.12cm,at end] {$\f B(Dx,\!Ex)$\!}(160.0000,962.3622);
  \path[draw=black,line join=miter,line cap=butt,line width=0.650pt] (94.9482,974.5086) .. controls (102.6256,982.1860) and (121.2350,982.3622) .. 
  node[above left=0.12cm,at end] {$E_x$\!}(160.0000,982.3622);
  \path[fill=black] (26.75,962.36218) node[circle, draw, line width=0.65pt, minimum width=5mm, fill=white, inner sep=0.25mm] (text5442) {$j$
     };
  \path[fill=black] (93,973.11218) node[circle, draw, line width=0.65pt, minimum width=5mm, fill=white, inner sep=0.25mm] (text5446) {$f$
   };
\end{tikzpicture}\rlap{ .}}}
\end{equation*}
These two operations are mutually inverse and respect composition, whence we have for each $\f A$ and $\f B$ a fully faithful functor $\f B(1, \thg) \colon \cat{Cat}_\kappa(\C)(\f A,\f B) \to \cat{Mod}_\kappa(\C\lo)(\f A,\f B)$.

It is clear that the image $\f B(1, 1_\f B)$ of the identity $\C$-functor $1_{\f B} \colon \f B \to \f B$ is isomorphic, via unit constraints in $\C\lo$, to the identity module on $\f B$. Moreover, if we have $\C$-functors $D\colon \f A \to \f B$ and $E\colon \f B \to \f C$, then there is a split coequaliser diagram
\vskip0.3\baselineskip
\begin{equation*}
  \xymatrix@-.5pc{
    \textstyle \sum_{u,v}^{\phantom{\f B}} \scriptstyle
    \f C(w,E u) \circ E_u \circ \f B(u,v) \circ  \f B(v,D x) \circ D_x
    \ar@<1mm>[r] \ar@<-1mm>[r]
    &
    \textstyle\sum_{u}^{\phantom{\f B}} \scriptstyle
    \f C(w,E u) \circ E_u \circ \f B(u,D x) \circ D_x
    \ar[r] \ar@/_24pt/[l]
    &
    \scriptstyle
    \f C(w,E D x) \circ E_{D x} \circ D_x
    \ar@/_24pt/[l]
  }\rlap{ ,}
\end{equation*}
wherein both leftwards-pointing arrows factor through the coprojection labeled by $u\coloneqq D x$ via the morphism $j\colon 1_{\epsilon(D x)} \to \f B(D x,D x)$.
It follows that $\f C(1,E D)$, whose $1$-cell components are the target of this split coequaliser, is isomorphic to the composite $\f C(1,E) \circ \f B(1,D)$ in $\cat{Mod}_\kappa(\C\lo)$. We have now described the action on homs and the  nullary and binary functoriality constraints of a functor $J_\kappa \colon \cat{Cat}_\kappa(\C) \to \cat{Mod}_\kappa(\C\lo)$; the corresponding coherence axioms are straightforward to verify. Clearly, $J_\kappa$ is locally fully faithful and bijective on objects, and, since $\C\lo$ is locally $\kappa$-cocomplete, so also is $\cat{Mod}_\kappa(\C\lo)$; and so we have described a locally $\kappa$-cocomplete equipment, which we denote by $\fcat{Mod}_\kappa(\C)$.

Given a tight morphism $f \colon v \to w$, we have the $\C$-functor $\hat f \colon \hat v \to \hat w$ inducing the $\hat v$-$\hat w$-module $\hat w(1, \hat f)$ with $\hat w(1, \hat f)(\star, \star) = 1_w \circ f$. On the other hand, the loose morphism underlying $f$ induces directly the isomorphic $\hat v$-$\hat w$-module $\hat f$ with $\hat f(\star, \star) = f$, and so we have a morphism of equipments $\hat{(\thg)} \colon \C \to \fcat{Mod}_\kappa(\C)$ as displayed in the square
\begin{equation*}
\cd{
 \C\ti \ar[d]_{J_\C} \ar[r]^-{\hat{(\thg)}} \twocong{dr} &
 \cat{Cat}_\kappa(\C) \ar[d]^{J_\kappa} \\
 \C\lo \ar[r]_-{\hat{(\thg)}} &
 \cat{Mod}_\kappa(\C_\lambda)\rlap{ .}
}
\end{equation*}

\subsection{Locally $\kappa$-cocomplete equipments via enrichment}
\label{subsec:equip-enrichment}

We shall now show that the equipment morphism $\hat{(\thg)} \colon \C \to \fcat{Mod}_\kappa(\C)$ exhibits $\fcat{Mod}_\kappa(\C)$ as the free completion of $\C$ under a suitable class of enriched bicolimits. In this section, we construct the base for enrichment: a monoidal bicategory $\F_\kappa$ such that $\F_\kappa$-bicategories are locally $\kappa$-cocomplete equipments. 

We begin by applying Theorem~\ref{thm:comma} to the forgetful functor $R\colon \cat{Colim}_\kappa \to \cat{Cat}$ from Section~\ref{subsec:loc-cocplt}.
This yields a complete and cocomplete monoidal bicategory $(\cat{Cat}\downarrow R)$, whose objects $A$ are functors $j_A\colon A\ti\to A\lo$ such that $A\lo$ is $\kappa$-cocomplete.
Now let $\F_\kappa$ denote the full subcategory of $(\cat{Cat}\downarrow R)$ on the fully faithful functors.
Such functors are the right class of the (essentially surjective on objects, fully faithful) factorization system on $\cat{Cat}$; thus $\F_\kappa$ is a reflective sub-bicategory of $(\cat{Cat}\downarrow R)$ and hence complete and cocomplete.
Bilimits in $\F_\kappa$ are constructed as in $(\cat{Cat}\downarrow R)$, which is to say that we take bilimits of the tight and loose parts separately (with limits in $\cat{Colim}_\kappa$, as usual, created by $R$).
Bicolimits in $\F_\kappa$ are computed by taking bicolimits of the tight and loose parts separately, then applying the reflector, which is the (essentially surjective, fully faithful) factorization.

Now the class $^{\perp}\F_\kappa$ of morphisms sent to equivalences by the reflector $(\cat{Cat} \downarrow R) \to \F_\kappa$ comprises all morphisms
\[\vcenter{\xymatrix{
    A\ti\ar[r]^{f\ti}\ar[d]_j \ar@{}[dr]|\cong &
    B\ti\ar[d]^j\\
    A\lo\ar[r]_{f\lo} &
    B\lo
  }}\]
for which $f\lo$ is an equivalence in $\cat{Colim}_\kappa$, and for any $b\in B\ti$ there exists $a\in A\ti$ such that $f\lo(j(a)) \cong j(b)$.
$^{\perp}\F_\kappa$ is closed under the monoidal product of $(\cat{Cat}\downarrow R)$, and so by Theorem~\ref{thm:reflect}, $\F_\kappa$ inherits a monoidal structure from $(\cat{Cat} \downarrow R)$ which is both symmetric and closed.
Its right hom can be computed as the pullback:
\[\xymatrix@C=3pc{
[B,C]\ti \ar[r]^{j_{[B,C]}} \ar[dd] & [B,C]\lo \ar@{=}[d] \\
& [B\lo,C\lo]_\kappa \ar[d]^{[j_B,1]} \\
[B\ti,C\ti] \ar[r]_{[1,j_C]} & [B\ti,C\lo]\rlap{ .}
}\]
Here $[B\lo,C\lo]_\kappa$ denotes the right hom of $\cat{Colim}_\kappa$, i.e.\ the category of $\kappa$-cocontinuous functors $B\lo \to C\lo$.

We have forgetful functors $U\ti \colon \F_\kappa \to \cat{Cat}$ and $U\lo \colon \F_\kappa \to \cat{Colim}_\kappa$, of which $U\ti$ is lax monoidal and $U\lo$ is strong monoidal.
By Theorems~\ref{thm:comma} and~\ref{thm:reflect}, to give an $\F_\kappa$-bicategory is equally to give a bicategory $\C\ti$, a locally $\kappa$-cocomplete bicategory $\C\lo$, and an identity-on-objects and locally fully faithful functor $J_\C\colon \C\ti \to \C\lo$; in other words, a locally $\kappa$-cocomplete equipment. Similarly, to give an $\F_\kappa$-functor $\C \to \D$ is equally to give a morphism of equipments whose loose part is locally $\kappa$-cocontinuous; and so on.

We may repeat the above arguments replacing $R \colon \cat{Colim}_\kappa \to \cat{Cat}$ by the identity functor $\cat{Cat} \to \cat{Cat}$; on doing so, we obtain a monoidal bicategory $\F$ such that $\F$-bicategories are equipments subject to no local cocompleteness requirements. 
This is just the bicategorical version of the $\F$-categories of~\cite{ls:limlax}; our bicategory $\F$ is not identical to their 2-category $\F$, but is chosen to make our $\F$-bicategories related to their $\F$-categories in the same way that ordinary bicategories are related to strict 2-categories.
There is an evident lax monoidal forgetful functor $U\colon \F_\kappa \to \F$, such that $U$ computes the underlying $\F$-bicategory of an  $\F_\kappa$-bicategory; and in fact, this $U$ forms part of a monoidal biadjunction.
\begin{Prop}\label{prop:free-kappa-equipment}
  The forgetful functor $U\colon \F_\kappa \to \F$ admits a strong monoidal left adjoint $H$.
\end{Prop}
\begin{proof}
  We define $H$ applied to a fully faithful functor $A\ti\to A\lo$ to be the composite $A\ti \to A\lo \xrightarrow{\eta} R L A\lo$, where $L\dashv R$ denotes the adjunction $\cat{Cat} \rightleftarrows \cat{Colim}_\kappa$ as usual.
  Since $RLA\lo$ can be identified, up to equivalence, with the closure of $A\lo$ in its presheaf category under the appropriate colimits, the adjunction unit $\eta$ is fully faithful and thus $H$ lands in $\F_\kappa$.
  The adjunction $H \dashv U$ is easy to check.

  Finally, given $A\ti\to A\lo$ and $B\ti\to B\lo$, we have $L A\lo \otimes L B\lo \cong L(A\lo \times B\lo)$ since $L$ is strong monoidal.
  Thus $A\ti\times B\ti \to R (L A\lo \otimes L B\lo)$ is already fully faithful, and hence is $H A \otimes H B$, which is therefore equivalent to $H (A\times B)$.
\end{proof}
Since $\F$ and $\F_\kappa$ are closed and complete and $U\colon \F_\kappa \to \F$ creates bilimits, it follows from Theorem~\ref{thm:adjoints} that any $\F$-weight $W$ gives rise to an $\F_\kappa$-weight $H W$, with the property that an $H W$-weighted colimit in an $\F_\kappa$-bicategory $\C$ is nothing more or less than a $W$-weighted colimit in the underlying $\F$-bicategory of $\C$. We shall use this fact shortly.

\subsection{Tight collages}
\label{subsec:tightcollages}
We will now  describe the kinds of $\F_\kappa$-enriched colimits under which the locally $\kappa$-cocomplete equipment $\C$ is to be completed to obtain $\fcat{Mod}_\kappa(\C)$. As before, we first give an elementary presentation of these colimits. Given a locally $\kappa$-cocomplete equipment and a $\kappa$-ary $\C$-category $\f A$, a \emph{tight collage} of $\f A$ is 
given by an object $v \in \C$ that birepresents both functors
\[ \cat{Cat}_\kappa(\C)({\f{A}},\hat{(\thg)}) \colon \C\ti \to \cat{Cat} \qquad \text{and} \qquad
\cat{Mod}_\kappa(\C\lo)({\f{A}},\hat{(\thg)}) \colon \C\lo \to \cat{Cat}\]
in a compatible manner. By this, we mean that there are given universal elements $T\ti \in \cat{Cat}_\kappa(\C)({\f{A}},\hat{v})$ and $T\lo \in \cat{Mod}_\kappa(\C\lo)({\f{A}},\hat{v})$---thus in particular, $T\lo$ exhibits $v$ as the collage of $\f A$ in $\C\lo$---such that $T\ti$ is a tightening of $T\lo$.
Observe that a tight collage $(T\ti, T\lo)$ is determined to within unique isomorphism by $T\lo$ for which there exists some such $T\ti$.
More precisely, a collage $T\lo$ for $\f A$ forms part of a tight collage if and only if:
\begin{enumerate}[(i)]
\item $T\lo$ admits a tightening;
\item $T\lo$ \emph{detects tightness}: which is to say that, given $f \colon v \to w$ in $\C\lo$ and a diagram
\begin{equation*}
\cd{
 \f A \ar[d]_{T\lo} \ar[r]^U \twocong[0.33]{dr} & \hat w \\
 \hat v \ar[ur]_{\hat f} & {}
}
\end{equation*}
in $\cat{Mod}_\kappa(\C)$, if $U$ admits a tightening, then so does $f$.
\end{enumerate}

From the definitions of $\C$-functor and $\C$-module, we see that tightenings of an $\f A$-$\hat v$-module $T$ are the same as tightenings for each coprojection $T(\star, x) \colon \epsilon x \to v$ in $\C\lo$.
Hence, arguing as in Section~\ref{subsec:collages}, a tight collage for $\f A$ is equally a \emph{tight lax colimit} for $\f A \colon \nabla_O \to \C\lo$, amounting to a universal oplax cocone $T \colon \f A \Rightarrow \Delta v$ in $\C\lo$ whose components come equipped with tightenings and  jointly detect tightness.  It follows that the notion of tight collage, which was defined only for locally $\kappa$-cocomplete equipments, in fact makes sense for any equipment.

The special cases are again worth mentioning:

\begin{Ex}
  If $\ob{\f{A}}$ is a singleton, so that ${\f{A}}$ is a (loose) monad as in Example~\ref{ex:kleisli}, then a tight collage of ${\f{A}}$ is determined by a (loose) Kleisli object whose coprojection admits a tightening and detects tightness.
\end{Ex}

\begin{Ex}
  If $O$ is a $\kappa$-small set and ${\f{A}}$ is constructed from an $O$-indexed set of objects in $\C$, as in Example~\ref{ex:coproducts}, then a tight collage of ${\f{A}}$ is simply a coproduct in $\C\ti$ that is preserved by $J_\C \colon \C\ti\to\C\lo$.
\end{Ex}

\begin{Rk}
  Let $\C$ be a locally $\kappa$-cocomplete bicategory, and let $\tilde{\C}$ be the locally $\kappa$-cocomplete equipment with $\tilde{\C}\lo = \C$ and $\tilde{\C}\ti$ its locally full sub-bicategory consisting exactly of the maps.
  Then the first halves of Theorem~\ref{thm:collages2}(\ref{item:c4}) and~(\ref{item:c5}) say exactly that all collages in $\C$ underlie tight collages in $\tilde{\C}$.
  Conversely, suppose that $\C$ is a map equipment (recall that this means that every tight morphism is a map).
  Then a tight collage in $\C$ automatically satisfies the first half of Theorem~\ref{thm:collages2}(\ref{item:c4}) and an $\F$-analogue of the first half of~(\ref{item:c5}).
  In~\cite{wood:proarrows-ii}, the existence of (finite) coproducts and Kleisli objects with these properties were taken as additional axioms for proarrow equipments.
\end{Rk}

Inspecting the proof of our Theorem~\ref{thm:collages}(\ref{item:c1}) in~\cite[Proposition 2.2(a)]{street:cauchy-enr}, we see that if it is carried out using tight collages, then the resulting lax colimits will also be tight.
This proves:
\begin{Thm}\label{thm:tightcollages}
  If $\C$ is a locally $\kappa$-cocomplete equipment, then the following are equivalent.
  \begin{enumerate}
  \item $\C$ admits tight lax colimits of lax functors into $\C\lo$ with $\kappa$-small domain.\label{item:c1t}
  \item $\C$ admits tight $\kappa$-ary collages.\label{item:c2t}
  \item $\C$ admits tight Kleisli objects and tight $\kappa$-small coproducts.\label{item:c3t}
  \end{enumerate}
\end{Thm}

\subsection{Tight collages as $\F_\kappa$-colimits}

We now exhibit tight collages as $\F_\kappa$-weighted colimits. We do so by first exhibiting them as $\F$-weighted colimits and then transporting across the monoidal biadjunction $\F \rightleftarrows \F_\kappa$. To do that,  we need to examine the notion of $\F$-weighted colimit in more detail.

The same argument that identifies an $\F$-bicategory $\B$ with a functor $J_\B \colon \B\ti \to \B\lo$ shows that a right $\B$-module $W$ can be identified with a right $\B\ti$-module $W\ti$, a right $\B\lo$-module $W\lo$, and a pointwise fully faithful $\B\ti$-module morphism $J_W \colon W\ti \to W\lo(J_\B)$.
If $\C$ is another $\F$-bicategory, $F\colon \B\to \C$ an $\F$-functor, and $v\in\ob\C$, then a cylinder $\phi \colon W \to \C(F, v)$ consists of transformations
$\phi\ti\colon W\ti \to \C\ti(F\ti,v)$ and $\phi\lo\colon W\lo \to \C\lo(F\lo,v)$ together with an invertible modification
\begin{equation*}
\cd{
W\ti \ar[d]_{J_W} \ar[r]^-{\phi\ti} \dtwocell{dr}{J_\phi} &
\C\ti(F\ti, v) \ar[d]^{\C\lo(J_F, 1) \circ J_\C} \\
W\lo(J_\B) \ar[r]_-{\phi\lo} &
\C\lo(F\lo J_\B, v)\rlap{ .}
}
\end{equation*}
Since $J_\C$ is locally fully faithful, such a cylinder is determined up to unique isomorphism by $\phi\lo$ for which there exists some such $\phi\ti$.
Thus, to give $\phi$ is equally to give a transformation $\phi\lo\colon W\lo \to \C\lo(F\lo,v)$ such that for every $a \in W\ti(x)$, the morphism $\phi(a) \colon F(x) \to v$ comes equipped with a chosen tightening. Thus, arguing as in Proposition~3.6 of~\cite{ls:limlax}, we conclude that:
\begin{Prop}\label{prop:Fcolimits}
  For $\B,\C,W,F,v$ as above, a $W$-weighted colimit of $F$ is given by a $W\lo$-weighted colimit $\phi\colon W\lo \to \C\lo(F\lo,v)$ of $F\lo$ such that
 the morphisms
\[
\{\phi(a) \colon F(x) \to v \,\mid\, x \in \ob \B \text{ and } a \in W\ti(x)\}
\] 
come equipped with tightenings and jointly detect tightness.
\end{Prop}

Now given $O$ a $\kappa$-small set, let $O^\diamond$ be as in Section~\ref{subsec:collages}, and let $\O$ be the inchordate equipment on $O^\diamond$.
In $\O$, only identity morphisms are tight, and so any functor $O^\diamond \to \C\lo$ preserves tightness, and as such, admits an essentially-unique extension to an $\F$-functor $\O \to \C$.
Let $V_O$ be the right $O^\diamond$-module from Section~\ref{subsec:collages} such that $V_O$-weighted colimits are collages.
For each $x\in O$ there is a distinguished ``generating'' element $p_x \in V_O(x)$, such that in a $V_O$-weighted cylinder the images of the $p_x$'s are the coprojections.
From these, we obtain a morphism of right $\O_\tau$-modules
\begin{equation*}
\spn{p_x} \colon \sum_{x \in O} \O_\tau(\thg, x) \to J_\O^\ast(V_O)\rlap{ ,}
\end{equation*}
and taking the second half of the pointwise (bijective on objects, fully faithful) factorisation of this, we obtain a right $\O$-module $\tilde V_O$. By Proposition~\ref{prop:Fcolimits}, a $\tilde V_O$-weighted colimit of $\O \to \C$ is precisely a tight lax colimit of the corresponding $\f A \colon \nabla_O \to \C\lo$.

Letting $\tilde{W}_O = H (\tilde{V}_O)$, where $H$ is the free-$\kappa$-cocomplete-equipment functor from Proposition~\ref{prop:free-kappa-equipment}, it follows from Theorems~\ref{thm:tightcollages} and~\ref{thm:adjoints}
that a locally $\kappa$-cocomplete equipment admits tight collages of $\kappa$-ary categories if and only if, when regarded as a $\F_\kappa$-bicategory, it is $\tilde \Phi_\kappa$-cocomplete for the class of weights $\tilde\Phi_\kappa = \{\tilde W_O \mid \text{$O$ a $\kappa$-small set}\}$. 

Our aim now is to show that $\fcat{Mod}_\kappa(\C) \simeq \tilde\Phi_\kappa(\C)$. First we need an $\F$-analogue of Section~\ref{subsec:coll-cocplt}.

\subsection{Tight collages in $\F_\kappa$}
\label{sec:tcoll-Fk}
  Since $\F_\kappa$ is a right closed monoidal bicategory, we can regard it as enriched over itself, hence as a locally $\kappa$-cocontinuous equipment.
  A loose morphism in $\F_\kappa$ from $A\ti\to A\lo$ to $B\ti\to B\lo$ is a $\kappa$-cocontinuous functor $A\lo\to B\lo$, while a tight morphism is a pseudo-commutative square
  \[\vcenter{\xymatrix@-.5pc{
      A\ti\ar[r]\ar[d] \ar@{}[dr]|\cong &
      B\ti\ar[d]\\
      A\lo\ar[r] &
      B\lo
    }}\]
  in which the lower arrow is $\kappa$-cocontinuous.
Regarding $\F_\kappa$ as an equipment in this way, let ${\f{A}}$ be an $\F_\kappa$-category.
Thus, it consists of objects $x$ with extents $\epsilon x \in\F_\kappa$, together with $\kappa$-cocontinuous functors ${\f{A}}(x,y)\colon (\epsilon y)\lo \to (\epsilon x)\lo$ making its loose part ${\f{A}}\lo$ into a $\cat{Colim}_\kappa$-category as in Section~\ref{subsec:coll-cocplt}.
Then its tight collage $\Vert{\f{A}}\Vert$ is, in particular, a collage in $(\F_\kappa)\lo$.
Since $(\F_\kappa)\lo$ is equivalent to $\cat{Colim}_\kappa$, the loose part $\Vert{\f{A}}\Vert\lo$ of $\Vert{\f{A}}\Vert$ is just the ordinary collage $|{\f{A}}\lo|$ as described in~\ref{subsec:coll-cocplt}.

We claim that the rest of $\Vert{\f{A}}\Vert$ admits the following description.
The set of objects of $\Vert{\f{A}}\Vert\ti$ is the disjoint union of the objects of the categories $(\epsilon x)\ti$.
A morphism in $\Vert{\f{A}}\Vert\ti$ from $\omega\in(\epsilon x)\ti$ to $\varpi\in(\epsilon y)\ti$ is a morphism $\omega \to {\f{A}}(x,y)(\varpi)$ in $(\epsilon x)\ti$.
When ${\f{A}}$ has one object, this is just the ordinary Kleisli category of a ($\kappa$-cocontinuous) monad (i.e.\ the Kleisli object in $\cat{Cat}$, as opposed to the Kleisli object in $\cat{Colim}_\kappa$, which as we have seen is the ordinary Eilenberg--Moore category).
In the general case, this describes precisely the collage of ${\f{A}}$ regarded as a $\cat{Cat}$-category (i.e.\ forgetting the fact that the functors ${\f{A}}(x,y)$ are $\kappa$-cocontinuous).
In particular, we have coprojections $(\epsilon x)\ti \to \Vert{\f{A}}\Vert\ti$ for all $x\in\ob{\f{A}}$.

The functor $j_{\Vert{\f{A}}\Vert}\colon \Vert{\f{A}}\Vert\ti \to \Vert{\f{A}}\Vert\lo$ sends an object $\omega\in(\epsilon x)\ti$ to the tuple $(\xi_z \coloneqq {\f{A}}(z,x)(\omega))_{z}$, and sends a morphism $\omega \to {\f{A}}(x,y)(\varpi)$ to the tuple of composites
\[j_{\Vert{\f{A}}\Vert}(\omega)_z = {\f{A}}(z,x)(\omega) \longrightarrow {\f{A}}(z,x)({\f{A}}(x,y)(\varpi)) \longrightarrow {\f{A}}(z,y)(\varpi) = j_{\Vert{\f{A}}\Vert}(\varpi)_z.\]
When ${\f{A}}$ has one object, this is just the ordinary inclusion of the Kleisli category of a ($\kappa$-cocontinuous) monad into its Eilenberg--Moore category, which is well-known to be fully faithful; the general case is fully faithful for the same reason.
The desired universal property for $\Vert{\f{A}}\Vert$ follows by combining the universal property of $\Vert{\f{A}}\Vert\ti$, as a collage in $\cat{Cat}$, with that of $\Vert{\f{A}}\Vert\lo$, as a collage in $\cat{Colim}_\kappa$.

\begin{Rk}
 By Theorem~\ref{thm:collages2}, a tight collage of a $\C$-category $\f A$ may be equivalently characterized as a functor $R\colon \f A \to\hat{v}$ whose underlying module is an equivalence (i.e.\ $R$ is a ``Morita equivalence'') and which ``detects representability of modules''.
  In particular, if $\C$ has tight $\kappa$-ary collages, then the inclusion $\hat{(-)}\colon \C \to \fcat{Mod}_\kappa(\C)$ of $\F$-bicategories has a left $\F$-biadjoint.
  Once we show that $\fcat{Mod}_\kappa(\C)$ has tight $\kappa$-ary collages, the converse will follow, since a left $\F$-biadjoint must preserve $\F$-colimits such as tight collages.
\end{Rk}

The following lemma says that the construction $\fcat{Mod}_\kappa$ preserves map equipments.
\begin{Lemma}\label{lemma:fmod-preserves-equip}
  If $D\colon \f A \to \f B$ is a $\C$-functor such that each morphism $D_x$ has a loose right adjoint, then $\f B(1,D)$ has a loose right adjoint in $\fcat{Mod}_\kappa(\C)$.
  In particular, if $\C$ is a map equipment, then so too is $\fcat{Mod}_\kappa(\C)$.
\end{Lemma}
\begin{proof}
  Write $D_x^*$ for the right adjoint of $D_x$.
  We define an $\f A$-$\f B$ module $\f B(D,1)$, whose 1-cell components are $D_x^* \circ \f B(D x,u)$.
  The rest of the structure is completely dual to that of $\f B(1,D)$, and the adjunction $\f B(1,D) \dashv \f B(D,1)$ is easy to check.
\end{proof}

\begin{Thm}\label{thm:free-tight-collages}
  For any $\F_\kappa$-bicategory $\C$, the free cocompletion of $\C$ under tight $\kappa$-ary collages is equivalent to $\fcat{Mod}_\kappa(\C)$.
\end{Thm}
\begin{proof}
  As in Theorem~\ref{thm:free-collages}, the desired free cocompletion is the closure of $\C$ in $\M\C$ under tight $\kappa$-ary collages, which contains the full sub-$\F$-bicategory of $\M\C$ whose objects are tight collages of the image of some $\kappa$-ary $\C$-category.
  We will show that the latter is equivalent to $\fcat{Mod}_\kappa(\C)$ and that it is closed in $\M\C$ under tight $\kappa$-ary collages, and hence coincides with the desired closure.

  We proceed exactly as in Theorem~\ref{thm:free-collages}, considering the equivalent $\F_\kappa$-bicategory whose objects are literally the $\kappa$-ary $\C$-categories.
  For ${\f{A}}$ a $\kappa$-ary $\C$-category, the tight collage $\Vert Y{\f{A}}\Vert$ of its image in $\M\C$ can be defined by
  \[ \Vert Y{\f{A}}\Vert(c) \coloneqq \Vert\C(c,{\f{A}})\Vert. \]
  Here $\C(c,{\f{A}})$ is an $\F_\kappa$-category that enhances Example~\ref{ex:represented-Catcat} in an obvious way, while $\Vert\thg\Vert$ denotes its collage in $\F_\kappa$.

  Now suppose ${\f{A}}$ and ${\f{B}}$ are two $\kappa$-ary $\C$-categories.
  By the universal property of collages, a loose morphism $\Vert Y{\f{A}}\Vert  \to \Vert Y{\f{B}}\Vert$ in $\M\C$ is determined by a lax cocone of loose morphisms under $Y{\f{A}}$ with vertex $\Vert Y{\f{B}}\Vert$, i.e.\ by a collection of loose right $\C$-module morphisms $Y(\epsilon x) \to \Vert Y{\f{B}}\Vert$, for $x\in\ob{\f{A}}$, together with associative module transformations.
  By the Yoneda lemma, a loose right $\C$-module morphism $Y(\epsilon x) \to \Vert Y{\f{B}}\Vert$ is uniquely determined by a loose object of $\Vert Y{\f{B}}\Vert (\epsilon x) = \Vert\C(\epsilon x,{\f{B}})\Vert$.
  Since we saw in Section~\ref{sec:tcoll-Fk} that the loose parts of tight collages in $\F_\kappa$ are simply collages in $\cat{Colim}_\kappa$, the proof of Theorem~\ref{thm:free-collages} carries all the way through to show that the bicategory of loose morphisms is equivalent to $\cat{Mod}_\kappa(\C\lo)$.

  On the other hand, a \emph{tight} morphism $\Vert Y{\f{A}}\Vert  \to \Vert Y{\f{B}}\Vert$ is determined by the same data such that each right $\C$-module morphism $Y(\epsilon x) \to \Vert Y{\f{B}}\Vert$ is tight, or equivalently that the corresponding object of $\Vert\C(\epsilon x,{\f{B}})\Vert$ is tight.
  By Section~\ref{sec:tcoll-Fk} and the definition of $\C(\epsilon x, \f B)$, a tight object of $\Vert\C(\epsilon x,{\f{B}})\Vert$ consists of an object of $\f B$---call it, say, $D x$---together with a tight morphism $\epsilon x \to \epsilon(D x)$---call it, say, $D_x$.

  The additional data of a tight morphism $\Vert Y{\f{A}}\Vert  \to \Vert Y{\f{B}}\Vert$ consists of, for each $x,y\in\ob\f A$, a morphism from $(D x,D_x)$ to $(D y,D_y)$ in $\Vert\C(\epsilon x,{\f{B}})\Vert\ti$.
  By Section~\ref{sec:tcoll-Fk} and the definition of $\C(\epsilon x, \f B)$, such a morphism consists of a 2-cell $D_x \to \f B(D x, E x) \circ E_x$ in $\C$.
  The associativity of these morphisms reduces exactly to the axioms of a $\C$-functor.
  (This is the same argument used in~\cite{ls:ftm2} for the case $\kappa=\{1\}$.)
  
  Finally, the description of $j_{\Vert\C(\epsilon x,{\f{B}})\Vert}$ in Section~\ref{sec:tcoll-Fk} identifies it exactly with the operation $D \mapsto \f B(1,D)$ described in Section~\ref{sec:equip-Ccats}.
  Thus, a morphism $\Vert Y{\f{A}}\Vert  \to \Vert Y{\f{B}}\Vert$ in $\M\C$ admits a tightening exactly when the corresponding $\f A$-$\f B$-module is representable.
  This is sufficient to identify the ``one-step'' closure of $\C$ in $\M\C$ under tight $\kappa$-ary collages with $\fcat{Mod}_\kappa(\C)$.

  It remains to show that $\fcat{Mod}_\kappa(\C)$ is closed in $\M\C$ under tight $\kappa$-ary collages, so that this ``one-step'' closure is in fact the full closure. This is proved in the  following lemma.
\end{proof}

\begin{Lemma}
  Let $\C$ be an $\F_\kappa$-bicategory.
  Let $\f B$ be a $\kappa$-ary $\C$-category and $S\colon \f B \to \hat{v}$ a tight collage.
  For each $x\in\ob\f B$, let $\f A_x$ be a $\kappa$-ary $\C$-category and $R_x\colon \f A_x \to \widehat{\epsilon x}$ a tight collage.
  Then there is a $\kappa$-ary $\C$-category $\f A$ with object set $\sum_{x\in\ob\f B} \ob \f A_x$, and with extents induced from those of each $\f A_x$, whose tight collage is $\hat{v}$.
\end{Lemma}
\begin{proof}
  Note that since $R_x$ is a (tight) collage, each tight coprojection $R_z \colon \epsilon z \to \epsilon x$, for $z\in\ob\f A_x$, is a map in $\C\lo$.
  We denote its (loose) right adjoint by $R_z^*$.
  Now the objects and extents of $\f A$ are given; for $z\in\ob\f A_x$ and $w\in\ob \f A_y$, we define
  \[ \f A(z,w) \coloneqq R_z^* \circ \f B(x,y) \circ R_w. \]
  Then $\f A$ is a $\kappa$-ary $\C$-category, with structure morphisms induced by those of $\f B$ and the units and counits of the adjunctions $R_z\dashv R_z^*$.
  The adjunct morphisms to the definition of $\f A(z,w)$:
  \[ R_z \circ \f A(z,w) \longrightarrow \f B(x,y) \circ R_w \]
  make $R$ into a $\C$-functor $\f A \to \f B$.
  Since each $R_z$ is a map, by Lemma~\ref{lemma:fmod-preserves-equip}, the $\f B$-$\f A$-module $\f B(1,R)$ has a right adjoint $\f B(R,1)$.
  The usual split-coequalizer arguments show that $\f B(1,R)$ and $\f B(R,1)$ are inverse equivalences in $\fcat{Mod}_\kappa(\C)$.
  However, by Theorem~\ref{thm:collages2}, $S$ (or equivalently $\hat{v}(1,S)$) is also an equivalence from $\f B$ to $\hat{v}$.
  Thus, the composite functor $S R$ is equivalence from $\f A$ to $\hat{v}$, and hence is a loose collage.

  It remains to prove the tight part of the universal property.
  Suppose $g\colon v \to w$ is a loose morphism in $\C$.
  Since $S$ is a tight collage, $g$ is tight if and only if $g S_x$ is tight for each $x\in\ob\f B$.
  But since $R_x$ is a tight collage, $g S_x$ is tight if and only if $g S_x R_z$ is tight for each $z\in\ob\f A_x$.
  Thus, $g$ is tight if and only if each $g S_x R_z$ is tight.
\end{proof}

\begin{Rk}
  Suppose that we wanted to let $\f A(x,y)$ denote a morphism in $\C$ from $\epsilon x$ to $\epsilon y$, rather than the reverse, as we have done.
  Then in order for $\C$-transformations to reduce to their usual meaning in examples, a $\C$-transformation $D\Rightarrow E$ would have to involve a 2-cell $E_x \to \f B(D x, E x) \circ D_x$ rather than $D_x \to \f B(D x, E x) \circ E_x$.
  However, there seems no way to reverse the order of $D$ and $E$ in the 2-cells of a free cocompletion.
  Thus, in order for Theorem~\ref{thm:free-tight-collages} to hold, we are essentially forced to take $\f A(x,y)$ to be a morphism from $\epsilon y$ to $\epsilon x$.
\end{Rk}

\begin{Rk}
  By Theorem~\ref{thm:tightcollages}, $\fcat{Mod}_\kappa(\C)$ also has all tight lax colimits of lax functors into $\C\lo$ with $\kappa$-small domain.
  Thus, as in Remark~\ref{remark:laxcolim-completion}, we can conclude that it is also the free cocompletion of $\C$ under such tight lax colimits.
\end{Rk}

\begin{Rk}\label{remark:fmatmod}
  As in Remark~\ref{remark:matmod}, we have $\fcat{Mod}_\kappa(\C) \simeq \fcat{Mod}_{\{1\}}(\fcat{Matr}_\kappa(\C))$, where $\fcat{Matr}_\kappa(\C)$ is the free cocompletion of $\C$ under tight $\kappa$-small coproducts.
  Its loose bicategory is $\cat{Matr}_\kappa(\C_\lambda)$ as in Remark~\ref{remark:matmod}, while a tight morphism from $\epsilon_1\colon O_1\to \ob\C$ to $\epsilon_2\colon O_2 \to \ob\C$ consists of a function $D\colon O_1 \to O_2$ and tight morphisms $D_x\colon \epsilon_1 x \to \epsilon_2(D x)$.
\end{Rk}

\begin{Rk}\label{remark:mealy}
  There are various other similar monoidal bicategories for which we could consider analogous free cocompletions.
  For instance, if we used $\F$-bicategories---that is, mere equipments without any local cocompleteness---then the free cocompletion of an equipment $\C$ under tight collages would consist of $\C$-categories, $\C$-functors, ``$\C$-Mealy morphisms'', and Mealy cells, as in Example~\ref{ex:mealy}.
  And if we used $\cat{Cat}$-bicategories---that is, ordinary bicategories---then the free cocompletion under collages would consist of $\C$-categories, $\C$-Mealy morphisms, and Mealy cells.

  The latter, in the case $\kappa=\{1\}$, is precisely what was observed in~\cite{ls:ftm2}.
  Thus, from the starting point of~\cite{ls:ftm2}, the main observations of the second part of the current paper are that (1) by using equipments instead of bicategories, we can recover functors exactly as the tight morphisms in the free cocompletion, and (2) by including local cocompleteness in the enrichment, we can recover the more common and useful notion of modules, instead of Mealy morphisms, as the (loose) morphisms in the free cocompletion.
\end{Rk}

\end{document}